\documentclass[]{interact}
\usepackage{amsmath}
\usepackage{amsfonts}
\usepackage{amssymb}
\usepackage{enumerate}
\usepackage{mathrsfs}

\newtheorem{theorem}{Theorem}[section]

\newtheorem{lemma}[theorem]{Lemma}

\theoremstyle{definition}
\newtheorem{definition}[theorem]{Definition}

\newtheorem{remark}[theorem]{Remark}

\newtheorem{assumption}{Assumption}[subsection]
\newtheorem{ilemma}[assumption]{Lemma}

\numberwithin{equation}{section}

\begin{document}

\title{Absolute Logic: an alternative framework for formal logic}

\author{
\name{Mauro Avon\thanks{CONTACT Mauro Avon. Email: avonma@gmail.com}}
\affil{orcid: 0000-0003-4368-4759\\
Spilimbergo, Italy}
}

\maketitle

\begin{abstract}
This paper introduces an alternative framework for formal logic, termed ``Absolute Logic", designed to overcome the intrinsic limitations and structural relativity of standard First-Order Logic (FOL). While FOL relies on external Tarskian structures to assign meaning and restricts quantification to a single, predetermined domain, the proposed system establishes an invariant, self-contained semantics where every symbol possesses an intrinsic meaning. By unifying the traditional syntactic dichotomy between terms and formulas into a single concept of ``expression" and evaluating them relative to formalized variable-expression contexts, the system closely mirrors the natural, cumulative nature of human mathematical deduction. Furthermore, we address the critical balance between expressive power and constructive utility by integrating recursion-theoretic constraints, demonstrating how computability theory acts as a necessary bound to preserve foundational validity. The consistency and structural properties of the resulting language are formally established, providing a novel perspective on non-hierarchical, absolute logical systems.
\end{abstract}


\section{Introduction}\label{Ch:intro}

This paper outlines a system or approach to mathematical logic which is different from the standard
one. By `the standard approach to logic' I mean the one presented in chapter~2 of Enderton's book~\cite{Enderton} and there named `First-Order Logic'.  The same approach is also outlined in chapter 2 of
Mendelson's book~\cite{Mendelson}, where it is named `Quantification Theory'.\\

An online article by W. Ewald~\cite{Ewald}, in the Stanford Encyclopedia of Philosophy, describes the process that led to the establishment of first-order logic as the standard system of mathematical logic. However, the conclusion is that there are no clear reasons why this occurred.\\

\begin{quotation}
How did first-order logic come to be regarded as a privileged logical system—that is, as (in some sense) the “correct” logic for investigations in foundations of mathematics? That question, too, is highly complicated. Even after the Gödel results were widely understood, logicians continued to work in type theory, and it took years before first-order logic attained canonical status. The transition was gradual, and cannot be given a specific date.\\
\end{quotation}

First-order logic has been around for many decades, but to date no absolute evidence has been
found that first-order logic is the best possible logic system. In this regard I may quote a stronger
statement at the beginning of Jos\`{e} Ferreir\'{o}s' paper `The road to modern logic – an interpretation'
(\cite{Ferreiros}).\\

\begin{quotation}
It will be my contention that, contrary to a frequent assumption (at least among philosophers), First-Order Logic is \emph{not}
a `natural unity', i.e. a system the scope and limits of which could be justified solely by rational argument.\\
\end{quotation}

Honestly, in my opinion, the approach to logic I am going to propose seems to be a `natural unity' much more
than first-order logic is. The basic idea behind this system is indeed to build a logical system that is as natural, general, and absolute as possible, and to have a faithful model of the human deductive process, as far as possible.\\

The proposed `system' seems `natural' enough to me in many respects, but I can't say for sure it's a truly general and absolute approach, or the only valid approach to logic. In fact, for instance, I believe that a logical system must satisfy some computability requirements. Although computability theory was born in the 1930s, therefore after mathematical logic and the formalization of first-order logic by Hilbert and Ackermann, when formalizing a logical system it is not possible to ignore basic concepts inherent in computability theory. I suspect that this very requirement could be an obstacle to the possibility of obtaining a general and absolute logical system, or a unique approach to logic.\\

However, the alternative framework proposed here intentionally dispenses with certain architectural features of first-order logic that introduce artificial limitations or relativity. Chief among these are the rigid constraints on the `order' of expressions—which restrict reasoning to a single, predetermined domain—and the reliance on external `structures' to assign meaning.\\

Let's first discuss these two features.\\ 

In first-order logic variables range over individuals, but in mathematics there are statements in
which both quantifiers over individuals and quantifiers over sets of individuals occur. One simple
example is the following condition:\\

for each subset X of $\mathbb{N}$ and for each x $\in \mathbb{N}$ we have x $\in$ X or x $\notin$ X .\\

We will explicitly show in section~\ref{Ch:expMixed} that this assertion can easily be expressed in our language. Another example is the condition in which we state that every bounded, non empty set of real
numbers has a supremum. Formalisms which are better suited than first-order logic to express such conditions are second-order logic and type theory, but these systems have a certain level of complexity and are based on different types of variable. In our system we can express the conditions we mentioned above, and
we absolutely don't need different types of variables, the set to which the quantifier refers is
explicitly written in the expression, this ultimately makes things easier and allows a more general
approach. If we read the statement of a theorem in a mathematics book, usually in this statement
some variables are introduced, and when introducing them often the set in which they are varying is
explicitly specified, so from this point of view our approach is consistent with the actual processes
of mathematics. \\

Our logic is not a first-order or second-order or n-order logic, it doesn't involve types, so from this point of view it is an `absolute' type of logic.\\

Let's examine how our system behaves when giving a meaning and possibly a truth value to
expressions. Standard logic doesn't plainly associate meanings and truth values to formulas. It
introduces some related notion as the concepts of `structure' (defined in section 2.2 of Enderton's
book), truth in a structure, validity, satisfiability. Within first-order logic a structure is used, first of
all, to define the collection of things to which a quantifier refers to. Moreover, some symbols such
as connectives and quantifiers have a fixed meaning, while for other symbols the meaning is given
by the structure. Notions such as validity and satisfiability reveal a question-based approach: `what happens when we change the meaning of some symbols?' Although this may be an interesting perspective, this is not our approach, understanding what happens when we change the meaning of the symbols does not have a primary interest for us, although it's quite obvious that we'll also try to enunciate some results that are valid regardless of the meaning of the symbols. In this regard, if we had this perspective, in the first place it would have to be discussed if there are anyway symbols (e.g. connectives, quantifiers and others too) whose meaning cannot change. \\

Consequently, every symbol in the proposed system possesses an intrinsic, invariant meaning, thereby rendering the Tarskian notion of external structures obsolete. From this semantic perspective as well, the system qualifies as an `absolute' logic, since the truth-value of a sentence is determined exclusively by its constituent symbols.\\

We now list other features of our system, pointing out the differences and improvements with respect
to standard logic.\\

Unlike first-order logic, which maintains a strict syntactic dichotomy between terms and formulas, our approach unifies these notions into a single, comprehensive concept of `expression'. Furthermore, every expression is evaluated relative to a specific `context', formalised as a sequence of variable-expression pairs.\\

Our approach requires to build all at the same time, contexts, expressions, states and meanings.
We'll call sentences those expressions which are related to an empty context and whose meaning is
true or false. The meaning of a sentence depends solely on the meaning of the symbols it contains, it
doesn't depend on external `structures'.\\

In first-order logic we have terms and formulas and we cannot apply a predicate to one or more
formulas, and it seems this can be a limitation. With our system we can apply predicates to formulas.\\

Our deductive system seeks to provide a good model of human mathematical deductive process.
The concept of proof we'll feature is probably the most simple and intuitive that comes to mind, we
try to anticipate some of it.\\

If we define $S$ as the set of sentences then an axiom is a subset of $S$, an n-ary rule is a subset of $S^{n+1}$. If $\varphi$ is a sentence then a proof of $\varphi$ is a sequence $(\psi_1, \dots ,\psi_m )$ of sentences such that 

\begin{itemize}
\item there exists an axiom $A$ such that $\psi_1 \in A$ ;
\item if $m>1$ then for each $j=2 \dots m$ one of the following holds
\begin{itemize}
\item there exists an axiom $A$ such that $\psi_j \in A$ ,
\item there exists an n-ary rule $R$ and $i_1, \dots , i_n < j$ such that $(\psi_{i_1} , \dots , \psi_{i_n} , \psi_j ) \in R$;
\end{itemize}

\item $\psi_m = \varphi$ .\\

\end{itemize} 

As regards the soundness of the system, it is proved at the beginning of section~\ref{Ch:dedSysandProofs}. Consistency, proved in paragraph~\ref{sec:consistency}, is a direct consequence of soundness. We discuss (in paragraph~\ref{sec:completeness}) on the completeness of our deductive systems.\\

We have examined the main features of the system. If the reader will ask what is the basic idea
behind a system of this type, in agreement with what I said earlier I could say that the principle is to try to provide something like a general, absolute and unifying approach to logic and a faithful model of human
mathematical deductive process.\\

This statement about our system of course is not a mathematical statement, so I cannot give a
mathematical proof of it. I'm not even sure that I have truly and fully achieved the declared objectives and that they are fully achievable. A key aspect in this regard is the computability requirements that a logical system must satisfy, and in this version of the manuscript we pay due attention to these requirements.\\

On the other hand, logic exists with the specific primary purpose of being
a model to human deduction. In general, suppose we want to provide a mathematical model of some
process or reality. The fairness of the model can be judged much more through experience than
through mathematics. In fact, mathematics always has to do with models and not directly with
reality.\\

This paper's purpose is to present an approach to logic, but clearly we cannot provide here all
possible explanations and comparisons in any way related to the approach itself.
The author believes that this paper provides a fairly comprehensive presentation of the approach in
question, this introduction includes significant elements of explanation, justification and
comparison with the standard approach to logic. Other material in this regard is presented in the
subsequent parts.\\

Further investigations on this approach will be conducted, in the future, if and when possible, by the
author and/or other people. If any claim of this introduction would seem inappropriate, the author is
ready to reconsider and possibly fix it. In any case he believes the most important part of this paper
is not in the introduction, but in the subsequent sections.\\

The paper is quite long, but the time required to get an idea of the content is not very high. In fact, the author has
chosen to include all the proofs, but quite often they aren't difficult proofs. In addition, the most
complex part is perhaps definition~\ref{D:expr-new} which has a certain complexity, but at a first reading it is not necessary to take care of all the details.

\section{Changes from previous version}\label{Ch:changes}

Here we describe the main changes of the paper with respect to the previous version.\\

We decided to remove the set-builder notation from the system. It doesn't appear to be a fundamental part, and by eliminating it, we've simplified the system. This also makes the system a bit more similar to first-order logic, and comparisons may be easier.\\

We've made some changes to the classes of functional transformations included in the system, attempting to introduce constructs characterized by a high degree of generality and uniformity.\\

We also describe the important changes introduced in the previous version.\\

First of all, we have introduced computability constraints in the definition of the system. The process by which we generate expressions in our language is an inductive process. At each step we must ensure that the set of the new expressions related to a certain context is a recursive set. This ensures that the global set of expressions (related to a context) is a recursively enumerable set and so are the set of sentences etc.. We also introduced the constraint that axioms and rules must be r.e. sets, which seems reasonable.\\

Besides this we also added a new example of deduction.

\section{The language of our logic system}\label{Ch:lang}

In this section we want to define the language, which is the entity that underlies our logic system. The language is actually made up of various elements including some sets of symbols.\\

First we need a set of symbols $\mathcal{V}$. $\mathcal{V}$ members are also called ‘variables’ and just play the role of variables in the construction of our expressions (this implies that $\mathcal{V}$ members have no meaning
associated). We assume  $\mathcal{V}$ is a finite or countable set.\\

In addition we need another set of symbols $\mathcal{C}$. $\mathcal{C}$ members are also called `constants' and have a meaning. For each $c \in \mathcal{C}$ we denote by $\#(c)$ the meaning of $c$. We assume $\mathcal{C}$ is a finite set.\\

Let $f$ be a member of $\mathcal{C}$. Being $f$ endowed with meaning, $f$ is always an expression of our language. However, the meaning of $f$ could also be a function. In this case $f$ can also play the role of an `operator' in the construction of expressions that are more complex than the simple constant $f$.\\

Not all the operators that we need, however, are identifiable as functions. Think to the logical connectives (logical negation, logical implication, quantifiers, etc..), but also to the membership predicate `$\in$' and to the equality predicate `$=$'. The meaning of these operators cannot be mapped to a precise mathematical object, therefore these operators won't have a precise meaning in our language, but we'll need to give meaning to the application of the operator to objects, where the operator is applicable.\\

In mathematics and in the real world objects can have properties, such as having a certain color, or
being true, or being false. A property is therefore something that can be assigned to a single object, to no object, or to more than one object. For example, with reference to color, one or more objects are red or
have the property `to be of red color'. But more generally one or more objects have a color. Suppose
we denote, for objects $x$ that have a color, the color of $x$ with $C(x)$. So we can say that $C$ is a
property applicable to a class of objects. On the same object class we can indicate with $R(x)$ the
condition `$x$ has the red color'. $R$ is in turn a property applicable to a class of objects, with the
characteristic that for all $x$ $R(x)$ is true or false. A property with this additional feature can be called
a `predicate'.\\

The class of objects to which a property may be assigned may be called the domain of the property.
The members of that domain may be individual objects or sequences of objects, for example, if $x$ is
an object and $X$ is a set, the condition `$x \in X$' involves two objects, and then the domain of the membership property consists of the ordered pairs $(x,X)$, where $x$ is an object and $X$ is a set.\\
Generally we are dealing with properties such that the objects of their domain are all individual
objects, or all ordered pairs. Theoretically there may also be properties such that the objects of their domain
are sequences of more than two items or even the number of items in sequence may be different in
different elements of the domain.\\

As mentioned above the concept of `property' is similar to the concept of function, but in mathematics there are properties that are not functions. For example, the condition `$x \in X$' just introduced can be applied to an arbitrary object and an arbitrary set, so the `membership property' does not have a well-defined domain and cannot be considered a function in a strict sense.\\

So, in order to build our language, we need another set of symbols $\mathcal{F}$, where each $f$ in $\mathcal{F}$ represents a property $P_f$. Symbols in $\mathcal{F}$ are also called operators or `property symbols'. We assume $\mathcal{F}$ is a finite set. We will not assign a meaning to operators, because a property cannot be mapped to a consistent mathematical object (function or other). However, for each $f$ 

\smallskip

\begin{itemize}
\item we need to determine a condition $A_f(x_1, \dots ,x_n)$ that given a positive integer $n$ and $x_1, \dots , x_n$ arbitrary objects  indicates if $P_f$ is applicable to $x_1, \dots , x_n$. The condition $A_f(x_1, \dots ,x_n)$ does not have to be decidable in an absolute sense, but it must be so when it is used in the process by which we construct our expressions;
\item for each positive integer $n$ and $x_1, \dots , x_n$ arbitrary objects such that $A_f(x_1, \dots ,x_n)$ holds we must be able to calculate the value of $P_f(x_1, \dots ,x_n)$. This doesn't mean that $P_f$ must be a computable function in a strict sense, but we must be able to know the value of $P_f(x_1, \dots ,x_n)$ when this calculation is required in the construction of our expressions.\\
\end{itemize}

We immediately explain these concepts by specifying what are the most important operators that we may include in our language, providing for each of them the conditions $A_f(x_1, \dots ,x_n)$ and $P_f(x_1, \dots ,x_n)$ (in general $P_f(x_1, \dots ,x_n)$ is a generic value, but in these cases it is a condition, i.e. its value can be true or false).

\smallskip

\begin{itemize}
\item Logical conjunction: it's the symbol $\wedge$ and we have \\
for $n \ne 2$  $A_{\wedge}(x_1, \dots ,x_n)$ is false ,\\
$A_{\wedge}(x_1, x_2) =$ ( $x_1$ is true or $x_1$ is false ) and ( $x_2$ is true or $x_2$ is false ), \\
$P_{\wedge}(x_1, x_2) =$ both $x_1$ and $x_2$ are true ;

\medskip

\item Logical disjunction: it's the symbol $\vee$ and we have \\
for $n \ne 2$  $A_{\vee}(x_1, \dots ,x_n)$ is false ,\\
$A_{\vee}(x_1, x_2) =$ ( $x_1$ is true or $x_1$ is false ) and ( $x_2$ is true or $x_2$ is false ), \\
$P_{\vee}(x_1, x_2) =$ at least one between $x_1$ and $x_2$ is true ;

\medskip

\item Logical implication: it's the symbol $\to$ and we have \\
for $n \ne 2$  $A_{\to}(x_1, \dots ,x_n)$ is false ,\\
$A_{\to}(x_1, x_2) =$ ( $x_1$ is true or $x_1$ is false ) and ( $x_2$ is true or $x_2$ is false ), \\
$P_{\to}(x_1, x_2) =$ $x_1$ is false or $x_2$ is true ;

\medskip

\item Double logical implication: it's the symbol $\leftrightarrow$ and we have \\
for $n \ne 2$  $A_{\leftrightarrow}(x_1, \dots ,x_n)$ is false ,\\
$A_{\leftrightarrow}(x_1, x_2) =$ ( $x_1$ is true or $x_1$ is false ) and ( $x_2$ is true or $x_2$ is false ), \\
$P_{\leftrightarrow}(x_1, x_2) =$ $P_{\to}(x_1, x_2)$ and  $P_{\to}(x_2, x_1)$ ;

\medskip

\item Logical negation: it's the symbol $\neg$ and we have \\
for $n > 1$  $A_{\neg}(x_1, \dots ,x_n)$ is false ,\\
$A_{\neg}(x_1)$ is true, \\
$P_{\neg}(x_1) =$ $x_1$ is false ;

\medskip

\item Membership predicate: it's the symbol $\in$ and we have \\
for $n \ne 2$  $A_{\in}(x_1, \dots ,x_n)$ is false ,\\
$A_{\in}(x_1, x_2) =$ $x_2$ is a set, \\
$P_{\in}(x_1, x_2) =$ $x_1$ is a member of $x_2$ ;

\medskip

\item Equality predicate: it's the symbol $=$ and we have \\
for $n \ne 2$  $A_{=}(x_1, \dots ,x_n)$ is false ,\\
$A_{=}(x_1, x_2)$ is true, \\
$P_{=}(x_1, x_2) =$ $x_1$ is equal to $x_2$ .\\

\end{itemize}

In principle we can think and use also other operators. First of all we could think to the universal and existential quantifier as operators, the universal quantifier could be represented by the symbol `$\forall$' and we would require the following:

\begin{itemize}
\item for $n > 1$  $A_{\forall}(x_1, \dots ,x_n)$ is false,\\
$A_{\forall}(x_1) = $ $x_1$ is a set and for each $x$ in $x_1$ ($x$ is true or $x$ is false), \\
$P_{\forall}(x_1) =$ for each $x$ in $x_1$ ($x$ is true).
\end{itemize}

Similarly, the existential quantifier would be represented by `$\exists$'; however, to avoid excessive complexity of the system we prefer not to use those operators. In our system of course we will still be able to perform universal and existential quantification and we will still use to this end the symbols `$\forall$' and `$\exists$'.\\

Operations between sets such as union or intersection could still be represented through an operator, etc.. In any case, we must choose our operators in such a way as to ensure computability in the construction of our expressions, and for this reason we must impose limits on the choice of operators. For example, set operators of the type just mentioned will not be used.\\

Our set $\mathcal{F}$ will typically be contained in the set $\{ \neg, \wedge, \vee, \to, \leftrightarrow, \in, = \}$, where each of the just mentioned symbols has been defined above. However, we want to have a more general approach than the one in which the operators are explicitly indicated, so we will also allow other types of operators, as long as they fall into one of the following categories.\\

The first admitted category of operators is the category of the symbols $f$ such that

\begin{itemize}
\item for $n \ne 2$  $A_f(x_1, \dots ,x_n)$ is false,
\item $A_f(x_1, x_2) =$ ( $x_1$ is true or $x_1$ is false ) and ( $x_2$ is true or $x_2$ is false ),
\item $P_f(x_1, x_2)$ is true or false.
\end{itemize}

\medskip

Since for $n \ne 2$  $A_f(x_1, \dots ,x_n)$ is false, we say the symbols in this category have a multiplicity of 2.\\

All of the symbols $\wedge, \vee, \to, \leftrightarrow$ fall within this category.\\

Another admitted category of operators is the category of the symbols $f$ such that

\begin{itemize}
\item for $n > 1$  $A_f(x_1, \dots ,x_n)$ is false,
\item $A_f(x_1)$ is true,
\item $P_f(x_1)$ is true or false.
\end{itemize}

\medskip

Since for $n > 1$  $A_f(x_1, \dots ,x_n)$ is false, we say the symbols in this category have a multiplicity of 1.\\

The symbol $\neg$ falls within this category.\\

Another admitted category of operators is the category of the symbols $f$ such that

\begin{itemize}
\item for $n \ne 2$  $A_f(x_1, \dots ,x_n)$ is false,
\item $A_f(x_1, x_2) =$ $x_2$ is a set,
\item $P_f(x_1, x_2)$ is true or false.
\end{itemize}

\medskip

Clearly the symbols in this category have a multiplicity of 2.\\

The symbol $\in$ falls within this category.\\

Finally, another admitted category of operators is the category of the symbols $f$ such that

\begin{itemize}
\item for $n \ne 2$  $A_f(x_1, \dots ,x_n)$ is false,
\item $A_f(x_1, x_2)$ is true,
\item $P_f(x_1, x_2)$ is true or false.
\end{itemize}

\medskip

Clearly the symbols in this category have a multiplicity of 2.\\

The symbol $=$ falls within this category.\\

We require that all the symbols in $\mathcal{F}$ fall within one of the mentioned categories, and so they must have a multiplicity of 1 or 2.\\

As we said above the quantifier symbols `$\forall$' and `$\exists$' are included in our system, but not with the role of operator. Besides the quantifier symbols, we need parentheses and commas to avoid ambiguity in the reading of our expressions; to this end we use the following symbols: left parenthesis `(', right parenthesis `)', comma `,' and colon `:'. We can indicate this further set of symbols with $\mathcal{Z}$.\\

To avoid ambiguity in reading our expressions we require that the sets $\mathcal{V}$, $\mathcal{C}$, $\mathcal{F}$ and $\mathcal{Z}$ are disjoint. It is also requested that no symbol correspond to a concatenation of more symbols of the language. More generally, given $\alpha_1, \dots ,\alpha_n$ and $\beta_1, \dots ,\beta_m$ symbols of our language, and using the symbol `$\|$' to indicate the concatenation of characters and strings, we assume that equality of the two concatenations $\alpha_1 \| \dots \|\alpha_n$ and $\beta_1 \| \dots \|\beta_m$ is achieved when and only when $m = n$ and for each $i = 1 \dots n \ \alpha_i = \beta_i$. We also specify that by `string' we mean a concatenation of symbols of our language.\\

While the set $\mathcal{Z}$ will be always the same, the sets $\mathcal{V}$, $\mathcal{C}$, $\mathcal{F}$ may change according to what is the language that we describe. If we think to our language as a language as defined in formal language theory, once we have chosen $\mathcal{V}$, $\mathcal{C}$ and $\mathcal{F}$, the alphabet $\Sigma$ of our language is given by $\Sigma = \mathcal{V} \cup \mathcal{C} \cup \mathcal{F} \cup \mathcal{Z}$.\\

Another variable element that we add to our language is made by a finite number of non-empty sets $D_1, \dots , D_p$ such that:

\begin{itemize}
\item for each $i, j = 1 \dots p$ such that $i \ne j$ $D_i \ne D_j$ and $D_i \cap D_j = \emptyset$;
\item for each $i = 1 \dots p$ and for each $x \in D_i$ $x$ is not a set;
\item for each $i = 1 \dots p$ and for each $x \in D_i$ $x$ is not true and $x$ is not false.
\end{itemize}

Note that these additional sets may not be required, in which case $p=0$.\\

A notion that we will soon use in the continuation is the notion of power set. Given a set $A$ we'll indicate with $\mathcal{P}(A)$ the set of the subsets of $A$, but in our definition the empty set will not be a member of $\mathcal{P}(A)$, so $\mathcal{P}(A)$ for us is the set of the non empty subsets of $A$.\\

We also define $\mathcal{P}^q(A)$ for any positive integer $q$. Of course $\mathcal{P}^1(A) = \mathcal{P}(A)$ by definition, and given a positive integer $q$ $\mathcal{P}^{q+1}(A) = \mathcal{P}(\mathcal{P}^q(A))$.\\

A specific language of our logic system is described by its variable elements which are the sets $\mathcal{V}$, $\mathcal{C}$, $\mathcal{F}$, the function $\#$ which associates a meaning to every element of $\mathcal{C}$ and in addition the (potentially empty) set of sets $\{ D_1, \dots , D_p \}$. Therefore our language is identified by the 5-tuple $(\mathcal{V}, \mathcal{F}, \mathcal{C}, \#, \{D_1, \dots, D_p \})$. Since the `meaning' of an operator is not a mathematical object, operators must be seen as symbols which are tightly coupled with their meaning.\\

We also need to set some constraints on our constants, which must not refer to the empty set or to a set which has the empty set as a member and so on. In order to do that, we want to define formally some predicates that we'll soon use in the continuation. We actually define the following predicates.\\

$Set_1(x)$ = $x$ is a set.\\

$Event_1(x)$ = $x$ is true or $x$ is false.\\

Given a positive integer $q$

\begin{itemize}
\item $Set_{q+1}(x)$ = $x$ is a set and for each $u \in x$ $Set_q(u)$; 
\item $Event_{q+1}(x)$ = $x$ is a set and for each $u \in x$ $Event_q(u)$. 
\end{itemize}

\medskip

If $Set_1(x)$ holds we define\\

$NotEmpty_1(x) = (x \ne \emptyset)$.\\

Given a positive integer $q$, if $Set_{q+1}(x)$ holds we define\\

$NotEmpty_{q+1}(x)$ = $NotEmpty_1(x)$ and for each $u \in x$ $NotEmpty_q(u)$.\\

The constraints we want to put on our constants can now be stated as follows: for each $c \in \mathcal{C}$, for each $i = 1 \dots p$ and for any positive integer $q$ we must be able to decide all of the following conditions

\begin{itemize}
\item $Set_q( \#(c))$;
\item $Event_q( \#(c))$;
\item $\#(c) \in D_i$;
\item $\#(c) \in \mathcal{P}^q(D_i)$;
\item if ($Set_q( \#(c))$) then ($NotEmpty_q( \#(c))$).
\end{itemize}

Moreover, the last condition must be decided as true.

\subsection{Other definitions and results}

Before we can describe the process of constructing expressions we still need to introduce some notation. In fact in that process we'll use the notion of `context' and the notion of `state'. Context and states have a similar form, here we define a notion of state-like pair and related results that will apply to states, but similar definitions and results will be given for contexts.\\

We define $\mathcal{D} = \{\emptyset\} \cup \{ \{1, \dots , m \} \vert \ m \text{ is a positive integer} \} $. 

\medskip

Suppose $x$ is a function whose domain $dom(x)$ belongs to $\mathcal{D}$. Suppose $C \in \mathcal{D}$ is such that $C \subseteq dom(x)$. Then we define $x_{/C}$ as a function whose domain is $C$ and such that for each $j \in C \ x_{/C}(j) = x(j)$ .

\medskip

Suppose $x$ and $\varphi$ are two functions with the same domain $D$, and $D \in \mathcal{D}$. Then we say that $(x,\varphi)$ is a `\emph{state-like pair}'. 

\medskip

Given a state-like pair $k = (x,\varphi)$ the domain of $x$ will be also called the \emph{domain of k}. Therefore $dom(k) = dom(x) = dom(\varphi)$.

\medskip

Furthermore $dom(k) \in \mathcal{D}$ and given $C \in  \mathcal{D}$ such that $C \subseteq dom(k)$ we can define \linebreak $k_{/C} = (x_{/C}, \varphi_{/C})$. Clearly $k_{/C}$ is a state-like pair.

\medskip

We define $\mathcal{R}(k) = \{ k_{/C} \vert \ C \in  \mathcal{D}, C \subseteq dom(k)   \}$.

\medskip

Given another state-like pair $h$ we write $h \sqsubseteq k$ if and only if $h \in \mathcal{R}(k) $ .

\medskip

Suppose $h \in \mathcal{R}(k)$, then there exists $C \in \mathcal{D}$ such that $C \subseteq dom(k)$, $h=k_{/C}=(x_{/C}, \varphi_{/C})$. Therefore $dom(h) = C$ and $k_{/dom(h)} = k_{/C} = h$.

\medskip

Suppose $h \in \mathcal{R}(k)$ and $g \in \mathcal{R}(h)$. This means there exist $C \in \mathcal{D}$ such that $C \subseteq dom(k), \ h = k_{/C}$, and there exist $D \in \mathcal{D}$ such that $D \subseteq dom(h), \ g = h_{/D}$. So $D \subseteq dom(h) = C \subseteq dom(k)$, $g = (k_{/C})_{/D} = (x_{/C}, \varphi_{/C})_{/D} = (x_{/D}, \varphi_{/D}) = k_{/D}$. Therefore $g \in \mathcal{R}(k)$.
\\

Suppose $k = (x,\varphi)$ is a state-like pair whose domain is $D$. Suppose $(y,\psi)$ is an ordered pair. Then we can define the `addition' of $(y,\psi)$ to $k$.\\
Suppose $D = \{1, \dots , m\}$, then we define $D' = \{1, \dots , m+1\}$. We define $x'$ as a function whose domain is $D'$ such that for each $\alpha = 1 \dots m \ x'(\alpha)=x(\alpha)$, and $x'(m+1) = y$. We define $\varphi'$ as a function whose domain is $D'$ such that for each $\alpha = 1 \dots m \ \varphi'(\alpha)=\varphi(\alpha)$, $\varphi'(m+1) = \psi$. Then we define $k + (y,\psi) = (x',\varphi')$. Obviously $(k + (y,\psi))_{/\{1, \dots, m \}} = k$, so $k \in \mathcal{R}(k + (y,\psi))$.\\
If $D = \emptyset$ then clearly $D' = \{1\}$. We define $x'$ as a function whose domain is $D'$ such that $x'(1) = y$. We define $\varphi'$ as a function whose domain is $D'$ such that $\varphi'(1) = \psi$. Then we define $k + (y,\psi) = (x',\varphi')$. Obviously $(k + (y,\psi))_{/\emptyset} = \emptyset = k$, so $k \in \mathcal{R}(k + (y,\psi))$.\\
In both cases $k + (y,\psi)$ is a state-like pair and $k \in \mathcal{R}(k + (y,\psi))$, which implies $dom(k) \subseteq dom(k + (y,\psi))$.\\
We have also seen that $(k + (y,\psi))_{/dom(k)} = (k + (y,\psi))_{/D} = k$.

\medskip

We also define $\epsilon = (\emptyset,\emptyset)$, so $\epsilon$ is a state-like pair.

\medskip

Given a state-like pair $k = (x, \varphi)$ we define $var(k)$ as the image of the function $x$. In other words if $k = \epsilon$ then $x = \emptyset$, so $var(k) = \emptyset$, otherwise $x$ has a domain $\{1, \dots ,m\}$ and $var(k) = \{x_i | i = 1 \dots m \}$.\\

Clearly, if we assume that $k + (y,\psi) = (x',\varphi')$, we can easily see that 
\[ var(k + (y,\psi)) = \{ x_i' | i \in dom(x_i') \} = \{ x_i | i \in dom(x_i) \} \cup \{ y \} = var(k) \cup \{ y \}  . \]

\bigskip

In the next lemma we prove that, when a state-like pair is obtained as $k + (y,\psi)$, then $k$, $y$, and $\psi$ are univocally determined.\\

\begin{lemma}\label{L:add-slp-univ}
Suppose $k_1 = (x_1,\varphi_1)$ is a state-like pair whose domain is $D_1$, and $(y_1,\psi_1)$ is an ordered pair. Suppose $k_2 = (x_2,\varphi_2)$ is a state-like pair whose domain is $D_2$, and $(y_2,\psi_2)$ is an ordered pair. Finally suppose $k_1 + (y_1,\psi_1) = k_2 + (y_2,\psi_2)$. Under these assumptions we can prove that $k_1 = k_2, y_1=y_2, \psi_1=\psi_2$. 
\end{lemma}

\begin{proof}

We define $h = k_1 + (y_1,\psi_1) = k_2 + (y_2,\psi_2)$. Since $h = k_1 + (y_1,\psi_1)$ we can have two possibilities:

\begin{itemize}
\item $D_1 = \emptyset, D_1' = \{1 \}$ and there exist two functions $x_1'$ and $\varphi_1'$ whose domain is $D_1'$ such that $h = (x_1',\varphi_1')$ ;
\item there exists a positive integer $m_1$ such that $D_1 = \{1, \dots , m_1\}, D_1' = \{1, \dots , m_1+1\}$ and there exist two functions $x_1'$ and $\varphi_1'$ whose domain is $D_1'$ such that $h = (x_1',\varphi_1')$.
\end{itemize}

Similarly, since $h = k_2 + (y_2,\psi_2)$ we can have two possibilities:

\begin{itemize}
\item $D_2 = \emptyset, D_2' = \{1 \}$ and there exist two functions $x_2'$ and $\varphi_2'$ whose domain is $D_2'$ such that $h = (x_2',\varphi_2')$ ;
\item there exists a positive integer $m_2$ such that $D_2 = \{1, \dots , m_2\}, D_2' = \{1, \dots , m_2+1\}$ and there exist two functions $x_2'$ and $\varphi_2'$ whose domain is $D_2'$ such that $h = (x_2',\varphi_2')$.
\end{itemize}

It follows that $(x_1',\varphi_1') = h = (x_2',\varphi_2')$, so $x_1' = x_2'$ and $\varphi_1' = \varphi_2'$, and $D_1' = D_2'$.

\medskip

Suppose $D_1 = \emptyset$. This implies that $D_2' = D_1' = \{1 \}$, thus $D_2 = \emptyset$.\\
In this case $k_1 = \epsilon = k_2$, $y_1 = x_1'(1) = x_2'(1) = y_2$, $\psi_1 = \varphi_1'(1) = \varphi_2'(1) = \psi_2$ .

\medskip

Suppose there exists a positive integer $m_1$ such that $D_1 = \{1, \dots , m_1\}$. This implies that $D_2' = D_1' = \{1, \dots , m_1+1\}$, thus $D_2 = \{1, \dots , m_1\}$.\\
In this case for each $\alpha = 1 \dots m_1$ $x_1(\alpha) = x_1'(\alpha) = x_2'(\alpha) = x_2(\alpha)$, $\varphi_1(\alpha) = \varphi_1'(\alpha) = \varphi_2'(\alpha) = \varphi_2(\alpha)$ . So $k_1 = (x_1,\varphi_1) = (x_2,\varphi_2) = k_2$; and moreover $y_1 = x_1'(m_1+1) = x_2'(m_1+1) = y_2$, $\psi_1 = \varphi_1'(m_1+1) = \varphi_2'(m_1+1) = \psi_2$ .
\end{proof}

\bigskip

Other useful results are the following.

\medskip

\begin{lemma}\label{L:useful-slp-0}
Suppose $h= (x, \varphi), \ k = (z, \psi)$ are state-like pairs such that $h \in \mathcal{R}(k)$. Then, for each $j \in dom(h)$ 
$x_j = z_j$ and $\varphi_j = \psi_j$.
\end{lemma}

\begin{proof}

There exists $C \in \mathcal{D}$ such that $C \subseteq dom(k)$, $h=k_{/C}=(z_{/C}, \psi_{/C})$. Therefore $x = z_{/C}$ and $\varphi = \psi_{/C}$. For each $j \in dom(h) = C$ $x_j = z_j$ and $\varphi_j = \psi_j$.
\end{proof}

\medskip

\begin{lemma}\label{L:useful-slp}
Suppose $h= (x, \varphi), \ k = (z, \psi)$ are state-like pairs such that $h \in \mathcal{R}(k)$ and for each $i, j \in dom(k)$ $i \ne j \to z_i \ne z_j$. Then, for each $i \in dom(k), \ j \in dom(h)$ $z_i = x_j \to \psi_i = \varphi_j$.
\end{lemma}

\begin{proof}
Let $i \in dom(k), \ j \in dom(h)$ and $z_i = x_j$. Clearly $j \in dom(k), \ x_j = z_j$, thus $z_i = z_j$, $i = j$, $\varphi_j = \psi_j = \psi_i$.
\end{proof}

\medskip

\begin{lemma}\label{L:useful-slp-ext}
Suppose $k = (x, \varphi)$ and $h = (y, \psi)$ are state-like pairs such that for each $i \in dom(k)$, $j \in dom(h)$ $x_i = y_j \to \varphi_i = \psi_j$. Suppose $(u, \theta)$ is an ordered pair and $u \notin var(k)$, $u \notin var(h)$. Let $k' = k + (u, \theta)$ and $h' = h + (u, \theta)$. Let also $k' = (x', \varphi')$ and $h' = (y', \psi')$, then for each $i \in dom(k')$, $j \in dom(h')$ $x'_i = y'_j \to \varphi'_i = \psi'_j$.
\end{lemma}

\begin{proof}
Let $i \in dom(k')$, $j \in dom(h')$ such that $x'_i = y'_j$. 

\medskip

Suppose $i \in dom(k)$. If $j \notin dom(h)$ then $x'_i = x_i \in var(k)$, $y'_j = u \notin var(k)$ so $x'_i \ne y'_j$. So $j \in dom(h)$ and $\varphi'_i = \varphi_i = \psi_j = \psi'_j$.

\medskip

Suppose $i \notin dom(k)$. If $j \in dom(h)$ then $x'_i = u \notin var(h)$ and $y'_j = y_j \in var(h)$, so $x'_i \ne y'_j$.
Then obviously also $j \notin dom(h)$ and $\varphi'_i = \theta = \psi'_j$.
\end{proof}

\bigskip

\begin{lemma}\label{L:useful-sub-slp}
Suppose $k = (x, \varphi)$ and $h = (y, \vartheta)$ are state-like pairs such that for each $i \in dom(k)$, $j \in dom(h)$ $x_i = y_j \to \varphi_i = \vartheta_j$. Suppose $\kappa = (z, \phi) \sqsubseteq k$ and $g = (w, \theta) \sqsubseteq h$. Then for each $i \in dom(\kappa)$, $j \in dom(g)$ $z_i = w_j \to \phi_i = \theta_j$.
\end{lemma}

\begin{proof}

\smallskip

There exists $C \in \mathcal{D}$ such that $C \subseteq dom(k)$, $\kappa = k_{/C} = (x_{/C}, \varphi_{/C})$. Therefore $dom(\kappa) = C \subseteq dom(k)$.

\medskip

Similarly there exists $D \in \mathcal{D}$ such that $D \subseteq dom(h)$, $g = h_{/D} = (y_{/D}, \vartheta_{/D})$. Therefore $dom(g) = D \subseteq dom(h)$.

\medskip

Let $i \in dom(\kappa)$, $j \in dom(g)$, $z_i = w_j$, then $i \in dom(k)$, $j \in dom(h)$,
\[ x_i = (x_{/C})_i = z_i = w_j = (y_{/D})_j = y_j \ . \]

Then 
\[ \phi_i = (\varphi_{/C})_i = \varphi_i = \vartheta_j = (\vartheta_{/D})_j = \theta_j \ . \]
\end{proof}

\bigskip

\begin{lemma}\label{L:useful1-slp}
Suppose $h = (x, \varphi)$ is a state-like pair, $(y, \phi)$ is an ordered pair and define $k = h + (y, \phi)$. Suppose $g \in \mathcal{R}(k)$ is such that $g \ne k$. Then $g \in \mathcal{R}(h)$.
\end{lemma}

\begin{proof}

\smallskip

Let $D = dom(h)$. 

\medskip

Suppose $m$ is a positive integer and $D = \{ 1, \dots , m \}$. Then $k = (x', \varphi')$ has a domain \linebreak $\{ 1, \dots , m+1 \}$. Moreover there exists $C \in \mathcal{D}$ such that $C \subseteq \{ 1, \dots , m+1 \}$ and $g = k_{/C}$. Since $g \ne k$ we must have $C \subseteq \{ 1, \dots , m \}$. We have
\[ g = k_{/C} = (x'_{/C}, \varphi'_{/C}) = ((x'_{/D})_{/C} , (\varphi'_{/D})_{/C} ) = ( x_{/C}, \varphi_{/C} ) = h_{/C} \ . \]

\medskip

Now suppose $D = \emptyset$. Then $k = (x', \varphi')$ has a domain $\{ 1 \}$. Moreover there exists $C \in \mathcal{D}$ such that $C \subseteq \{ 1 \}$ and $g = k_{/C}$. Since $g \ne k$ we must have $C = \emptyset$ and $g = (\emptyset, \emptyset) = h$.

\medskip

In both cases $g \in \mathcal{R}(h)$.
\end{proof}

\bigskip

\begin{lemma}
Let $x$ be a function such that $dom(x) \in \mathcal{D}$, let $C, D \in \mathcal{D}$ such that \linebreak $C \subseteq D \subseteq dom(x)$. Then we can define $x_{/C}$ and $(x_{/D})_{/C}$, and we have $(x_{/D})_{/C} = x_{/C}$.
\end{lemma}

\begin{proof}
Define $y = x_{/D}$, we have $dom(y) = D$ and for each $j \in D$ $y(j) = x(j)$. Moreover $dom(y_{/C}) = C = dom(x_{/C})$ and for each $j \in dom(C)$ $y_{/C}(j) = y(j) = x(j) = x_{/C}(j)$.
\end{proof}

\bigskip

\begin{lemma}\label{L:useful5-slp}
Let $k = (x, \varphi)$ be a state-like pair, let $C, D \in \mathcal{D}$ such that $C \subseteq D \subseteq dom(k)$. Then we can define $k_{/C}$ and $(k_{/D})_{/C}$, and we have $(k_{/D})_{/C} = k_{/C}$.
\end{lemma}

\begin{proof}
\[ (k_{/D})_{/C} = (x_{/D}, \varphi_{/D})_{/C} = ( (x_{/D})_{/C}, (\varphi_{/D})_{/C}  ) = (x_{/C}, \varphi_{/C}) = k_{/C} . \]
\end{proof}

\bigskip

\begin{lemma}
Let $g, h, k$ be state-like pairs, let $g \sqsubseteq h, \ h \sqsubseteq k$. Then $g \sqsubseteq k$.
\end{lemma}

\begin{proof}

\smallskip

There exists $C \in \mathcal{D}$ such that $C \subseteq dom(h)$, $g = h_{/C}$.\\
There exists $D \in \mathcal{D}$ such that $D \subseteq dom(k)$, $h = k_{/D}$.

\medskip

This implies that $C \subseteq dom(h) = D$, so $g = h_{/C} = (k_{/D})_{/C} = k_{/C}$.

\medskip

Since $C \subseteq dom(k)$, $g \sqsubseteq k$.
\end{proof}

\bigskip

\begin{lemma}\label{L:useful3-slp}
Let $g, h$ and $k = (x, \varphi)$ be state-like pairs such that $g, h \in \mathcal{R}(k), \linebreak dom(g) \subseteq dom(h)$. Then $g \in \mathcal{R}(h)$.
\end{lemma}

\begin{proof}
There exists $C \in \mathcal{D}$ such that $C \subseteq dom(k)$, $g = k_{/C}$. And there exists $D \in \mathcal{D}$ such that $D \subseteq dom(k)$, $h = k_{/D}$. It results $C = dom(g) \subseteq dom(h) = D$. Then, clearly 
\[  g = (x, \varphi)_{/C} = (x_{/C}, \varphi_{/C}) = ( (x_{/D})_{/C}, (\varphi_{/D})_{/C} ) = (x_{/D}, \varphi_{/D})_{/C} = h_{/C}  \ .\]
\end{proof}

\bigskip

\begin{lemma}\label{L:useful4-slp}
Suppose $h = (x, \varphi)$ is a state-like pair, $(y, \phi)$ is an ordered pair and define $k = h + (y, \phi)$. Then $k_{/dom(h)} = h$.
\end{lemma}

\begin{proof}

\smallskip

Let $D = dom(h)$ and $k = (x', \varphi')$. Then $k_{/dom(h)} = (x'_{/D}, \varphi'_{/D}) = (x, \varphi) = h$.

\end{proof}

\bigskip

We also require some notation concerning generic strings, this notation will be useful when applied to
our expressions, which are non-empty strings. If $t$ is a string we can indicate with $\ell(t)$ $t$'s length, i.e. the number of characters in $t$. If $\ell(t) > 0$ then for each $\alpha \in \{1, \dots , \ell(t) \}$ at position $\alpha$ within $t$ there is a character, this symbol will be indicated with $t[\alpha]$. We call `depth of $\alpha$ within $t$' (briefly $d(t,\alpha)$) the number which is obtained by subtracting the number of right round brackets `)' that occur in $t$ before position $\alpha$ from the number of left round brackets `(' that occur in $t$ before position $\alpha$.\\

The following lemma will be useful later within proofs of unique readability.

\medskip

\begin{lemma}\label{L:depth-sum}
Let $\vartheta$, $\varphi$, $\eta$ be strings with $\ell(\vartheta)>0$, $\ell(\varphi)>0$, and let $t = \vartheta \| \varphi \| \eta$; let also $\alpha \in \{1, \dots , \ell(\varphi) \}$. The following result clearly holds: 
\[
d(t,\ell(\vartheta)+\alpha) =  d(t, \ell(\vartheta)+1) + d(\varphi,\alpha) .
\]
\end{lemma}

\begin{proof}
In order to prove this we provide some further simple definition.\\

Given a string $t$ such that $\ell(t) > 0$ and given $\alpha \in \{1, \dots , \ell(t) \}$ we define $L(t, \alpha)$ as the number of left round brackets `(' that occur in $t$ before position $\alpha$. We define $R(t, \alpha)$ as the number of right round brackets `)' that occur in $t$ before position $\alpha$. Given $\beta \in \{1, \dots , \ell(t) \}$ such that $\beta \leqslant \alpha$ we define $L(t, \beta, \alpha)$ as the number of left round brackets `(' that occur in $t$ starting from position $\beta$ and before position $\alpha$. We define $R(t, \beta, \alpha)$ as the number of right round brackets `)' that occur in $t$ starting from position $\beta$ and before position $\alpha$.\\

With such definitions we can observe what follows:
\small
\begin{align*}
d(t,\ell(\vartheta)+\alpha) &=  L(t,\ell(\vartheta)+\alpha) - R(t,\ell(\vartheta)+\alpha)\\
&= L(t,\ell(\vartheta)+1) + L(t,\ell(\vartheta)+1, \ell(\vartheta)+\alpha) - R(t,\ell(\vartheta)+1) - R(t,\ell(\vartheta)+1, \ell(\vartheta)+\alpha)\\
&= L(t,\ell(\vartheta)+1) + L(\varphi, \alpha) - R(t,\ell(\vartheta)+1) - R(\varphi, \alpha)\\
&= d(t,\ell(\vartheta)+1) + d(\varphi, \alpha)
\end{align*}

\normalsize

\end{proof}

\medskip

Before we describe the process of constructing expressions for our language we must also prove some useful lemmas related to the predicates we have defined above.\\

\begin{lemma}\label{L:useful-on-predicates-1}
Given $i = 1 \dots p$ and a positive integer $q$, for each $x \in \mathcal{P}^q(D_i)$ we have

\begin{itemize}
\item $Set_q(x)$,
\item for each $r > q$ $\neg Set_r(x)$
\end{itemize}

\end{lemma}

\begin{proof}

We proceed by induction on $q$.\\

Let $q = 1$. We assume $x \in \mathcal{P}(D_i)$, then clearly $Set_1(x)$.\\

Given a positive integer $r$, we assume $Set_{r+1}(x)$ and try to derive a contradiction. Since $x \ne \emptyset$ we can take $z \in x$, we have $Set_r(z)$ and $z \in D_i$, this actually is a contradiction. So we have proved $\neg(Set_{r+1}(x))$.\\

Let now $q$ be a positive integer and assume for each $x \in \mathcal{P}^q(D_i)$ we have

\begin{itemize}
\item $Set_q(x)$,
\item for each $r > q$ $\neg Set_r(x)$.
\end{itemize}

\smallskip

We want to show that for each $x \in \mathcal{P}^{q+1}(D_i)$ we have

\begin{itemize}
\item $Set_{q+1}(x)$,
\item for each $r > q+1$ $\neg Set_r(x)$.
\end{itemize}

\smallskip

We have that $x \ne \emptyset$, for each $z \in x$ $z \in \mathcal{P}^q(D_i)$, so for each $z \in x$ $Set_q(z)$. Therefore $Set_{q+1}(x)$.\\

Given a positive integer $r > q+1$, we assume $Set_r(x)$ and try to derive a contradiction. Let $z \in x$, we have $z \in \mathcal{P}^q(D_i)$ and $Set_{r-1}(z)$. Since $r-1 > q$ $\neg Set_{r-1}(z)$ should hold, and we have derived a contradiction.
\end{proof}

\bigskip

\begin{lemma}\label{L:useful-on-predicates-2}
Given $i, j = 1 \dots p$, $q, r$ positive integers such that $q \ne r$ $\mathcal{P}^q(D_i) \cap \mathcal{P}^r(D_j) = \emptyset$.

\end{lemma}

\begin{proof}
Let's suppose, absurdly, $x \in \mathcal{P}^q(D_i) \cap \mathcal{P}^r(D_j)$. Suppose $q < r$.\\

Using lemma~\ref{L:useful-on-predicates-1} we have both $Set_r(x)$ and $\neg Set_r(x)$. Therefore we must have $\mathcal{P}^q(D_i) \cap \mathcal{P}^r(D_j) = \emptyset$.\\

In the case $q > r$ we can apply the same type of reasoning.
\end{proof}

\bigskip

\begin{lemma}\label{L:useful-on-predicates-3}
Given $i, j = 1 \dots p$ such that $i \ne j$ and a positive integer $q$ we have $\mathcal{P}^q(D_i) \cap \mathcal{P}^q(D_j) = \emptyset$.
\end{lemma}

\begin{proof}

We proceed by induction on $q$.\\

Let $q = 1$. Assume $\mathcal{P}(D_i) \cap \mathcal{P}(D_j) \ne \emptyset$ and let $x \in \mathcal{P}(D_i) \cap \mathcal{P}(D_j)$.\\

We have $x \ne \emptyset$, $x \subseteq D_i$, $x \subseteq D_j$, so $x \subseteq D_i \cap D_j$, and $D_i \cap D_j \ne \emptyset$, against our assumptions.\\

In order to perform the inductive step, let $q$ be a positive integer, we assume $\mathcal{P}^q(D_i) \cap \mathcal{P}^q(D_j) = \emptyset$ and we try to show $\mathcal{P}^{q+1}(D_i) \cap \mathcal{P}^{q+1}(D_j) = \emptyset$.\\

We assume, absurdly, $x \in \mathcal{P}^{q+1}(D_i) \cap \mathcal{P}^{q+1}(D_j)$. We have $x \ne \emptyset$, $x \subseteq \mathcal{P}^q(D_i)$, $x \subseteq \mathcal{P}^q(D_j)$. So $x \subseteq \mathcal{P}^q(D_i) \cap \mathcal{P}^q(D_j)$ and $\mathcal{P}^q(D_i) \cap \mathcal{P}^q(D_j) \ne \emptyset$ against our assumptions.
\end{proof}

\bigskip

\begin{lemma}\label{L:useful-on-predicates-4}
Given $i = 1 \dots p$ and a positive integer $q$ for each $x \in \mathcal{P}^q(D_i)$ and $r \leqslant q$ we have $Set_r(x)$ and $NotEmpty_r(x)$.
\end{lemma}

\begin{proof}

We proceed by induction on $q$.\\

Let $q = 1$ and let $x \in \mathcal{P}(D_i)$. Clearly $Set_1(x)$ and $NotEmpty_1(x)$.\\

In order to perform the inductive step, let $q$ be a positive integer, we assume for each $x \in \mathcal{P}^q(D_i)$ and $r \leqslant q$ we have $Set_r(x)$ and $NotEmpty_r(x)$.\\

Let now $x \in \mathcal{P}^{q+1}(D_i)$, we want to show that for each $r \leqslant q+1$ we have $Set_r(x)$ and $NotEmpty_r(x)$.\\

Clearly $Set_1(x)$ and $NotEmpty_1(x)$ both hold true.\\

Given $r \leqslant q+1$ such that $r > 1$ we want to prove $Set_r(x)$ and $NotEmpty_r(x)$.\\

In order to prove this we just need to prove that for each $u \in x$ $Set_{r-1}(u)$ and $NotEmpty_{r-1}(u)$.\\

Since $x \in \mathcal{P}^{q+1}(D_i)$ then $x \subseteq \mathcal{P}^q(D_i)$ and for each $u \in x$ $u \in \mathcal{P}^q(D_i)$. Since $r-1 \leqslant q$ we have indeed $Set_{r-1}(u)$ and $NotEmpty_{r-1}(u)$.
\end{proof}

\bigskip

\begin{lemma}\label{L:useful-on-predicates-5}
Given $i = 1 \dots p$ and a positive integer $q$ for each $x \in \mathcal{P}^q(D_i)$ and $r \leqslant q+1$ we have $\neg Event_r(x)$.
\end{lemma}

\begin{proof}
We proceed by induction on $q$.\\

Let $q = 1$. Let $x \in \mathcal{P}(D_i)$, $x$ is a set and I think we can assume $\neg Event_1(x)$.\\
Moreover $x \subseteq D_i$ so for each $u \in x$ $u \in D_i$, for each $u \in x$ $\neg Event_1(x)$, so it is false that for each $u \in x$ $Event_1(x)$, and it follows that $\neg(Event_2(x))$.\\

For the inductive step, let $q$ be a positive integer and we assume for each $x \in \mathcal{P}^q(D_i)$ and $r \leqslant q+1$ we have $\neg Event_r(x)$. Let $x \in  \mathcal{P}^{q+1}(D_i)$, let $r \leqslant q+2$, we want to show that $\neg Event_r(x)$ holds.\\

If $r = 1$ since $x$ is a set we can assume $\neg Event_1(x)$ holds.\\ 

If $r > 1$ we have $x \subseteq \mathcal{P}^q(D_i)$, so for each $u \in x$ $u \in \mathcal{P}^q(D_i)$, and then for each $u \in x$ $\neg Event_{r-1}(u)$. So it is false that for each $u \in x$ $Event_{r-1}(u)$, and $Event_r(x)$ is false.
\end{proof}

\bigskip

\begin{lemma}\label{L:useful-on-predicates-6}
For each positive integer $q$ and for each $x$ if $Event_{q+1}(x)$ then $Set_q(x)$.
\end{lemma}

\begin{proof}
For the initial step of the proof, if $Event_2(x)$ then $Set_1(x)$.\\

For the inductive step, let $q \geqslant 1$ and $Event_{q+2}(x)$, then $x$ is a set and for each $u \in x$ $Event_{q+1}(u)$. It then follows that for each $u \in x$ $Set_q(u)$, and then $Set_{q+1}(x)$.
\end{proof}

\bigskip

\begin{lemma}\label{L:useful-on-predicates-7}
Given $i = 1 \dots p$ and a positive integer $q$ for each $x \in \mathcal{P}^q(D_i)$ and $r > q+1$ we have $\neg Event_r(x)$.
\end{lemma}

\begin{proof}
Let $x \in \mathcal{P}^q(D_i)$ and let $r > q+1$. Since $r-1 > q$ by lemma~\ref{L:useful-on-predicates-1} $\neg Set_{r-1}(x)$ and by lemma~\ref{L:useful-on-predicates-6} $\neg Event_r(x)$.
\end{proof}

\bigskip

\begin{lemma}\label{L:useful-on-predicates-8}
For each positive integer $q$ and for each $x$ if $Set_{q+1}(x)$ then $Set_q(x)$.
\end{lemma}

\begin{proof}
For the initial step of the proof, if $Set_2(x)$ then $Set_1(x)$.\\

For the inductive step, let $q \geqslant 1$ and $Set_{q+2}(x)$, then $x$ is a set and for each $u \in x$ $Set_{q+1}(u)$. It then follows that for each $u \in x$ $Set_q(u)$, and then $Set_{q+1}(x)$.
\end{proof}

\section{Computability theory}\label{Ch:computability}

The proposed logical system aims to satisfy every computability requirement that is desirable. To verify these requirements, we briefly recall some foundational concepts of computability theory.\\

We use Cutland's book~\cite{Cutland} as the main reference for this. The book defines the concept of computable function: given a set $A$ of natural numbers and a function $f: A \to \mathbb{N}$, we say that f is \emph{computable} when it is URM-computable. We will not define here the concept of URM-computability, the reader can find the definition in the mentioned book.\\

As suggested by the book we use the symbol $\mathscr{C}$ to indicate the set of the computable functions from a subset of $\mathbb{N}$ to $\mathbb{N}$ (also called the `partial functions' from $\mathbb{N}$ to $\mathbb{N}$).\\

The book also provides many alternative definitions of the notion of effective computability and affirms that `the remarkable result of investigation by many researchers is the following: Each of the above proposals for a characterisation of the notion of effective computability gives rise to the same class of functions, the class that we have denoted with $\mathscr{C}$'.\\

Finally the book also states the famous `Church's thesis' in the following terms: `The intuitively and informally defined class of effectively computable partial functions coincides exactly with the class $\mathscr{C}$ of URM-computable functions'.\\

If $A$ is a subset of $\mathbb{N}$ we can define the \emph{characteristic function} of $A$ as the function $c_A$ given by: if $x \in A$ $c_A(x) = 1$; if $x \notin A$ $c_A(x) = 0$. Then $A$ is said to be \emph{recursive} if $c_A$ is computable.\\

If $A$ is a subset of $\mathbb{N}$ we can define the \emph{semi-characteristic function} of $A$ as the function $s_A$ given by: if $x \in A$ $s_A(x) = 1$; if $x \notin A$ $s_A(x)$ is undefined. Then $A$ is said to be \emph{recursively enumerable (r.e.)} if $s_A$ is computable.\\

A recursive set is obviously also recursively enumerable.\\

Given a subset $A$ of $\mathbb{N}$ the following statements are equivalent:

\begin{itemize}
\item $A$ is r.e.;
\item $A = \emptyset$ or $A$ is the range of a total computable function;
\item $A$ is the range of a partial computable function.
\end{itemize}

\medskip

Please refer to Cutland's book for the proof of the equivalence.\\

We now state and prove a theorem required for our framework, which does not explicitly appear in~\cite{Cutland}.\\

\begin{theorem}\label{Th:union-of-re}
Let $A$ be a r.e. subset of $\mathbb{N}$, let $f$ be a function defined on $A$ such that for each $x \in A$ $f(x)$ is a r.e. subset of $\mathbb{N}$. Then $\bigcup_{x \in A} f(x)$ is r.e..
\end{theorem}

\begin{proof}

There exists a partial computable function $\xi$ such that $A = ran(\xi)$. For each $x \in A$ there also exists a partial function $\chi_x$ such that $f(x) = ran(\chi_x)$.\\

Let's consider the function $\pi : \mathbb{N}^2 \to \mathbb{N}$ (named the Cantor's pairing function) defined by 
\[ \pi(x, y) = (x+y)(x+y+1)/2 + y \ . \] 

This function is a bijection and the inverse function $\zeta : \mathbb{N} \to \mathbb{N}^2$ is a computable function itself (cfr. Wikipedia `https://en.wikipedia.org/wiki/Pairing{\_}function').\\

Let's now consider a function $\phi$ defined over $\mathbb{N}$ such that $\phi(z)$ is calculated as follows: we first calculate $\zeta(z) = (z_1, z_2)$, then we calculate $\xi(z_1)$, if it terminates $\xi(z_1) \in A$ and we can set $\phi(z) = \chi_{\xi(z_1)}(z_2)$.\\

The function $\phi$ is a partial computable function and we can show that $\bigcup_{x \in A} f(x) = ran(\phi)$.\\

Given $y \in \bigcup_{x \in A} f(x)$ we will prove that $y \in ran(\phi)$. In fact there exists $x \in A$ such that $y \in f(x)$. There exists $z_1 \in \mathbb{N}$ such that $x = \xi(z_1)$, and there exists $z_2 \in \mathbb{N}$ such that $y = \chi_x(z_2) = \chi_{\xi(z_1)}(z_2)$. There exists $z \in \mathbb{N}$ such that $\zeta(z) = (z_1, z_2)$ and therefore $y = \phi(z)$.\\

Vice versa given $y \in ran(\phi)$ we want to show that $y \in \bigcup_{x \in A} f(x)$. There exists $z \in \mathbb{N}$ such that $y = \phi(z)$. If we set $(z_1, z_2) = \zeta(z)$ then $y = \phi(z) = \chi_{\xi(z_1)}(z_2)$. We have that $\xi(z_1) \in A$ and $y = \chi_{\xi(z_1)}(z_2) \in f(\xi(z_1))$.

\end{proof}

\bigskip

Our reference book also explains how to apply the definition of computability and the related ones to a domain $D$ which is different from $\mathbb{N}$. This requires the availability of a coding.\\

A \emph{coding} of a domain $D$ of objects is and explicit and effective injection $\alpha: D \to \mathbb{N}$.\\

We can actually assume that the range of $\alpha$ is $\mathbb{N}$ (and in this case $\alpha$ is a bijection) or at least that $ran(\alpha)$ is recursive.\\

A partial function $f: D \to D$ is coded by the function $f^* = \alpha \circ f \circ \alpha^{-1}$, so $f^*$ is a partial function $\mathbb{N} \to \mathbb{N}$. We say that $f$ is \emph{computable} if and only if $f^*$ is computable.\\

Given a set $A \subseteq D$ we can define $A^* = \{ \alpha(d) | d \in A \}$. We say that $A$ is \emph{recursive} iff $A^*$ is recursive, and that  $A$ is \emph{recursively enumerable} iff $A^*$ is recursively enumerable.\\

Given $A \subseteq D$ we can define the \emph{characteristic function} of $A$ as the function $c_A$ whose domain is $D$ given by: if $x \in A$ $c_A(x) = 1$; if $x \notin A$ $c_A(x) = 0$. We can also define the \emph{semi-characteristic function} of $A$ as the function $s_A$ whose domain is A, such that for each $x \in A$ $s_A(x) = 1$.\\

In relation to the former definitions, we can prove the following lemma.\\ 

\begin{lemma}\label{L:recursive-re-sets}
Let $A \subseteq D$, then

\begin{itemize}
\item $A$ is recursive if and only if $c_A$ is computable;
\item $A$ is r.e. if and only if $s_A$ is computable.
\end{itemize}

\end{lemma}

\begin{proof}

First of all we notice that given $n \in \mathbb{N}$

\begin{itemize} 
\item if $n \in A^*$ then $\alpha^{-1}(n) \in A$;
\item if $\alpha^{-1}(n) \in A$ then $n \in A^*$.
\end{itemize}

\smallskip

We also notice that given $x \in D$, if $x \notin A$ then $\alpha(x) \notin A^*$. In fact if $\alpha(x) \in A^*$ then there exists $y \in A$ such that $\alpha(x) = \alpha(y)$, but since $x \ne y$ and $\alpha$ is injective we cannot have $\alpha(x) = \alpha(y)$.\\ 

Let's assume $A$ is recursive, we want to show that $c_A$ is computable.\\

We know that $c_{A^*}$ is computable. Given $x \in D$

\begin{itemize}
\item if $x \in A$ then $\alpha(x) \in A^*$ $c_A(x) = 1 = c_{A^*}(\alpha(x))$;
\item if $x \notin A$ then $\alpha(x) \notin A^*$ $c_A(x) = 0 = c_{A^*}(\alpha(x))$.
\end{itemize}

\smallskip

Therefore in every case $c_A(x) = c_{A^*}(\alpha(x))$, and then $c_A$ is computable.\\

Vice versa we now assume $c_A$ is computable and we want to show that $A$ is recursive.\\

Given $n \in \mathbb{N}$,

\begin{itemize}
\item if $n \notin ran(\alpha)$ we have $c_{ran(\alpha)}(n) = 0$, so $c_{A^*}(n) = 0 = c_{ran(\alpha}(n)$.
\item if $n \in ran(\alpha)$ we have $c_{ran(\alpha)}(n) = 1$ and
	\begin{itemize}
	\item if $n \in A^*$ then $\alpha^{-1}(n) \in A$, $c_{A^*}(n) = 1 = c_A(\alpha^{-1}(n))$;
	\item if $n \notin A^*$ then $\alpha^{-1}(n) \notin A$, $c_{A^*}(n) = 0 = c_A(\alpha^{-1}(n))$.
	\end{itemize}
\end{itemize} 

\smallskip

Clearly we can compute $c_{A^*}(n)$ as follows:\\

If $c_{ran(\alpha)}(n) = 0$ then $c_{A^*}(n) = 0$;\\

if $c_{ran(\alpha)}(n) = 1$ then $c_{A^*}(n) = c_A(\alpha^{-1}(n))$.\\

Let's assume $A$ is r.e., we want to show that $s_A$ is computable.\\

Given $x \in D$,

\begin{itemize}
\item if $x \in A$ then $\alpha(x) \in A^*$, $x \in dom(s_A)$, $\alpha(x) \in dom(s_{A^*})$ $s_A(x) = 1 = s_{A^*}(\alpha(x))$;
\item if $x \notin A$ then $\alpha(x) \notin A^*$ $x \notin dom(s_A)$, $\alpha(x) \notin dom(s_{A^*})$, $s_A(x)$ and $s_{A^*}(\alpha(x))$ are both divergent.
\end{itemize}

\smallskip

Therefore $s_A(x)$ can be calculated by $s_{A^*}(\alpha(x))$, and $s_A$ is computable.\\

Vice versa we now assume $s_A$ is computable and we want to show that $A$ is r.e..\\

Given $n \in \mathbb{N}$,

\begin{itemize}
\item if $n \notin ran(\alpha)$ we have $n \notin A^* = dom(s_{A^*})$, therefore $s_{A^*}(n)$ is divergent, and $s_A(\alpha^{-1}(n))$ is divergent too.
\item if $n \in ran(\alpha)$ we have $c_{ran(\alpha)}(n) = 1$ and
	\begin{itemize}
	\item if $n \in A^*$ then $\alpha^{-1}(n) \in A$, $s_{A^*}(n) = 1 = s_A(\alpha^{-1}(n))$;
	\item if $n \notin A^*$ then $\alpha^{-1}(n) \notin A$, $s_{A^*}(n)$ and $s_A(\alpha^{-1}(n))$ are both divergent.
	\end{itemize}
\end{itemize} 

\smallskip

Therefore in all cases $s_{A^*}(n)$ can be calculated as $s_A(\alpha^{-1}(n))$, and so $s_{A^*}$ is computable and $A$ is r.e..

\end{proof}

\bigskip

In the theorem~\ref{Th:union-of-re} above we proved that a r.e. union of r.e. sets is still a r.e. set. This theorem was stated for subsets of $\mathbb{N}$, and we will generalize it to generic domains.\\

\begin{theorem}\label{Th:gen-union-of-re}
Let $D_1$ and $D_2$ be two `domains' to which we can apply the notions of computability using two codings $\alpha_1: D_1 \to \mathbb{N}$ and $\alpha_2: D_2 \to \mathbb{N}$.
Let $A$ be a r.e. subset of $D_1$, let $f$ be a function defined on $A$ such that for each $x \in A$ $f(x)$ is a r.e. subset of $D_2$. Then $\bigcup_{x \in A} f(x)$ is a r.e. subset of $D_2$.
\end{theorem}

\begin{proof}

We call $W$ the set $\bigcup_{x \in A} f(x) \subseteq D_2$. We'll prove that $W$ is r.e. by proving that $W^* = \{ \alpha_2(y) | y \in W \}$ is r.e.. Let's assume that actually $W^* = \bigcup_{z \in A^*} f(\alpha_1^{-1}(z))^*$.\\

If the equality we have assumed holds, then we can consider that $A^*$ is a r.e. subset of $\mathbb{N}$, and that for each $z \in A^*$ $\alpha_1^{-1}(z) \in A$, $f(\alpha_1^{-1}(z))$ is a r.e. subset of $D_2$, $f(\alpha_1^{-1}(z))^*$ is a r.e. subset of $\mathbb{N}$. Therefore, if the equality holds, we have proved that $W^*$ is r.e. and our proof is finished.\\

Let's then show that $W^* = \bigcup_{z \in A^*} f(\alpha_1^{-1}(z))^*$ actually holds.\\

Let $w \in W^*$ then there exists $y \in W$: $w = \alpha_2(y)$, and there exists $x \in A$: $y \in f(x)$. Let $z = \alpha_1(x) \in A^*$, then $y \in f(\alpha_1^{-1}(z))$ and $w = \alpha_2(y) \in f(\alpha_1^{-1}(z))^*$. So we can confirm that $w \in  \bigcup_{z \in A^*} f(\alpha_1^{-1}(z))^*$.\\

Conversely let $w \in  \bigcup_{z \in A^*} f(\alpha_1^{-1}(z))^*$ and we want to prove that $w \in W^*$. There exists $z \in A^*$ such that $w \in f(\alpha_1^{-1}(z))^*$. Let $x = \alpha_1^{-1}(z) \in A$ then $w \in f(x)^*$. Let $y = \alpha_2^{-1}(w) \in f(x)$. We have $y \in W$ and so $w = \alpha_2(y) \in W^*$.

\end{proof}

\bigskip

Typically we will be dealing with a finite or countable alphabet $\Sigma$, and the domain to which we will have to apply the concepts of computability will be the set $\Sigma^*$ of all the empty or finite strings with characters in the mentioned alphabet. But we may also need to apply those concepts e.g. to $(\Sigma^*)^n$. So let us examine some sets to which we can actually apply computability notions, in order to be able to apply such concepts wherever we need them.\\ 

First of all we consider the set $\mathbb{N}^2$. There is a coding $\pi: \mathbb{N}^2 \to \mathbb{N}$ and this coding is the Cantor pairing function defined by
\[ \pi(k_1, k_2) = (k_1 + k_2)(k_1 + k_2 + 1)/2 + k_2 \ . \]

So, obviously, we can apply computability notions to $\mathbb{N}^2$, and actually we are able to apply them also to $\mathbb{N}^n$ for an arbitrary integer $n > 2$. In fact if we assume $\pi_n$ is a coding $\mathbb{N}^n \to \mathbb{N}$, with $n \geqslant 2$, then we can define a function $\pi_{n+1}: \mathbb{N}^{n+1} \to \mathbb{N}$ as follows:
\[ \pi_{n+1}(x_1, \dots , x_n, x_{n+1}) = \pi( \pi_n(x_1, \dots , x_n), x_{n+1} ) \ . \]

And this function is actually a coding $\mathbb{N}^{n+1} \to \mathbb{N}$.\\

At this point given $n$ domains $D_1, \dots , D_n$ such that for each $i = 1 \dots n$ there exists a coding $\alpha_i: D_i \to \mathbb{N}$ we can build a coding $\alpha: D_1 \times \dots \times D_n \to \mathbb{N}$. Our coding will be defined as follows:
\[ \alpha(d_1, \dots , d_n) = \pi_n( \alpha_1(d_1), \dots , \alpha_n(d_n) ) \ . \]

We said earlier that typically we will be dealing with a finite or countable alphabet $\Sigma$, and the domain to which we will have to apply the concepts of computability will be the set $\Sigma^*$. With respect to this, we notice that $\mathbb{N}$ can be itself considered as an alphabet, so we first try to find a coding $\mathbb{N}^* \to \mathbb{N}$.\\

Here we notice that $\mathbb{N}^* = \{ \epsilon \} \cup \bigcup_{i \geqslant 1} \mathbb{N}^i$.\\

We have seen that for each $i \geqslant 2$ $\pi_i$ is a coding $\mathbb{N}^i \to \mathbb{N}$, so we can create a coding $\gamma: \bigcup_{i \geqslant 1} \mathbb{N}^i \to \mathbb{N}^2$ as follows:

\begin{itemize}
\item for each $x \in \mathbb{N}$ $\gamma(x) = (0,x)$;
\item for each $i > 1$, $(x_1, \dots , x_i) \in \mathbb{N}^i$ $\gamma(x_1, \dots , x_i) = (i-1, \pi_i(x_1, \dots , x_i))$.
\end{itemize}

\smallskip

We now want to create a coding $\alpha: \mathbb{N}^* \to \mathbb{N}$. We define our coding $\alpha$ as follows:

\begin{itemize}
\item $\alpha(\epsilon) = 0$;
\item for each $x \in \bigcup_{i \geqslant 1} \mathbb{N}^i$ $\alpha(x) = \pi_2(\gamma(x)) + 1$.
\end{itemize}

\smallskip

Given a finite or countable alphabet $\Sigma$ we now want to define a coding $\Sigma^* \to \mathbb{N}$. Of course there exists a coding $\sigma: \Sigma \to \mathbb{N}$. We first want to define a coding $\delta: \Sigma^* \to \mathbb{N}^*$, and, since $\Sigma^* = \{ \epsilon \} \cup \bigcup_{i \geqslant 1} \Sigma^i$, we can define it as follows. 

\begin{itemize}
\item $\delta(\epsilon) = \epsilon$;
\item for each $i \geqslant 1$ $(x_1, \dots , x_i) \in \Sigma^i$ $\delta(x_1, \dots , x_i) = (\sigma(x_1), \dots , \sigma(x_i)) \in \mathbb{N}^i$.
\end{itemize}

\smallskip

At this point if $\alpha$ is a coding $\mathbb{N}^* \to \mathbb{N}$ then $\gamma = \alpha \circ \delta$ is a coding  $\Sigma^* \to \mathbb{N}$.\\

We can notice that if $\Sigma$ is finite then $\sigma$ is not surjective, and so also $\delta$ and $\gamma$ are not surjective. Nevertheless, $ran(\gamma)$ remains recursive since given $x \in \mathbb{N}$ we can decide whether $x \in ran(\gamma)$. In order to do this we can calculate $\alpha^{-1}(x) \in N^*$, and here we can determine if $\alpha^{-1}(x) \in ran(\delta)$, if this is true since $x = \alpha(\alpha^{-1}(x))$ then $x \in ran(\gamma)$. If on the contrary $\alpha^{-1}(x) \notin ran(\delta)$ then $x \notin ran(\gamma)$. In fact if $x \in ran(\gamma)$ then exists $y \in \Sigma^*$ such that $x = \gamma(y) = \alpha(\delta(y))$, so $\alpha(\alpha^{-1}(x)) = \alpha(\delta(y))$, and $\alpha^{-1}(x) = \delta(y) \in ran(\delta)$.\\

Once we have a coding for $\Sigma^*$ we have it also for $(\Sigma^*)^n$, where $n$ is a positive integer, and if $\Gamma$ is another alphabet we have a coding for $(\Sigma^*)^n \times (\Gamma^*)^m$, where $m$ is another positive integer.\\

Given $n$ domains $D_1, \dots , D_n$ and given $A_1 \subseteq D_1, \dots , A_n \subseteq D_n$ if $A_1, \dots , A_n$ are r.e. then $A_1 \times \dots \times A_n$ is also r.e.. In fact given $(x_1, \dots , x_n) \in D_1 \times \dots \times D_n$ we can compute $s_{A_1 \times \dots \times A_n}(x_1, \dots , x_n)$ as follows: for each $i = 1 \dots n$ we calculate $s_{A_i}(x_i)$ and if we obtain a result for each $i$ then we emit the result $1$.\\

Given a domain $D$ and a coding $\alpha: D \to \mathbb{N}$ and given a function $f: \mathbb{N} \to D$ we can say that $f$ is computable when $\alpha \circ f: \mathbb{N} \to \mathbb{N}$ is computable.\\

We can prove the following lemma:

\begin{lemma}
Given a set $A \subseteq D$ $A$ is r.e. if and only if $A = \emptyset$ or there exists a total computable function $f: N \to D$ such that $A = ran(f)$.
\end{lemma}

\begin{proof}

Let $A^* = \{ \alpha(d) | d \in A \}$.\\

If $A$ is r.e. then $A^*$ is r.e. and so $A^* = \emptyset$ or there exists a total computable function $g: \mathbb{N} \to \mathbb{N}$ such that $A^* = ran(g)$. If $A^* = \emptyset$ then clearly $A = \emptyset$, otherwise since $A^* \subseteq ran(\alpha)$ we can define $f = \alpha^{-1} \circ g$, $f$ is a function $\mathbb{N} \to D$ and $ran(f) = A$.\\

In fact if $d \in ran(f)$ then there exists $x \in \mathbb{N}$: $d = \alpha^{-1}(g(x))$. We know that $g(x) \in A^*$ and so $d = \alpha^{-1}(g(x)) \in A$. Conversely if $d \in A$ then $\alpha(d) \in A^*$ and there exists $x \in \mathbb{N}$: $g(x) = \alpha(d)$, $d = \alpha^{-1}(g(x)) = f(x) \in ran(f)$.\\

Conversely if $A= \emptyset$ then $A^* = \emptyset$, $A^*$ is r.e. and $A$ is r.e.. If there exists a total computable function $f: \mathbb{N} \to D$ such that $A = ran(f)$ then $g = \alpha \circ f: \mathbb{N} \to \mathbb{N}$ is computable and $A^* = ran(g)$.\\

In fact if $y \in A^*$ then there exists $d \in A$: $y = \alpha(d)$ and there exists $x \in \mathbb{N}$: $d = f(x)$, so $y = \alpha(f(x)) \in ran(g)$. Conversely if $y \in ran(g)$ then there exists $x \in \mathbb{N}$: $y = \alpha(f(x))$, and since $f(x) \in A$ then $y \in A^*$.\\

So in the latest case too $A^*$ is r.e. and $A$ is r.e.. 
\end{proof}

Let $D_1$ and $D_2$ be two `domains' to which we can apply the notions of computability using two codings $\alpha_1: D_1 \to \mathbb{N}$ and $\alpha_2: D_2 \to \mathbb{N}$. Using codings, we can define the notion of `computable function' also for a function $f: D_1 \to D_2$. We say that $f$ is computable if and only if $ \alpha_2 \circ f \circ \alpha_1^{-1} : \mathbb{N} \to \mathbb{N}$ is computable. We can notice that since the domain of $\alpha_1^{-1}$ could be a proper subset of $\mathbb{N}$ the mentioned function could actually be a partial function $\mathbb{N} \to \mathbb{N}$.\\

Using the just introduced notion, we can prove the following lemma:

\begin{lemma}\label{L:fre-is-re}
Let $A$ be a r.e. subset of $D_1$ and let $f: D_1 \to D_2$ be a computable function. If we define $B = \{ f(a) | \, a \in A \}$, then $B \subseteq D_2$ is r.e..
\end{lemma}

\begin{proof}
If $A = \emptyset$ then $B = \emptyset$ is r.e..\\

Otherwise there exists a computable function $g: \mathbb{N} \to D_1$ such that $A = ran(g)$. We have $f \circ g: \mathbb{N} \to D_2$ and $B = ran(f \circ g)$.\\

In fact if $d_2 \in B$ then there exists $a \in A$ such that $d_2 = f(a)$ and there exists $x \in \mathbb{N}$ such that $a = g(x)$, so $d_2 = f(g(x))$. Conversely if $d_2 \in ran(f \circ g)$ then there exists $x \in \mathbb{N}$ such that $d_2 = f(g(x))$, so $g(x) \in A$ and $d_2 = f(g(x)) \in B$. 
\end{proof}

\section{Premise: description of contexts}

We  want to be able to show that a certain set of contexts is recursive or recursively enumerable. To this end we are going to define contexts as strings. Given a finite or countable alphabet $\Sigma$ which contains a finite or countable set of variables $\mathcal{V}$, the symbols `:' and `,' and doesn't contain the symbols `$<$' and `$>$' we can define an alphabet $\Gamma = \Sigma \cup \{<,>\}$.\\

Let $\epsilon$ be the string `$<>$' and let $\Theta( \Sigma, \mathcal{V})$ (henceforth $\Theta$) be the set 
\[ \{ \epsilon \} \cup \{ <<x_1: \varphi_1> \dots <x_m: \varphi_m>> | \, x_1, \dots , x_m \in \mathcal{V}, \, \varphi_1 \dots \varphi_m \in \Sigma^*  \} \, . \]

\smallskip

\begin{lemma}\label{unique-read-contexts}
Let $k \in \Theta - \{ \epsilon \}$, let $m$ positive integer, $x_1, \dots , x_m \in \mathcal{V}, \, \varphi_1, \dots , \varphi_m \in \Sigma^*$ such that $k = <<x_1: \varphi_1> \dots <x_m: \varphi_m>>$.  Let also $p$ positive integer, $y_1, \dots , y_p \in \mathcal{V}, \, \psi_1, \dots , \psi_p \in \Sigma^*$ such that $k = <<y_1: \psi_1> \dots <y_p: \psi_p>>$. Then $p = m$, for each $i = 1 \dots m$ $y_i = x_i$ and $\psi_i = \varphi_i$.
\end{lemma}

\begin{proof}

\smallskip

First of all, as in former parts of the paper, if $t$ is a string we will indicate with $\ell(t)$ $t$'s length, i.e. the number of characters in $t$.\\

Also, given $\alpha = 1 \dots \ell(t)$ let $t[\alpha]$ indicate the character with position $\alpha$ inside $t$, and given $\alpha, \beta = 1 \dots \ell(t)$ with $\alpha \leqslant \beta$ let $t[\alpha, \beta]$ be the substring of $t$ which begins at character $\alpha$ and ends at character $\beta$.\\

For each $i = 1 \dots m$ let $\varphi_i = k[\alpha_i, \mu_i]$. For each $j = 1 \dots p$ let $\psi_j = k[\beta_j, \nu_j]$.\\

Clearly $y_1 = k[3] = x_1$.\\

If $\nu_1 > \mu_1$ then $\nu_1 \geqslant \mu_1 + 1 > \alpha_1 = 5 = \beta_1$ so `$>$' $ = k[\mu_1 + 1]$ is a character in $\psi_1$: this cannot be true, so $\nu_1 > \mu_1$ is false and similarly $\mu_1 > \nu_1$ is false, so $\nu_1 = \mu_1$ and it follows that $\psi_1 = \varphi_1$.\\

If $m = 1$ then $k[\mu_1 + 2] = $ `$>$'. In this case if $p > 1$ then $k[\nu_1 + 2] = $ `$<$'. Therefore it must be $p = 1$ and our proof for $m=1$ is finished.\\

Let's consider the case $m > 1$. In this case given $i = 1 \dots m - 1$ we assume we have proved that for each $j = 1 \dots i$ $p \geqslant j$ $y_j = x_j$, $\psi_j = \varphi_j$ and that $\nu_i = \mu_i$. We want to show that $p \geqslant i+1$, $y_{i+1} = x_{i+1}$, $\psi_{i+1} = \varphi_{i+1}$, $\nu_{i+1} = \mu_{i+1}$.\\

Since $i < m$ $k[\mu_i + 2] = $ `$<$', therefore also $k[\nu_i + 2] = $ `$<$' and this implies that $p \geqslant i+1$.\\

We also notice that $y_{i+1} = k[\nu_i + 3] = k[\mu_i + 3] = x_{i+1}$.\\

We also notice that $\beta_{i+1} = \nu_i + 5 = \mu_i + 5 = \alpha_{i+1}$.\\

If $\nu_{i+1} > \mu_{i+1}$ then $\nu_{i+1} \geqslant \mu_{i+1} + 1 > \alpha_{i+1} = \beta_{i+1}$ so `$>$' $ = k[\mu_{i+1} + 1]$ is a character in $\psi_{i+1}$: this cannot be true, so $\nu_{i+1} > \mu_{i+1}$ is false and similarly $\mu_{i+1} > \nu_{i+1}$ is false, so $\nu_{i+1} = \mu_{i+1}$ and it follows that $\psi_{i+1} = \varphi_{i+1}$.\\

We have now proved that for each $i = 1 \dots m$ $p \geqslant i$ $y_i = x_i$, $\psi_i = \varphi_i$ and $\nu_i = \mu_i$.\\

We have also $\ell(k) = \mu_m + 2$. If $p > m$ then $\ell(k) > \nu_m + 2 = \mu_m + 2$, and this is a contradiction, therefore $p = m$.

\end{proof}

\smallskip

Given $k \in \Theta$ we define $dom(k)$ (i.e. the domain of $k$) as follows. 

\begin{itemize}
\item if $k = \epsilon$ then $dom(k) = \emptyset$,
\item if $k \ne \epsilon$ then $k = <<x_1: \varphi_1> \dots <x_m: \varphi_m>>$ and we define $dom(k) = \{ 1, \dots , m \} \, . $
\end{itemize}

\medskip

We define $\mathcal{D} = \{\emptyset\} \cup \{ \{1, \dots , m \} \vert \ m \text{ is a positive integer} \} $. \\

Of course, given $k \in \Theta$, $dom(k) \in \mathcal{D}$.\\

Given $k \in \Theta$ and $C \in \mathcal{D}$ such that $C \subseteq dom(k)$ we can define $k_{/C}$, i.e. the `restriction' of $k$ to the domain $C$, as follows:

\begin{itemize}
\item if $k = \epsilon$ or $C = \emptyset$ then $k_{/C} = \epsilon \in \Theta$ (so $dom(k_{/C}) = \emptyset = C$),
\item if $k \ne \epsilon$ and $C \ne \emptyset$ then $k = <<x_1: \varphi_1> \dots <x_m: \varphi_m>>$ and $C = \{ 1, \dots , p \}$ where $1 \leqslant p \leqslant m$, we define $k_{/C} = <<x_1: \varphi_1> \dots <x_p: \varphi_p>>  \in \Theta$ (so $dom(k_{/C}) = \{ 1, \dots , p \} = C$).
\end{itemize}

\medskip

We also define $\mathcal{R}(k) = \{ k_{/C} \vert \ C \in  \mathcal{D}, C \subseteq dom(k)   \}$.

\medskip

Given another $h \in \Theta$ we write $h \sqsubseteq k$ if and only if $h \in \mathcal{R}(k) $ .

\medskip

Suppose $h \in \mathcal{R}(k)$, then there exists $C \in \mathcal{D}$ such that $C \subseteq dom(k)$, $h=k_{/C}$. As we have seen in this case $dom(h) = C$ and $k_{/dom(h)} = k_{/C} = h$.\\

Given $k \in \Theta$ we define $var(k)$ as follows. 

\begin{itemize}
\item if $k = \epsilon$ then $var(k) = \emptyset$,
\item if $k \ne \epsilon$ then $k = <<x_1: \varphi_1> \dots <x_m: \varphi_m>>$ and we define $var(k) = \{ x_1, \dots , x_m \} \, . $
\end{itemize}

\medskip

We are now going to define the `addition' of a new element to a context string.\\

\begin{definition}
Let $h \in \Theta, \, x \in \mathcal{V}, \, \varphi \in \Sigma^*$, we define $h + <x, \varphi>$ as follows:

\begin{itemize}
\item if $h = \epsilon$ then $h + <x, \varphi> = <<x, \varphi>> \in \Theta$;
\item if $h \ne \epsilon$ then let $m$ positive integer $x_1, \dots , x_m \in \mathcal{V}$, $\varphi_1, \dots , \varphi_m \in \Sigma^*$ such that $h = <<x_1: \varphi_1> \dots <x_m: \varphi_m>>$, we define\\
$h + <x, \varphi> = <<x_1: \varphi_1> \dots <x_m: \varphi_m><x, \varphi>> \in \Theta$. 
\end{itemize}

\end{definition}

\bigskip

\begin{lemma}\label{L:addition-properties}
Let $h \in \Theta, \, x \in \mathcal{V}, \, \varphi \in \Sigma^*$, let $k = h + <x, \varphi>$. Then the following hold true:

\begin{itemize}
\item $dom(h) \subseteq dom(k)$,
\item $h = k_{/dom(h)}$,
\item $h \in \mathcal{R}(k)$,
\item $var(k) = var(h) \cup \{ x \}$.
\end{itemize}

\end{lemma}

\begin{proof}
It is obvious by the definition of $k$ that $dom(h) \subseteq dom(k)$.\\

If $h = \epsilon$ then $dom(h) = \emptyset$ and $k_{/dom(h)} = \epsilon = h$.\\

if $h \ne \epsilon$ then let $h = <<x_1: \varphi_1> \dots <x_m: \varphi_m>>$, this implies\\ 
$k = <<x_1: \varphi_1> \dots <x_m: \varphi_m><x, \varphi>>$ and clearly $k_{/dom(h)} = h$.\\

It also follows that $h \in \mathcal{R}(k)$.\\

If $h = \epsilon$ then $var(k) = \{x \} = var(h) \cup \{ x \}$.\\

If $h \ne \epsilon$ and $h = <<x_1: \varphi_1> \dots <x_m: \varphi_m>>$ then 
\[ var(k) = \{ x_1, \dots , x_m \} \cup \{ x \} = var(h) \cup \{ x \} \ . \]
\end{proof}

\bigskip

\begin{lemma}\label{L:add-context-univ}
Given $k \in \Theta - \{\epsilon \}$ there exist $h \in \Theta$, $x \in \mathcal{V}$, $\varphi \in \Sigma^*$ such that 
$k = h + <x, \varphi>$. Moreover $h$, $x$ and $\varphi$ are univocally determined.

\end{lemma}

\begin{proof}

Let $k \in \Theta - \{\epsilon \}$, there exist a positive integer $m$, $x_1, \dots , x_m \in \mathcal{V}, \, \varphi_1, \dots , \varphi_m \in \Sigma^*$ such that $k = <<x_1: \varphi_1> \dots <x_m: \varphi_m>>$.\\

If $m = 1$ then $k = <<x_1: \varphi_1>> = \epsilon + <x_1, \varphi_1>$.\\

If $m > 1$ then $k = <<x_1: \varphi_1> \dots <x_{m-1}: \varphi_{m-1}><x_m: \varphi_m>>$, and so\\
$k = <<x_1: \varphi_1> \dots <x_{m-1}: \varphi_{m-1}>> + <x_m: \varphi_m>$.\\

We have seen there exist $h \in \Theta$, $x \in \mathcal{V}$, $\varphi \in \Sigma^*$ such that 
$k = h + <x, \varphi>$. Suppose there also exist $g \in \Theta$, $y \in \mathcal{V}$, $\psi \in \Sigma^*$ such that 
$k = g + <y, \psi>$.\\

Suppose $h = \epsilon$ and $g \ne \epsilon$, then there exist a positive integer $m$ and $y_1, \dots , y_m \in \mathcal{V}$, $\psi_1, \dots , \psi_m \in \Sigma^*$ such that $g = <<y_1: \psi_1> \dots <y_m: \psi_m>>$. It follows that $k = <<x, \varphi>>$ and $k = <<y_1: \psi_1> \dots <y_m: \psi_m><y, \psi>>$, hence $1 = m+1$, which is false. Therefore $h = \epsilon$ and $g \ne \epsilon$ is false and similarly $h \ne \epsilon$ and $g = \epsilon$ is false.\\

Let's consider the case where $h = g = \epsilon$. In this case $<<x, \varphi>> = k = <<y, \psi>>$, so $y = x$ and $\psi = \varphi$.\\

Finally we consider the case where both $h \ne \epsilon$ and $g \ne \epsilon$. There exist $m$ positive integer $x_1, \dots , x_m \in \mathcal{V}$, $\varphi_1, \dots , \varphi_m \in \Sigma^*$ such that $h = <<x_1: \varphi_1> \dots <x_m: \varphi_m>>$. There also exist $p$ positive integer, $y_1, \dots , y_p \in \mathcal{V}$, $\psi_1, \dots , \psi_p \in \Sigma^*$ such that $g = <<y_1: \psi_1> \dots <y_p: \psi_p>>$. It follows that 

\begin{itemize}
\item $k = <<x_1: \varphi_1> \dots <x_m: \varphi_m><x, \varphi>>$
\item $k = <<y_1: \psi_1> \dots <y_p: \psi_p><y, \psi>>$
\end{itemize}

Hence $p+1 = m+1$, $p = m$, for each $i = 1 \dots m$ $y_i = x_i$, $y = x$, $\psi = \varphi$. Finally $g = h$ also  holds. 

\end{proof}

\begin{lemma}\label{L:useful-contexts-1}
Let $h \in \Theta, \, x \in \mathcal{V}, \, \varphi \in \Sigma^*$, $k = h + <y, \varphi>$. Suppose $g \in \mathcal{R}(k)$ is such that $g \ne k$. Then $g \in \mathcal{R}(h)$.
\end{lemma}

\begin{proof}
Let $D = dom(h)$.\\

We first consider the case where $h \ne \epsilon$. In this case there exists a positive integer $m$ such that $D = \{ 1, \dots , m \}$, and clearly $dom(k) = \{ 1, \dots , m+1 \}$. Since $g \in  \mathcal{R}(k)$ there exists $C \in \mathcal{D}$ such that $C \subseteq \{ 1, \dots , m+1 \}$ and $g = k_{/C}$. Since $g \ne k$ we must have $C \subseteq \{ 1, \dots , m \}$. We have
\[ g = k_{/C} = (k_{/D})_{/C} = h_{/C} \, . \]

\medskip

Let's now consider the case where $h = \epsilon$. In this case $D = \emptyset$ and $dom(k) = \{ 1 \}$. Moreover there exists $C \in \mathcal{D}$ such that $C \subseteq \{ 1 \}$ and $g = k_{/C}$. Since $g \ne k$ we must have $C = \emptyset$ and $g = \epsilon = h$.\\

In both cases $g \in \mathcal{R}(h)$, of course.
\end{proof}

\medskip

\begin{lemma}\label{L:useful-context-2}
Let $k = <<x_1, \varphi_1> \dots <x_m, \varphi_m>> \in \Theta - \{ \epsilon \}$, let $h \in \mathcal{R}(k)$. If $h \ne \epsilon$ then there exists $p = 1 \dots m$ such that $h = <<x_1, \varphi_1> \dots <x_p, \varphi_p>>$.
\end{lemma}

\begin{proof}

If $h \in \mathcal{R}(k)$ then there exists $C \in \mathcal{D}$ such that $C \subseteq dom(k)$, $h=k_{/C}$. If $C = \emptyset$ then $h = \epsilon$, so $C \ne \emptyset$, since $dom(k) = \{ 1, \dots , m \}$ there exists $p = 1 \dots m$ such that $C = \{ 1, \dots , p \}$ and $h = <<x_1, \varphi_1> \dots <x_p, \varphi_p>>$.
\end{proof}

\medskip

\begin{lemma}\label{L:Theta-recursive}
$\Theta$ is recursive.
\end{lemma}

\begin{proof}
Let $k \in \Gamma^*$ and let's see how we can decide wheter $k \in \Theta$.\\

If $k$ has less than two characters then clearly $k \notin \Theta$. If $k$ has exactly two characters then: if $k = \ \texttt{`<>'}$ then $k \in \Theta$ otherwise $k \notin \Theta$.\\

We still have to examine the case where $k$ has more than two characters. Recall that if $t$ is a string we can indicate with $\ell(t)$ $t$'s length, i.e. the number of characters in $t$. If $\ell(t) > 0$ then for each $\alpha \in \{1, \dots , \ell(t) \}$ at position $\alpha$ within $t$ there is a character, this symbol can be indicated with $t[\alpha]$.\\

If $k[1] \ne \texttt{`<'}$ or $k[\ell(k)] \ne \texttt{`>'}$ then $k \notin \Theta$.\\

Otherwise we can assume $k[1] = \texttt{`<'}$ and $k[\ell(k)] = \texttt{`>'}$.\\

If $\ell[k] = 3$ then $k \notin \Theta$, so hereafter $\ell[k] \geqslant 4$.\\

If $k[2] \ne \texttt{`<'}$ or $k[\ell(k)-1] \ne \texttt{`>'}$ then $k \notin \Theta$.\\

Otherwise we can assume $k[2] = \texttt{`<'}$ and $k[\ell(k)-1] = \texttt{`>'}$.\\

At this point we can try to create two arrays of indexes $\alpha$ and $\beta$ with the same size $p \geqslant 1$.\\

We set $\alpha[1] = 2$, and let $\beta[1]$ be the first index $i$: $3 \leqslant i \leqslant \ell(k)-1$ such that $k[i] = \texttt{`>'}$.\\

If $\beta[1] = \ell(k)-1$ then we decide that $p = 1$ and our arrays $\alpha$ and $\beta$ are defined.\\

Otherwise $\beta[1] < \ell(k)-1$. Here if $k[\beta[1] + 1] \ne \texttt{`<'}$ then $k \notin \Theta$ and we have finished, so let's assume $k[\beta[1] + 1] = \texttt{`<'}$. We also assume $p \geqslant 2$ and let $\alpha[2] = \beta[1] + 1$ and $\beta[2]$ be the first index $i$: $\beta[1] + 1 < i \leqslant \ell(k)-1$ such that $k[i] = \texttt{`>'}$.\\

If $\beta[2] = \ell(k)-1$ then we decide that $p = 2$ and our arrays $\alpha$ and $\beta$ are defined.\\

Otherwise $\beta[2] < \ell(k)-1$ and we can continue as above.\\

Suppose at a certain step $j \geqslant 2$ we have defined $\alpha[j] = \beta[j-1] + 1$ and $\beta[j]$ as the first index $i$: $\beta[j-1] + 1 < i \leqslant \ell(k)-1$ such that $k[i] = \texttt{`>'}$.\\

If $\beta[j] = \ell(k)-1$ then we decide that $p = j$ and our arrays $\alpha$ and $\beta$ are defined.\\

Otherwise $\beta[j] < \ell(k)-1$. Here if $k[\beta[j] + 1] \ne \texttt{`<'}$ then $k \notin \Theta$ and we have finished, so let's assume $k[\beta[j] + 1] = \texttt{`<'}$. We also assume $p \geqslant j+1$ and let $\alpha[j+1] = \beta[j] + 1$ and $\beta[j+1]$ be the first index $i$: $\beta[j] + 1 < i \leqslant \ell(k)-1$ such that $k[i] = \texttt{`>'}$.\\

If $\beta[j+1] = \ell(k)-1$ then we decide that $p = j+1$ and our arrays $\alpha$ and $\beta$ are defined.\\

This process goes on like this until it terminates, it can terminate because we have established that $k \notin \Theta$, or otherwise because at a certain index $h$ we have $\beta[h] = \ell(k)-1$. If we don't reach the condition $k \notin \Theta$, we will certainly sooner or later reach $\beta[h] = \ell(k)-1$. In fact at each step j $\beta[j+1] > \beta[j]$, and our string $k$ has a finite number of characters.\\

So, our process can decide that $k \notin \Theta$, otherwise it will determine the size $p$ of our arrays and the arrays will themselves be defined at the end of the process. In this case for each $i = 1 \dots p$ we'll need to determine if the substring $\varphi$ of $k$ between the characters $\alpha[i]$ and $\beta[i]$, which we can indicate for instance with $k[\alpha[i], \beta[i]]$, has the form we need in order to decide that $k \in \Theta$.\\

The string $\varphi$ needs to be of the form $<x:\psi>$ where $x \in \mathcal{V}$, $\psi \in \Sigma^*$. How can we decide if this is the case?\\

Of course $\varphi$ has a finite length and for each $j = 1 \dots \ell(\varphi)$ we must be able to take a decision over $\varphi[j]$. The decisions are the following:

\begin{itemize}
\item $\varphi[1] = \texttt{`<'}$,
\item $\varphi[2] \in \mathcal{V}$,
\item $\varphi[3] = \texttt{`:'}$,
\item $\varphi[\ell(\varphi)] = \texttt{`>'}$,
\item if $\ell(\varphi) > 4$ for each $j = 4 \dots \ell(\varphi) - 1$ $\varphi[j] \in \Sigma$. 
\end{itemize}

How can we decide that $\varphi[2] \in \mathcal{V}$? Since in our alphabet there are only a finite number of characters that are not variables, we just need to verify that $\varphi[2]$ is not one of those characters.\\

How can we decide that $\varphi[j] \in \Sigma$? Here we just need to verify that $\varphi[j]$ is not $\texttt{`<'}$ or $\texttt{`>'}$. The other decisions are more than trivial, so we are able to decide if $\varphi$ has the desired form.\\

If for each $i = 1 \dots p$ $k[\alpha[i], \beta[i]]$ has the desired form we decide that $k \in \Theta$, otherwise $k \notin \Theta$.
\end{proof}

\section{Building the expressions of our system}

We can now describe the process of constructing expressions for our language $\mathcal{L}$. This is an
inductive process in which not only we build expressions, but also we associate them with meaning,
and in parallel also define the fundamental concept of `context'. This process will be identified as
`Definition~\ref{D:expr-new}' although actually it is a process in which we give the definitions and prove
properties which are needed in order to set up those definitions.\\

Within this definition we will define the expressions of our language. Such expressions are finite sequences of characters of the alphabet $\Sigma = \mathcal{V} \cup \mathcal{C} \cup \mathcal{F} \cup \mathcal{Z}$. In other words they are members of $\Sigma^*$.\\

Since this is a complex definition, we will first try to provide an informal idea of the entities we'll
define in it. The definition is by induction on positive integers, we now introduce the sets and concepts
we'll define for a generic positive integer n (this first listing is not the true definition, it's just to introduce the concepts, to enable the reader to understand their role).

\medskip

$K(n)$ is the set of `contexts' at step n. If we define $\Gamma = \Sigma \cup \{ <, > \}$, contexts will be defined as members of $\Gamma^*$, and they will be strings of the form $<<x_1,\varphi_1>, \dots , <x_m,\varphi_m >>$, where $x_1, \dots , x_m \in \mathcal{V}$ and $\varphi_1, \dots , \varphi_m$ are expressions. The string $<>$ which we'll also name $\epsilon$ is also a possibile context, and when we use the symbol $\epsilon$ with respect to a context we actually mean $<>$.

\medskip

For each $k \in K(n)$ $\Xi(k)$ is the set of `states' bound to context $k$. If $n>1$ and $k \in K(n-1)$ then $\Xi(k)$ has already been defined at step $n-1$ or formerly, otherwise it will be defined at step $n$.

\medskip

If $k = <<x_1,\varphi_1>, \dots , <x_m,\varphi_m >>$ is a context, a state on $k$ is a state-like pair $\sigma = (x,s)$ where of course x is the function which associates $x_i$ to each $i = 1 \dots m$ and (roughly speaking) for each $i = 1 \dots m$ $s_i$ is a member of the meaning of the corresponding expression $\varphi_i$ .

\medskip

For each $k \in K(n)$ $E(n,k)$ is the set of expressions bound to step $n$ and context $k$. And here it is important to underline that \textbf{we  need to ensure that $E(n,k)$ (as a subset of $\Sigma^*$) is a recursive set}.

\medskip

$E(n)$ is the union of $E(n,k)$ for $k \in K(n)$ (this will not be explicitly recalled on each iteration in the definition).

\medskip

For each $k \in K(n)$, $t \in E(n,k)$, $\sigma \in \Xi(k)$ we'll define $\#(k,t,\sigma)$ which stands for `the meaning of $t$ bound to $k$ and $\sigma$'. 

\medskip

The following set $E_s(n,k)$ should be defined in the same way at each step, we put here its definition, to avoid to repeat that definition each time. For each $k \in K(n)$ we define 
\[ E_s(n,k) = \{t | t \in E(n,k), \forall \sigma \in \Xi(k) \, \#(k,t,\sigma) \text{ is a set} \} . \]

\subsection{Definition process}

This section contains only definition~\ref{D:expr-new}. This definition is an inductive definition process within which we have assumptions, lemmas etc.. Symbols like $\qed$ within this definition are not intended to terminate the definition, they just terminate an assumption or lemma etc. which is internal to the definition. 

\begin{definition}\label{D:expr-new}
We are now ready to begin the actual definition process, so we perform the simple initial step of our inductive process.\\

We define $K(1) = \{\epsilon\}$, $\Xi(\epsilon) = \{\epsilon\}$, $E(1,\epsilon) = \mathcal{C}$.\\

Clearly when we define $K(1)$ with $\epsilon$ we mean the string $<>$, while when defining $\Xi(\epsilon) = \{\epsilon\}$ the $\epsilon$ on the left side is $<>$ and the $\epsilon$ on the right side is $(\emptyset, \emptyset)$.\\

For each $t \in E(1,\epsilon)$ we define $\#(\epsilon,t,\epsilon) = \#(t)$.\\

The inductive step is a bit more complex. Suppose all our definitions have been given at step $n$ and
let's proceed with step $n+1$. In this inductive step we'll need some assumptions which will be
identified with a title like `Assumption~\ref{D:expr-new}.x'. Each assumption is a statement that must be valid at
step $1$, we suppose is valid at step $n$ and needs to be proved true at step $n+1$ at the end of our
definition process.\\

The first assumptions we need are the following.\\

\begin{assumption}\label{A:Kn-sub-Theta}
$K(n) \subseteq \Theta$. \hfill $\qed$
\end{assumption}

\begin{assumption}\label{A:Kn-recursive}
$K(n)$ is recursive and $\epsilon \in K(n)$. \hfill $\qed$
\end{assumption}

\begin{assumption}\label{A:Xik-notempty}
For each $k \in K(n)$ $\Xi(k) \ne \emptyset$. \hfill $\qed$
\end{assumption}

\begin{assumption}\label{A:k-m-in-k-n-new}
If $n > 1$ then for each $m < n$ $K(m) \subseteq K(n)$. \hfill $\qed$
\end{assumption}

\begin{assumption}\label{A:Ekn-sub-Closure}
For each $k \in K(n)$ $E(n,k) \subseteq \Sigma^*$. \hfill $\qed$
\end{assumption}

\begin{assumption}\label{A:Ekn-recursive}
For each $k \in K(n)$ $E(n,k)$ is recursive. \hfill $\qed$
\end{assumption}

\begin{assumption}\label{A:simple-k-sigma-slp-new}
For each $k \in K(n)$ $k \in \Theta$ and for each $\sigma \in \varXi(k)$ $\sigma$ is a state-like pair and $dom(\sigma) = dom(k)$. \hfill $\qed$
\end{assumption}

\begin{assumption}\label{A:k-is-epsilon-or-new}
For each $k \in K(n)$ $k = \epsilon$ and $\varXi(k) = \{\epsilon\}$ or\\
($n > 1$ and there exist $m<n$, $h \in K(m)$, $\phi \in E_s(m,h)$, $y \in (\mathcal{V}-var(h))$ such that $k = h + <y, \phi>$, $\Xi(k) = \{ \sigma + (y,s) | \, \sigma \in \Xi(h), s \in \#(h,\phi,\sigma) \} )$. \hfill $\qed$
\end{assumption}

\begin{assumption}\label{A:rest-of-state-is-state}
If $n > 1$ then for each $k \in K(n): k \ne \epsilon$, $\sigma \in \Xi(k)$, $h \in \mathcal{R}(k): h \ne k$, there exists $m < n$ such that $h \in K(m)$ and it results $\sigma_{/dom(h)} \in \Xi(h)$. \hfill $\qed$
\end{assumption}

\begin{assumption}\label{A:Ekn-decidable-predicates}
For each $k \in K(n)$, $\varphi \in E(n,k)$, $q$ positive integer, $i = 1 \dots p$ we must be able to decide all of the following conditions:

\begin{itemize}
\item for each $\sigma \in \Xi(k)$ $Set_q( \#(k, \varphi, \sigma))$;
\item for each $\sigma \in \Xi(k)$ $Event_q( \#(k, \varphi, \sigma))$;
\item for each $\sigma \in \Xi(k)$ $\#(k, \varphi, \sigma) \in D_i$;
\item for each $\sigma \in \Xi(k)$ $\#(k, \varphi, \sigma) \in \mathcal{P}^q(D_i)$;
\item if (for each $\sigma \in \Xi(k)$ $Set_q( \#(k, \varphi, \sigma))$) then\\
(for each $\sigma \in \Xi(k)$ $NotEmpty_q( \#(k, \varphi, \sigma))$).
\end{itemize}

\medskip

Moreover, for the last condition, we must be able to decide it is true. \hfill $\qed$
\end{assumption}

Clearly assumption~\ref{A:Ekn-decidable-predicates} is valid with $n = 1$, in fact in this case $k = \epsilon$, $E(n,k) = E(1,\epsilon) = \mathcal{C}$, $\Xi(k) = \{ \epsilon \}$, so $\varphi \in \mathcal{C}$ and $\#(k, \varphi, \sigma) = \#( \varphi )$, and the conditions are the following:

\begin{itemize}
\item $Set_q( \#(\varphi))$;
\item $Event_q( \#(\varphi))$;
\item $\#(\varphi) \in D_i$;
\item $\#(\varphi) \in \mathcal{P}^q(D_i)$;
\item if ($Set_q( \#(\varphi))$) then ($NotEmpty_q( \#(\varphi))$).
\end{itemize}

We have in fact assumed to be able to decide all of these conditions, and to be able to decide as true the last of these conditions.\\

We can go on with the inductive step and define
\begin{align} 
K(n)^+ &= \{ h + <y,\phi> | \, h \in K(n), \phi \in E_s(n,h), y \in (\mathcal{V}-var(h)) \} - K(n) \ , \notag  \\
K(n+1) &= K(n) \cup K(n)^+ \  \notag .
\end{align}

\medskip

Let $k \in K(n)^+$. Then there exist $h \in K(n), \phi \in E_s(n,h), y \in (\mathcal{V}-var(h))$ such that $k = h + <y,\phi>$. By lemma~\ref{L:add-context-univ} we know that $h, \phi, y$ are \emph{univocally determined}.\\

We can assume that $\Xi(k)$ is defined for $k \in K(n)$, and we need to define this for $k \in K(n+1) - K(n)$, i.e. for $k \in K(n)^+$. If $k \in K(n)^+$ there exist $h \in K(n), \phi \in E_s(n,h), y \in (\mathcal{V}-var(h))$ such that $k = h + <y,\phi>$; and $h, \phi, y$ are univocally determined. So we can define 
\[
\Xi(k) = \{ \sigma + (y,s) | \, \sigma \in \Xi(h), s \in \#(h,\phi,\sigma) \} \, .
\]

\medskip

A consequence of lemma~\ref{L:add-slp-univ} is the following: for each $k \in K(n)^+$ and $\sigma + (y,s)$ in $\Xi(k)$, $\sigma$, $y$ and $s$ are \emph{univocally determined}.\\

To ensure the unique readability of our expressions we need the following assumption (which is clearly satisfied for $n=1$).

\begin{assumption}\label{A:expr-cont-depth}
For each $t \in E(n)$
\begin{itemize}
\item $t[\ell(t)] \ne \text{`('}$ ;
\item if $t[\ell(t)] = \text{`)'}$ then $d(t, \ell(t)) = 1$, else $d(t, \ell(t)) = 0$ ;
\item for each $\alpha \in \{1, \dots, \ell(t) \}$ if $(t[\alpha]=\text{`:'}) \vee (t[\alpha]=\text{`,'}) \vee (t[\alpha]=\text{`)'}) $ then $d(t, \alpha) \geqslant 1$.

\end{itemize}
\end{assumption}

\begin{flushright}
$\qed$\\
\end{flushright}

\medskip

We immediately prove the following.

\begin{proof}[Proof of~\ref{A:Kn-sub-Theta}]
Given that $K(n) \subseteq \Theta$ we have to show that $K(n+1) \subseteq \Theta$.\\

Let $k \in K(n+1)$, if $k \in K(n)$ then $k \in \Theta$, else there exist $h \in K(n) \subseteq \Theta$, $\phi \in E_s(n,h) \subseteq \Sigma^*$, $y \in (\mathcal{V}-var(h))$ such that $k = h + <y,\phi> \in \Theta$.
\end{proof}

\medskip

\begin{proof}[Proof of~\ref{A:Kn-recursive}]
We have to show that $K(n+1)$ is recursive and $\epsilon \in K(n+1)$.\\

We have assumed by inductive hypothesis that $K(n)$ is recursive and that $\epsilon \in K(n)$.\\

First of all it is obvious that $\epsilon \in K(n+1)$ because $K(n) \subseteq K(n+1)$.\\

Let $k \in \Gamma^*$ and let's try to decide whether $k \in K(n+1)$. We can decide whether $k \in K(n)$, if this holds then $k \in K(n+1)$, otherwise we know $k \ne \epsilon$.\\

At this point, given $\Theta$ is recursive, we can decide if $k \in \Theta$, if $k \notin \Theta$ using lemma~\ref{A:Kn-sub-Theta} we can decide that $k \notin K(n+1)$.\\

If $k \in \Theta$ then $k \in \Theta - \{ \epsilon \}$, so there exist $h \in \Theta$, $y \in \mathcal{V}$, $\psi \in \Sigma^*$ such that $k = h + <y, \psi>$, we know how to calculate  $h$, $y$, $\psi$, and they are univocally determined.\\

Now consider the following conditions

\begin{itemize}
\item $h \in K(n)$,
\item $y \in \mathcal{V} - var(h)$,
\item $\psi \in E_s(n,h)$.
\end{itemize}

\medskip

If all these conditions hold, then $k \in K(n+1)$, else (knowing that $k \notin K(n)$) $k \notin K(n+1)$.\\

All of the mentioned conditions are decidable. In fact $K(n)$ is recursive and so we can decide whether $h \in K(n)$. Moreover $E(n,h)$ is recursive, and given $\psi \in E(n,h)$ the condition `for each $\rho \in \Xi(h)$ $Set_1(\#(h, \psi, \rho))$' is decidable. Therefore $\psi \in E_s(n,h)$ is decidable.\\

As regards the condition $y \in \mathcal{V} - var(h)$, we know that $y$ is a variable, so if it doesn't belong to $var(h)$ this means it belongs to $\mathcal{V} - var(h)$, and so we can also decide this condition.\\

Therefore we have proved that $K(n+1)$ is recursive. 
\end{proof}

\medskip

\begin{proof}[Proof of~\ref{A:Xik-notempty}]

Let $k \in K(n+1)$, we have to show $\Xi(k) \ne \emptyset$.\\

If $k \in K(n)$ then $\Xi(k) \ne \emptyset$, else there exist $h \in K(n)$, $\phi \in E_s(n,h)$, $y \in (\mathcal{V}-var(h))$ such that $k = h + <y, \phi>$, and $\Xi(k) = \{ \sigma + (y,s) | \, \sigma \in \Xi(h), s \in \#(h,\phi,\sigma) \}$.\\

By the inductive hypothesis $\Xi(h) \ne \emptyset$, let's then take $\sigma \in \Xi(h)$, then $Set_1(\#(h, \phi, \sigma))$ and $NotEmpty_1(\#(h, \phi, \sigma))$. If we take $s \in \#(h,\phi,\sigma)$ then $\sigma + (y,s) \in \Xi(k)$ .
\end{proof}

\bigskip

It is time to define $E(n+1,k)$, for each $k$ in $K(n+1)$. Then for each $t$ in $E(n+1,k)$ and $\sigma$ in $\Xi(k)$ we need to define $\#(k,t,\sigma)$. We begin to do this by defining some new sets of expressions bound to context $k$, and for the expressions in each new set we define the proposed value of $\#(k,t,\sigma)$.\\

For each $k = h + <y,\phi> \in K(n)^+$ we define

\[
E_a(n+1,k) = \{ y \}.
\]

\medskip

Clearly $E_a(n+1,k) \subseteq \Sigma^*$ and $E_a(n+1,k)$ recursive.

\medskip

For each $t \in E_a(n+1,k), \  \sigma = \rho + (y,s) \in \Xi(k)$ we define:

\[
\#(k,t,\sigma)_{(n+1,k,a)} = s.
\]

\bigskip

We notice that $\epsilon \in K(n)$ and define $E_b(n+1,\epsilon) = \emptyset$.

\medskip

For each $k = h + <y,\phi> \in K(n) - \{ \epsilon \}$ we define

\[
E_b(n+1,k) = \{ t | t \in E(n,h), t\notin E(n,k) \}.
\]

\medskip

Clearly $E_b(n+1,k) \subseteq E(n,h) \subseteq \Sigma^*$ and $E_b(n+1,k)$ recursive.

\medskip

For each $t \in E_b(n+1,k), \  \sigma = \rho + (y,s) \in \Xi(k)$ we define the proposed value of $\#(k,t,\sigma)$:

\[
\#(k,t,\sigma)_{(n+1,k,b)} = \#(h,t,\rho).
\]

\bigskip

Given $k \in K(n)$ and a constant $c \in \mathcal{C}$ we can define the following set 
\[ H_c(n+1, k) = \{ (c)( \varphi_1, \dots , \varphi_m) | \, \varphi_1, \dots , \varphi_m \in E(n,k) \} \, . \]

and we can prove it is recursive using some auxiliary lemma.\\

\begin{ilemma}\label{IL:Hc-recursive-aux0}
Let $m$ positive integer, $m > 1$, let $\psi_1, \dots , \psi_m \in E(n,k)$, $\varphi = (c)( \psi_1, \dots , \psi_m)$. Let $q_1, \dots, q_{m-1}$  the positions of the explicit occurrences of `,' in the representation $(c)( \psi_1, \dots , \psi_m)$ of $\varphi$. Then for each $j = 1 \dots m-1$ $d(\varphi, q_j) = 1$.
\end{ilemma}

\begin{proof}
Let $\vartheta = \varphi[1,4]$ be the substring of $\varphi$ which begins at character $1$ and ends at character $4$. We first want to prove that $d(\varphi, q_1) = 1$. We have that 
\[ d(\varphi,q_1 - 1) = d(\varphi, \ell(\vartheta) + \ell(\psi_1)) = d(\varphi, \ell(\vartheta) + 1) + d( \psi_1, \ell(\psi_1)) = 1 + d( \psi_1, \ell(\psi_1)) \ . \]

\medskip

If $\varphi[q_1 - 1] = \psi_1[\ell(\psi_1)] = \text{`)'}$ then $d(\varphi, q_1) = d(\varphi,q_1 - 1) - 1 = d(\psi_1, \ell(\psi_1)) = 1$ .\\
Else $\varphi[q_1 - 1] = \psi_1[\ell(\psi_1)] \notin \{ \text{`('}, \text{`)'} \}$ so $d(\varphi, q_1) = d(\varphi, q_1 - 1) = 1 + d(\psi_1, \ell(\psi_1)) = 1$.\\

If $m > 2$ we want to prove that for each $j = 2 \dots m-1$ $d(\varphi, q_j) = 1$. Let $j = 2 \dots m-1$, we can assume $d(\varphi, q_{j-1}) = 1$. Let $\eta = \varphi[1, q_{j-1} ]$, then
\[ d(\varphi,q_j - 1) = d(\varphi, \ell(\eta) + \ell(\psi_j)) = d(\varphi, q_{j-1} + 1) + d( \psi_j, \ell(\psi_j)) = 1 + d( \psi_j, \ell(\psi_j)) \ . \]

\medskip

If $\varphi[q_j - 1] = \psi_j[\ell(\psi_j)] = \text{`)'}$ then $d(\varphi, q_j) = d(\varphi,q_j - 1) - 1 = d(\psi_j, \ell(\psi_j)) = 1$ .\\
Else $\varphi[q_j - 1] = \psi_j[\ell(\psi_j)] \notin \{ \text{`('}, \text{`)'} \}$ so $d(\varphi, q_j) = d(\varphi, q_j - 1) = 1 + d(\psi_j, \ell(\psi_j)) = 1$.
\end{proof}

\medskip

\begin{ilemma}\label{IL:Hc-recursive-aux1}
Let $\psi \in \Sigma^*$ and let $\varphi = (c)(\psi) \in \Sigma^*$. Suppose for each $r$ positive integer such that $4 < r < \ell(\varphi)$ and $\varphi[r]$ = `,' we have $d(\varphi, r) \ne 1$. Then $\varphi \in H_c(n+1, k)$ if and only if $\psi \in E(n,k)$.
\end{ilemma}

\begin{proof}
It is obvious that if $\psi \in E(n,k)$ then $\varphi \in H_c(n+1, k)$.\\

Conversely, if $\varphi \in H_c(n+1, k)$ then there exist a positive integer $m$ and $\psi_1, \dots , \psi_m \in E(n,k)$ such that $\varphi = (c)( \psi_1, \dots , \psi_m)$.\\

If $m > 1$ then let $r$ be the first explicit occurrence of `,' in $(c)( \psi_1, \dots , \psi_m)$. By lemma~\ref{IL:Hc-recursive-aux0} we have that $d(\varphi, r) = 1$. Anyway we required that $d(\varphi, r) \ne 1$, so we have a contradiction and it cannot be $m>1$.\\

It follows that $m = 1$, then $(c)(\psi) = \varphi = (c)( \psi_1)$ and so $\psi = \psi_1 \in E(n,k)$.
\end{proof}

\medskip

\begin{ilemma}\label{IL:Hc-recursive-aux2}
Let $\psi \in \Sigma^*$ and let $\varphi = (c)(\psi) \in \Sigma^*$. Consider the set of the positive integers $r$ such that $4 < r < \ell(\varphi)$, $\varphi[r]$ = `,' and $d(\varphi, r) = 1$. Assume this set is not empty and let's name $r_1, \dots, r_h$ its members (in increasing order).\\ 
Let's also define the following. Let $\epsilon$ be the empty string over the alphabet $\Sigma$.\\

If $r_1 = 5$ then let $\psi_1 = \epsilon$ else $r_1 > 5$ and let $\psi_1 = \varphi[5, r_1 - 1]$.\\

If $h > 1$ then for each $i = 1 \dots h-1$: 
if $r_{i+1} = r_i + 1$ then let $\psi_{i+1} = \epsilon$ else $r_{i+1} > r_i + 1$ and let $\psi_{i+1} = \varphi[r_i + 1, r_{i+1} - 1]$.\\

Finally if $r_h = \ell(\varphi) - 1$ we define $\psi_{h+1} = \epsilon$, else $r_h < \ell(\varphi) - 1$ and let $\psi_{h+1} = \varphi[r_h + 1, \ell(\varphi) - 1]$ .\\

With these definitions we have $\varphi = (c)( \psi_1, \dots , \psi_{h+1})$ and $\varphi \in H_c(n+1, k)$ if and only if for each $i = 1 \dots h+1$ $\psi_i \in E(n,k)$.
\end{ilemma}

\begin{proof}
Clearly if for each $i = 1 \dots h+1$ $\psi_i \in E(n,k)$ then $\varphi \in H_c(n+1, k)$.\\

Conversely, if $\varphi \in H_c(n+1, k)$ then there exist a positive integer $m$ and $\chi_1, \dots , \chi_m \in E(n,k)$ such that $\varphi = (c)( \chi_1, \dots , \chi_m)$.\\ 

If $m = 1$ then $\varphi = (c)( \chi_1 )$, we have $d(\varphi, r_1) = 1$. We have also $4 < r_1 < \ell(\varphi)$ $\chi_1[r_1 - 4] = \varphi[r_1]$ = `,'. Moreover let $\vartheta = \varphi[1,4]$ be the substring of $\varphi$ which begins at character $1$ and ends at character $4$. We have that
\[ 1 = d(\varphi, r_1) = d(\varphi, 4 + r_1 - 4) = d( \varphi, \ell(\theta) + r_1 - 4) = d( \varphi, \ell(\theta) + 1) + d(\chi_1, r_1 - 4) \ . \]

It follows that $1 = 1 + d(\chi_1, r_1 - 4)$ and so $d(\chi_1, r_1 - 4) = 0$. This contradicts assumption~\ref{A:expr-cont-depth} and therefore we cannot have $m = 1$.\\

Since $m > 1$ we can indicate with $q_1, \dots, q_{m-1}$  the positions of the explicit occurrences of `,' in the representation $(c)( \chi_1, \dots , \chi_m)$ of $\varphi$.\\

By lemma~\ref{IL:Hc-recursive-aux0} we have that for each $j = 1 \dots m-1$ $d(\varphi, q_j) = 1$, therefore $\{ q_1, \dots, q_{m-1} \} \subseteq \{ r_1, \dots, r_h \}$.\\

Suppose there exists $i = 1 \dots h$ such that $r_i \notin \{ q_1, \dots, q_{m-1} \}$. In this case one of these conditions will occur: 

\begin{itemize}
\item $r_i < q_1$,
\item $r_i > q_{m-1}$,
\item $m-1 > 1$ and there exists $j = 1 \dots m-2$ such that $q_j < r_i < q_{j+1}$.
\end{itemize}

If $r_i < q_1$ then $4 < r_i$ also holds, $\chi_1 = \varphi[5, q_1 - 1]$, $\ell(\chi_1) = q_1 - 1 - 5 + 1 =  q_1 - 5$, for each $\alpha = 1 \dots q_1 - 5$ $\chi_1[\alpha] = \varphi[4 + \alpha]$. So $r_i - 4 \geqslant 1$, $r_i - 4 < q_1 - 4$ and then $r_i - 4 \leqslant q_1 - 5 = \ell(\chi_1)$. Then also $\chi_1[r_i - 4] = \varphi[r_i]$ = `,' . If we define $\vartheta = \varphi[1,4]$ and observe that $r_i = 4 + (r_i - 4) = \ell(\vartheta) + r_i - 4$ then we can also observe that $1 = d(\varphi, r_i) = d(\varphi, \ell(\vartheta) + 1) + d(\chi_1, r_i - 4) = 1 + d(\chi_1, r_i - 4)$, so $d(\chi_1, r_i - 4) = 0$. This contradicts assumption~\ref{A:expr-cont-depth} and therefore we cannot have  $r_i < q_1$.\\

If $r_i > q_{m-1}$ then $r_i < \ell(\varphi)$ also holds, $\chi_m = \varphi[q_{m-1} + 1, \ell(\varphi) - 1]$, $\ell(\chi_m) = \ell(\varphi) - 1 - (q_{m-1} + 1) + 1 = \ell(\varphi) - q_{m-1} - 1$. For each $\alpha = 1 \dots \ell(\varphi) - q_{m-1} - 1$ $\chi_m[\alpha] = \varphi[q_{m-1} + \alpha]$. So $r_i - q_{m-1} \geqslant 1$, $r_i - q_{m-1} < \ell(\varphi) - q_{m-1}$ and then $r_i - q_{m-1} \leqslant \ell(\varphi) - q_{m-1} - 1 = \ell(\chi_m)$. Then also $\chi_m[r_i - q_{m-1}] = \varphi[r_i]$ = `,'. If we define $\vartheta = \varphi[1,q_{m-1}]$ and observe that $r_i = q_{m-1} + (r_i - q_{m-1}) = \ell(\vartheta) + r_i - q_{m-1}$ then we can also observe that $1 = d(\varphi, r_i) = d(\varphi, \ell(\vartheta) + 1) + d(\chi_m, r_i - q_{m-1}) = 1 + d(\chi_m, r_i - q_{m-1})$. Therefore $d(\chi_m, r_i - q_{m-1}) = 0$. This contradicts assumption~\ref{A:expr-cont-depth} and therefore we cannot have $r_i > q_{m-1}$.\\

Finally assume $m-1 > 1$ and there exists $j = 1 \dots m-2$ such that $q_j < r_i < q_{j+1}$. In this case $\chi_{j+1} = \varphi[q_j + 1, q_{j+1} - 1]$, $\ell(\chi_{j+1}) = q_{j+1} - 1 - (q_j + 1) + 1 = q_{j+1} - q_j - 1$. For each $\alpha = 1 \dots q_{j+1} - q_j - 1$ $\chi_{j+1}[\alpha] = \varphi[q_j + \alpha]$. So $r_i - q_j \geqslant 1$, $r_i - q_j < q_{j+1} - q_j$ and then $r_i - q_j \leqslant q_{j+1} - q_j - 1 = \ell(\chi_{j+1})$. Then also $\chi_{j+1}[r_i - q_j] = \varphi[r_i]$ = `,'. If we define $\vartheta = \varphi[1,q_j]$ and observe that $r_i = q_j + (r_i - q_j) = \ell(\vartheta) + r_i - q_j$ then we can also observe that $1 = d(\varphi, r_i) = d(\varphi, \ell(\vartheta) + 1) + d(\chi_{j+1}, r_i - q_j) = 1 + d(\chi_{j+1}, r_i - q_j)$. Therefore $d(\chi_{j+1}, r_i - q_j) = 0$. This contradicts assumption~\ref{A:expr-cont-depth} and therefore we cannot have that $m-1 > 1$ and there exists $j = 1 \dots m-2$ such that $q_j < r_i < q_{j+1}$.\\

So we have to conclude that $\{ q_1, \dots, q_{m-1} \} = \{ r_1, \dots, r_h \}$. This means that $h+1 = m$ and for each $i = 1 \dots h+1$ $\psi_i = \chi_i \in E(n,k)$.
\end{proof}

\medskip

\begin{ilemma}\label{IL:Hc-recursive}
Given $k \in K(n)$ and $c \in \mathcal{C}$ $H_c(n+1, k)$ is recursive.
\end{ilemma}

\begin{proof}

Let $\varphi \in \Sigma^*$. If $\varphi$ doesn't begin with the four characters $(c)($ or doesn't end with the character $)$ then $\varphi \notin H_c(n+1, k)$.\\

Then assume we are in the case $\varphi = (c)(\psi)$ where $\psi \in \Sigma^*$. Consider the set of the positive integers $r$ such that $4 < r < \ell(\varphi)$, $\varphi[r]$ = `,' and $d(\varphi, r) = 1$.\\

If the mentioned set is empty then $\varphi \in H_c(n+1, k)$ if and only if $\psi \in E(n,k)$.\\

If the mentioned set is not empty then let's name $r_1, \dots, r_h$ its members (in increasing order). Let's also provide some other definitions.\\

If $r_1 = 5$ then let $\psi_1 = \epsilon$ else $r_1 > 5$ and let $\psi_1 = \varphi[5, r_1 - 1]$.\\

If $h > 1$ then for each $i = 1 \dots h-1$: 
if $r_{i+1} = r_i + 1$ then let $\psi_{i+1} = \epsilon$ else $r_{i+1} > r_i + 1$ and let $\psi_{i+1} = \varphi[r_i + 1, r_{i+1} - 1]$.\\

Finally if $r_h = \ell(\varphi) - 1$ we define $\psi_{h+1} = \epsilon$, else $r_h < \ell(\varphi) - 1$ and let $\psi_{h+1} = \varphi[r_h + 1, \ell(\varphi) - 1]$ .\\

With these definitions we have $\varphi = (c)( \psi_1, \dots , \psi_{h+1})$ and $\varphi \in H_c(n+1, k)$ if and only if for each $i = 1 \dots h+1$ $\psi_i \in E(n,k)$.
\end{proof}

\medskip

\begin{ilemma}\label{IL:Hc-aux3}
Let $k \in K(n)$, $c \in \mathcal{C}$. There exists an algorithm that given $\varphi \in \Sigma^*$
\begin{itemize}
\item determines if $\varphi \in H_c(n+1, k)$,
\item if $\varphi \in H_c(n+1, k)$ it also identifies a positive integer $m$ and $\psi_1, \dots , \psi_m \in E(n,k)$ such that $\varphi = (c)( \psi_1, \dots , \psi_m)$.
\end{itemize}
\end{ilemma}

\begin{proof}
See the proof of lemma~\ref{IL:Hc-recursive}.
\end{proof}

\medskip

\begin{ilemma}\label{IL:Hc-aux4}
Given $k \in K(n)$, $c \in \mathcal{C}$, $\varphi \in H_c(n+1, k)$ there exist $m$ positive integer, $\chi_1, \dots , \chi_m \in E(n,k)$ such that $\varphi = (c)( \chi_1, \dots , \chi_m)$ and $m$ and $\chi_1, \dots , \chi_m$ are univocally determined.
\end{ilemma}

\begin{proof}
It is obvious by the definition of $H_c(n+1, k)$ there exist $m$ positive integer, $\chi_1, \dots , \chi_m \in E(n,k)$ such that $\varphi = (c)( \chi_1, \dots , \chi_m)$.\\

Suppose there are also $p$ positive integer and $\varphi_1, \dots , \varphi_p$ such that $\varphi = (c)( \varphi_1, \dots , \varphi_p)$. Of course we want to show that $p = m$ and for each $i = 1 \dots m$ $\varphi_i = \chi_i$.\\

To this end we consider there exists $\psi \in \Sigma^*$ such that $\varphi = (c)( \psi )$. Consider the set of the positive integers $r$ such that $4 < r < \ell(\varphi)$, $\varphi[r]$ = `,' and $d(\varphi, r) = 1$.\\

Suppose the mentioned set is empty. In this case if $m > 1$ then let $r$ be the first explicit occurrence of `,' in $(c)( \chi_1, \dots , \chi_m)$. Clearly we would have $d(\varphi, r) = 1$, so it cannot be $m>1$. Similarly it cannot be $p>1$, so $m = 1 = p$ and $\varphi_1 = \psi = \chi_1$.\\

Now assume this set is not empty and let's name $r_1, \dots, r_h$ its members (in increasing order). Let's also define the following.\\

If $r_1 = 5$ then let $\psi_1 = \epsilon$ else $r_1 > 5$ and let $\psi_1 = \varphi[5, r_1 - 1]$.\\

If $h > 1$ then for each $i = 1 \dots h-1$: 
if $r_{i+1} = r_i + 1$ then let $\psi_{i+1} = \epsilon$ else $r_{i+1} > r_i + 1$ and let $\psi_{i+1} = \varphi[r_i + 1, r_{i+1} - 1]$.\\

Finally if $r_h = \ell(\varphi) - 1$ we define $\psi_{h+1} = \epsilon$, else $r_h < \ell(\varphi) - 1$ and let $\psi_{h+1} = \varphi[r_h + 1, \ell(\varphi) - 1]$ .\\

We have seen in lemma~\ref{IL:Hc-recursive-aux2} that we cannot have $m = 1$ and that since $m > 1$ we can indicate with $q_1, \dots, q_{m-1}$  the positions of the explicit occurrences of `,' in the representation $(c)( \chi_1, \dots , \chi_m)$ of $\varphi$.\\

For each $j = 1 \dots m-1$ $d(\varphi, q_j) = 1$, therefore $\{ q_1, \dots, q_{m-1} \} \subseteq \{ r_1, \dots, r_h \}$. In the mentioned lemma we have seen that actually $\{ q_1, \dots, q_{m-1} \} = \{ r_1, \dots, r_h \}$. This means that $h+1 = m$ and for each $i = 1 \dots h+1$ $\chi_i = \psi_i \in E(n,k)$.\\

Similarly we obtain that $h+1 = p$ and for each $i = 1 \dots h+1$ $\varphi_i = \psi_i \in E(n,k)$.\\

Therefore finally $p = h+1 = m$ and for each $i = 1 \dots m$ $\varphi_i = \psi_i = \chi_i$.
\end{proof}

\bigskip

Given a constant $c \in \mathcal{C}$ if $\#(c)$ is a particular type of function then for each $k \in K(n)$ we can define a set of expressions related to $c$ and $k$, and we'll call $E^c(n+1,k)$ this set of expressions.\\

Let's examine the categories of functions to which we refer.\\

Let $m$ be a positive integer, let $i: \{ 1, \dots , m \} \to \{ 1, \dots , p \}$, let $q$ be another function with the same domain $\{ 1, \dots , m \}$ such that for each $j = 1 \dots m$ $q(j) \geqslant 0$. For each $j = 1 \dots m$ we also define $M_j$ as follows: 
\begin{itemize}
\item if $q(j) = 0$ then $M_j = D_{i(j)}$,
\item if $q(j) > 0$ then $M_j = \mathcal{P}^{q(j)}(D_{i(j)})$
\end{itemize}

\medskip

Our costant $c$ can have as meaning three types of function, in every case the domain of $\#(c)$ is $M_1 \times \dots \times M_m$, but we have a different codomain in the three cases, which are the following:

\begin{itemize}
\item let $\alpha \in \{ 1, \dots , p \}$ and let $\#(c) : M_1 \times \dots \times M_m \to D_{\alpha}$;
\item let $\alpha \in \{ 1, \dots , p \}$, $r$ positive integer and let $\#(c) : M_1 \times \dots \times M_m \to \mathcal{P}^r(D_{\alpha})$;
\item let $\#(c)$ be a function over $M_1 \times \dots \times M_m$ such that for each $(d_1, \dots , d_m) \in M_1 \times \dots \times M_m$ $\#(c)(d_1, \dots , d_m)$ is true or false.
\end{itemize}

Let's call $\mathcal{C}'$ the set of all constants $c$ which have as meaning one of this three types of function.\\

In all of the three cases we define $E^c(n+1, k)$ as the set of the strings $(c)( \varphi_1, \dots , \varphi_m) \in H_c(n+1,k)$ such that:

\begin{itemize}

\item $\varphi_1, \dots , \varphi_m \in E(n,k)$;

\item for each $j = 1 \dots m$, $\sigma \in \Xi(k)$ $\#(k, \varphi_j, \sigma) \in M_j$;

\item $(c)( \varphi_1, \dots , \varphi_m) \notin E(n,k)$;

\item $(c)( \varphi_1, \dots , \varphi_m) \notin E_b(n+1,k)$.
\end{itemize}

\medskip

The set $E^c(n+1,k)$ is recursive since given $\psi \in \Sigma^*$ we can determine if $\psi \in H_c(n+1,k)$ and if so we can identify a positive integer $u$ and $\varphi_1, \dots , \varphi_u \in E(n,k)$ such that $\psi = (c)( \varphi_1, \dots , \varphi_u)$. As we have seen  $u$ and $\varphi_1, \dots , \varphi_u$ are univocally determined, so if $u \ne m$ then $\psi \notin E^c(n+1,k)$. If $u = m$ then, for each $j = 1 \dots m$, we can decide if for each $\sigma \in \Xi(k)$ $\#(k, \varphi_j, \sigma) \in M_j$, and we can also decide if the following conditions hold: 

\begin{itemize}
\item $(c)( \varphi_1, \dots , \varphi_m) \notin E(n,k)$,
\item $(c)( \varphi_1, \dots , \varphi_m) \notin E_b(n+1,k)$.
\end{itemize}

\medskip

For each $t = (c)(\varphi_1, \dots , \varphi_m) \in E^c(n+1,k)$ we define
\[
\#(k,t,\sigma)_{(n+1,k,<c>)} = \#(c)( \#(k, \varphi_1, \sigma), \dots , \#(k, \varphi_m, \sigma) ).
\]

\bigskip

Given $k \in K(n)$ and $f \in \mathcal{F}$ we can define the set $H_f(n+1, k)$ as follows. If $f$ has multiplicity $1$ then
\[ H_f(n+1, k) = \{ f( \varphi_1 ) | \, \varphi_1 \in E(n,k) \} \, . \]

If $f$ has multiplicity 2 then
\[ H_f(n+1, k) = \{ f( \varphi_1, \varphi_2) | \, \varphi_1, \varphi_2 \in E(n,k) \} \, . \]

We can prove $H_f(n+1, k)$ is recursive using some auxiliary lemma.\\

\begin{ilemma}\label{IL:Hf-recursive-aux1}
Let $f \in \mathcal{F}$ and assume $f$ has multiplicity 1. Let $\psi \in \Sigma^*$ and let $\varphi = f(\psi) \in \Sigma^*$. Then $\varphi \in H_f(n+1, k)$ if and only if $\psi \in E(n,k)$.
\end{ilemma}

\begin{proof}
It is obvious that if $\psi \in E(n,k)$ then $\varphi \in H_f(n+1, k)$.\\

Conversely, if $\varphi \in H_f(n+1, k)$ then there exists $\chi \in E(n,k)$ such that $\varphi = f( \chi )$.\\

Therefore $\psi = \chi \in E(n,k)$.
\end{proof}

\medskip

\begin{ilemma}\label{IL:Hf-recursive-aux2}
Let $f \in \mathcal{F}$ and assume $f$ has multiplicity 2. Let $\psi \in \Sigma^*$ and let $\varphi = f(\psi) \in \Sigma^*$.\\

Consider the set of the positive integers $r$ such that $2 < r < \ell(\varphi)$, $\varphi[r]$ = `,' and $d(\varphi, r) = 1$. If this set has just one member $r_1$ then we can define the following.\\

If $r_1 = 3$ then $\psi_1 = \epsilon$, else $\psi_1 = \varphi[3, r_1 - 1]$.\\

If $r_1 = \ell(\varphi) - 1$ then $\psi_2 = \epsilon$ else $\psi_2 = \varphi[r_1 + 1, \ell(\varphi) - 1]$.\\

With these definitions we have that $\varphi \in H_f(n+1, k)$ if and only if 
\begin{itemize}
\item the set of the positive integers $r$ such that $2 < r < \ell(\varphi)$, $\varphi[r]$ = `,' and $d(\varphi, r) = 1$ has just one member $r_1$,
\item $\psi_1, \psi_2 \in E(n,k)$.
\end{itemize}
\end{ilemma}

\begin{proof}

If the two conditions 
\begin{itemize}
\item the set of the positive integers $r$ such that $2 < r < \ell(\varphi)$, $\varphi[r]$ = `,' and $d(\varphi, r) = 1$ has just one member $r_1$,
\item $\psi_1, \psi_2 \in E(n,k)$.
\end{itemize}

both hold then clearly $\varphi = f(\psi_1, \psi_2) \in H_f(n+1, k)$.\\

Conversely if $\varphi \in H_f(n+1, k)$ then there exist $\chi_1, \chi_2 \in E(n, k)$ such that $\varphi = f(\chi_1, \chi_2)$. Let's call $q_1$ the position of the explicit occurrence of `,' in the representation $f(\chi_1, \chi_2)$ of $\varphi$.\\

Let $\vartheta = \varphi[1,2]$ be the substring of $\varphi$ which begins at character $1$ and ends at character $2$. We first want to prove that $d(\varphi, q_1) = 1$. We have that 
\[ d(\varphi,q_1 - 1) = d(\varphi, \ell(\vartheta) + \ell(\chi_1)) = d(\varphi, \ell(\vartheta) + 1) + d( \chi_1, \ell(\chi_1)) = 1 + d( \chi_1, \ell(\chi_1)) \ . \]

\medskip

If $\varphi[q_1 - 1] = \chi_1[\ell(\chi_1)] = \text{`)'}$ then $d(\varphi, q_1) = d(\varphi,q_1 - 1) - 1 = d(\chi_1, \ell(\chi_1)) = 1$ .\\
Else $\varphi[q_1 - 1] = \chi_1[\ell(\chi_1)] \notin \{ \text{`('}, \text{`)'} \}$ so $d(\varphi, q_1) = d(\varphi, q_1 - 1) = 1 + d(\chi_1, \ell(\chi_1)) = 1$.\\

In both cases $d(\varphi, q_1) = 1$ and $q_1$ is a member of the set of the positive integers $r$ such that $2 < r < \ell(\varphi)$, $\varphi[r]$ = `,' and $d(\varphi, r) = 1$.\\

Let's then call $r_1, \dots , r_h$ the members of the set of the positive integers $r$ such that $2 < r < \ell(\varphi)$, $\varphi[r]$ = `,' and $d(\varphi, r) = 1$. We have already seen that $q_1 \in \{ r_1, \dots , r_h \}$. Suppose $h > 1$ and there exists $i = 1 \dots h$ such that $r_i \ne q_1$. In this case one of the following conditions will occur:

\begin{itemize}
\item $r_i < q_1$,
\item $r_i > q_1$.
\end{itemize}

\medskip

If $r_i < q_1$ then $2 < r_i$ also holds, $\chi_1 = \varphi[3, q_1 - 1]$, $\ell(\chi_1) = q_1 - 1 - 2 =  q_1 - 3$, for each $\alpha = 1 \dots q_1 - 3$ $\chi_1[\alpha] = \varphi[2 + \alpha]$. So $r_i - 2 \geqslant 1$, $r_i - 2 < q_1 - 2$ and then $r_i - 2 \leqslant q_1 - 3 = \ell(\chi_1)$. Then also $\chi_1[r_i - 2] = \varphi[r_i]$ = `,'.  If we define $\vartheta = \varphi[1,2]$ and observe that $r_i = 2 + (r_i - 2) = \ell(\vartheta) + r_i - 2$ then we can also observe that $1 = d(\varphi, r_i) = d(\varphi, \ell(\vartheta) + 1) + d(\chi_1, r_i - 2) = 1 + d(\chi_1, r_i - 2)$, so $d(\chi_1, r_i - 2) = 0$. This contradicts assumption~\ref{A:expr-cont-depth} and therefore we cannot have  $r_i < q_1$.\\

If $r_i > q_1$ then $r_i < \ell(\varphi)$ also holds, $\chi_2 = \varphi[q_1 + 1, \ell(\varphi) - 1]$, $\ell(\chi_2) = \ell(\varphi) - 1 - q_1$. For each $\alpha = 1 \dots \ell(\varphi) - 1 - q_1$ $\chi_2[\alpha] = \varphi[q_1 + \alpha]$. So $r_i - q_1 \geqslant 1$, $r_i - q_1 < \ell(\varphi) - q_1$ and then $r_i - q_1 \leqslant \ell(\varphi) - 1 - q_1 = \ell(\chi_2)$. Then also $\chi_2[r_i - q_1] = \varphi[r_i]$ = `,'.\\

Moreover if we define $\vartheta = \varphi[1, q_1]$ then $\varphi$ is the concatenation of $\vartheta$, $\chi_2$ and $)$. Then $d(\varphi, r_i) = d(\varphi, q_1 + (r_i - q_1)) = d(\varphi, \ell(\vartheta) + (r_i - q_1)) = d(\varphi, \ell(\vartheta) + 1) + d(\chi_2, (r_i - q_1))$. It follows that $d(\varphi, r_i) = d(\varphi, q_1) + d(\chi_2, (r_i - q_1))$, and so $d(\chi_2, (r_i - q_1)) = 0$. This contradicts assumption~\ref{A:expr-cont-depth} and therefore we cannot have $r_i > q_1$.\\

So we have to conclude that $h = 1$, $r_1 = q_1$, $\psi_i = \chi_i \in E(n,k)$.
\end{proof}

\medskip

\begin{ilemma}\label{IL:Hf-recursive-mult2}
Let $f \in \mathcal{F}$ and assume $f$ has multiplicity 2. Then $H_f(n+1, k)$ is recursive. 
\end{ilemma}

\begin{proof}
Let $\varphi \in \Sigma^*$ and let's see how we decide if $\varphi \in H_f(n+1, k)$.\\ 

If $\varphi$ doesn't begin with the characters $f($ or doesn't end with the character $)$ then $\varphi \notin H_f(n+1, k)$.\\

Then assume we are in the case $\varphi = f(\psi)$ where $\psi \in \Sigma^*$. Consider the set of the positive integers $r$ such that $2 < r < \ell(\varphi)$, $\varphi[r]$ = `,' and $d(\varphi, r) = 1$.\\

If the mentioned set is empty or has not exactly one member then $\varphi \notin H_f(n+1, k)$.\\

If this set has just one member $r_1$ then we can define the following.\\

If $r_1 = 3$ then $\psi_1 = \epsilon$, else $\psi_1 = \varphi[3, r_1 - 1]$.\\

If $r_1 = \ell(\varphi) - 1$ then $\psi_2 = \epsilon$ else $\psi_2 = \varphi[r_1 + 1, \ell(\varphi) - 1]$.\\

If $\psi_1, \psi_2 \in E(n,k)$ then we can decide $\varphi \in H_f(n+1, k)$, otherwise 
$\varphi \notin H_f(n+1, k)$.
\end{proof}

\medskip

\begin{ilemma}\label{IL:Hf-aux3}
Let $f \in \mathcal{F}$ and assume $f$ has multiplicity 2. There exists an algorithm that given $\varphi \in \Sigma^*$
\begin{itemize}
\item determines if $\varphi \in H_f(n+1, k)$,
\item if $\varphi \in H_f(n+1, k)$ it also identifies $\psi_1, \psi_2 \in E(n,k)$ such that $\varphi = f( \psi_1, \psi_2)$.
\end{itemize}
\end{ilemma}

\begin{proof}
See the proof of lemma~\ref{IL:Hf-recursive-mult2}.
\end{proof}

\medskip

\begin{ilemma}\label{IL:Hf-aux4}
Let $f \in \mathcal{F}$ and assume $f$ has multiplicity 2. Given $\varphi \in H_f(n+1, k)$ there exist $\chi_1, \chi_2 \in E(n,k)$ such that $\varphi = f( \chi_1, \chi_2)$ and $\chi_1, \chi_2$ are univocally determined.
\end{ilemma}

\begin{proof}

It is obvious by the definition of $H_f(n+1, k)$ that there exist $\chi_1, \chi_2 \in E(n,k)$ such that $\varphi = f( \chi_1, \chi_2)$. We have also seen in lemma~\ref{IL:Hf-recursive-aux2} that the set of the positive integers $r$ such that $2 < r < \ell(\varphi)$, $\varphi[r]$ = `,' and $d(\varphi, r) = 1$ has just one member $r_1$.\\

By the same lemma we can define the following.\\

If $r_1 = 3$ then $\psi_1 = \epsilon$, else $\psi_1 = \varphi[3, r_1 - 1]$.\\

If $r_1 = \ell(\varphi) - 1$ then $\psi_2 = \epsilon$ else $\psi_2 = \varphi[r_1 + 1, \ell(\varphi) - 1]$.\\

And we can see in the lemma that $\psi_1 = \chi_1, \psi_2 = \chi_2$.\\

We can assume there also exist $\phi_1, \phi_2 \in E(n,k)$ such that $\varphi = f( \phi_1, \phi_2)$. Clearly we can apply lemma~\ref{IL:Hf-recursive-aux2} also in this case and obtain $\psi_1 = \phi_1, \, \psi_2 = \phi_2$.\\

It obviously follows that $\phi_1 = \chi_1, \, \phi_2 = \chi_2$.

\end{proof}

\medskip

\begin{ilemma}\label{IL:Hf-recursive-mult1}
Let $f \in \mathcal{F}$ and assume $f$ has multiplicity 1. Then $H_f(n+1, k)$ is recursive. 
\end{ilemma}

\begin{proof}
Let $\varphi \in \Sigma^*$ and let's see how we decide if $\varphi \in H_f(n+1, k)$.\\ 

If $\varphi$ doesn't begin with the characters $f($ or doesn't end with the character $)$ then $\varphi \notin H_f(n+1, k)$.\\

Then assume we are in the case $\varphi = f(\psi)$ where $\psi \in \Sigma^*$.\\

In this case using lemma~\ref{IL:Hf-recursive-aux1} if $\psi \in E(n,k)$ we'll decide that $\varphi \in H_f(n+1, k)$, otherwise we'll decide that $\varphi \notin H_f(n+1, k)$.
\end{proof}

\medskip

\begin{ilemma}
Let $f \in \mathcal{F}$ and assume $f$ has multiplicity 1. There exists an algorithm that given $\varphi \in \Sigma^*$
\begin{itemize}
\item determines if $\varphi \in H_f(n+1, k)$,
\item if $\varphi \in H_f(n+1, k)$ it also identifies $\psi \in E(n,k)$ such that $\varphi = f( \psi )$.
\end{itemize}
\end{ilemma}

\begin{proof}
See the proof of lemma~\ref{IL:Hf-recursive-mult1}.
\end{proof}

\medskip

\begin{ilemma}
Let $f \in \mathcal{F}$ and assume $f$ has multiplicity 1. Given $\varphi \in H_f(n+1, k)$ there exists $\chi \in E(n,k)$ such that $\varphi = f( \chi )$ and $\chi$ is univocally determined.
\end{ilemma}

\begin{proof}

It is obvious by the definition of $H_f(n+1, k)$ that there exists $\chi \in E(n,k)$ such that $\varphi = f( \chi )$.\\

We can also assume there exists $\phi \in E(n,k)$ such that $\varphi = f( \phi )$, then obviously $\phi = \chi$.
\end{proof}

\medskip

For each $k \in K(n)$ and $f \in \mathcal{F}$ if $f$ has multiplicity 2 we define $E^f(n+1,k)$ as the set of the strings $f( \varphi_1, \varphi_2) \in H_f(n+1,k)$ such that:

\begin{itemize}
\item $\varphi_1, \varphi_2 \in E(n,k)$;

\item for each $\sigma \in \Xi(k)$ $A_f( \#(k,\varphi_1, \sigma), \#(k,\varphi_2, \sigma) )$ is true;

\item $f( \varphi_1, \varphi_2) \notin E(n,k)$;

\item $f( \varphi_1, \varphi_2) \notin E_b(n+1,k)$.
\end{itemize}

\medskip

For instance, this means that if $f$ is the `logical conjunction' symbol `$\wedge$' and it belongs to $\mathcal{F}$, $\varphi_1$, $\varphi_2$ belong to $E(n,k)$, for each $\sigma \in \Xi(k)$ both $\#(k,\varphi_1, \sigma)$ and $\#(k,\varphi_2, \sigma)$ are true or false, $\wedge(\varphi_1, \varphi_2) \notin E(n,k)$, $\wedge(\varphi_1, \varphi_2) \notin E_b(n+1,k)$ then $\wedge(\varphi_1, \varphi_2)$ belongs to $E^f(n+1,k)$.\\

We now show that $E^f(n+1,k)$ is recursive. Given $\varphi \in \Sigma^*$ we can determine if $\varphi \in H_f(n+1, k)$. Clearly if $\varphi \notin H_f(n+1, k)$ then $\varphi \notin E^f(n+1, k)$. If $\varphi \in H_f(n+1, k)$ then we can identify $\psi_1, \psi_2 \in E(n,k)$ such that $\varphi = f( \psi_1, \psi_2)$. We have seen that $\psi_1, \psi_2$ are univocally determined.\\

For $f$ with multiplicity 2 $A_f( \#(k,\varphi_1, \sigma), \#(k,\varphi_2, \sigma) )$ can be one of the following

\begin{itemize}
\item $Event_1(\#(k,\varphi_1, \sigma))$ and $Event_1(\#(k,\varphi_2, \sigma))$,
\item $Set_1(\#(k,\varphi_2, \sigma))$,
\item `something which is true' (e.g. $1 = 1$) .
\end{itemize}

In every mentioned case the condition `for each $\sigma \in \Xi(k)$ $A_f( \#(k,\varphi_1, \sigma), \#(k,\varphi_2, \sigma) )$' is decidable, and we can also decide if the following conditions hold

\begin{itemize}
\item $f( \varphi_1, \varphi_2) \notin E(n,k)$,
\item $f( \varphi_1, \varphi_2) \notin E_b(n+1,k)$.
\end{itemize}

\bigskip

For each $f$ with multiplicity 2, $t = f(\varphi_1, \varphi_2) \in E^f(n+1,k)$ we define
\[
\#(k,t,\sigma)_{(n+1,k,<f>)} = P_f( \#(k, \varphi_1, \sigma), \#(k, \varphi_2, \sigma) ).
\]

\bigskip

If f has multiplicity 1 we define $E^f(n+1,k)$ as the set of the strings $f( \varphi_1 ) \in H_f(n+1,k)$ such that:

\begin{itemize}
\item $\varphi_1 \in E(n,k)$;

\item for each $\sigma \in \Xi(k)$ $A_f( \#(k,\varphi_1, \sigma) )$ is true;

\item $f( \varphi_1 ) \notin E(n,k)$.

\item $f( \varphi_1 ) \notin E_b(n+1,k)$.
\end{itemize}

\medskip

We now show that $E^f(n+1,k)$ is recursive. Given $\varphi \in \Sigma^*$ we can determine if $\varphi \in H_f(n+1, k)$. Clearly if $\varphi \notin H_f(n+1, k)$ then $\varphi \notin E^f(n+1, k)$. If $\varphi \in H_f(n+1, k)$ then we can identify $\psi \in E(n,k)$ such that $\varphi = f( \psi )$. We have seen that $\psi$ is univocally determined.\\

For $f$ with multiplicity 1 $A_f( \#(k,\varphi_1, \sigma) )$ is the following:

\begin{itemize}
\item `something which is true' (e.g. $1 = 1$).
\end{itemize}

Therefore the condition `for each $\sigma \in \Xi(k)$ $A_f( \#(k,\varphi_1, \sigma)$' is trivially decidable, and we can also decide if the following conditions hold

\begin{itemize}
\item $f( \varphi_1) \notin E(n,k)$,
\item $f( \varphi_1) \notin E_b(n+1,k)$.
\end{itemize}

\bigskip

For each $f$ with multiplicity 1, $t = f(\varphi_1) \in E^f(n+1,k)$ we define
\[
\#(k,t,\sigma)_{(n+1,k,<f>)} = P_f( \#(k, \varphi_1, \sigma) ).
\]

\bigskip

In our language we have the quantifier symbols `$\forall$' and `$\exists$'. Let's now see how we can use them.\\

Given $Q \in \{ \forall, \exists \}$ and $k \in K(n)$ we can define the set $H_Q(n+1, k)$ as the set of the strings $Q(x:\varphi,\phi)$ such that 

\begin{itemize}
\item $\varphi \in E_s(n,k)$,
\item $x \in \mathcal{V} - var(k)$,
\item if we define $k' = k + <x, \varphi>$ then $k' \in K(n)$ and $\phi \in S(n,k')$ .
\end{itemize}

\medskip

We can prove $H_Q(n+1, k)$ is recursive using some auxiliary lemma.\\

\begin{ilemma}\label{IL:HQ-recursive-aux1}
Let $\psi \in \Sigma^*$, $x \in \mathcal{V}$ and let $\chi = Q(x:\psi) \in \Sigma^*$.\\

Consider the set of the positive integers $r$ such that $4 < r < \ell(\chi)$, $\chi[r]$ = `,' and $d(\chi, r) = 1$. If this set has just one member $r_1$ then we can define the following.\\

If $r_1 = 5$ then $\varphi = \epsilon$, else $\varphi = \chi[5, r_1 - 1]$.\\

If $r_1 = \ell(\chi) - 1$ then $\phi = \epsilon$ else $\phi = \chi[r_1 + 1, \ell(\chi) - 1]$.\\

With these definitions we have that $\chi \in H_Q(n+1, k)$ if and only if 
\begin{itemize}
\item $x \in \mathcal{V} - var(k)$,
\item the set of the positive integers $r$ such that $4 < r < \ell(\chi)$, $\chi[r]$ = `,' and $d(\chi, r) = 1$ has just one member $r_1$,
\item $\varphi \in E_s(n,k)$.
\item if we define $k' = k + <x, \varphi>$ then $k' \in K(n)$ and $\phi \in S(n,k')$ .
\end{itemize}
\end{ilemma}

\begin{proof}

Assume the following conditions hold:
\begin{itemize}
\item $x \in \mathcal{V} - var(k)$,
\item the set of the positive integers $r$ such that $4 < r < \ell(\chi)$, $\chi[r]$ = `,' and $d(\chi, r) = 1$ has just one member $r_1$,
\item $\varphi \in E_s(n,k)$.
\item if we define $k' = k + <x, \varphi>$ then $k' \in K(n)$ and $\phi \in S(n,k')$ .
\end{itemize}

Then $\chi = Q(x:\varphi,\phi) \in H_Q(n+1, k)$.\\

Conversely if $\chi \in H_Q(n+1, k)$ then there exist $\vartheta \in E_s(n,k)$, $y \in \mathcal{V} - var(k)$ such that if we define $\kappa = k + <y, \vartheta>$ then $\kappa \in K(n)$ and there also exists $\theta \in S(n,\kappa)$ such that $\chi = Q(y:\vartheta, \theta)$.\\

Let's call $q_1$ the position of the explicit occurrence of `,' in the representation $Q(y:\vartheta, \theta)$ of $\chi$. We first want to prove that $d(\chi, q_1) = 1$.\\

Let $\eta = \chi[1,4]$ be the substring of $\chi$ which begins at character $1$ and ends at character $4$. We have that 
\[ d(\chi,q_1 - 1) = d(\chi, \ell(\eta) + \ell(\vartheta)) = d(\chi, \ell(\eta) + 1) + d( \vartheta, \ell(\vartheta)) = 1 + d( \vartheta, \ell(\vartheta)) \ . \]

\medskip

If $\chi[q_1 - 1] = \vartheta[\ell(\vartheta)] = \text{`)'}$ then $d(\chi, q_1) = d(\chi,q_1 - 1) - 1 = d(\vartheta, \ell(\vartheta)) = 1$ .\\
Else $\chi[q_1 - 1] = \vartheta[\ell(\vartheta)] \notin \{ \text{`('}, \text{`)'} \}$ so $d(\chi, q_1) = d(\chi, q_1 - 1) = 1 + d(\vartheta, \ell(\vartheta)) = 1$.\\

In both cases $d(\chi, q_1) = 1$ and $q_1$ is a member of the set of the positive integers $r$ such that $4 < r < \ell(\chi)$, $\chi[r]$ = `,' and $d(\chi, r) = 1$.\\

Let's then call $r_1, \dots , r_h$ the members of the set of the positive integers $r$ such that $4 < r < \ell(\chi)$, $\chi[r]$ = `,' and $d(\chi, r) = 1$. We have already seen that $q_1 \in \{ r_1, \dots , r_h \}$. Suppose $h > 1$ and there exists $i = 1 \dots h$ such that $r_i \ne q_1$. In this case one of the following conditions will occur:

\begin{itemize}
\item $r_i < q_1$,
\item $r_i > q_1$.
\end{itemize}

If $r_i < q_1$ then $4 < r_i$ also holds, $\vartheta = \chi[5, q_1 - 1]$, $\ell(\theta) = q_1 - 1 - 4 =  q_1 - 5$, for each $\alpha = 1 \dots q_1 - 5$ $\vartheta[\alpha] = \chi[4 + \alpha]$. So $r_i - 4 \geqslant 1$, $r_i - 4 < q_1 - 4$ and then $r_i - 4 \leqslant q_1 - 5 = \ell(\vartheta)$. Then also $\vartheta[r_i - 4] = \chi[r_i]$ = `,'.  If we define $\eta = \chi[1,4]$ and observe that $r_i = 4 + (r_i - 4) = \ell(\eta) + r_i - 4$ then we can also observe that $1 = d(\chi, r_i) = d(\chi, \ell(\eta) + 1) + d(\vartheta, r_i - 4) = 1 + d(\vartheta, r_i - 4)$, so $d(\vartheta, r_i - 4) = 0$. This contradicts assumption~\ref{A:expr-cont-depth} and therefore we cannot have  $r_i < q_1$.\\

If $r_i > q_1$ then $r_i < \ell(\chi)$ also holds, $\theta = \chi[q_1 + 1, \ell(\chi) - 1]$, $\ell(\theta) = \ell(\chi) - 1 - q_1$. For each $\alpha = 1 \dots \ell(\chi) - 1 - q_1$ $\theta[\alpha] = \chi[q_1 + \alpha]$. So $r_i - q_1 \geqslant 1$, $r_i - q_1 < \ell(\chi) - q_1$ and then $r_i - q_1 \leqslant \ell(\chi) - 1 - q_1 = \ell(\theta)$. Then also $\theta[r_i - q_1] = \chi[r_i]$ = `,'.\\

Moreover if we define $\eta = \chi[1, q_1]$ then $\chi$ is the concatenation of $\eta$, $\theta$ and $)$. Then $d(\chi, r_i) = d(\chi, q_1 + (r_i - q_1)) = d(\chi, \ell(\eta) + (r_i - q_1)) = d(\chi, \ell(\eta) + 1) + d(\theta, (r_i - q_1))$. It follows that $1 = d(\chi, r_i) = d(\chi, q_1) + d(\theta, (r_i - q_1)) = 1 + d(\theta, (r_i - q_1))$, and so $d(\theta, (r_i - q_1)) = 0$. This contradicts assumption~\ref{A:expr-cont-depth} and therefore we cannot have $r_i > q_1$.\\

So we have to conclude that $h = 1$, $r_1 = q_1$, $\varphi = \vartheta \in E_s(n,k)$, $x = y \in \mathcal{V} - var(k)$, $k' = k + <x, \varphi> = k + <y, \vartheta> = \kappa \in K(n)$, $\phi = \theta \in S(n, k')$.
\end{proof}

\medskip

\begin{ilemma}\label{IL:HQ-recursive}
$H_Q(n+1, k)$ is recursive. 
\end{ilemma}

\begin{proof}
Let $\chi \in \Sigma^*$ and let's see how we decide if $\chi \in H_Q(n+1, k)$.\\ 

If $\chi$ has not at least $5$ characters then $\chi \notin H_Q(n+1, k)$.\\

If $\chi[1] \ne Q$ then $\chi \notin H_Q(n+1, k)$.\\

If $\chi[2] \ne \ $`(' then $\chi \notin H_Q(n+1, k)$.\\

If $\chi[3] \notin \mathcal{V}$ (and $\Sigma - \mathcal{V}$ is a finite set) then $\chi \notin H_Q(n+1, k)$.\\

If $\chi[4] \ne \ $`:' then $\chi \notin H_Q(n+1, k)$.\\

Let's now assume we are in the case where $\chi = Q(x:\psi)$ with $x \in \mathcal{V}$ and $\psi \in \Sigma^*$.\\

If $x \in var(k)$ then $\chi \notin H_Q(n+1, k)$.\\

Consider the set of the positive integers $r$ such that $4 < r < \ell(\chi)$, $\chi[r]$ = `,' and $d(\chi, r) = 1$.\\

If this set is empty or has more than one member then $\chi \notin H_Q(n+1, k)$.\\

If this set has just one member $r_1$ then we can define the following.\\

If $r_1 = 5$ then $\varphi = \epsilon$, else $\varphi = \chi[5, r_1 - 1]$.\\

If $r_1 = \ell(\chi) - 1$ then $\phi = \epsilon$ else $\phi = \chi[r_1 + 1, \ell(\chi) - 1]$.\\

if $\varphi \notin E_s(n,k)$ then $\chi \notin H_Q(n+1, k)$.\\

Let now $k' = k + <x, \varphi>$, if $k' \notin K(n)$ then $\chi \notin H_Q(n+1, k)$.\\

If $\phi \notin S(n,k')$ then $\chi \notin H_Q(n+1, k)$.\\

At this point, if we have not decided that $\chi \notin H_Q(n+1, k)$, then $\chi \in H_Q(n+1, k)$.\\
\end{proof}

\bigskip

\begin{ilemma}\label{IL:HQ-aux2}
There exists an algorithm that given $\chi \in \Sigma^*$
\begin{itemize}
\item determines if $\chi \in H_Q(n+1, k)$,
\item if $\chi \in H_Q(n+1, k)$ 
\begin{itemize}
\item it also identifies $x \in \mathcal{V} - var(k)$, $\varphi \in E_s(n,k)$ such that if we define  $k' = k + <x, \varphi>$ then $k' \in K(n)$; 
\item it also identifies $\phi \in S(n,k')$ such that $\chi = Q(x:\varphi,\phi)$.
\end{itemize}
\end{itemize}
\end{ilemma}

\begin{proof}
See the proof of lemma~\ref{IL:HQ-recursive}.
\end{proof}

\bigskip

\begin{ilemma}\label{IL:HQ-aux3}
Given $\chi \in H_Q(n+1, k)$ there exist $x \in \mathcal{V} - var(k)$, $\varphi \in E_s(n,k)$ such that if we define  $h = k + <x, \varphi>$ then $h \in K(n)$ and there exists $\phi \in S(n,h)$ such that $\chi = Q(x:\varphi,\phi)$. Moreover $x$, $\varphi$ and $\phi$ are univocally determined. 
\end{ilemma}

\begin{proof}
It is obvious by the definition of $H_Q(n+1, k)$ that there exist $\vartheta \in E_s(n,k)$, $x \in \mathcal{V} - var(k)$ such that if we define $h = k + <x, \vartheta>$ then $h \in K(n)$ and there also exists $\theta \in S(n,h)$ such that $\chi = Q(x:\vartheta, \theta)$.\\

We have seen in lemma \ref{IL:HQ-recursive-aux1} that we have the following 
\begin{itemize}
\item $x \in \mathcal{V} - var(k)$,
\item the set of the positive integers $r$ such that $4 < r < \ell(\chi)$, $\chi[r]$ = `,' and $d(\chi, r) = 1$ has just one member $r_1$,
\end{itemize}

\smallskip

Moreover let's define he following:\\

If $r_1 = 5$ then $\varphi = \epsilon$, else $\varphi = \chi[5, r_1 - 1]$.\\

If $r_1 = \ell(\chi) - 1$ then $\phi = \epsilon$ else $\phi = \chi[r_1 + 1, \ell(\chi) - 1]$.\\

With these definitions we have also 

\begin{itemize}
\item $\varphi \in E_s(n,k)$.
\item if we define $k' = k + <x, \varphi>$ then $k' \in K(n)$ and $\phi \in S(n,k')$ .
\end{itemize}

\smallskip

We have also seen in the mentioned lemma that $\vartheta = \varphi$ and $\theta = \phi$.\\

We now suppose there also exist $\zeta \in E_s(n,k)$, $z \in \mathcal{V} - var(k)$ such that if we define $g = k + <x, \zeta>$ then $g \in K(n)$ and there also exists $\eta \in S(n,g)$ such that $\chi = Q(z:\zeta, \eta)$.\\

Cleary it must be $z = x$ and by lemma \ref{IL:HQ-recursive-aux1} we have the same facts that we have listed before, and finally we have $\zeta = \varphi = \vartheta$ and $\eta = \phi = \theta$.
\end{proof}

\bigskip

Given $Q \in \{ \forall, \exists \}$ and $k \in K(n)$ we can define the set $E_Q(n+1, k)$ as the set of the strings $\chi = Q(x:\varphi,\phi) \in H_Q(n+1, k)$ such that 

\begin{itemize}
\item $\chi \notin E(n,k)$;
\item $\chi \notin E_b(n+1,k)$.
\end{itemize}

\bigskip

Clearly $E_Q(n+1, k)$ is recursive, in fact given $\chi \in \Sigma^*$ we can determine if $\chi \in H_Q(n+1, k)$ or not, if $\chi \notin H_Q(n+1, k)$ then $\chi \notin E_Q(n+1, k)$. If $\chi \in H_Q(n+1, k)$ then we can also determine if $\chi \in E(n,k)$ or not and if $\chi \in E_b(n+1,k)$ or not.

\bigskip

Given $\chi = Q(x:\varphi,\phi) \in E_Q(n+1, k)$ and $\sigma \in \Xi(k)$ we define $\#(k, \chi, \sigma)_{(n+1,k,Q)}$ as follows. First of all let $k' = k + <x, \varphi>$, then: if $Q = \forall$:
\[ \#(k, \chi, \sigma)_{(n+1,k,Q)} = \text{for each} \ \sigma' \in \Xi(k'): \sigma \sqsubseteq \sigma' \ \#(k', \phi, \sigma') \ . \]

Else if $Q = \exists$:
\[ \#(k, \chi, \sigma)_{(n+1,k,Q)} = \text{exists} \ \sigma' \in \Xi(k'): \sigma \sqsubseteq \sigma' \ \#(k', \phi, \sigma') \ . \]

\medskip

Here it is important to notice that in the expression 
\[\text{for each} \ \sigma' \in \Xi(k'): \sigma \sqsubseteq \sigma' \ \#(k', \phi, \sigma') \]
we are just applying a predicate to a set. The set actually is the set of $\#(k', \phi, \sigma')$ such that $\sigma' \in \Xi(k')$, $\sigma \sqsubseteq \sigma'$. The predicate just applies a logical conjunction to the elements of the set. If we provisionally represent this predicate with $U$ and our set with $A$ then our expression could be written $U(A)$. Given another set $B$ if we can apply $U$ to $B$ then we can consider the expression $U(B)$. Now it is important to notice that if $A = B$ then obviously $U(A) = U(B)$.\\

The same observation can be made for the expression
\[ \text{exists} \ \sigma' \in \Xi(k'): \sigma \sqsubseteq \sigma' \ \#(k', \phi, \sigma') \ . \]
Here we are also applying a predicate to a set, we can name the predicate $E$ and the set is the same we have already considered, so let's provisionally name it $A$. So our expression in this case can be written $E(A)$ and if $B$ is another set and we can apply $E$ to $B$ then we can consider the expression $E(B)$. Now it is important to notice that if $A = B$ then obviously $E(A) = E(B)$.

\bigskip

We have terminated the definition of the `new sets' (of expressions bound to context $k$) and the related work, we are now ready to define $E(n+1,k)$ for $k \in K(n+1)$.\\

We recall we defined $\mathcal{C}'$ as the set of the constants $c \in \mathcal{C}$ for which, given $k \in K(n)$, we can define $E^c(n+1, k)$.\\

If $k \in K(n)^+$ we have defined $E_a(n+1,k)$, we also define 
\begin{itemize}
\item $E(n+1,k) = E_a(n+1,k)$.
\end{itemize}

\medskip

If $k \in K(n)$ we have defined $E_b(n+1,k)$, $E^c(n+1,k)$ (for each $c \in \mathcal{C}'$), $E^f(n+1,k)$ (for each $f \in \mathcal{F}$), $E_\forall(n+1,k)$, $E_\exists(n+1,k)$, and we can also define a set of tuples $\mathcal{H}(n+1,k)$ where each of the tuples represents one of the sets we have defined, and there is also a tuple which stands for $E(n,k)$. The formal definition is the following
\small
\[ \mathcal{H}(n+1,k) = \{ (n,k), (n+1,k, b), (n+1,k, \forall), (n+1,k, \exists) \} \cup \{ (n+1,k, c) | c \in \mathcal{C}' \} \cup \{ (n+1,k,f) | f \in \mathcal{F} \}  . \]
\normalsize

We also define a function $\mathcal{E}$ over $\mathcal{H}(n+1,k)$ as follows
\begin{itemize}
\item $\mathcal{E}(n,k) = E(n,k)$,
\item $\mathcal{E}(n+1,k,b) = E_b(n+1,k)$,
\item $\mathcal{E}(n+1,k,Q) = E_Q(n+1,k)$,
\item for each $c \in \mathcal{C}'$ $\mathcal{E}(n+1,k,c) = E^c(n+1,k)$,
\item for each $f \in \mathcal{F}$ $\mathcal{E}(n+1,k,f) = E^f(n+1,k)$.
\end{itemize}

\medskip

We can now define
\[ E(n+1,k) = \bigcup_{u \in \mathcal{H}(n+1,k)} \mathcal{E}(u) \ . \]

\begin{ilemma}
If $k \in K(n)$, given $u,v \in \mathcal{H}(n+1,k)$ such that $u \ne v$, $\mathcal{E}(u) \cap \mathcal{E}(v) = \emptyset$.
\end{ilemma}

\begin{proof}
It is obvious by definition that $E_b(n+1,k) \cap E(n,k) = \emptyset$.\\

It is also obvious that $E_Q(n+1,k) \cap E(n,k) = \emptyset$ and $E_Q(n+1,k) \cap E_b(n+1,k) = \emptyset$.\\

It is also obvious that $E_\exists(n+1,k) \cap E_\forall(n+1,k) = \emptyset$.\\

Given $c \in \mathcal{C}'$
\begin{itemize}
\item $E^c(n+1,k) \cap E(n,k) = \emptyset$,
\item $E^c(n+1,k) \cap E_b(n+1,k) = \emptyset$,
\item $E^c(n+1,k) \cap E_Q(n+1,k) = \emptyset$.
\end{itemize}

Given $c_1, c_2 \in \mathcal{C}'$, with $c_2 \ne c_1$, $E^{c_1}(n+1,k) \cap E^{c_2}(n+1,k) = \emptyset$.\\

Given $f \in \mathcal{F}$
\begin{itemize}
\item $E^f(n+1,k) \cap E(n,k) = \emptyset$,
\item $E^f(n+1,k) \cap E_b(n+1,k) = \emptyset$,
\item $E^f(n+1,k) \cap E_Q(n+1,k) = \emptyset$.
\end{itemize}

Given $f \in \mathcal{F}$, $c \in \mathcal{C}'$ $E^f(n+1,k) \cap E^c(n+1,k) = \emptyset$.\\

Given $f_1, f_2 \in \mathcal{F}$, with $f_2 \ne f_1$, $E^{f_1}(n+1,k) \cap E^{f_2}(n+1,k) = \emptyset$.\\
\end{proof}

\bigskip

For every $k \in K(n+1)$, $t \in E(n+1,k)$ and $\sigma \in \Xi(k)$ we need that $\#(k,t,\sigma)$ is defined.\\

If  $k \in K(n)^+$ we just need to define $\#(k,t,\sigma)$ for each $t \in E_a(n+1,k)$. Obviously we define $\#(k,t,\sigma) = \#(k,t,\sigma)_{(n+1,k,a)}$.\\

If $k \in K(n)$, how do we define $\#(k,t,\sigma)$ for each $t \in E(n+1,k)$?\\

Given $t \in E(n,k)$ $\#(k,t,\sigma)$ is already defined and we don't need to redefine it.\\

Given $t \in E_b(n+1,k)$ we define $\#(k,t,\sigma) = \#(k,t,\sigma)_{(n+1,k,b)}$.\\

Given $t \in E_Q(n+1,k)$ we define $\#(k,t,\sigma) = \#(k,t,\sigma)_{(n+1,k,Q)}$.\\

Given $c \in \mathcal{C}'$, $t \in E^c(n+1,k)$ we define $\#(k,t,\sigma) = \#(k,t,\sigma)_{(n+1,k,<c>)}$.\\

Given $f \in \mathcal{F}$, $t \in E^f(n+1,k)$ we define $\#(k,t,\sigma) = \#(k,t,\sigma)_{(n+1,k,<f>)}$.\\

Notice that if $k \in K(n)^+$ we have not defined $E_b(n+1,k), \ E_Q(n+1,k)$\, given $c \in \mathcal{C}'$ we have not defined  $E^c(n+1,k)$ and given $f \in \mathcal{F}$ we have not defined $E^f(n+1,k)$. We can conventionally define all of these sets as the empty set.\\

Also notice that if $k \in K(n)$ we have not defined $E_a(n+1,k)$ and we can conventionally define it as the empty set.\\

In the last part of our definition we need to prove that all the assumptions we have made at step $n$ are true at step $n+1$.\\

\begin{proof}[Proof of ~\ref{A:k-m-in-k-n-new}]

\smallskip

Let $m < n+1$. If $m = n$ then clearly $K(m) = K(n) \subseteq K(n+1)$. Else $m < n$ so $K(m) \subseteq K(n) \subseteq K(n+1)$. 

\end{proof}

\bigskip

\begin{proof}[Proof of~\ref{A:Ekn-sub-Closure} and~\ref{A:Ekn-recursive} ]

Let $k \in K(n+1)$, if $k \in K(n)^+$ then $E(n+1,k) = E_a(n+1,k) \subseteq \Sigma^*$.\\

If $k \in K(n)$ then $E(n+1,k) = \bigcup_{u \in \mathcal{H}(n+1,k)} \mathcal{E}(u)$. Then in order to prove that $E(n+1,k) \subseteq \Sigma^*$ we just need to prove that for each $u \in \mathcal{H}(n+1,k)$ $\mathcal{E}(u) \subseteq \Sigma^*$.\\

We actually have the following:\\

\begin{itemize}
\item $E(n,k) \subseteq \Sigma^*$,
\item $E_b(n+1,k) \subseteq \Sigma^*$,
\item $E_Q(n+1,k) \subseteq \Sigma^*$,
\item for each $c \in \mathcal{C}'$ $E^c(n+1,k) \subseteq \Sigma^*$,
\item for each $f \in \mathcal{F}$ $E^f(n+1,k) \subseteq \Sigma^*$.
\end{itemize}

\medskip

Let's now see how we prove that $E(n+1,k)$ is recursive.\\

Let $t \in \Sigma^*$, we have to decide if $t \in E(n+1,k)$. First we can decide if $k \in K(n)$, if this is false then we just need to decide if $t \in E_a(n+1,k)$.\\

If instead $k \in K(n)$ holds true, we check the following conditions

\begin{itemize}
\item $t \in E(n,k)$,
\item $t \in E_b(n+1,k)$,
\item $t \in E_\forall(n+1,k)$,
\item $t \in E_\exists(n+1,k)$,
\item the condition $t \in E^c(n+1,k)$ (for each $c \in \mathcal{C}'$),
\item the condition $t \in E^f(n+1,k)$ (for each $f \in \mathcal{F}$).
\end{itemize}

If at least one of the conditions is true then we can decide $t \in E(n+1,k)$, otherwise $t \notin E(n+1,k)$.

\end{proof}

\bigskip

\begin{proof}[Proof of~\ref{A:simple-k-sigma-slp-new}]

\smallskip

We have to show that for each $k \in K(n+1)$ $k \in \Theta$ and for each $\sigma \in \varXi(k)$ $\sigma$ is a state-like pair and $dom(\sigma) = dom(k)$.\\

If $k \in K(n)$ this is clearly true because it is precisely our assumption.\\

If $k \in K(n)^+$ then there exist $h \in K(n), \phi \in E_s(n,h), z \in (\mathcal{V}-var(h))$ such that $k = h + <z,\phi>$ and $\Xi(k) = \{ \rho + (z,s) | \, \rho \in \Xi(h), s \in \#(h,\phi,\rho) \}$.\\

For each $\sigma \in \Xi(k)$ $\sigma = \rho + (z,s)$ with $\rho \in \Xi(h), s \in \#(h,\phi,\rho)$, so $\sigma$ is a state-like pair.\\

Moreover, we can assume $dom(h) = dom(\rho) = \emptyset$ or $dom(h) = dom(\rho) = \{ 1, \dots , m \}$ for a positive  integer $m$. In the first case $dom(\sigma) = \{ 1 \} = dom(k)$, else
\[ dom(\sigma) = dom(\rho) \cup \{m+1\} = dom(h) \cup \{m+1\} = dom(k) \ . \]
\end{proof}

\bigskip

\begin{proof}[Proof of~\ref{A:k-is-epsilon-or-new}]

\smallskip

We have to show that for each $k \in K(n+1)$ $k = \epsilon$ and $\varXi(k) = \{\epsilon\}$ or\\
(there exist $m<n+1$, $h \in K(m)$, $\phi \in E_s(m,h)$, $y \in (\mathcal{V}-var(h))$ such that $k = h + <y, \phi>$, $\Xi(k) = \{ \sigma + (y,s) | \, \sigma \in \Xi(h), s \in \#(h,\phi,\sigma) \}$ ).

\medskip

If $k \in K(n)$ by the inductive hypothesis $k = \epsilon$ and $\varXi(k) = \{\epsilon\}$ or\\
($n > 1$ and there exist $m<n<n+1$, $h \in K(m)$, $\phi \in E_s(m,h)$, $y \in (\mathcal{V}-var(h))$ such that $k = h + <y, \phi>$, $\Xi(k) = \{ \sigma + (y,s) | \, \sigma \in \Xi(h), s \in \#(h,\phi,\sigma) \} )$.

\medskip

Otherwise $k \in K(n)^+$ so there exist $h \in K(n), \phi \in E_s(n,h), y \in (\mathcal{V}-var(h))$ such that $k = h + <y,\phi>$, $\Xi(k) = \{ \sigma + (y,s) | \, \sigma \in \Xi(h), s \in \#(h,\phi,\sigma) \}$.
\end{proof}

\bigskip

\begin{proof}[Proof of~\ref{A:rest-of-state-is-state}]
We have to show that for each $k \in K(n+1): k \ne \epsilon$, $\sigma \in \Xi(k)$, $h \in \mathcal{R}(k): h \ne k$, there exists $m < n+1$ such that $h \in K(m)$ and it results $\sigma_{/dom(h)} \in \Xi(h)$.\\

We first consider the case where $n+1 = 2$. Here we have to show that for each $k \in K(2): k \ne \epsilon$, $\sigma \in \Xi(k)$, $h \in \mathcal{R}(k): h \ne k$, $h \in K(1)$ and it results $\sigma_{/dom(h)} \in \Xi(h)$.\\

Let $k \in K(2): k \ne \epsilon$, $\sigma \in \Xi(k)$, $h \in \mathcal{R}(k): h \ne k$. Clearly $k \in K(1)^+$, so there exist $g \in K(1)$, $\phi \in E_s(1,g)$, $y \in \mathcal{V} - var(g)$ such that $k = g + <y, \phi>$. By lemma~\ref{L:useful-contexts-1} we obtain that $h \in \mathcal{R}(g)$. Since $g = \epsilon$ then also $h = \epsilon \in K(1)$ , so $\sigma_{/dom(h)} = \sigma_{/\emptyset} = \epsilon \in \Xi(\epsilon) = \Xi(h)$.\\

Let's now examine the case where $n+1 > 2$. Let $k \in K(n+1): k \ne \epsilon$, let $\sigma \in \Xi(k)$, $h \in \mathcal{R}(k): h \ne k$, we have to show there exists $m < n+1$ such that $h \in K(m)$ and it results $\sigma_{/dom(h)} \in \Xi(h)$.\\

As we have just proved in relation to assumption~\ref{A:k-is-epsilon-or-new}, there exist $m<n+1$, $g \in K(m)$, $\phi \in E_s(m,g)$, $y \in (\mathcal{V}-var(g))$ such that $k = g + <y, \phi>$, $\Xi(k) = \{ \rho + (y,s) | \, \rho \in \Xi(g), s \in \#(g,\phi,\rho) \}$ ).\\

This implies there exist $\rho \in \Xi(g), s \in \#(g,\phi,\rho)$ such that $\sigma = \rho + (y,s)$. By assumption~\ref{A:simple-k-sigma-slp-new} and lemma~\ref{L:useful4-slp} we have that $\sigma_{/dom(g)} = \sigma_{/dom(\rho)} = \rho$.\\ 

\medskip

If $h = g$ then $\sigma_{/dom(h)} = \sigma_{/dom(g)} = \rho \in \Xi(h)$.\\

Otherwise we have $h \ne g$. Since $k = g + <y,\phi>$, $h \in \mathcal{R}(k)$, $h \ne k$ by lemma~\ref{L:useful-contexts-1} we have that $h \in \mathcal{R}(g)$. If $g = \epsilon$ we would have $h = \epsilon = g$, so $g \ne \epsilon$. This implies that $m \geqslant 2$. By our inductive hypothesis we obtain there exists $q < m \leqslant n$ such that $h \in K(q)$ and $\rho_{/dom(h)} \in \Xi(h)$. Now by lemma~\ref{L:useful5-slp}
\[ \sigma_{/dom(h)} = (\sigma_{/dom(g)})_{/dom(h)} =  \rho_{/dom(h)} \in \Xi(h) . \]
\end{proof}

\bigskip

\begin{proof}[Proof of~\ref{A:Ekn-decidable-predicates}]
Given $k \in K(n+1)$, $\varphi \in E(n+1,k)$, $q$ positive integer, $i = 1 \dots p$ we must show we are able to decide all of the following conditions:
\begin{itemize}
\item for each $\sigma \in \Xi(k)$ $Set_q( \#(k, \varphi, \sigma))$;
\item for each $\sigma \in \Xi(k)$ $Event_q( \#(k, \varphi, \sigma))$;
\item for each $\sigma \in \Xi(k)$ $\#(k, \varphi, \sigma) \in D_i$;
\item for each $\sigma \in \Xi(k)$ $\#(k, \varphi, \sigma) \in \mathcal{P}^q(D_i)$;
\item if (for each $\sigma \in \Xi(k)$ $Set_q( \#(k, \varphi, \sigma))$) then\\
(for each $\sigma \in \Xi(k)$ $NotEmpty_q( \#(k, \varphi, \sigma))$).
\end{itemize}

\medskip

And we also need to verify that the last condition holds true.\\

We have seen that if $k \in K(n)^+$ $E(n+1,k) = E_a(n+1,k)$, and if $k \in K(n)$ 
\scriptsize
\begin{align*}
E(n+1,k) &= \bigcup_{u \in \mathcal{H}(n+1,k)} \mathcal{E}(u)\\
&= E(n,k) \cup E_b(n+1,k) \cup E_\forall(n+1,k) \cup E_\exists(n+1,k) \cup \bigcup_{c \in \mathcal{C}'} E^c(n+1,k) \cup \bigcup_{f \in \mathcal{F}} E^f(n+1,k) \ .
\end{align*}
\normalsize

\medskip

Let now $k \in K(n)^+$ and let's try to show how we can decide the conditions for $\varphi \in E_a(n+1,k)$.\\

Since $k \in K(n)^+$  there exist $h \in K(n), \phi \in E_s(n,h), y \in (\mathcal{V}-var(h))$ such that $k = h + <y,\phi>$, $\Xi(k) = \{ \rho + (y,s) | \, \rho \in \Xi(h), s \in \#(h,\phi,\rho) \}$. Moreover $h$, $y$ and $\phi$ are clearly identifiable within $k$ and $E_a(n+1, k) = \{y \}$.\\

Given $\varphi \in E_a(n+1, k)$, $\sigma = \rho + (y, s) \in \Xi(k)$ $\#(k, \varphi , \sigma) = s \in \#(h,\phi,\rho)$.\\

We first consider the condition `for each $\sigma \in \Xi(k)$ $Set_q(\#(k, \varphi, \sigma))$' (where $q$ is a positive integer).\\

We consider that $\phi \in E(n,h)$ and by inductive hypothesis we can decide whether `for each $\rho \in \Xi(h)$ $Set_{q+1}( \#(h, \phi, \rho))$'.\\

If we decide this is true then for each $\sigma \in \Xi(k)$ there exist $\rho \in \Xi(h)$, $s \in \#(h,\phi,\rho)$ such that $\sigma = \rho + (y,s)$, $\#(k, \varphi, \sigma) = s \in \#(h,\phi,\rho)$, and since $Set_{q+1}( \#(h, \phi, \rho))$ we have $Set_q(\#(k, \varphi, \sigma))$.\\

So if we decide `for each $\rho \in \Xi(h)$ $Set_{q+1}( \#(h, \phi, \rho))$' is true we can use this to decide `for each $\sigma \in \Xi(k)$ $Set_q(\#(k, \varphi, \sigma))$' is true.\\

If instead we decide `for each $\rho \in \Xi(h)$ $Set_{q+1}( \#(h, \phi, \rho))$' is false this means there exists $\rho \in \Xi(h)$ such that $\neg (Set_{q+1}( \#(h, \phi, \rho)))$. Since $\phi \in E_s(n,h)$ we have that $\#(h, \phi, \rho)$ is a set and so there exists $s \in \#(h, \phi, \rho)$ such that $\neg(Set_q(s))$. If we set $\sigma = \rho + (y,s)$ then $\sigma \in \Xi(k)$ and $\#(k, \varphi, \sigma) = s$ so $\neg(Set_q(\#(k, \varphi, \sigma))$.\\

So if we decide `for each $\rho \in \Xi(h)$ $Set_{q+1}( \#(h, \phi, \rho))$' is false we can use this to decide `for each $\sigma \in \Xi(k)$ $Set_q(\#(k, \varphi, \sigma))$' is false too.\\

We now want to prove the following:\\

if (for each $\sigma \in \Xi(k)$ $Set_q( \#(k, \varphi, \sigma))$) then\\
(for each $\sigma \in \Xi(k)$ $NotEmpty_q( \#(k, \varphi, \sigma))$).\\

We assume (for each $\sigma \in \Xi(k)$ $Set_q( \#(k, \varphi, \sigma))$), clearly this implies\\ 
(for each $\rho \in \Xi(h)$ $Set_{q+1}( \#(h, \phi, \rho))$), which (by inductive hypothesis) implies\\
(for each $\rho \in \Xi(h)$ $NotEmpty_{q+1}( \#(h, \phi, \rho))$).\\

We can then consider that for each $\sigma \in \Xi(k)$ there exist $\rho \in \Xi(h)$, $s \in \#(h,\phi,\rho)$ such that $\sigma = \rho + (y,s)$, $\#(k, \varphi, \sigma) = s \in \#(h,\phi,\rho)$. Since $NotEmpty_{q+1}( \#(h, \phi, \rho))$ holds then $NotEmpty_q( \#(k, \varphi, \sigma))$ holds too.\\

Given $i = 1 \dots p$ we now want to consider the predicate `for each $\sigma \in \Xi(k)$ $\#(k, \varphi, \sigma) \in D_i$'.\\

By the inductive hypothesis we are able to decide the predicate `for each $\rho \in \Xi(h)$ $\#(h,\phi,\rho) \in \mathcal{P}(D_i)$'.\\

If we decide the last condition is true then as seen above for each $\sigma \in \Xi(k)$ there exist $\rho \in \Xi(h)$, $s \in \#(h,\phi,\rho)$ such that $\sigma = \rho + (y,s)$, $\#(k, \varphi, \sigma) = s \in \#(h,\phi,\rho)  \in \mathcal{P}(D_i)$, therefore $\#(k, \varphi, \sigma) \in D_i$.\\

If instead we decide the mentioned condition is false, then there exists $\rho \in \Xi(h)$: $\#(h,\phi,\rho) \notin \mathcal{P}(D_i)$. Since $\#(h,\phi,\rho)$ is a set and it is not empty, this means there exists $s \in \#(h,\phi,\rho)$: $s \notin D_i$. If we set $\sigma = \rho + (y,s)$ then $\sigma \in \Xi(k)$ and $\#(k, \varphi, \sigma) = s \notin D_i$, so there exists $\sigma \in \Xi(k)$: $\#(k, \varphi, \sigma) \notin D_i$.\\

Given $i = 1 \dots p$ and a positive integer $q$ we now want to consider the predicate `for each $\sigma \in \Xi(k)$ $\#(k, \varphi, \sigma) \in \mathcal{P}^q(D_i)$'.\\

By the inductive hypothesis we are able to decide the predicate `for each $\rho \in \Xi(h)$ $\#(h,\phi,\rho) \in \mathcal{P}^{q+1}(D_i)$'.\\

If we decide the last condition is true then as seen above for each $\sigma \in \Xi(k)$ there exist $\rho \in \Xi(h)$, $s \in \#(h,\phi,\rho)$ such that $\sigma = \rho + (y,s)$, $\#(k, \varphi, \sigma) = s \in \#(h,\phi,\rho)  \in \mathcal{P}^{q+1}(D_i)$, therefore $\#(k, \varphi, \sigma) \in \#(h,\phi,\rho) \subseteq  \mathcal{P}^q(D_i)$.\\

If instead we decide the mentioned condition is false, then there exists $\rho \in \Xi(h)$: $\#(h,\phi,\rho) \notin \mathcal{P}^{q+1}(D_i)$. Since $\#(h,\phi,\rho)$ is a set and it is not empty, this means there exists $s \in \#(h,\phi,\rho)$: $s \notin \mathcal{P}^q(D_i)$. If we set $\sigma = \rho + (y,s)$ then $\sigma \in \Xi(k)$ and $\#(k, \varphi, \sigma) = s \notin \mathcal{P}^q(D_i)$, so there exists $\sigma \in \Xi(k)$: $\#(k, \varphi, \sigma) \notin \mathcal{P}^q(D_i)$.\\

Given a positive integer $q$ we now want to consider the predicate `for each $\sigma \in \Xi(k)$ $Event_q(\#(k, \varphi, \sigma))$'.\\

By the inductive hypothesis we are able to decide the predicate `for each $\rho \in \Xi(h)$ $Event_{q+1}(\#(h, \phi, \rho))$'.\\

If we decide the last condition is true then as seen above for each $\sigma \in \Xi(k)$ there exists $\rho \in \Xi(h)$ such that $\#(k, \varphi, \sigma) \in \#(h,\phi,\rho)$. Since $Event_{q+1}(\#(h, \phi, \rho))$ we have $Event_q(\#(k, \varphi, \sigma))$.\\

If instead we decide the mentioned condition is false, then there exists $\rho \in \Xi(h)$: $\neg(Event_{q+1}(\#(h, \phi, \rho)))$. This implies there exists $s \in \#(h,\phi,\rho)$: $\neg(Event_q(s))$. If we set $\sigma = \rho + (y,s)$ then $\sigma \in \Xi(k)$ and $\#(k, \varphi, \sigma) = s$, so $\neg(Event_q(\#(k, \varphi, \sigma)))$. This means there exists $\sigma \in \Xi(k)$: $\neg(Event_q(\#(k, \varphi, \sigma)))$.\\

Let's now move to the second step of our proof, where we expect to prove that if $\mathbf{k \in K(n)}$ then for each $\varphi \in E(n+1,k)$ we can decide all of our conditions.\\

We recall that
\scriptsize
\begin{align*}
E(n+1,k) = E(n,k) \cup E_b(n+1,k) \cup E_\forall(n+1,k) \cup E_\exists(n+1,k) \cup \bigcup_{c \in \mathcal{C}'} E^c(n+1,k) \cup \bigcup_{f \in \mathcal{F}} E^f(n+1,k) \ .
\end{align*}
\normalsize

By the inductive hypothesis (i.e. what we assumed true at level $n$) we can assume that we can decide each of our conditions for $\mathbf{\varphi \in E(n,k)}$.\\

Let's now try to prove we can decide all of our conditions for each $\mathbf{\varphi \in E_b(n+1, k)}$.\\

If $k = \epsilon$ we have $E_b(n+1, k) = \emptyset$, so in this case there is nothing to decide.\\

We'll then consider the case $k \ne \epsilon$. Here by our assumption~\ref{A:k-is-epsilon-or-new} $n > 1$ and there exist $m<n$, $h \in K(m)$, $\phi \in E_s(m,h)$, $y \in (\mathcal{V}-var(h))$ such that $k = h + <y, \phi>$, $\Xi(k) = \{ \rho + (y,s) | \, \rho \in \Xi(h), s \in \#(h,\phi,\rho) \} )$.\\

For each $\varphi \in E_b(n+1,k), \  \sigma = \rho + (y,s) \in \Xi(k)$ $\#(k,\varphi,\sigma) = \#(h,\varphi,\rho)$.\\

We first consider the condition `for each $\sigma \in \Xi(k)$ $Set_q(\#(k, \varphi, \sigma))$' (where $q$ is a positive integer).\\

By the inductive hypothesis we can decide if the following condition holds: `for each $\rho \in \Xi(h)$ $Set_q(\#(h, \varphi, \rho))$'.\\

If the just mentioned condition holds we can consider that for each $\sigma \in \Xi(k)$ there exist $\rho \in \Xi(h)$, $s \in \#(h,\phi,\rho)$ such that $\sigma = \rho + (y,s)$ and $\#(k,\varphi,\sigma) = \#(h,\varphi,\rho)$. Since $Set_q(\#(h, \varphi, \rho))$ then $Set_q(\#(k, \varphi, \sigma))$ holds too.\\

If the mentioned condition is decided as false then there exists $\rho \in \Xi(h)$: $\neg(Set_q(\#(h, \varphi, \rho)))$. We have that for each $\delta \in \Xi(h)$ $Set_1(\#(h, \phi, \delta))$, so for each $\delta \in \Xi(h)$ $NotEmpty_1(\#(h, \phi, \delta))$. So let $s \in \#(h, \phi, \rho)$ and let $\sigma = \rho + (y,s)$, then $\sigma \in \Xi(k)$ and $\#(k,\varphi,\sigma) = \#(h,\varphi,\rho)$ and so $\neg(Set_q(\#(k, \varphi, \sigma)))$.\\

We now want to prove the following:\\

if (for each $\sigma \in \Xi(k)$ $Set_q( \#(k, \varphi, \sigma))$) then\\
(for each $\sigma \in \Xi(k)$ $NotEmpty_q( \#(k, \varphi, \sigma))$).\\

We assume (for each $\sigma \in \Xi(k)$ $Set_q( \#(k, \varphi, \sigma))$), clearly this implies\\ 
(for each $\rho \in \Xi(h)$ $Set_q( \#(h, \phi, \rho))$), which (by inductive hypothesis) implies\\
(for each $\rho \in \Xi(h)$ $NotEmpty_q( \#(h, \phi, \rho))$).\\

We can then consider that for each $\sigma \in \Xi(k)$ there exist $\rho \in \Xi(h)$, $s \in \#(h,\phi,\rho)$ such that $\sigma = \rho + (y,s)$, $\#(k, \varphi, \sigma) = \#(h,\phi,\rho)$. Since $NotEmpty_q( \#(h, \phi, \rho))$ holds then $NotEmpty_q( \#(k, \varphi, \sigma))$ holds too.\\

Given $i = 1 \dots p$ we now want to consider the predicate `for each $\sigma \in \Xi(k)$ $\#(k, \varphi, \sigma) \in D_i$'.\\

By the inductive hypothesis we can decide the condition `for each $\rho \in \Xi(h)$ $\#(h, \varphi, \rho) \in D_i$'.\\

If the just mentioned condition holds then we consider that for each $\sigma \in \Xi(k)$ there exist $\rho \in \Xi(h)$, $s \in \#(h,\phi,\rho)$ such that $\sigma = \rho + (y,s)$ and $\#(k,\varphi,\sigma) = \#(h,\varphi,\rho)$. Therefore clearly $\#(k, \varphi, \sigma) \in D_i$.\\

If on the contrary the mentioned condition is decided as false then there exists $\rho \in \Xi(h)$: $\#(h, \varphi, \rho) \notin D_i$. We have that for each $\delta \in \Xi(h)$ $Set_1(\#(h, \phi, \delta))$, so for each $\delta \in \Xi(h)$ $NotEmpty_1(\#(h, \phi, \delta))$. So let $s \in \#(h, \phi, \rho)$ and let $\sigma = \rho + (y,s)$, then $\sigma \in \Xi(k)$ and $\#(k,\varphi,\sigma) = \#(h,\varphi,\rho) \notin D_i$.\\

Given $i = 1 \dots p$ and a positive integer $q$ we now want to consider the predicate `for each $\sigma \in \Xi(k)$ $\#(k, \varphi, \sigma) \in \mathcal{P}^q(D_i)$'.\\

By the inductive hypothesis we are able to decide the predicate `for each $\rho \in \Xi(h)$ $\#(h,\varphi,\rho) \in \mathcal{P}^q(D_i)$'.\\

If the just mentioned condition holds then we consider that for each $\sigma \in \Xi(k)$ there exist $\rho \in \Xi(h)$, $s \in \#(h,\phi,\rho)$ such that $\sigma = \rho + (y,s)$ and $\#(k,\varphi,\sigma) = \#(h,\varphi,\rho)$. Therefore clearly $\#(k, \varphi, \sigma) \in \mathcal{P}^q(D_i)$.\\

If on the contrary the mentioned condition is decided as false then there exists $\rho \in \Xi(h)$: $\#(h, \varphi, \rho) \notin \mathcal{P}^q(D_i)$. We have that for each $\delta \in \Xi(h)$ $Set_1(\#(h, \phi, \delta))$, so for each $\delta \in \Xi(h)$ $NotEmpty_1(\#(h, \phi, \delta))$. So let $s \in \#(h, \phi, \rho)$ and let $\sigma = \rho + (y,s)$, then $\sigma \in \Xi(k)$ and $\#(k,\varphi,\sigma) = \#(h,\varphi,\rho) \notin \mathcal{P}^q(D_i)$.\\

Given a positive integer $q$ we now want to consider the condition `for each $\sigma \in \Xi(k)$ $Event_q(\#(k, \varphi, \sigma))$'.\\

By the inductive hypothesis we are able to decide the condition `for each $\rho \in \Xi(h)$ $Event_q(\#(h, \varphi, \rho))$'.\\

If the just mentioned condition holds then we consider that for each $\sigma \in \Xi(k)$ there exist $\rho \in \Xi(h)$, $s \in \#(h,\phi,\rho)$ such that $\sigma = \rho + (y,s)$ and $\#(k,\varphi,\sigma) = \#(h,\varphi,\rho)$. Therefore clearly $Event_q(\#(k, \varphi, \sigma))$.\\

If on the contrary the mentioned condition is decided as false then there exists $\rho \in \Xi(h)$: $\neg(Event_q(\#(h, \varphi, \rho)))$. We have that for each $\delta \in \Xi(h)$ $Set_1(\#(h, \phi, \delta))$, so for each $\delta \in \Xi(h)$ $NotEmpty_1(\#(h, \phi, \delta))$. So let $s \in \#(h, \phi, \rho)$ and let $\sigma = \rho + (y,s)$, then $\sigma \in \Xi(k)$ and $\#(k,\varphi,\sigma) = \#(h,\varphi,\rho)$, so $\neg(Event_q(\#(k, \varphi, \sigma)))$.\\

Let's now try to prove we can decide all of our conditions for each $\mathbf{\chi \in E_Q(n+1, k)}$.\\

We recall that for every $\chi = Q(x:\varphi,\phi) \in E_Q(n+1, k)$ and $\sigma \in \Xi(k)$ we defined $\#(k, \chi, \sigma)$ as follows.\\

Let $k' = k + <x, \varphi>$, then: if $Q = \forall$:
\[ \#(k, \chi, \sigma) = \text{for each} \ \sigma' \in \Xi(k'): \sigma \sqsubseteq \sigma'  \ \#(k', \phi, \sigma') \ . \]

Else if $Q = \exists$:
\[ \#(k, \chi, \sigma) = \text{exists} \ \sigma' \in \Xi(k'): \sigma \sqsubseteq \sigma'  \ \#(k', \phi, \sigma') \ . \]
 
We first consider the condition `for each $\sigma \in \Xi(k)$ $Set_q(\#(k, \chi, \sigma))$' (where $q$ is a positive integer).\\

For each $\sigma \in \Xi(k)$ $\neg(Set_1(\#(k, \chi, \sigma))$, so $\neg(Set_q(\#(k, \chi, \sigma))$ and since $\Xi(k) \ne \emptyset$ we can decide the condition `for each $\sigma \in \Xi(k)$ $Set_q(\#(k, \chi, \sigma))$' is false.\\

Given $i = 1 \dots p$ we must be able to decide the condition `for each $\sigma \in \Xi(k)$ $\#(k, \chi, \sigma) \in D_i$'.\\

Clearly for each $\sigma \in \Xi(k)$ $\#(k, \chi, \sigma)$ is true or false and so $\#(k, \chi, \sigma) \notin D_i$. Since $\Xi(k) \ne \emptyset$ the condition `for each $\sigma \in \Xi(k)$ $\#(k, \chi, \sigma) \in D_i$' is false.\\

Given $i = 1 \dots p$ and a positive integer $q$ we must be able to decide the condition `for each $\sigma \in \Xi(k)$ $\#(k, \chi, \sigma) \in \mathcal{P}^q(D_i)$'.\\

For each $\sigma \in \Xi(k)$ if it was $\#(k, \chi, \sigma) \in \mathcal{P}^q(D_i)$ then it would be $Set_1(\#(k, \chi, \sigma))$, but we have seen that $Set_1(\#(k, \chi, \sigma))$ is false, so $\#(k, \chi, \sigma) \in \mathcal{P}^q(D_i)$ is also false. Since $\Xi(k) \ne \emptyset$ the condition `for each $\sigma \in \Xi(k)$ $\#(k, \chi, \sigma) \in  \mathcal{P}^q(D_i)$' is false.\\

We must be able to decide the condition `for each $\sigma \in \Xi(k)$ $Event_1( \#(k, \chi, \sigma))$'. This condition is obviously true.\\

Given a positive integer $q > 1$ we must be able to decide the condition `for each $\sigma \in \Xi(k)$ $Event_q( \#(k, \chi, \sigma))$'.\\

For each $\sigma \in \Xi(k)$ if it was $Event_q( \#(k, \chi, \sigma))$ then it would be $Set_1(\#(k, \chi, \sigma))$, but we have seen that $Set_1(\#(k, \chi, \sigma))$ is false, so $Event_q( \#(k, \chi, \sigma))$ is also false. Since $\Xi(k) \ne \emptyset$ the condition `for each $\sigma \in \Xi(k)$ $Event_q( \#(k, \chi, \sigma))$' is false.\\

Given $c \in \mathcal{C}'$, let's now try to prove we can decide all of our conditions for each $\mathbf{t \in E^c(n+1, k)}$.\\

For each $t = (c)(\varphi_1, \dots , \varphi_m) \in E^c(n+1,k)$ we defined
\[
\#(k,t,\sigma) = \#(c)( \#(k, \varphi_1, \sigma), \dots , \#(k, \varphi_m, \sigma) ).
\]

Our costant $c$ can have as meaning three types of function, in every case the domain of $\#(c)$ is $M_1 \times \dots \times M_m$, but we have a different codomain in the three cases, which are the following:

\begin{itemize}
\item $\alpha \in \{ 1, \dots , p \}$ and $\#(c) : M_1 \times \dots \times M_m \to D_{\alpha}$;
\item $\alpha \in \{ 1, \dots , p \}$, $r$ positive integer and $\#(c) : M_1 \times \dots \times M_m \to \mathcal{P}^r(D_{\alpha})$;
\item $\#(c)$ is a function over $M_1 \times \dots \times M_m$ such that for each $(d_1, \dots , d_m) \in M_1 \times \dots \times M_m$ $\#(c)(d_1, \dots , d_m)$ is true or false.
\end{itemize}

Let's consider the first of those cases, clearly in this case for each $\sigma \in \Xi(k)$ $\#(k,t,\sigma) \in D_{\alpha}$.\\

Let's consider $\beta = 1 \dots p$ such that $\beta \ne \alpha$ and we want to decide the condition `for each $\sigma \in \Xi(k)$ $\#(k,t,\sigma) \in D_{\beta}$'.\\

We have assumed that $D_\alpha \cap D_\beta = \emptyset$, so for each $\sigma \in \Xi(k)$ $\#(k,t,\sigma) \notin D_{\beta}$. Since $\Xi(k) \ne \emptyset$ the condition `for each $\sigma \in \Xi(k)$ $\#(k,t,\sigma) \in D_{\beta}$' is false.\\

Given $\beta = 1 \dots p$ and a positive integer $q$ we want to decide the condition `for each $\sigma \in \Xi(k)$ $\#(k,t,\sigma) \in \mathcal{P}^q(D_{\beta})$'.\\

For each $\sigma \in \Xi(k)$ $\#(k,t,\sigma) \in D_{\alpha}$, so $\#(k,t,\sigma)$ is not a set, so $\#(k,t,\sigma) \notin \mathcal{P}^q(D_{\beta})$. Since $\Xi(k) \ne \emptyset$ the condition `for each $\sigma \in \Xi(k)$ $\#(k, t, \sigma) \in  \mathcal{P}^q(D_\beta)$' is false.\\

We then consider the condition `for each $\sigma \in \Xi(k)$ $Set_q(\#(k, t, \sigma))$' (where $q$ is a positive integer).\\

For each $\sigma \in \Xi(k)$ $\#(k,t,\sigma) \in D_{\alpha}$, so $\#(k,t,\sigma)$ is not a set, so $Set_q(\#(k, t, \sigma))$ is false. Since $\Xi(k) \ne \emptyset$ the condition `for each $\sigma \in \Xi(k)$ $Set_q(\#(k, t, \sigma))$' is false.\\

We must also be able to decide the condition `for each $\sigma \in \Xi(k)$ $Event_1( \#(k, t, \sigma))$'.\\

For each $\sigma \in \Xi(k)$ $\#(k,t,\sigma) \in D_{\alpha}$, so $Event_1( \#(k, t, \sigma))$ is false. Since $\Xi(k) \ne \emptyset$ the condition `for each $\sigma \in \Xi(k)$ $Event_1( \#(k, t, \sigma))$' is false.\\

Give a positive integer $q > 1$ we must be able to decide the condition `for each $\sigma \in \Xi(k)$ $Event_q( \#(k, t, \sigma))$'.\\

For each $\sigma \in \Xi(k)$ $\#(k,t,\sigma) \in D_{\alpha}$, so $Set_1(\#(k,t,\sigma))$ is false, so $Event_q( \#(k, t, \sigma))$ is false. Since $\Xi(k) \ne \emptyset$ the condition `for each $\sigma \in \Xi(k)$ $Event_q( \#(k, t, \sigma))$' is false.\\

\medskip

Let's now consider the case where $\alpha \in \{ 1, \dots , p \}$, $r$ positive integer and $\#(c) : M_1 \times \dots \times M_m \to \mathcal{P}^r(D_{\alpha})$. Clearly in this case for each $\sigma \in \Xi(k)$ $\#(k,t,\sigma) \in \mathcal{P}^r(D_{\alpha})$.\\

Given $\beta = 1 \dots p$ we want to decide the condition `for each $\sigma \in \Xi(k)$ $\#(k,t,\sigma) \in D_{\beta}$'.\\

For each $\sigma \in \Xi(k)$ $\#(k,t,\sigma) \in \mathcal{P}^r(D_{\alpha})$, so $\#(k,t,\sigma)$ is a set and $\#(k,t,\sigma) \notin D_{\beta}$. Since $\Xi(k) \ne \emptyset$ the condition `for each $\sigma \in \Xi(k)$ $\#(k,t,\sigma) \in D_{\beta}$' is false.\\

As we have seen the condition `for each $\sigma \in \Xi(k)$ $\#(k,t,\sigma) \in \mathcal{P}^r(D_{\alpha})$' is true.\\

Let's consider a positive integer $s$ and $\beta = 1 \dots p$ such that $s \ne r$ or $\beta \ne \alpha$, we want to decide the condition `for each $\sigma \in \Xi(k)$ $\#(k,t,\sigma) \in \mathcal{P}^s(D_{\beta})$'.\\

We first consider the case $s \ne r$. Using lemma~\ref{L:useful-on-predicates-2} we have that $\mathcal{P}^s(D_{\beta}) \cap \mathcal{P}^r(D_{\alpha}) = \emptyset$, so for each $\sigma \in \Xi(k)$ $\#(k,t,\sigma) \notin \mathcal{P}^s(D_{\beta})$. Since $\Xi(k) \ne \emptyset$ the condition `for each $\sigma \in \Xi(k)$ $\#(k,t,\sigma) \in \mathcal{P}^s(D_{\beta})$' is false.\\

Let now $s = r$ and let $\beta \ne \alpha$, using lemma~\ref{L:useful-on-predicates-3} $\mathcal{P}^s(D_{\beta}) \cap \mathcal{P}^r(D_{\alpha}) = \emptyset$ and again for each $\sigma \in \Xi(k)$ $\#(k,t,\sigma) \notin \mathcal{P}^s(D_{\beta})$, and the condition `for each $\sigma \in \Xi(k)$ $\#(k,t,\sigma) \in \mathcal{P}^s(D_{\beta})$' is false.\\

Given a positive integer $s$ we want to decide the condition `for each $\sigma \in \Xi(k)$ $Set_s(\#(k,t,\sigma))$', and we also want to verify that the following condition `if (for each $\sigma \in \Xi(k)$ $Set_s(\#(k,t,\sigma))$) then (for each $\sigma \in \Xi(k)$ $NotEmpty_s(\#(k,t,\sigma)))$' is true.\\

We first consider the case $s \leqslant r$. For each $\sigma \in \Xi(k)$ $\#(k,t,\sigma) \in \mathcal{P}^r(D_{\alpha})$, so by lemma~\ref{L:useful-on-predicates-4} we have $Set_s(\#(k,t,\sigma))$ and $NotEmpty_s(\#(k,t,\sigma))$.\\

Let's now consider the case $s > r$. Here for each $\sigma \in \Xi(k)$ $\#(k,t,\sigma) \in \mathcal{P}^r(D_{\alpha})$, so by lemma~\ref{L:useful-on-predicates-1}  $\neg(Set_s(\#(k,t,\sigma)))$, so the condition `for each $\sigma \in \Xi(k)$ $Set_s(\#(k,t,\sigma))$' is false.\\

Given a positive integer $s$ we want to decide the condition `for each $\sigma \in \Xi(k)$ $Event_s(\#(k,t,\sigma))$'.\\

We can still consider that for each $\sigma \in \Xi(k)$ $\#(k,t,\sigma) \in \mathcal{P}^r(D_{\alpha})$, and by lemmas~\ref{L:useful-on-predicates-5} and~\ref{L:useful-on-predicates-7} we obtain that $\neg(Event_s(\#(k,t,\sigma)))$. This implies that the condition `for each $\sigma \in \Xi(k)$ $Event_s(\#(k,t,\sigma))$' is false.\\

\medskip

Finally let's consider the case where $\#(c)$ is a function over $M_1 \times \dots \times M_m$ such that for each $(d_1, \dots , d_m) \in M_1 \times \dots \times M_m$ $\#(c)(d_1, \dots , d_m)$ is true or false. Clearly in this case for each $\sigma \in \Xi(k)$ $\#(k,t,\sigma)$ is true or false.\\

Given a positive integer $q$ we want to decide the condition `for each $\sigma \in \Xi(k)$ $Set_q(\#(k,t,\sigma))$'.\\

For each $\sigma \in \Xi(k)$ $\#(k,t,\sigma)$ is true or false, so $Set_1(\#(k,t,\sigma))$ is false and $Set_q(\#(k,t,\sigma))$ is false. Since $\Xi(k) \ne \emptyset$ we can decide the condition `for each $\sigma \in \Xi(k)$ $Set_q(\#(k,t,\sigma))$' is false.\\

Given $\beta = 1 \dots p$ we want to decide the condition `for each $\sigma \in \Xi(k)$ $\#(k,t,\sigma) \in D_{\beta}$'.\\

For each $\sigma \in \Xi(k)$ $\#(k,t,\sigma)$ is true or false, so $\#(k,t,\sigma) \notin D_{\beta}$. Therefore the condition `for each $\sigma \in \Xi(k)$ $\#(k,t,\sigma) \in D_{\beta}$' is false.\\

Let's consider a positive integer $q$ and $\beta = 1 \dots p$, we want to decide the condition `for each $\sigma \in \Xi(k)$ $\#(k,t,\sigma) \in \mathcal{P}^q(D_{\beta})$'.\\

For each $\sigma \in \Xi(k)$ $\#(k,t,\sigma)$ is true or false, so $Set_1(\#(k,t,\sigma))$ is false, so $\#(k,t,\sigma) \notin \mathcal{P}^q(D_{\beta})$'. Therefore the condition `for each $\sigma \in \Xi(k)$ $\#(k,t,\sigma) \in \mathcal{P}^q(D_{\beta})$' is false.\\

We can also decide the condition `for each $\sigma \in \Xi(k)$ $Event_1(\#(k,t,\sigma))$'. Clearly the condition is true.\\

Given a positive integer $q$ we also want to decide the condition `for each $\sigma \in \Xi(k)$ $Event_q(\#(k,t,\sigma))$'.\\

For each $\sigma \in \Xi(k)$ $\#(k,t,\sigma)$ is true or false, so $Set_1(\#(k,t,\sigma))$ is false, so $Event_q(\#(k,t,\sigma))$ is false. Therefore the condition `for each $\sigma \in \Xi(k)$ $Event_q(\#(k,t,\sigma))$' is false.\\

Given $f \in \mathcal{F}$, let's now try to prove we can decide all of our conditions for each $\mathbf{t \in E^f(n+1, k)}$.\\

We first consider the case where $f$ has multiplicity 2.\\

By our definitions, for each $t = f(\varphi_1, \varphi_2) \in E^f(n+1,k)$,
\[
\#(k,t,\sigma) = P_f( \#(k, \varphi_1, \sigma), \#(k, \varphi_2, \sigma) ).
\]

Clearly for each $t \in E^f(n+1,k)$ and $\sigma \in \Xi(k)$ $\#(k,t,\sigma)$ is true or false.\\

Let's also consider the case where $f$ has multiplicity 1.\\

By our definitions, for each $t = f(\varphi_1) \in E^f(n+1,k)$,
\[
\#(k,t,\sigma) = P_f( \#(k, \varphi_1, \sigma) ).
\]

It is also true in this case that for each $t \in E^f(n+1,k)$ and $\sigma \in \Xi(k)$ $\#(k,t,\sigma)$ is true or false, and using this very property we can decide each of our conditions.\\

Given $\alpha = 1 \dots p$ we must be able to decide the condition `for each $\sigma \in \Xi(k)$ $\#(k, t, \sigma) \in D_{\alpha}$'.\\

Since for each $\sigma \in \Xi(k)$ $\#(k,t,\sigma)$ is true or false, then for each $\sigma \in \Xi(k)$ $\#(k, t, \sigma) \notin D_{\alpha}$, and so then condition `for each $\sigma \in \Xi(k)$ $\#(k, t, \sigma) \in D_{\alpha}$' is false.\\

Given $\alpha = 1 \dots p$ and a positive integer $q$ we must be able to decide the condition `for each $\sigma \in \Xi(k)$ $\#(k, t, \sigma) \in \mathcal{P}^q(D_\alpha)$'.\\

Given $\sigma \in \Xi(k)$ $Event_1(\#(k,t,\sigma))$ and by lemma~\ref{L:useful-on-predicates-5} this implies $\#(k, t, \sigma) \notin \mathcal{P}^q(D_\alpha)$. Therefore the condition `for each $\sigma \in \Xi(k)$ $\#(k, t, \sigma) \in \mathcal{P}^q(D_\alpha)$' is false.\\

Given a positive integer $r$ we must be able to decide the condition `for each $\sigma \in \Xi(k)$ $Set_r(\#(k, t, \sigma))$', and when this condition is decided as true we must also be able to decide that for each $\sigma \in \Xi(k)$ $NotEmpty_r(\#(k, t, \sigma))$.\\

For each $\sigma \in \Xi(k)$ $Event_1(\#(k,t,\sigma))$ so $\neg Set_1(\#(k, t, \sigma))$ and then also $\neg Set_r(\#(k, t, \sigma))$. Therefore the condition `for each $\sigma \in \Xi(k)$ $Set_r(\#(k, t, \sigma))$' is false.\\

Given a positive integer $r$ we must be able to decide the condition `for each $\sigma \in \Xi(k)$ $Event_r(\#(k, t, \sigma))$'.\\

Clearly the condition is true for $r = 1$, while for $r > 1$ given $\sigma \in \Xi(k)$ $\neg Set_1(\#(k, t, \sigma))$ and so $\neg Event_r(\#(k, t, \sigma))$, so the condition is false for $r > 1$.\\
\end{proof}

\bigskip

\begin{proof}[Proof of~\ref{A:expr-cont-depth}]

\smallskip

We need to prove that for each $k \in K(n+1), \ t \in E(n+1,k)$
\begin{itemize}
\item $t[\ell(t)] \ne \text{`('}$ ;
\item if $t[\ell(t)] = \text{`)'}$ then $d(t, \ell(t)) = 1$, else $d(t, \ell(t)) = 0$ ;
\item for each $\alpha \in \{1, \dots, \ell(t) \}$ if $(t[\alpha]=\text{`:'}) \vee (t[\alpha]=\text{`,'}) \vee (t[\alpha]=\text{`)'}) $ then $d(t, \alpha) \geqslant 1$.
\end{itemize}

\medskip

We have seen that if $k \in K(n)^+$ $E(n+1,k) = E_a(n+1,k)$, and if $k \in K(n)$ 
\scriptsize
\begin{align*}
E(n+1,k) &= \bigcup_{u \in \mathcal{H}(n+1,k)} \mathcal{E}(u)\\
&= E(n,k) \cup E_b(n+1,k) \cup E_\forall(n+1,k) \cup E_\exists(n+1,k) \cup \bigcup_{c \in \mathcal{C}'} E^c(n+1,k) \cup \bigcup_{f \in \mathcal{F}} E^f(n+1,k) \ .
\end{align*}
\normalsize

\medskip

Let $k \in K(n)^+$ and $t \in \mathbf{E_a(n+1,k)}$. There exist $h \in K(n)$, $\phi \in E_s(n,h)$, $y \in \mathcal{V} - var(h)$ such that $k = h + <y,\phi>$. We also have $t = y$, so $t$ has just one character, $t[1]$ differs from `(', `:', `,', `)' and $d(t,\ell(t)) = 0$.

\medskip

Let $k \in K(n)$ and $\mathbf{t \in E(n,k)}$, this means that $t \in E(n)$. In this case we just need to apply assumption~\ref{A:expr-cont-depth}.

\medskip

Let $k = h + <y,\phi> \in K(n) - \{ \epsilon \}$ and $t \in \mathbf{E_b(n+1,k)}$. We have $h \in K(n)$, $t \in E(n,h)$, so we can apply assumption~\ref{A:expr-cont-depth} to finish.\\

Let $k \in K(n)$ and $t \in \mathbf{E_Q(n+1,k)}$. This implies $t \in H_Q(n+1,k)$, where $H_Q(n+1,k)$ is the set of the strings $Q(x:\varphi,\phi)$ such that 

\begin{itemize}
\item $\varphi \in E_s(n,k)$,
\item $x \in \mathcal{V} - var(k)$,
\item if we define $k' = k + <x, \varphi>$ then $k' \in K(n)$ and $\phi \in S(n,k')$ .
\end{itemize}

Given $t \in H_Q(n+1,k)$ clearly it has the form $Q(x:\psi)$, with $\psi \in \Sigma^*$, so we can apply lemma~\ref{IL:HQ-recursive-aux1}. It follows that
\begin{itemize}
\item $x \in \mathcal{V} - var(k)$,
\item the set of the positive integers $r$ such that $4 < r < \ell(t)$, $t[r]$ = `,' and $d(t, r) = 1$ has just one member $r_1$.
\end{itemize}

Then we can define the following.\\

If $r_1 = 5$ then $\varphi = \epsilon$, else $\varphi = t[5, r_1 - 1]$.\\
If $r_1 = \ell(t) - 1$ then $\phi = \epsilon$ else $\phi = t[r_1 + 1, \ell(t) - 1]$.\\

And we have also
\begin{itemize}
\item $\varphi \in E_s(n,k)$,
\item if we define $k' = k + <x, \varphi>$ then $k' \in K(n)$ and $\phi \in S(n,k')$,
\item $t = Q(x:\varphi,\phi)$.
\end{itemize}

\medskip

We have
\begin{align*}
d(t, \ell(t) - 1) &= d(t, r_1 + \ell(\phi)) =\\
&= d(t, r_1 + 1) + d(\phi, \ell(\phi)) = 1 + d(\phi, \ell(\phi)).
\end{align*}

\medskip

If $t[\ell(t) - 1] = \phi[\ell(\phi)] = \text{`)'}$ then $d(t, \ell(t)) = d(t, \ell(t) - 1) - 1 = d(\phi, \ell(\phi)) = 1$.\\
Else $t[\ell(t) - 1] = \phi[\ell(\phi)] \notin \{ \text{`('}, \text{`)'} \}$ so $d(t, \ell(t)) = d(t, \ell(t) - 1) = 1 + d(\phi, \ell(\phi)) = 1$.\\

Let's now examine the facts we have to prove. It is true that $t[\ell(t)] \ne \text{`('}$. It's also true that $t[\ell(t)] = \text{`)'}$ and $d(t, \ell(t)) = 1$.

\medskip

Now let $\alpha \in \{ 1, \dots , \ell(t) \}$ and ( $t[\alpha]=\text{`:'}$ or $t[\alpha]=\text{`,'}$ or $t[\alpha]=\text{`)'}$ ).

\medskip

If $\alpha \in \{4, r_1, \ell(t) \}$ we have already shown that $d(t, \alpha) = 1$. Otherwise there are these alternative possibilities:

\begin{enumerate}[a.]
\item $4 < \alpha < r_1$,
\item $r_1 < \alpha < \ell(t)$.
\end{enumerate}

\medskip

In the situation a. we have
\begin{align} 
&4 < \alpha < r_1 \ , \notag  \\
&0 < \alpha - 4 < r_1 - 4  \ , \notag  \\
&1 \leqslant \alpha - 4 \leqslant r_1 - 4 - 1 = \ell(\varphi)  \ , \notag  \\
&\varphi[\alpha - 4 ] = t[\alpha]  \  \notag ,
\end{align}
\begin{align*}
d(t, \alpha) &= d( t, 4 + (\alpha - 4) ) = d(t, 4 + 1) + d(\varphi, \alpha - 4) =\\
&= 1 + d(\varphi, \alpha - 4) \geqslant 2.
\end{align*}

\medskip

In the situation b. we have
\begin{align} 
&r_1 < \alpha < \ell(t) \ , \notag  \\
&0 < \alpha - r_1 < \ell(t) - r_1 \ , \notag  \\
&1 \leqslant \alpha - r_1 \leqslant \ell(t) - r_1 - 1 = \ell(\phi) \ , \notag  \\
&\phi[\alpha - r_1] = t[\alpha] \  \notag ,
\end{align}
\begin{align*}
d(t, \alpha) &= d( t, r_1 + (\alpha - r_1) ) = d(t, r_1 + 1) + d(\phi, \alpha - r_1) =\\
&= 1 + d(\phi, \alpha - r_1) \geqslant 2.
\end{align*}

\bigskip

Let $k \in K(n)$, $c \in \mathcal{C}'$ and $t \in \mathbf{E^c(n+1,k)}$. Then $t \in H_c(n+1,k)$, so there exist $\varphi_1, \dots , \varphi_m$ in $E(n,k)$ such that $t = (c)( \varphi_1, \dots , \varphi_m)$.\\

In this representation of $t$ we see `explicit occurrences' of the symbols `(' ,  `)' and `,'. There are explicit occurrences of `,' only when $m>1$. The first explicit occurrence of `)' is in position $3$, and the second explicit occurrence of `)' is clearly in position $\ell(t)$. If $m>1$ we indicate with $q_1, \dots, q_{m-1}$  the positions of the explicit occurrences of `,'. By lemma~\ref{IL:Hc-recursive-aux0} we have that for each $j = 1 \dots m-1$ $d(t, q_j) = 1$.\\

Moreover we have $d(t,2) = 1$ and also $d(t,3) = 1$, $d(t,4) = 0$, $d(t,5) = 1$.\\

We now want to show that $d(t, \ell(t)) = 1$.

\medskip

If $m=1$ then 
\begin{align*}
d(t,\ell(t) - 1) &= d(t, 4 + \ell(\varphi_1) ) = d(t,4+1) + d(\varphi_1, \ell(\varphi_1)) = 1 + d(\varphi_1, \ell(\varphi_1)).
\end{align*}

\medskip

If $m>1$ then 
\[
d(t,\ell(t) - 1) = d(t, q_{m-1} + \ell(\varphi_m) ) = d(t, q_{m-1} + 1) + d(\varphi_m, \ell(\varphi_m)) = 1 + d(\varphi_m, \ell(\varphi_m)).
\]

\medskip

If $t[\ell(t) - 1] = \varphi_m[\ell(\varphi_m)] = \text{`)'}$ then $d(t, \ell(t)) = d(t, \ell(t) - 1) - 1 = d(\varphi_m, \ell(\varphi_m)) = 1$.\\
Else $t[\ell(t) - 1] = \varphi_m[\ell(\varphi_m)] \notin \{ \text{`('}, \text{`)'} \}$ so $d(t, \ell(t)) = d(t, \ell(t) - 1) = 1 + d(\varphi_m, \ell(\varphi_m)) = 1$.

\medskip

Let's now examine the facts we have to prove. It is true that $t[\ell(t)] \ne \text{`('}$. It's also true that $t[\ell(t)] = \text{`)'}$ and $d(t, \ell(t)) = 1$.

\medskip

Now let $\alpha \in \{ 1, \dots , \ell(t) \}$ and ( $t[\alpha]=\text{`:'}$ or $t[\alpha]=\text{`,'}$ or $t[\alpha]=\text{`)'}$ ). This implies $\alpha \notin \{1,2,4 \}$. 

\medskip

If $\alpha \in \{3, q_1, \dots , q_{m-1}, \ell(t) \}$ we have already shown that $d(t, \alpha) = 1$. Otherwise there are these alternative possibilities:

\begin{enumerate}[a.]
\item $(m=1) \wedge (\alpha > 4) \wedge (\alpha < \ell(t))$,
\item $(m > 1) \wedge (\alpha > 4) \wedge (\alpha < q_1)$,
\item $(m > 2) \wedge (\exists i = 1 \dots m-2: (\alpha > q_i) \wedge (\alpha < q_{i+1}))$,
\item $(m > 1) \wedge (\alpha > q_{m-1}) \wedge (\alpha < \ell(t))$.
\end{enumerate}

In the situation a. we have
\begin{align} 
&4 < \alpha < \ell(t) \ , \notag  \\
&0 < \alpha - 4 < \ell(t) - 4 \ , \notag  \\
&1 \leqslant \alpha - 4 \leqslant \ell(t) - 5 = \ell(\varphi_1) \ , \notag  \\
&\varphi_1[\alpha - 4] = t[\alpha] \  \notag ,
\end{align}
\begin{align*}
d(t, \alpha) &= d( t, 4 + (\alpha - 4) ) = d(t, 4+1) + d(\varphi_1, \alpha - 4) =\\
&= 1 + d(\varphi_1, \alpha - 4) \geqslant 2.
\end{align*}

\medskip

In the situation b. we have
\begin{align} 
&4 < \alpha < q_1 \ , \notag  \\
&0 < \alpha - 4 < q_1 - 4 \ , \notag  \\
&1 \leqslant \alpha - 4 \leqslant q_1 - 5 = \ell(\varphi_1) \ , \notag  \\
&\varphi_1[\alpha - 4] = t[\alpha] \  \notag ,
\end{align}
\begin{align*}
d(t, \alpha) &= d( t, 4 + (\alpha - 4) ) = d(t, 4+1) + d(\varphi_1, \alpha - 4) =\\
&= 1 + d(\varphi_1, \alpha - 4) \geqslant 2.
\end{align*}

\medskip

In the situation c. we have
\begin{align} 
&q_i < \alpha < q_{i+1} \ , \notag  \\
&0 < \alpha - q_i < q_{i+1} - q_i \ , \notag  \\
&1 \leqslant \alpha - q_i \leqslant q_{i+1} - q_i - 1 = \ell(\varphi_{i+1}) \ , \notag  \\
&\varphi_{i+1}[\alpha - q_i] = t[\alpha] \  \notag ,
\end{align}
\begin{align*}
d(t, \alpha) &= d( t, q_i + (\alpha - q_i) ) = d(t, q_i + 1) + d(\varphi_{i+1}, \alpha - q_i) =\\
&= 1 + d(\varphi_{i+1}, \alpha - q_i) \geqslant 2.
\end{align*}

\medskip

In the situation d. we have
\begin{align} 
&q_{m-1} < \alpha < \ell(t) \ , \notag  \\
&0 < \alpha - q_{m-1} < \ell(t) - q_{m-1} \ , \notag  \\
&1 \leqslant \alpha - q_{m-1} \leqslant \ell(t) - q_{m-1} - 1 = \ell(\varphi_m) \ , \notag  \\
&\varphi_m[\alpha - q_{m-1}] = t[\alpha] \  \notag ,
\end{align}
\begin{align*}
d(t, \alpha) &= d( t, q_{m-1} + (\alpha - q_{m-1}) ) = d(t, q_{m-1} + 1) + d(\varphi_m, \alpha - q_{m-1}) =\\
&= 1 + d(\varphi_m, \alpha - q_{m-1}) \geqslant 2.
\end{align*}

\bigskip

Let $k \in K(n)$, $f \in \mathcal{F}$ and $t \in \mathbf{E^f(n+1,k)}$. Then $t \in H_f(n+1,k)$, so if $f$ has multiplicity $2$ there exist $\varphi_1, \varphi_2 \in E(n,k)$ such that $t = f(\varphi_1, \varphi_2)$, if $f$ has multiplicity $1$ there exists $\varphi_1 \in E(n,k)$ such that $t = f(\varphi_1)$.\\

We first consider the case where $f$ has multiplicity $1$. Here we first want to show that $d(t, \ell(t)) = 1$.\\

We have
\begin{align*}
d(t,\ell(t) - 1) &= d(t, 2 + \ell(\varphi_1) ) = d(t,2+1) + d(\varphi_1, \ell(\varphi_1)) = 1 + d(\varphi_1, \ell(\varphi_1)).
\end{align*}

\medskip

If $t[\ell(t) - 1] = \varphi_1[\ell(\varphi_1)] = \text{`)'}$ then $d(t, \ell(t)) = d(t, \ell(t) - 1) - 1 = d(\varphi_1, \ell(\varphi_1)) = 1$.\\
Else $t[\ell(t) - 1] = \varphi_1[\ell(\varphi_1)] \notin \{ \text{`('}, \text{`)'} \}$ so $d(t, \ell(t)) = d(t, \ell(t) - 1) = 1 + d(\varphi_1, \ell(\varphi_1)) = 1$.\\

Let's now examine the facts we have to prove. It is true that $t[\ell(t)] \ne \text{`('}$. It's also true that $t[\ell(t)] = \text{`)'}$ and $d(t, \ell(t)) = 1$.\\

Now let $\alpha \in \{ 1, \dots , \ell(t) \}$ and ( $t[\alpha]=\text{`:'}$ or $t[\alpha]=\text{`,'}$ or $t[\alpha]=\text{`)'}$ ). This implies $\alpha \notin \{1,2 \}$.\\

If $\alpha = \ell(t)$ we have already shown that $d(t, \alpha) = 1$. Otherwise clearly $2 < \alpha < \ell(t)$ and
\begin{align} 
&0 < \alpha - 2 < \ell(t) - 2 \ , \notag  \\
&1 \leqslant \alpha - 2 \leqslant \ell(t) - 3 = \ell(\varphi_1) \ , \notag  \\
&\varphi_1[\alpha - 2] = t[\alpha] \  \notag ,
\end{align}
\begin{align*}
d(t, \alpha) &= d( t, 2 + (\alpha - 2) ) = d(t, 2+1) + d(\varphi_1, \alpha - 2) =\\
&= 1 + d(\varphi_1, \alpha - 2) \geqslant 2.
\end{align*}

\bigskip

Let's then consider the case where $f$ has multiplicity $2$. Here we indicate with $q_1$ the position of the explicit occurrence of `$,$' within $t$, where $t = f(\varphi_1, \varphi_2)$. First of all we want to prove that $d(t, q_1) = 1$. To this end we consider that 
\[
d(t, q_1 - 1) = d(t, 2+\ell(\varphi_1)) = d(t, 2+1) + d(\varphi_1, \ell(\varphi_1)) = 1 + d(\varphi_1, \ell(\varphi_1)) .
\]

\medskip

If $t[q_1 - 1] = \varphi_1[\ell(\varphi_1)] = \text{`)'}$ then $d(t,q_1) = d(t,q_1 - 1) - 1 = d(\varphi_1, \ell(\varphi_1)) = 1$ .\\
Else $t[q_1 - 1] = \varphi_1[\ell(\varphi_1)] \notin \{ \text{`('}, \text{`)'} \}$ so $d(t,q_1) = d(t,q_1 - 1) = 1 + d(\varphi_1, \ell(\varphi_1)) = 1$.\\

We then want to show that $d(t, \ell(t)) = 1$. We have
\[
d(t,\ell(t) - 1) = d(t, q_1 + \ell(\varphi_2) ) = d(t, q_1 + 1) + d(\varphi_2, \ell(\varphi_2)) = 1 + d(\varphi_2, \ell(\varphi_2)).
\]

\medskip

If $t[\ell(t) - 1] = \varphi_2[\ell(\varphi_2)] = \text{`)'}$ then $d(t, \ell(t)) = d(t, \ell(t) - 1) - 1 = d(\varphi_2, \ell(\varphi_2)) = 1$.\\
Else $t[\ell(t) - 1] = \varphi_2[\ell(\varphi_2)] \notin \{ \text{`('}, \text{`)'} \}$ so $d(t, \ell(t)) = d(t, \ell(t) - 1) = 1 + d(\varphi_2, \ell(\varphi_2)) = 1$.\\

Let's now examine the facts we have to prove. It is true that $t[\ell(t)] \ne \text{`('}$. It's also true that $t[\ell(t)] = \text{`)'}$ and $d(t, \ell(t)) = 1$.\\

Now let $\alpha \in \{ 1, \dots , \ell(t) \}$ and ( $t[\alpha]=\text{`:'}$ or $t[\alpha]=\text{`,'}$ or $t[\alpha]=\text{`)'}$ ). This implies $\alpha \notin \{1,2 \}$.\\

If $\alpha \in \{ q_1, \ell(t) \}$ we have already shown that $d(t, \alpha) = 1$. Otherwise there are these alternative possibilities:

\begin{enumerate}[a.]
\item $(\alpha > 2) \wedge (\alpha < q_1)$,
\item $(\alpha > q_1) \wedge (\alpha < \ell(t))$.
\end{enumerate} 

\medskip

In the situation a. we have 
\begin{align} 
&2 < \alpha < q_1 \ , \notag  \\
&0 < \alpha - 2 < q_1 - 2 \ , \notag  \\
&1 \leqslant \alpha - 2 \leqslant q_1 - 3 = \ell(\varphi_1) \ , \notag  \\
&\varphi_1[\alpha - 2] = t[\alpha] \  \notag ,
\end{align}
\begin{align*}
d(t, \alpha) &= d( t, 2 + (\alpha - 2) ) = d(t, 2+1) + d(\varphi_1, \alpha - 2) =\\
&= 1 + d(\varphi_1, \alpha - 2) \geqslant 2.
\end{align*}

\medskip

In the situation b. we have
\begin{align} 
&q_1 < \alpha < \ell(t) \ , \notag  \\
&0 < \alpha - q_1 < \ell(t) - q_1 \ , \notag  \\
&1 \leqslant \alpha - q_1 \leqslant \ell(t) - q_1 - 1 = \ell(\varphi_2) \ , \notag  \\
&\varphi_2[\alpha - q_1] = t[\alpha] \  \notag ,
\end{align}
\begin{align*}
d(t, \alpha) &= d( t, q_1 + (\alpha - q_1) ) = d(t, q_1 + 1) + d(\varphi_2, \alpha - q_1) =\\
&= 1 + d(\varphi_2, \alpha - q_1) \geqslant 2.
\end{align*}

\end{proof}

\end{definition}

\section{Deductive systems and proofs}\label{Ch:dedSysandProofs}

In this section we will define deductive systems and proofs and we will introduce other concepts and results related to our deductive methodology. Given a language $\mathcal{L}=(\mathcal{V}, \mathcal{F}, \mathcal{C}, \#, \{D_1, \dots, D_p \})$, we begin with some preliminary definitions.\\ 

Let $K = \bigcup_{n \geqslant 1} K(n)$.\\

For each $k \in K$ let 
\[ E(k) = \bigcup_{n \geqslant 1: k \in K(n)} E(n,k) \ , \]
\[ E_s(k) = \{ t | t \in E(k), \forall \sigma \in \Xi(k) \ \#(k,t,\sigma) \text{ is a set }  \} \ . \]

\bigskip

Let $E = \bigcup_{k \in K} E(k)$; $E$ is the set of all expressions in our language.\\

One expression $t \in E(k)$ is a `sentence with respect to $k$' when for each $\sigma \in \Xi(k)$ $\#(k,t,\sigma)$ is true or $\#(k,t,\sigma)$ is false.\\

We define $S(k) = \{ t | t \in E(k), t \text{ is a sentence with respect to } k \}$.\\

For each $t \in E(\epsilon)$ we define $\#(t) = \#(\epsilon, t , \epsilon)$.\\

A sentence with respect to $\epsilon$ will simply be called a `sentence'.\\

At this point we can define what is a proof in our language. To define this we need to define the notions of axiom and rule. We first notice that the symbols of our language belong to the four disjoint sets $\mathcal{V}$, $\mathcal{C}$, $\mathcal{F}$ and $\mathcal{Z}$. Let's call $\Sigma = \mathcal{V} \cup \mathcal{C} \cup \mathcal{F} \cup \mathcal{Z}$ the set of all the symbols (or alphabet) of our language and $\Sigma^*$ the set of all the empty or finite strings built with the symbols in $\Sigma$. Clearly given $k \in K$ $S(k) \subseteq E(k) \subseteq \Sigma^*$.\\

An \emph{axiom} is a set $A$ such that
\begin{itemize}
\item $A \subseteq S(\epsilon) \subseteq \Sigma^*$,
\item $A$ is r.e.,
\item for each $\varphi \in A$ $\#(\varphi)$ holds.
\end{itemize}

\medskip

The property `for each $\varphi \in A$ $\#(\varphi)$ holds' states that axiom $A$ is `sound'.\\

Given a positive integer $n$ we indicate with $S(\epsilon)^n$ the set of all $n$-tuples $(\varphi_1, \dots , \varphi_n)$ for $\varphi_1, \dots , \varphi_n \in S(\epsilon)$. An $n$-ary \emph{rule} is a set $R$ such that
\begin{itemize}
\item $R \subseteq S(\epsilon)^{n+1} \subseteq (\Sigma^*)^{n+1}$,
\item $R$ is r.e.,
\item for each $(\varphi_1, \dots , \varphi_n, \varphi) \in R$ if $\#(\varphi_1), \dots , \#(\varphi_n)$ hold then $\#(\varphi)$ holds.
\end{itemize}

\medskip

The property `for each $(\varphi_1, \dots , \varphi_n, \varphi) \in R$ if $\#(\varphi_1), \dots , \#(\varphi_n)$ hold then $\#(\varphi)$ holds' states that rule $R$ is `sound'.\\

\medskip

Both in the definition of axiom and rule we have included a requirement of soundness.\\

A deductive system is built on top of our language $\mathcal{L}$, and is identified by a pair $(\mathcal{A}, \mathcal{R})$ where $\mathcal{A}$ is a finite set of axioms in $\mathcal{L}$ and $\mathcal{R}$ is a finite set of rules in $\mathcal{L}$.\\

We require that the set of the axioms and the set of the rules are finite since we need to be able to list each of them on a piece of paper.\\

Given a language $\mathcal{L}$, $\mathcal{D} = (\mathcal{A}, \mathcal{R})$ deductive system in $\mathcal{L}$, $\varphi$, $\psi_1, \dots, \psi_m$ sentences in $\mathcal{L}$, we say that $(\psi_1, \dots, \psi_m)$ is a \emph{proof} of $\varphi$ in $\mathcal{D}$ if and only if
\begin{itemize}
\item there exists $A \in \mathcal{A}$ such that $\psi_1 \in A$;
\item if $m>1$ then for each $j = 2 \dots m$ one of the following holds
\begin{itemize}
\item there exists $A \in \mathcal{A}$ such that $\psi_j \in A$,
\item there exist an $n$-ary rule $R \in \mathcal{R}$ and $i_1, \dots , i_n < j$ such that\\
$(\psi_{i_1}, \dots, \psi_{i_n}, \psi_j) \in R$;
\end{itemize}
\item $\psi_m = \varphi$.
\end{itemize}

\medskip

Given $\mathcal{D} = (\mathcal{A}, \mathcal{R})$ deductive system in $\mathcal{L}$ and $\varphi$ sentence in $\mathcal{L}$ we say that $\varphi$ \emph{is derivable in} $\mathcal{D}$ and write $\vdash_{\mathcal{D}} \varphi$ if and only if there exist $\psi_1, \dots , \psi_m$ sentences in $\mathcal{L}$ such that $(\psi_1, \dots, \psi_m)$ is a proof of $\varphi$ in $\mathcal{D}$.\\

A deductive system $\mathcal{D} = (\mathcal{A}, \mathcal{R})$ is said to be \emph{sound} if and only if for each $\varphi$ sentence in $\mathcal{L}$ if $\vdash_{\mathcal{D}} \varphi$ then $\#(\varphi)$ holds. In the next lemma we easily prove that each of our systems is sound.\\

\begin{lemma}\label{L:ds-sound}
Let $\mathcal{D} = (\mathcal{A}, \mathcal{R})$ be a deductive system in $\mathcal{L}$. Then $\mathcal{D}$ is sound.
\end{lemma}

\begin{proof}

\smallskip

Let $\varphi$ be a sentence in $\mathcal{L}$. Suppose $\vdash_{\mathcal{D}} \varphi$. There exist $\psi_1, \dots , \psi_m$ sentences in $\mathcal{L}$ such that $(\psi_1, \dots, \psi_m)$ is a proof of $\varphi$ in $\mathcal{D}$. We can show that for each $j = 1 \dots m$ $\#(\psi_j)$ holds.

\medskip

There exists $A \in \mathcal{A}$ such that $\psi_1 \in A$, so $\#(\psi_1)$ holds.

\medskip

If $m>1$ suppose $j = 2 \dots m$. We assume for each $i = 1 \dots j-1$ $\#(\varphi_i)$ holds.

\medskip

If there exists $A \in \mathcal{A}$ such that $\psi_j \in A$ then $\#(\psi_j)$ holds. 

\medskip

Otherwise there exist an $n$-ary rule $R \in \mathcal{R}$ and $i_1, \dots , i_n < j$ such that 
\[ (\psi_{i_1}, \dots, \psi_{i_n}, \psi_j) \in R \ . \]

\smallskip

Since $\#(\psi_{i_1}), \dots, \#(\psi_{i_n})$ all hold, then $\#(\psi_j)$ also holds.
\end{proof}

\medskip

We now want to point out some recursivity requirements with respect to the sets that we defined above: $E(k)$, $S(k)$, $E_s(k)$. We will prove these sets are recursively enumerable.\\

For each $k \in K$ we defined $E(k) = \bigcup_{n \geqslant 1: k \in K(n)} E(n,k)$.\\

The set $\{ n | n \in \mathbb{N}, n \geqslant 1, k \in K(n) \}$ is r.e.. In fact if we call $n_0$ the least $n \in \mathbb{N}$ such that $k \in K(n)$ we have that the just mentioned set is actually $\{ n | n \in \mathbb{N}, n \geqslant n_0 \}$, that is a recursive and r.e. set. Since for each $n$ in the mentioned r.e. set $E(n,k)$ is r.e. then $E(k)$ is also r.e..\\

Given a positive integer $n$ and $k \in K(n)$, let's define the following sets:
\begin{align*}
S(n,k) &= \{ \varphi | \ \varphi \in E(n,k) \text{, for each } \sigma \in \Xi(k) \ \#(k,\varphi,\sigma) \text{ is true or } \#(k,\varphi,\sigma) \text{ is false } \};\\
E_s(n,k) &= \{ \varphi | \ \varphi \in E(n,k) \text{, for each } \sigma \in \Xi(k) \ \#(k,\varphi,\sigma) \text{ is a set } \};\\
E_{D_i}(n,k) &= \{ \varphi | \ \varphi \in E(n,k) \text{, for each } \sigma \in \Xi(k) \ \#(k,\varphi,\sigma) \in D_i \}.
\end{align*}

For each $\varphi \in E(n,k)$ we can decide the following conditions:
\begin{itemize}
\item for each $\sigma \in \Xi(k)$ $\#(k, \varphi, \sigma)$ is true or false;
\item for each $\sigma \in \Xi(k)$ $\#(k, \varphi, \sigma)$ is a set;
\item for each $\sigma \in \Xi(k)$ $\#(k,\varphi,\sigma) \in D_i$.
\end{itemize} 

\smallskip

Therefore, since $E(n,k)$ is recursive, $S(n,k)$, $E_s(n,k)$ and $E_{D_i}(n,k)$ are recursive too.\\

It is easy to verify that $S(k) = \bigcup_{n \geqslant 1: k \in K(n)} S(n,k)$, therefore $S(k)$ is r.e..\\

Similarly, it is easy to verify that $E_s(k) = \bigcup_{n \geqslant 1: k \in K(n)} E_s(n,k)$, therefore $E_s(k)$ is r.e..\\

Moreover we can define $E_{D_i}(k) = \{ \varphi | \ \varphi \in E(k) \text{, for each } \sigma \in \Xi(k) \ \#(k,\varphi,\sigma) \in D_i \}$\\

Then it is easy to verify that $E_{D_i}(k) = \bigcup_{n \geqslant 1: k \in K(n)} E_{D_i}(n,k)$, therefore $E_{D_i}(k)$ is r.e..\\

We now want to introduce some notation that will help us when we define our axioms and rules and when we want to prove that we are dealing with recursively enumerable sets.\\

\begin{definition}
Let $x \in \mathcal{V}, \ \varphi \in E$. We define \[ H[x: \varphi] = \varphi \in E_s(\epsilon) \ . \]

\smallskip

If the condition $H[x: \varphi]$ holds then we define $k[x:\varphi] = \epsilon + <x, \varphi>$. Clearly $k[x:\varphi] \in K$. In fact there exists $n$ positive integer such that $\epsilon \in K(n) \wedge \varphi \in E_s(n,\epsilon)$, $x \in \mathcal{V} - var(\epsilon)$, so $k[x:\varphi] = \epsilon + <x, \varphi> \in K(n) \cup K(n)^+ = K(n+1) \subseteq K$.

\medskip

Moreover $k[x:\varphi] = <<x, \varphi>>$ and $var(k[x:\varphi]) = \{ x \}$.

\medskip

Let $m$ be a positive integer. Let $x_1, \dots , x_{m+1} \in \mathcal{V}$, with $x_i \ne x_j$ for $i \ne j$. Let $\varphi_1, \dots , \varphi_{m+1} \in E$. We can assume to have defined $H[x_1:\varphi_1, \dots , x_m:\varphi_m ]$ and if this holds to have defined also $k[x_1:\varphi_1, \dots , x_m:\varphi_m ] \in K$, such that 
\begin{align*}
k[x_1:\varphi_1, \dots , x_m:\varphi_m ] &= <<x_1,\varphi_1> \dots  <x_m:\varphi_m>>\\
var(k[x_1:\varphi_1, \dots , x_m:\varphi_m ]) &= \{ x_1, \dots , x_m \} 
\end{align*}

\smallskip

We define
\begin{multline*} 
H[x_1:\varphi_1, \dots , x_{m+1}:\varphi_{m+1} ] = H[x_1:\varphi_1, \dots , x_m:\varphi_m ]\\
\wedge \varphi_{m+1} \in E_s(k[x_1:\varphi_1, \dots , x_m:\varphi_m ]) \ .
\end{multline*}

\smallskip

If $H[x_1:\varphi_1, \dots , x_{m+1}:\varphi_{m+1} ]$ then we define 
\[ k[x_1:\varphi_1, \dots , x_{m+1}:\varphi_{m+1} ] = k[x_1:\varphi_1, \dots , x_m:\varphi_m ] + <x_{m+1}, \varphi_{m+1}> \ . \]

\smallskip

Clearly $k[x_1:\varphi_1, \dots , x_{m+1}:\varphi_{m+1} ] \in K$. In fact there exists a positive integer $n$ such that $k[x_1:\varphi_1, \dots , x_m:\varphi_m ] \in K(n)$ and
$\varphi_{m+1} \in E_s(n,k[x_1:\varphi_1, \dots , x_m:\varphi_m ])$, $x_{m+1} \in \mathcal{V} - var(k[x_1:\varphi_1, \dots , x_m:\varphi_m ])$, so $k[x_1:\varphi_1, \dots , x_{m+1}:\varphi_{m+1} ] \in K(n) \cup K(n)^+ = K(n+1)$.

\medskip
 
Moreover
\begin{align*}
k[x_1:\varphi_1, \dots , x_{m+1}:\varphi_{m+1} ] &= <<x_1,\varphi_1> \dots  <x_{m+1}:\varphi_{m+1}>>\\
var(k[x_1:\varphi_1, \dots , x_{m+1}:\varphi_{m+1} ]) &= \{ x_1, \dots , x_{m+1} \} \ .
\end{align*}
\end{definition}

\bigskip

\begin{lemma}
Let m positive integer, $x_1, \dots , x_m \in \mathcal{V}$, with $x_i \ne x_j$ for $i \ne j$, $\varphi_1, \dots , \varphi_m \in E$. Then $H[x_1:\varphi_1, \dots , x_m:\varphi_m ]$ is defined and if $H[x_1:\varphi_1, \dots , x_m:\varphi_m ]$ holds then $k[x_1:\varphi_1, \dots , x_m:\varphi_m ]$ is also defined and belongs to $K$. Moreover
\[ var( k[x_1:\varphi_1, \dots , x_m:\varphi_m ] ) = \{ x_1, \dots , x_m \}  \ . \]
\end{lemma}

\begin{proof}
This is an obvious consequence of the previous definition and has been verified, by induction on $m$, in the definition itself.
\end{proof}

\medskip

\begin{remark}\label{R:H[]-k[]}
Let $m$ be a positive integer. Let $x_1, \dots , x_m \in \mathcal{V}$, with $x_i \ne x_j$ for $i \ne j$. Let $\varphi_1, \dots , \varphi_m \in E$ and assume $H[x_1:\varphi_1, \dots , x_m:\varphi_m ]$. In these assumptions we can easily see that for each $i = 1 \dots m$ $H[x_1:\varphi_1, \dots , x_i:\varphi_i ]$ holds and so 

\begin{itemize}
\item $k[x_1:\varphi_1, \dots , x_i:\varphi_i ] \ = \ <<x_1,\varphi_1> \dots <x_i:\varphi_i>> \in K$,
\item $var(k[x_1:\varphi_1, \dots , x_i:\varphi_i]) = \{ x_1, \dots, x_i \}$,
\item $dom(k[x_1:\varphi_1, \dots , x_i:\varphi_i]) = \{ 1, \dots , i \}$.
\end{itemize}

\medskip

In fact this is clearly true for $i = m$. Given $i = 1 \dots m-1$, if we suppose this is true for $i+1$, then we have $H[x_1:\varphi_1, \dots , x_i:\varphi_i ]$, and so the remaining facts also hold.

\medskip

In these assumptions we can define $k_0 = \epsilon$ and for each $i = 1 \dots m$ \linebreak $k_i = k[x_1:\varphi_1, \dots , x_i:\varphi_i]$. We have $k_0 \in K$, $var(k_0) = \emptyset$, for each $i = 1 \dots m$ $k_i \in K$, $var(k_i) = \{ x_1, \dots, x_i \}$ and $dom(k_i) = \{ 1, \dots , i \}$, $\varphi_i \in E_s(k_{i-1})$, $k_i = k_{i-1} + <x_i, \varphi_i>$.\\

Hereafter we'll often use this kind of simplified notation.
\end{remark}

\bigskip

\begin{definition}\label{D:gamma-simple}
Let $m$ be a positive integer. Let $\varphi_1, \dots , \varphi_m \in \Sigma^*$. Let $\psi_1, \dots , \psi_m \in \Sigma^*$. Let $\varphi \in \Sigma^*$. Define
\[ \gamma[\psi_m: \varphi_m, \varphi] = \forall(\psi_m:\varphi_m, \varphi ) \ .  \]

\medskip

If $m > 1$ for each $i = 2 \dots m$ suppose we have defined $\gamma[\psi_i: \varphi_i, \dots , \psi_m: \varphi_m, \varphi]$  and define
\[ \gamma[\psi_{i-1}: \varphi_{i-1}, \dots , \psi_m: \varphi_m, \varphi] = \forall( \psi_{i-1}: \varphi_{i-1}, \gamma[\psi_i: \varphi_i, \dots , \psi_m: \varphi_m, \varphi] ) \ . \]

\medskip

With this we have also defined $\gamma[\psi_1: \varphi_1, \dots , \psi_m: \varphi_m, \varphi]$.

\end{definition}

\begin{flushright}
$\qed$
\end{flushright}

\medskip

We can define a function $\chi$ on the domain $(\Sigma^*)^{2m} \times \Sigma^*$ such that given $(\psi_1, \varphi_1, \dots , \psi_m, \varphi_m) \in (\Sigma^*)^{2m}$ and $\varphi \in \Sigma^*$, 
\[ \chi((\psi_1, \varphi_1, \dots , \psi_m, \varphi_m), \varphi) = \gamma[\psi_1: \varphi_1, \dots , \psi_m: \varphi_m, \varphi] \ . \] Clearly this function is computable since the result can be obtained by simply concatenating the elements of the input with other symbols of our language.\\

We can also observe the following.

\begin{lemma}\label{L:gamma-expand-last}
Let $m$ be a positive integer, $m > 1$. Let $\varphi_1, \dots , \varphi_m \in \Sigma^*$. Let $\psi_1, \dots , \psi_m \in \Sigma^*$. Let $\varphi \in \Sigma^*$. Then
\[ \gamma[\psi_1: \varphi_1, \dots , \psi_m: \varphi_m, \varphi] = \gamma[\psi_1: \varphi_1, \dots , \psi_{m-1}: \varphi_{m-1}, \forall( \psi_m:\varphi_m, \varphi )] \ . \]
\end{lemma}

\begin{proof}
We want to prove that for each $i = 1 \dots m-1$ 
\[ \gamma[\psi_i: \varphi_i, \dots , \psi_m: \varphi_m, \varphi] = \gamma[\psi_i: \varphi_i, \dots , \psi_{m-1}: \varphi_{m-1}, \forall(\psi_m:\varphi_m, \varphi )] \ . \]

We start the proof at $m-1$ and we are then going backwards by induction to $1$. So
\begin{align*}
\gamma[\psi_{m-1}: \varphi_{m-1}, \psi_m: \varphi_m, \varphi] &= \forall( \psi_{m-1}: \varphi_{m-1}, \gamma[\psi_m: \varphi_m, \varphi] )\\
&= \forall( \psi_{m-1}: \varphi_{m-1}, \forall( \psi_m: \varphi_m, \varphi) )\\
&= \gamma[\psi_{m-1}: \varphi_{m-1}, \forall(\psi_m:\varphi_m, \varphi ) ] 
\end{align*}

If $m = 2$ our proof is finished, whilst if $m > 2$ given $i = 2 \dots m-1$ we can assume 
\[ \gamma[\psi_i: \varphi_i, \dots , \psi_m: \varphi_m, \varphi] = \gamma[\psi_i: \varphi_i, \dots , \psi_{m-1}: \varphi_{m-1}, \forall(\psi_m:\varphi_m, \varphi ) ] \ . \]

And in this case we have 
\begin{align*}
\gamma[\psi_{i-1}: \varphi_{i-1}, \dots , \psi_m: \varphi_m, \varphi] &= \forall( \psi_{i-1}: \varphi_{i-1}, \gamma[\psi_i: \varphi_i, \dots , \psi_m: \varphi_m, \varphi] ) \\
&= \forall( \psi_{i-1}: \varphi_{i-1}, \gamma[\psi_i: \varphi_i, \dots , \psi_{m-1}: \varphi_{m-1}, \forall( \psi_m:\varphi_m, \varphi )] )  \\
&= \gamma[\psi_{i-1}: \varphi_{i-1}, \dots , \psi_{m-1}: \varphi_{m-1}, \forall( \psi_m:\varphi_m, \varphi )] \ .
\end{align*}

\end{proof}

\medskip

Given a positive integer $m$ let's call $R_m$ the set 
\scriptsize
\[ \{ (x_1, \varphi_1, \dots , x_m, \varphi_m) | \  x_1, \dots , x_m \in \mathcal{V} \text{ with } x_i \ne x_j \text{ for } i \ne j, \ \varphi_1, \dots , \varphi_m \in E, H[x_1:\varphi_1, \dots , x_m:\varphi_m ] \} \ . \]
\normalsize

\medskip

Clearly given $(x_1, \varphi_1, \dots , x_m, \varphi_m) \in R_m$ $k[x_1:\varphi_1, \dots , x_m:\varphi_m ] \in K$.\\

Let's also define
\[ Q_m = \bigcup_{(x_1, \varphi_1, \dots , x_m, \varphi_m) \in R_m} \{(x_1, \varphi_1, \dots , x_m, \varphi_m)\} \times S(k[x_1:\varphi_1, \dots , x_m:\varphi_m ]) \ . \]

Actually $Q_m \subseteq (\Sigma^*)^{2m} \times \Sigma^*$. Our goal is now to show that $Q_m$ is r.e. in order to be able to show that the set $\{ \chi((x_1, \varphi_1, \dots , x_m, \varphi_m), \varphi) | \, ((x_1, \varphi_1, \dots , x_m, \varphi_m), \varphi) \in Q_m \}$ is r.e. itself.\\

As a first remark in this proof we can notice that our set of variables $\mathcal{V}$ is recursive. In fact given a string $\varphi \in \Sigma^*$ if $\varphi$ has not exactly one character then it doesn't belong to $\mathcal{V}$. If it has just one character then, since apart from the variables our alphabet has a finite number of symbols, we can decide if $\varphi \in \mathcal{V}$.\\

The first step in this proof is to show that $R_m$ is r.e., i.e. the following lemma:\\

\begin{lemma}
For each $m$ positive integer $R_m$ is r.e..
\end{lemma}

\begin{proof}

In the initial step of the proof we have to show that $R_1$ is r.e.. We have 
\begin{align*}
R_1 &= \{ (x_1, \varphi_1) | \, x_1 \in \mathcal{V}, \varphi_1 \in E, H[x_1: \varphi_1] \}\\
&= \{ (x_1, \varphi_1) | \, x_1 \in \mathcal{V}, \varphi_1 \in E_s(\epsilon) \}\\
&= \mathcal{V} \times E_s(\epsilon)
\end{align*}

and since both $\mathcal{V}$ and $E_s(\epsilon)$ are r.e. then $R_1$ is r.e..\\

Given a positive integer $m$ we assume $R_m$ is r.e. and want to show that $R_{m+1}$ is r.e..\\

Actually 
\scriptsize
\begin{align*}
R_{m+1} &= \{ (x_1, \varphi_1, \dots , x_{m+1}, \varphi_{m+1}) | \  x_1, \dots , x_{m+1} \in \mathcal{V} \text{ with } x_i \ne x_j \text{ for } i \ne j, \ \varphi_1, \dots , \varphi_{m+1} \in E, H[x_1:\varphi_1, \dots , x_{m+1}:\varphi_{m+1} ] \}\\
&= \{ (x_1, \varphi_1, \dots , x_{m+1}, \varphi_{m+1}) | \  x_1, \dots , x_{m+1} \in \mathcal{V} \text{ with } x_i \ne x_j \text{ for } i \ne j, \ \varphi_1, \dots , \varphi_{m+1} \in E,\\ 
& \ H[x_1:\varphi_1, \dots , x_m:\varphi_m ] \wedge \varphi_{m+1} \in E_s(k[x_1:\varphi_1, \dots , x_m:\varphi_m ]) \}\\
&= \{ (x_1, \varphi_1, \dots , x_{m+1}, \varphi_{m+1}) | \ (x_1, \varphi_1, \dots , x_m, \varphi_m) \in R_m, x_{m+1} \in \mathcal{V} - \{x_1, \dots , x_m \}, \varphi_{m+1} \in E_s(k[x_1:\varphi_1, \dots , x_m:\varphi_m ]) \} \, .\\
\end{align*}
\normalsize

Let's now consider the set 
\[ U_{m+1} = \bigcup_{(x_1, \varphi_1, \dots , x_m, \varphi_m) \in R_m} \{(x_1, \varphi_1, \dots , x_m, \varphi_m)\} \times (\mathcal{V} - \{x_1, \dots , x_m \}) \times E_s(k[x_1:\varphi_1, \dots , x_m:\varphi_m ]) \ . \]

Given $(x_1, \varphi_1, \dots , x_m, \varphi_m) \in R_m$ the sets $\{(x_1, \varphi_1, \dots , x_m, \varphi_m)\} \subseteq (\Sigma^*)^{2m}$, $(\mathcal{V} - \{x_1, \dots , x_m \}) \subseteq \Sigma^*$ and $E_s(k[x_1:\varphi_1, \dots , x_m:\varphi_m ]) \subseteq \Sigma^*$ are r.e., so the cartesian product 
\[ \{(x_1, \varphi_1, \dots , x_m, \varphi_m)\} \times (\mathcal{V} - \{x_1, \dots , x_m \}) \times E_s(k[x_1:\varphi_1, \dots , x_m:\varphi_m ]) \]
is also r.e., and $U_{m+1} \subseteq (\Sigma^*)^{2m} \times \Sigma^* \times \Sigma^*$ is r.e..\\

The set $R_{m+1}$ is a subset of $(\Sigma^*)^{2m+2}$ which is not (necessarily) the same of $(\Sigma^*)^{2m} \times \Sigma^* \times \Sigma^*$. In fact a member of $(\Sigma^*)^{2m+2}$ can be expressed as $(\psi_1, \varphi_1, \dots , \psi_{m+1}, \varphi_{m+1})$ and a member of $(\Sigma^*)^{2m} \times \Sigma^* \times \Sigma^*$ can  be expressed as $((\psi_1, \varphi_1, \dots , \psi_m, \varphi_m), \psi_{m+1}, \varphi_{m+1})$. Anyway we can easily map members of the first set to the ones of the second set and vice-versa. In fact we can define a function $\kappa$ over $(\Sigma^*)^{2m+2}$ such that $\kappa(\psi_1, \varphi_1, \dots , \psi_{m+1}, \varphi_{m+1}) = ((\psi_1, \varphi_1, \dots , \psi_m, \varphi_m), \psi_{m+1}, \varphi_{m+1})$, and the function $\kappa$ is computable.\\

Given $(\psi_1, \varphi_1, \dots , \psi_{m+1}, \varphi_{m+1}) \in (\Sigma^*)^{2m+2}$ if $(\psi_1, \varphi_1, \dots , \psi_{m+1}, \varphi_{m+1}) \in R_{m+1}$ then $\kappa(\psi_1, \varphi_1, \dots , \psi_{m+1}, \varphi_{m+1}) \in U_{m+1}$ and vice-versa if $\kappa(\psi_1, \varphi_1, \dots , \psi_{m+1}, \varphi_{m+1}) \in U_{m+1}$ then $(\psi_1, \varphi_1, \dots , \psi_{m+1}, \varphi_{m+1}) \in R_{m+1}$.\\

As we have seen $U_{m+1}$ is r.e. so its semi-characteristic function $s_U$ is computable. Let's now consider the  function $s_U \circ \kappa$ which is defined over $(\Sigma^*)^{2m+2}$. Given $(\psi_1, \varphi_1, \dots , \psi_{m+1}, \varphi_{m+1}) \in (\Sigma^*)^{2m+2}$ if $(\psi_1, \varphi_1, \dots , \psi_{m+1}, \varphi_{m+1}) \in R_{m+1}$ then $\kappa(\psi_1, \varphi_1, \dots , \psi_{m+1}, \varphi_{m+1}) \in U_{m+1}$ and $s_U(\kappa(\psi_1, \varphi_1, \dots , \psi_{m+1}, \varphi_{m+1}) = 1$. If $(\psi_1, \varphi_1, \dots , \psi_{m+1}, \varphi_{m+1}) \notin R_{m+1}$ then $\kappa(\psi_1, \varphi_1, \dots , \psi_{m+1}, \varphi_{m+1}) \notin U_{m+1}$ and $s_U(\kappa(\psi_1, \varphi_1, \dots , \psi_{m+1}, \varphi_{m+1})$ diverges. So $s_U \circ \kappa$ is actually the semicharacteristic function of $R_{m+1}$ and it is clearly a computable function. This proves that $R_{m+1}$ is r.e..
\end{proof}

Now given $(x_1, \varphi_1, \dots , x_m, \varphi_m) \in R_m$ both $\{ (x_1, \varphi_1, \dots , x_m, \varphi_m) \}$ and $S(k[x_1:\varphi_1, \dots , x_m:\varphi_m ])$ are r.e., so $\{ (x_1, \varphi_1, \dots , x_m, \varphi_m) \} \times S(k[x_1:\varphi_1, \dots , x_m:\varphi_m ])$ is r.e. too, and so $Q_m$ is a r.e. subset of $(\Sigma^*)^{2m} \times \Sigma^*$.\\

We can now recall that we have defined a computable function $\chi : (\Sigma^*)^{2m} \times \Sigma^* \to \Sigma^*$. Because of lemma~\ref{L:fre-is-re} we have that the set $\{ \chi((x_1, \varphi_1, \dots , x_m, \varphi_m), \varphi) | \, ((x_1, \varphi_1, \dots , x_m, \varphi_m), \varphi) \in Q_m \}$  is a r.e. subset of $\Sigma^*$.\\

And finally the set 
\[ \bigcup_{m \geqslant 1} \{ \chi((x_1, \varphi_1, \dots , x_m, \varphi_m), \varphi) | \, ((x_1, \varphi_1, \dots , x_m, \varphi_m), \varphi) \in Q_m \} \]

is itself a r.e. set. It seems this is not particularly significant to us because this set is not an axiom, but we'll see sets that are very similar to this one and that we can use as an axiom in our deductive system.\\

\section{Deductive methodology}\label{Ch:dedMet}

We now need to introduce some other fundamental notions and results relevant to our deductive methodology.\\

At the beginning of section~\ref{Ch:lang} we have introduced the logical connectives. In our deductions, expressions will make an extensive use of the logical connectives, so we assume that all of these symbols: $\neg, \wedge, \vee, \to, \leftrightarrow$ are in our set $\mathcal{F}$. For each of these operators $f$ $A_f(x_1, \dots ,x_n)$ and $P_f(x_1, \dots ,x_n)$ are defined as specified at the beginning of section~\ref{Ch:lang}.

\medskip

\begin{lemma}\label{L:k-in-km+}
For each $n$ positive integer such that $n \geqslant 2$, $k \in K(n)$: $k \ne \epsilon$ there exists $m < n$ such that $k \in K(m)^+$.
\end{lemma}

\begin{proof}

\smallskip

We prove this by induction on $n$. Clearly if $k \in K(2)$ and $k \ne \epsilon$ then $k \in K(1)^+$.

\medskip

Let $n \geqslant 2$, $k \in K(n+1)$:  $k \ne \epsilon$. Clearly if $k \in K(n)^+$ our proof is finished. Otherwise $k \in K(n)$ and in this case we can apply the inductive hypothesis.
\end{proof}

\medskip

\begin{lemma}\label{L:k-in-km+ext}
For each $n$ positive integer such that $n \geqslant 2$, $k \in K(n)$: $k \ne \epsilon$
\begin{itemize}
\item there exist $m < n$, $h \in K(m), \phi \in E_s(m,h), y \in (\mathcal{V}-var(h))$ such that \linebreak $k = h + <y,\phi>$, $\Xi(k) = \{ \sigma + (y,s) | \, \sigma \in \Xi(h), s \in \#(h,\phi,\sigma) \}$;
\item for each $g \in K(n), \theta \in E_s(n,g), z \in (\mathcal{V}-var(g))$ such that $k = g + <z,\theta>$ \linebreak $\Xi(k) = \{ \sigma + (z,s) | \, \sigma \in \Xi(g), s \in \#(g,\theta,\sigma) \}$.
\end{itemize}
\end{lemma}

\begin{proof}

\smallskip

The first part clearly follows from lemma~\ref{L:k-in-km+}. The second part holds because we have $g = h, z = y, \theta = \phi$.
\end{proof}

\medskip

\begin{lemma}\label{L:rest-of-context-is-context}
For each $n$ positive integer such that $n \geqslant 2$, $k \in K(n): k \ne \epsilon$, $h \in \mathcal{R}(k): h \ne k$ there exists $m < n$ such that $h \in K(m)$.
\end{lemma}

\begin{proof}

\smallskip

We prove this by induction on $n$. Let $k \in K(2)$: $k \ne \epsilon$, $h \in \mathcal{R}(k): h \ne k$. There exist $m < n$, $g \in K(m), \phi \in E_s(m,g), y \in (\mathcal{V}-var(g))$ such that \linebreak $k = g + (y,\phi)$. In this case $m = 1$, so $g = \epsilon$. By lemma~\ref{L:useful-contexts-1} we have $h \in \mathcal{R}(\epsilon)$ and so $h = \epsilon \in K(1)$.

\medskip

In order  to perform the inductive step, let $k \in K(n+1): k \ne \epsilon$, $h \in \mathcal{R}(k)$ such that $h \ne k$. There exist $m < n+1$, $g \in K(m), \phi \in E_s(m,g), y \in (\mathcal{V}-var(g))$ such that $k = g + (y,\phi)$. By lemma~\ref{L:useful-contexts-1} we have $h \in \mathcal{R}(g)$. If $h = g \in K(m)$ our proof is finished. Otherwise $h \ne g$ and $g \ne \epsilon$, we can apply our inductive hypothesis and obtain that there exists $q < m < n+1$ such that $h \in K(q)$.
\end{proof}

\medskip

\begin{lemma}\label{L:rest-of-state-is_state}
For each $n$ positive integer such that $n \geqslant 2$, $k \in K(n): k \ne \epsilon$, $\sigma \in \Xi(k)$, $h \in \mathcal{R}(k): h \ne k$, there exists $m < n$ such that $h \in K(m)$ and it results $\sigma_{/dom(h)} \in \Xi(h)$.
\end{lemma}

\begin{proof}

\smallskip

We prove this by induction on $n$. Let $k \in K(2)$: $k \ne \epsilon$, $\sigma \in \Xi(k)$, $h \in \mathcal{R}(k)$ such that $h \ne k$. Clearly $k \in K(1)^+$, so there exist $g \in K(1)$, $\phi \in E_s(1,g)$, $y \in \mathcal{V} - var(g)$ such that $k = g + <y, \phi>$. By lemma~\ref{L:useful-contexts-1} we obtain that $h \in \mathcal{R}(g)$. Since $g = \epsilon$ then also $h = \epsilon \in K(1)$ , so $\sigma_{/dom(h)} = \sigma_{/\emptyset} = \epsilon \in \Xi(\epsilon) = \Xi(h)$. 

\medskip

In order  to perform the inductive step, let $k \in K(n+1): k \ne \epsilon$, $\sigma \in \Xi(k)$, $h \in \mathcal{R}(k)$ such that $h \ne k$ . By lemma~\ref{L:k-in-km+} there exists $m \leqslant n$ such that $k \in K(m)^+$. Then there exist $g \in K(m), \phi \in E_s(m,g), y \in (\mathcal{V}-var(g))$ such that $k = g + <y,\phi>$. Moreover
\[
\Xi(k) = \{ \rho + (y,s) | \, \rho \in \Xi(g), s \in \#(g,\phi,\rho) \} \, .
\]

\smallskip

Therefore there exist $\rho \in \Xi(g), s \in \#(g,\phi,\rho)$ such that $\sigma = \rho + (y,s)$. By assumption~\ref{A:simple-k-sigma-slp-new} and lemma~\ref{L:useful4-slp} we have that $\sigma_{/dom(g)} = \sigma_{/dom(\rho)} = \rho$.

\medskip

If $h = g$ then $\sigma_{/dom(h)} = \sigma_{/dom(g)} = \rho \in \Xi(h)$. 

\medskip

Otherwise we have $h \ne g$. Since $k = g + <y,\phi>$, $h \in \mathcal{R}(k)$, $h \ne k$ by lemma~\ref{L:useful-contexts-1} we have that $h \in \mathcal{R}(g)$. If $g = \epsilon$ we would have $h = \epsilon = g$, so $g \ne \epsilon$. This implies that $m \geqslant 2$. By our inductive hypothesis we obtain there exists $q < m \leqslant n$ such that $h \in K(q)$ and $\rho_{/dom(h)} \in \Xi(h)$. Now 
\[ \sigma_{/dom(h)} = (\sigma_{/dom(g)})_{/dom(h)} =  \rho_{/dom(h)} \in \Xi(h) . \]
\end{proof}

\medskip

\begin{lemma}\label{L:no-var-repetition-context}
For each $n$ positive integer $k = <<x_1, \varphi_1> \dots <x_m, \varphi_m>> \in K(n) - \{ \epsilon \}$, for each $i, j = 1 \dots m$ $i \ne j \to x_i \ne x_j$.
\end{lemma}

\begin{proof}

\smallskip

Since $K(1) - \{ \epsilon \} = \emptyset$ the initial step is trivially verified.

\medskip

Let $n$ be a positive integer, let $k = <<x_1, \varphi_1> \dots <x_m, \varphi_m>> \in K(n+1) - \{ \epsilon \}$, we want to verify that
for each $i, j = 1 \dots m$ $i \ne j \to x_i \ne x_j$.

\medskip

If $k \in K(n)$ this is obvioulsy verified.

\medskip

Otherwise $k \in K(n)^+$, so there exist $h \in K(n), \phi \in E_s(n,h), y \in (\mathcal{V}-var(h))$ such that $k = h + <y,\phi>$.\\

If $h = \epsilon$ then $k = <<y, \phi>>$, this implies $m = 1$ and we have finished.\\

If $h \ne \epsilon$ then $h = <<y_1, \psi_1> \dots <y_p, \psi_p>>$ and\\
$k = <<y_1, \psi_1> \dots <y_p, \psi_p><y, \phi>>$.\\

Clearly this implies $m = p+1$. Given $i,j = 1 \dots m$ with $i \ne j$ if $i, j \leqslant p$ then $x_i = y_i \ne y_j = x_j$.  If $i \leqslant p$ and $j = m$ then $x_i = y_i \ne y = x_m = x_j$.
\end{proof}

\medskip

\begin{lemma}\label{L:same-first-element}
For each $n$ positive integer, $k \in K(n)$, $\sigma = (z, \xi) \in \Xi(k)$:
\begin{itemize}
\item if $k = \epsilon$ then $z = \emptyset$, $var(\sigma) = \emptyset = var(k)$;
\item if $k \ne \epsilon$, $k = <<x_1, \varphi_1> \dots <x_m, \varphi_m>>$ then $dom(z) = \{ 1, \dots , m \}$, for each $i = 1 \dots m$ $z_i = x_i$, $var(\sigma) = var(k)$.
\end{itemize}
\end{lemma}

\begin{proof}

\smallskip

The initial step is trivially verified.

\medskip

Let $n$ be a positive integer, let $k \in K(n+1)$, let $\sigma = (z, \xi) \in \Xi(k)$.  If $k \in K(n)$ then we can assume the result is valid.

\medskip

Otherwise  $k \in K(n)^+$, so there exist $h \in K(n), \phi \in E_s(n,h), y \in (\mathcal{V}-var(h))$ such that $k = h + <y,\phi>$ and 
\[ \Xi(k) = \{ \rho + (y,s) | \, \rho \in \Xi(h), s \in \#(h,\phi,\rho) \} \, . \]

There exists $\rho = (u, \nu) \in \Xi(h), s \in \#(h,\phi,\rho)$ such that $\sigma = \rho + (y,s)$.\\

If $h = \epsilon$ then $k = <<y, \phi>>$, so $m = 1$, $\rho = \epsilon$, $dom(z) = \{ 1 \}$, $z(1) = y$. Moreover $x_1 = y = z(1)$, $var(k) = \{ y \} = var(\sigma)$.\\

Otherwise let $h = <<y_1, \psi_1> \dots <y_p, \psi_p>>$, so\\
$k = <<y_1, \psi_1> \dots <y_p, \psi_p><y, \phi>>$.\\

Using our inductive hypothesis we can state that $dom(u) = \{ 1, \dots , p \}$ for each $i = 1 \dots p$ $y_i = u_i$, $var(\rho) = var(h)$.\\

It follows that $dom(z) = \{ 1, \dots , p+1 \} = \{ 1, \dots , m \}$.\\

For each $i = 1 \dots p$ $x_i = y_i = u_i = z_i$, moreover $x_{p+1} = y = z_{p+1}$.\\

It also follows that $var(\sigma) = var(k)$.
\end{proof}

\medskip

\begin{lemma}\label{L:no-var-repetition-state}
For each $n$ positive integer, $k \in K(n)$, $\sigma = (z, \xi) \in \Xi(k)$, for each $i, j \in dom(\sigma)$ $i \ne j \to z_i \ne z_j$.
\end{lemma}

\begin{proof}

\smallskip

Clearly in the case $k = \epsilon$ we have $\sigma = \epsilon$ and the result is trivially verified.\\

Now suppose $k \ne \epsilon$, $k = <<x_1, \varphi_1> \dots <x_m, \varphi_m>>$.\\

From lemma~\ref{L:no-var-repetition-context} it follows that for each $i, j = 1 \dots m$ $i \ne j \to x_i \ne x_j$. From lemma~\ref{L:same-first-element} $dom(z) = \{ 1, \dots , m \}$, for each $i = 1 \dots m$ $z_i = x_i$.\\

It follows that for each $i,j \in dom(\sigma)$ if $i \ne j$ then $z_i = x_i \ne x_j = z_j$.\\
\end{proof}

\medskip

\begin{lemma}\label{L:context-and-subcontext-satisfy}
For each $n$ positive integer, $k = <<x_1, \varphi_1> \dots <x_m, \varphi_m>>, h = <<y_1, \psi_1> \dots <y_q, \psi_q>> \in K(n) - \{ \epsilon \}$ if $h \sqsubseteq k$ then for each $i \in dom(k)$, $j \in dom(h)$ $x_i = y_j \to \varphi_i = \psi_j$.
\end{lemma}

\begin{proof}

\smallskip

From lemma~\ref{L:no-var-repetition-context} it follows that for each $i, j \in dom(k)$ $i \ne j \to x_i \ne x_j$. With this we can apply lemma~\ref{L:useful-context-2} and obtain that there exists $p = 1 \dots m$ such that $h = <<x_1, \varphi_1> \dots <x_p, \varphi_p>>$.\\

At this point for each $i \in dom(k)$, $j \in dom(h)$ $x_i = y_j$ implies $x_i = x_j$ so $i = j$ and $\varphi_i = \varphi_j = \psi_j$.
\end{proof}

\medskip

\begin{lemma}\label{L:state-and-substate-satisfy}
For each $n$ positive integer, $h, k \in K(n)$, $\sigma = (x, \eta) \in \Xi(k)$, $\rho = (y, \theta) \in \Xi(h)$, if $\rho \sqsubseteq \sigma$ then for each $i \in dom(\sigma)$, $j \in dom(\rho)$ $x_i = y_j \to \eta_i = \theta_j$.
\end{lemma}

\begin{proof}

\smallskip

From lemma~\ref{L:no-var-repetition-state} it follows that for each $i, j \in dom(\sigma)$ $i \ne j \to x_i \ne x_j$. With this we can apply lemma~\ref{L:useful-slp} and obtain that\\
for each $i \in dom(\sigma)$, $j \in dom(\rho)$ $x_i = y_j \to \eta_i = \theta_j$.
\end{proof}

\medskip

\begin{lemma}\label{L:exp-in-sets-abcd}
For each $n$ positive integer such that $n \geqslant 2$, $k \in K(n)$, $t \in E(n,k)$ such that $t \notin \mathcal{C}$ one of the following two alternatives holds: 
\begin{itemize}
\item $t \in E_a(n,k) \cup E_\forall(n,k) \cup E_\exists(n,k) \cup \bigcup_{c \in \mathcal{C}'} E^c(n,k) \cup \bigcup_{f \in \mathcal{F}} E^f(n,k)$;
\item $n > 2$ and there exist $m$ positive integer such that $2 \leqslant m < n$, $h \in K(m)$ such that $h \sqsubseteq k$, $t \in E_a(m,h) \cup E_\forall(m,h) \cup E_\exists(m,h) \cup \bigcup_{c \in \mathcal{C}'} E^c(m,h) \cup \bigcup_{f \in \mathcal{F}} E^f(m,h)$ and for each $\sigma \in \Xi(k)$  $\sigma_{/dom(h)} \in \Xi(h)$ and $\#(k,t,\sigma) = \#(h,t,\sigma_{/dom(h)})$.
\end{itemize}
\end{lemma}

\begin{proof}

\smallskip

Of course we begin with the case $n = 2$. Let $k \in K(2)$, $t \in E(2,k)$ such that $t \notin \mathcal{C}$. We have $K(2) = K(1) \cup K(1)^+$. 

\medskip

If $k \in K(1)^+$ we have $E(2,k) = E_a(2,k)$, so $t \in E_a(2,k)$.

\medskip

If $k \in K(1)$ we have 
\[
E(2,k) = E(1,k) \cup E_b(2,k) \cup E_\forall(2,k) \cup E_\exists(2,k) \cup \bigcup_{c \in \mathcal{C}'} E^c(2,k) \cup \bigcup_{f \in \mathcal{F}} E^f(2,k) \ .
\]

\smallskip

Since $k = \epsilon$ we have $E(1,k) = \mathcal{C}$, $E_b(2,k) = \emptyset$, so 
\[
E(2,k) = \mathcal{C} \cup E_\forall(2,k) \cup E_\exists(2,k) \cup \bigcup_{c \in \mathcal{C}'} E^c(2,k) \cup \bigcup_{f \in \mathcal{F}} E^f(2,k) \ .
\]

\smallskip

Therefore in the case we are discussing
\[ t \in E_\forall(2,k) \cup E_\exists(2,k) \cup \bigcup_{c \in \mathcal{C}'} E^c(2,k) \cup \bigcup_{f \in \mathcal{F}} E^f(2,k) \ . \]

\medskip

Let now $n \geqslant 2$ and we try to prove the result for $n+1$. So let $k \in K(n+1)$, $t \in E(n+1, k)$ such that $t \notin \mathcal{C}$. We have $K(n+1) = K(n) \cup K(n)^+$. 

\medskip

If $k \in K(n)^+$ we have $E(n+1,k) = E_a(n+1,k)$ so $t \in E_a(n+1,k)$.

\medskip

We now need to examine the case $k \in K(n)$. Here we have 
\scriptsize
\begin{align*}
E(n+1,k) = E(n,k) \cup E_b(n+1,k) \cup E_\forall(n+1,k) \cup E_\exists(n+1,k) \cup \bigcup_{c \in \mathcal{C}'} E^c(n+1,k) \cup \bigcup_{f \in \mathcal{F}} E^f(n+1,k) \ .
\end{align*}
\normalsize

\smallskip

If $t \in E_\forall(n+1,k) \cup E_\exists(n+1,k) \cup \bigcup_{c \in \mathcal{C}'} E^c(n+1,k) \cup \bigcup_{f \in \mathcal{F}} E^f(n+1,k)$ then our result is verified.

\medskip

If $t \in E(n,k)$ we can apply our inductive hypothesis, which entails two alternatives: 
\begin{itemize}
\item $t \in E_a(n,k) \cup E_\forall(n,k) \cup E_\exists(n,k) \cup \bigcup_{c \in \mathcal{C}'} E^c(n,k) \cup \bigcup_{f \in \mathcal{F}} E^f(n,k)$;
\item $n > 2$ and there exist $m$ positive integer such that $2 \leqslant m < n$, $h \in K(m)$ such that $h \sqsubseteq k$, $t \in E_a(m,h) \cup E_\forall(m,h) \cup E_\exists(m,h) \cup \bigcup_{c \in \mathcal{C}'} E^c(m,h) \cup \bigcup_{f \in \mathcal{F}} E^f(m,h)$ and for each $\sigma \in \Xi(k)$  $\sigma_{/dom(h)} \in \Xi(h)$ and $\#(k,t,\sigma) = \#(h,t,\sigma_{/dom(h)})$.
\end{itemize}

\smallskip

In the first case we observe that $2 \leqslant n < n+1$, $k \in K(n)$, $k \sqsubseteq k$. Moreover for each $\sigma \in \Xi(k)$ $\sigma_{/dom(k)} = \sigma \in \Xi(k)$ and $\#(k,t,\sigma) = \#
(k,t,\sigma_{/dom(k)})$.

\medskip

So in the first case our result is verified.

\medskip

Let's examine the second case. Here $2 \leqslant m < n < n+1$, $h \in K(m)$, $h \sqsubseteq k$, for each $\sigma \in \Xi(k)$  $\sigma_{/dom(h)} \in \Xi(h)$ and $\#(k,t,\sigma) = \#(h,t,\sigma_{/dom(h)})$. So everything is as expected and our result is verified in this case too.\\

We have still one case to examine, which is the case of $t \in E_b(n+1,k)$. Here we have $k \ne \epsilon$ so by assumption~\ref{A:k-is-epsilon-or-new} there exist $m<n$, $h \in K(m)$, $\phi \in E_s(m,h)$, $y \in (\mathcal{V}-var(h))$ such that $k = h + <y, \phi>$. Moreover by the definition of $E_b(n+1,k)$ we know that $t \in E(n,h)$. So we can apply our inductive hypothesis, which again leads to two alternatives:

\begin{itemize}
\item $t \in E_a(n,h) \cup E_\forall(n,h) \cup E_\exists(n,h) \cup \bigcup_{c \in \mathcal{C}'} E^c(n,h) \cup \bigcup_{f \in \mathcal{F}} E^f(n,h)$;
\item $n > 2$ and there exist $p$ positive integer such that $2 \leqslant p < n$, $g \in K(p)$ such that $g \sqsubseteq h$, $t \in E_a(p,g) \cup E_\forall(p,g) \cup E_\exists(p,g) \cup \bigcup_{c \in \mathcal{C}'} E^c(p,g) \cup \bigcup_{f \in \mathcal{F}} E^f(p,g)$ and for each $\rho \in \Xi(h)$  $\rho_{/dom(g)} \in \Xi(g)$ and $\#(h,t,\rho) = \#(g,t,\rho_{/dom(g)})$.
\end{itemize}

\smallskip

In the first case we observe that $2 \leqslant n < n+1$, $h \in K(n)$, $h \sqsubseteq k$, $t \in E_a(n,h) \cup E_\forall(n,h) \cup E_\exists(n,h) \cup \bigcup_{c \in \mathcal{C}'} E^c(n,h) \cup \bigcup_{f \in \mathcal{F}} E^f(n,h)$, moreover for each $\sigma = \rho + (y,s) \in \Xi(k)$ we have 
\begin{itemize}
\item $\#(k,t,\sigma) = \#(h,t,\rho)$,
\item $\sigma_{/dom(h)} = \sigma_{/dom(\rho)} = \rho$, 
\item therefore $\#(k,t,\sigma) = \#(h,t,\sigma_{/dom(h)})$.
\end{itemize}

\medskip

Let's examine the second case. Here $2 \leqslant p < n < n+1$, $g \in K(p)$, $g \sqsubseteq h \sqsubseteq k$, $t \in E_a(p,g) \cup E_\forall(p,g) \cup E_\exists(p,g) \cup \bigcup_{c \in \mathcal{C}'} E^c(p,g) \cup \bigcup_{f \in \mathcal{F}} E^f(p,g)$. Moreover for each $\sigma = \rho + (y,s) \in \Xi(k)$ we have 
\begin{itemize}
\item $\#(k,t,\sigma) = \#(h,t,\rho)$,
\item $\sigma_{/dom(h)} = \sigma_{/dom(\rho)} = \rho$, 
\item $\sigma_{/dom(g)} = (\sigma_{/dom(h)})_{/dom(g)} = \rho_{/dom(g)}$,
\item $\#(k,t,\sigma) = \#(h,t,\rho) = \#(g,t,\rho_{/dom(g)}) = \#(g,t,\sigma_{/dom(g)})$.
\end{itemize}
\end{proof}

\bigskip

\begin{lemma}\label{L:meaning-of-constants-new}
For each $n$ positive integer, $k \in K(n)$, $t \in E(n,k)$ if $t \in \mathcal{C}$ then for each $\sigma \in \Xi(k)$ $\#(k, t, \sigma) = \#(t)$.
\end{lemma}

\begin{proof}

\smallskip

Let's verify the result for $n=1$. Here $k = \epsilon$, for each $\sigma \in \Xi(\epsilon)$ $\sigma = \epsilon$ so $\#(k, t, \sigma) = \#(\epsilon, t, \epsilon) = \#(t)$.

\medskip

Now let's examine the inductive step. Given $k \in K(n+1)$, $t \in E(n+1,k)$ such that $t \in \mathcal{C}$ and $\sigma \in \Xi(k)$ we want to show that $\#(k, t, \sigma) = \#(t)$.

\medskip

If $k \in K(n)^+$ then $t \in E_a(n+1,k)$, but since $t \in \mathcal{C}$ this cannot happen, so $k \in K(n)^+$ cannot happen.

\medskip

Therefore $k \in K(n)$ and\\
$t \in E(n,k) \cup E_b(n+1,k) \cup E_e(n+1,k) \cup \bigcup_{c \in \mathcal{C}'} E^c(n+1,k) \cup \bigcup_{f \in \mathcal{F}} E^f(n+1,k)$.

\medskip

Since $t \in \mathcal{C}$ it follows that $t \in E(n,k) \cup E_b(n+1,k)$.

\medskip

If $t \in E(n,k)$ clearly $\#(k, t, \sigma) = \#(t)$ holds by the inductive hypothesis.

\medskip

If $t \in E_b(n+1,k)$ then we have $k \ne \epsilon$ so by assumption~\ref{A:k-is-epsilon-or-new} there exist $m<n$, $h \in K(m)$, $\phi \in E_s(m,h)$, $y \in (\mathcal{V}-var(h))$ such that $k = h + <y, \phi>$, $\Xi(k) = \{ \rho + (y,s) | \, \rho \in \Xi(h), s \in \#(h,\phi,\rho) \} )$. Moreover by the definition of $E_b(n+1,k)$ we know that $t \in E(n,h)$. 

\medskip

Clearly there exist $\rho \in \Xi(h), \, s \in \#(h,\phi,\rho)$ such that $\sigma = \rho + (y, s)$ and $\#(k,t,\sigma) = \#(h,t,\rho)$. By the inductive hypothesis $\#(h,t,\rho) = \#(t)$, so $\#(k,t,\sigma) = \#(t)$.
\end{proof}

\bigskip

\begin{lemma}\label{L:kappa-acca-sigma-ro-aux1}
Let $k, h \in K(n)$ such that $h = \epsilon$ or $k = \epsilon$ or ($h, k \ne \epsilon$ and $k = <<x_1, \varphi_1> \dots <x_m, \varphi_m>>, \ h = <<y_1, \psi_1> \dots <y_q, \psi_q>>$ and for each $i \in dom(k)$, $j \in dom(h)$ $x_i = y_j \to \varphi_i = \psi_j$). Let $u \in \mathcal{V} - var(k)$: $u \in \mathcal{V} - var(h)$ and $\vartheta \in E(n)$, let $k' = k + <u, \vartheta>$ and $h' = h + <u, \vartheta>$. Since $k', h' \ne \epsilon$ there exist $x'_1, \dots , x'_p \in \mathcal{V}$, $\varphi'_1, \dots , \varphi'_p \in \Sigma^*$ such that $k' =  <<x'_1, \varphi'_1> \dots <x'_p, \varphi'_p>>$, $y'_1, \dots , y'_r \in \mathcal{V}$, $\psi'_1, \dots , \psi'_r \in \Sigma^*$ such that $h' = <<y'_1, \psi'_1> \dots <y'_r, \psi'_r>>$ and for each $i \in dom(k')$, $j \in dom(h')$ $x'_i = y'_j \to \varphi'_i = \psi'_j$.
\end{lemma}

\begin{proof}

If both $k, h = \epsilon$ then $k' = <<u, \vartheta >> = h'$ and our result is verified.\\

If $k \ne \epsilon$ and $h = \epsilon$ then let $k = <<x_1, \varphi_1> \dots <x_m, \varphi_m>>$, clearly $k' = <<x_1, \varphi_1> \dots <x_m, \varphi_m><u,\vartheta>>$ and $h = <<u, \vartheta >>$. Here we see that for each $i \in dom(k')$, $j \in dom(h')$ $x'_i = y'_j$ implies $j = 1$, $y'_j = u$, $x'_i = u$, so $\varphi'_i = \vartheta = \psi'_j$.\\

Finally if both $h, k \ne \epsilon$, $k = <<x_1, \varphi_1> \dots <x_m, \varphi_m>>, \ h = <<y_1, \psi_1> \dots <y_q, \psi_q>>$ then $k' = <<x_1, \varphi_1> \dots <x_m, \varphi_m><u,\vartheta>>$ and $h' = <<y_1, \psi_1> \dots <y_q, \psi_q><u,\vartheta>>$. Given $i \in dom(k')$, $j \in dom(h')$ such that $x'_i = y'_j$ we have $i = 1 \dots m+1$, $j = 1 \dots q+1$.\\

If $i \leqslant m$ and $j \leqslant q$ then clearly $x_i =  x'_i = y'_j = y_j$ and $\varphi'_i = \varphi_i = \psi_j = \psi'_j$.\\

If $i = m+1$ then $x'_i = u$, so $y'_j = u$ and $j = q+1$, so $\varphi'_i  = \vartheta = \psi'_j$.
\end{proof}

\bigskip

\begin{lemma}\label{L:kappa-acca-sigma-ro-aux2}
Let $k, h \in K(n)$ such that $h = \epsilon$ or $k = \epsilon$ or ($h, k \ne \epsilon$ and $k = <<x_1, \varphi_1> \dots <x_m, \varphi_m>>, \ h = <<y_1, \psi_1> \dots <y_q, \psi_q>>$ and for each $i \in dom(k)$, $j \in dom(h)$ $x_i = y_j \to \varphi_i = \psi_j$). Let $\kappa \sqsubseteq k$ and $g \sqsubseteq h$ then $\kappa = \epsilon$ or $g = \epsilon$ or 
\begin{itemize}
\item $\kappa, g \ne \epsilon$ and so $h, k \ne \epsilon$,
\item there exist $p, r$ positive integers such that $p \leqslant m$, $r \leqslant q$, $\kappa = <<x_1, \varphi_1> \dots <x_p, \varphi_p>>$, $g = <<y_1, \psi_1> \dots <y_r, \psi_r>>$ and for each $i \in dom(\kappa)$, $j \in dom(g)$ $x_i = y_j \to \varphi_i = \psi_j$.
\end{itemize}
\end{lemma}

\begin{proof}

Clearly we can have $\kappa = \epsilon$ or $g = \epsilon$, otherwise we have $\kappa, g \ne \epsilon$, so also ($h, k \ne \epsilon$ and $k = <<x_1, \varphi_1> \dots <x_m, \varphi_m>>, \ h = <<y_1, \psi_1> \dots <y_q, \psi_q>>$ and for each $i \in dom(k)$, $j \in dom(h)$ $x_i = y_j \to \varphi_i = \psi_j$).\\

By lemma~\ref{L:useful-context-2} there exist $p, r$ positive integers such that $p \leqslant m$, $r \leqslant q$, $\kappa = <<x_1, \varphi_1> \dots <x_p, \varphi_p>>$, $g = <<y_1, \psi_1> \dots <y_r, \psi_r>>$.\\

Moreover $dom(\kappa) \subseteq dom(k)$ and $dom(g) \subseteq dom(h)$ so for each $i \in dom(\kappa)$, $j \in dom(g)$ $x_i = y_j \to \varphi_i = \psi_j$.
\end{proof}

\bigskip

\begin{lemma}\label{L:kappa-acca-sigma-ro}
Let $k, h \in K(n)$ such that $h = \epsilon$ or $k = \epsilon$ or ($h, k \ne \epsilon$ and $k = <<x_1, \varphi_1> \dots <x_m, \varphi_m>>, \ h = <<y_1, \psi_1> \dots <y_q, \psi_q>>$ and for each $i \in dom(k)$, $j \in dom(h)$ $x_i = y_j \to \varphi_i = \psi_j$). Let $t \in E(n,k) \cap E(n,h)$. Let $\sigma = (x,z) \in \Xi(k)$, $\rho = (y,r) \in \Xi(h)$ such that for each $i \in dom(\sigma)$, $j \in dom(\rho)$ $x_i = y_j \to z_i = r_j$. Then $\#(k,t,\sigma) = \#(h,t,\rho)$.
\end{lemma}

\begin{proof}

\smallskip

We prove this by induction on a positive integer $n$.\\

Let's verify the initial step. Here we have $k, h \in K(1)$. This implies $h = \epsilon = k$. We have $t \in E(1, \epsilon) = \mathcal{C}$. We have $\sigma = (x,z) \in \Xi(\epsilon)$, $\rho = (y,r) \in \Xi(\epsilon)$. Of course this implies $\sigma = \epsilon = \rho$. Then $\#(k,t,\sigma) = \#(\epsilon, t, \epsilon) = \#(h,t,\rho)$.

\medskip 

Let us see the inductive step, that is given a positive integer $n$ we assume the result is true for each $m \leqslant n$ and we try to prove it for $n + 1$. In other words what we are trying to prove is that for each $k, h \in K(n+1)$ such that one of the following conditions holds 
\begin{itemize}
\item $h = \epsilon$
\item $k = \epsilon$
\item $h, k \ne \epsilon$ and $k = <<x_1, \varphi_1> \dots <x_m, \varphi_m>>, \ h = <<y_1, \psi_1> \dots <y_q, \psi_q>>$ and for each $i \in dom(k)$, $j \in dom(h)$ $x_i = y_j \to \varphi_i = \psi_j$
\end{itemize}

and for each $t \in E(n+1,k) \cap E(n+1,h)$, $\sigma = (x,z) \in \Xi(k)$, $\rho = (y,r) \in \Xi(h)$ such that for each $i \in dom(\sigma)$, $j \in dom(\rho)$ $x_i = y_j \to z_i = r_j$ we have $\#(k,t,\sigma) = \#(h,t,\rho)$.

\medskip

If $t \in \mathcal{C}$ then by lemma~\ref{L:meaning-of-constants-new} $\#(k,t,\sigma) = \#(t) = \#(h,t,\rho)$.

\medskip

Otherwise since $k \in K(n+1)$ and $t \in E(n+1,k)$ we can apply lemma~\ref{L:exp-in-sets-abcd} and obtain these two following alternative possibilities:

\begin{itemize}
\item $t \in E_a(n+1,k) \cup E_\forall(n+1,k) \cup E_\exists(n+1,k) \cup \bigcup_{c \in \mathcal{C}'} E^c(n+1,k) \cup \bigcup_{f \in \mathcal{F}} E^f(n+1,k)$;
\item $n+1 > 2$ and there exist $\mu$ positive integer such that $2 \leqslant \mu < n+1$, $\kappa \in K(\mu)$ such that $\kappa \sqsubseteq k$, $t \in E_a(\mu,\kappa) \cup E_\forall(\mu,\kappa) \cup E_\exists(\mu,\kappa) \cup \bigcup_{c \in \mathcal{C}'} E^c(\mu,\kappa) \cup \bigcup_{f \in \mathcal{F}} E^f(\mu,\kappa)$ and for each $\delta \in \Xi(k)$  $\delta_{/dom(\kappa)} \in \Xi(\kappa)$ and $\#(k,t,\delta) = \#(\kappa,t,\delta_{/dom(\kappa)})$.
\end{itemize}

\smallskip

Since $h \in K(n+1)$ and $t \in E(n+1,h)$ we can also use lemma~\ref{L:exp-in-sets-abcd} to obtain these two other following alternative possibilities:

\begin{itemize}
\item $t \in E_a(n+1,h) \cup E_\forall(n+1,h) \cup E_\exists(n+1,h) \cup \bigcup_{c \in \mathcal{C}'} E^c(n+1,h) \cup \bigcup_{f \in \mathcal{F}} E^f(n+1,h)$;
\item $n+1 > 2$ and there exist $\nu$ positive integer such that $2 \leqslant \nu < n+1$, $g \in K(\nu)$ such that $g \sqsubseteq h$, $t \in E_a(\nu,g) \cup E_\forall(\nu,g) \cup E_\exists(\nu,g) \cup \bigcup_{c \in \mathcal{C}'} E^c(\nu,g) \cup \bigcup_{f \in \mathcal{F}} E^f(\nu,g)$ and for each $\delta \in \Xi(h)$  $\delta_{/dom(g)} \in \Xi(g)$ and $\#(h,t,\delta) = \#(g,t,\delta_{/dom(g)})$.
\end{itemize}

\smallskip

So we have three possible cases to examine. The first is
\begin{itemize}
\item $t \in E_a(n+1,k) \cup E_\forall(n+1,k) \cup E_\exists(n+1,k) \cup \bigcup_{c \in \mathcal{C}'} E^c(n+1,k) \cup \bigcup_{f \in \mathcal{F}} E^f(n+1,k)$ and
\item $t \in E_a(n+1,h) \cup E_\forall(n+1,h) \cup E_\exists(n+1,h) \cup \bigcup_{c \in \mathcal{C}'} E^c(n+1,h) \cup \bigcup_{f \in \mathcal{F}} E^f(n+1,h)$.
\end{itemize} 

\smallskip

The second case is
\begin{itemize}
\item $t \in E_a(n+1,k) \cup E_\forall(n+1,k) \cup E_\exists(n+1,k) \cup \bigcup_{c \in \mathcal{C}'} E^c(n+1,k) \cup \bigcup_{f \in \mathcal{F}} E^f(n+1,k)$ and
\item $n+1 > 2$ and there exist $\nu$ positive integer such that $2 \leqslant \nu < n+1$, $g \in K(\nu)$ such that $g \sqsubseteq h$, $t \in E_a(\nu,g) \cup E_\forall(\nu,g) \cup E_\exists(\nu,g) \cup \bigcup_{c \in \mathcal{C}'} E^c(\nu,g) \cup \bigcup_{f \in \mathcal{F}} E^f(\nu,g)$ and for each $\delta \in \Xi(h)$  $\delta_{/dom(g)} \in \Xi(g)$ and $\#(h,t,\delta) = \#(g,t,\delta_{/dom(g)})$.
\end{itemize}

\smallskip

Another case to examine would be the following
\begin{itemize}
\item $n+1 > 2$ and there exist $\mu$ positive integer such that $2 \leqslant \mu < n+1$, $\kappa \in K(\mu)$ such that $\kappa \sqsubseteq k$, $t \in E_a(\mu,\kappa) \cup E_\forall(\mu,\kappa) \cup E_\exists(\mu,\kappa) \cup \bigcup_{c \in \mathcal{C}'} E^c(\mu,\kappa) \cup \bigcup_{f \in \mathcal{F}} E^f(\mu,\kappa)$ and for each $\delta \in \Xi(k)$  $\delta_{/dom(\kappa)} \in \Xi(\kappa)$ and $\#(k,t,\delta) = \#(\kappa,t,\delta_{/dom(\kappa)})$ and
\item $t \in E_a(n+1,h) \cup E_\forall(n+1,h) \cup E_\exists(n+1,h) \cup \bigcup_{c \in \mathcal{C}'} E^c(n+1,h) \cup \bigcup_{f \in \mathcal{F}} E^f(n+1,h)$.
\end{itemize}

\smallskip

Anyway this case is practically equal to the second one, so we don't need to consider it. Finally the third case is the following.
\begin{itemize}
\item $n+1 > 2$ and there exist $\mu$ positive integer such that $2 \leqslant \mu < n+1$, $\kappa \in K(\mu)$ such that $\kappa \sqsubseteq k$, $t \in E_a(\mu,\kappa) \cup E_\forall(\mu,\kappa) \cup E_\exists(\mu,\kappa) \cup \bigcup_{c \in \mathcal{C}'} E^c(\mu,\kappa) \cup \bigcup_{f \in \mathcal{F}} E^f(\mu,\kappa)$ and for each $\delta \in \Xi(k)$  $\delta_{/dom(\kappa)} \in \Xi(\kappa)$ and $\#(k,t,\delta) = \#(\kappa,t,\delta_{/dom(\kappa)})$ and
\item $n+1 > 2$ and there exist $\nu$ positive integer such that $2 \leqslant \nu < n+1$, $g \in K(\nu)$ such that $g \sqsubseteq h$, $t \in E_a(\nu,g) \cup E_\forall(\nu,g) \cup E_\exists(\nu,g) \cup \bigcup_{c \in \mathcal{C}'} E^c(\nu,g) \cup \bigcup_{f \in \mathcal{F}} E^f(\nu,g)$ and for each $\delta \in \Xi(h)$  $\delta_{/dom(g)} \in \Xi(g)$ and $\#(h,t,\delta) = \#(g,t,\delta_{/dom(g)})$.
\end{itemize}

\bigskip

We now examine the three different cases we have distinguished. We start with the first one, where we have five different subcases:\\
$t \in E_a(n+1,k) \cup E_\forall(n+1,k) \cup E_\exists(n+1,k) \cup \bigcup_{c \in \mathcal{C}'} E^c(n+1,k) \cup \bigcup_{f \in \mathcal{F}} E^f(n+1,k)$ 

\medskip

We start with the subcase $t \in E_a(n+1,k)$. We must have $t \in E_a(n+1,h)$. 

\medskip

If $k \in K(n)$ then $E_a(n+1,k) = \emptyset$ so $k \in K(n)^+$ and there exist $\kappa \in 
K(n), \ \theta \in E_s(n,\kappa), \  u \in (\mathcal{V}-var(\kappa))$ such that $k = \kappa 
+ <u,\theta>$, $E_a(n+1,k) = \{ u \}$. Since $\sigma \in \Xi(k)$ there exist $\xi \in \Xi
(\kappa), \ s \in \#(\kappa, \theta, \xi)$ such that $\sigma = \xi + (u,s)$,\\
$\#(k,t,\sigma) = \#(k,t,\sigma)_{(n+1,k,a)} = s$.

\medskip

If $h \in K(n)$ then $E_a(n+1,h) = \emptyset$ so $h \in K(n)^+$ and there exist $\vartheta \in 
K(n), \ \mu \in E_s(n,\vartheta), \  v \in (\mathcal{V}-var(\vartheta))$ such that $h = \vartheta 
+ <v,\mu>$, $E_a(n+1,k) = \{ v \}$. Since $\rho \in \Xi(h)$ there exist $\zeta \in \Xi
(\vartheta), \ q \in \#(\vartheta, \mu, \zeta)$ such that $\rho = \zeta + (v,q)$,\\
$\#(h,t,\rho) = \#(h,t,\rho)_{(n+1,h,a)} = q$.

\medskip

Since $t \in E_a(n+1,k)$ we have $t = u$, since $t \in E_a(n+1,h)$ we have $t = v$, therefore 
$u = v$.

\medskip

There exists $i \in dom(\sigma)$ such that $u = x_i, \ s = z_i$, there exists $j \in dom
(\rho)$ such that $v = y_j, \ q = r_j$.

\medskip

Therefore $x_i = u = v = y_j$ and $\#(k,t,\sigma) = s = z_i = r_j = q = \#(h,t,\rho)$.

\begin{flushright}
$\diamond$
\end{flushright}

\bigskip

We now consider the subcase $t \in E_\forall(n+1,k)$ (which implies $k \in K(n)$). Clearly $t \notin E_a(n+1,h) \cup E_\exists(n+1,h) \cup \bigcup_{c \in \mathcal{C}'} E^c(n+1,h) \cup \bigcup_{f \in \mathcal{F}} E^f(n+1,h)$, so $t \in E_\forall(n+1,h)$ which implies $h \in K(n)$.\\

Then $t \in H_\forall(n+1,k)$, so there exist $\vartheta \in E_s(n,k)$, $u \in \mathcal{V} - var(k)$ such that if we define $k' = k + <u, \vartheta>$ then $k' \in K(n)$ and there also exists $\theta \in S(n,k')$ such that $t = \forall(u: \vartheta, \theta)$.\\

For each $\sigma \in \Xi(k)$
\[ \#(k, t, \sigma) = \text{for each} \ \sigma' \in \Xi(k'): \sigma \sqsubseteq \sigma' \ \#(k', \theta, \sigma') \ . \]

Moreover $t \in H_\forall(n+1,h)$, so there exist $\xi \in E_s(n,h)$, $v \in \mathcal{V} - var(h)$ such that if we define $h' = h + <v, \xi>$ then $h' \in K(n)$ and there also exists $\phi \in S(n,h')$ such that $t = \forall(v: \xi, \phi)$.\\

For each $\rho \in \Xi(h)$
\[ \#(h, t, \rho) = \text{for each} \ \rho' \in \Xi(h'): \rho \sqsubseteq \rho' \ \#(h', \phi, \rho') \ . \]

\smallskip

Given $t \in H_\forall(n+1,k)$, by lemma~\ref{IL:HQ-recursive-aux1} it follows that 

\begin{itemize}
\item $u \in \mathcal{V} - var(k)$,
\item the set of the positive integers $r$ such that $4 < r < \ell(t)$, $t[r]$ = `,' and $d(t, r) = 1$ has just one member $r_1$.
\end{itemize}

\smallskip

If we define the following:\\
if $r_1 = 5$ then $\eta = \epsilon$, else $\eta = t[5, r_1 - 1]$,\\
if $r_1 = \ell(t) - 1$ then $\zeta = \epsilon$ else $\zeta = t[r_1 + 1, \ell(t) - 1]$,\\

then we have also
\begin{itemize}
\item $\eta = \vartheta \in E_s(n,k)$.
\item if we define $\kappa = k + <u, \eta>$ then $\kappa = k' \in K(n)$ and $\zeta = \theta \in S(n,\kappa)$ .
\end{itemize}

\smallskip

Given $t \in H_\forall(n+1,h)$, by lemma~\ref{IL:HQ-recursive-aux1} it also follows that 

\begin{itemize}
\item $u = v \in \mathcal{V} - var(h)$,
\item $\eta = \xi \in E_s(n,h)$.
\item if we define $g = h + <v, \eta>$ then $g = h' \in K(n)$ and $\zeta = \phi \in S(n,g)$ .
\end{itemize}

\smallskip

Finally it also follows that $\xi = \vartheta$ and $\phi = \theta$.\\

We want now verify that $\#(k,t,\sigma)$ = $\#(h, t, \rho)$.\\

We have just seen that 
\[ \#(k, t, \sigma) = \text{for each} \ \sigma' \in \Xi(k'): \sigma \sqsubseteq \sigma' \ \#(k', \theta, \sigma') \ . \]

and 
\[ \#(h, t, \rho) = \text{for each} \ \rho' \in \Xi(h'): \rho \sqsubseteq \rho' \ \#(h', \phi, \rho') \ . \]

So $\#(k, t, \sigma)$ is obtained by applying a `for each' predicate to the set 
\[ \{ \#(k', \theta, \sigma') | \  \sigma' \in \Xi(k'), \sigma \sqsubseteq \sigma'  \} \ . \] 

Similarly $\#(h, t, \rho)$ is obtained by applying a `for each' predicate to the set 
\[ \{  \#(h', \phi, \rho') | \  \rho' \in \Xi(h'), \rho \sqsubseteq \rho'  \}  \ . \]

Clearly if we can show that 
\[ \{ \#(k', \theta, \sigma') | \  \sigma' \in \Xi(k'), \sigma \sqsubseteq \sigma'  \} = \{  \#(h', \phi, \rho') | \  \rho' \in \Xi(h'), \rho \sqsubseteq \rho'  \} \ . \]

then we have shown that $\#(k,t,\sigma)$ = $\#(h, t, \rho)$.\\

Given that $\phi = \theta$, in order to show the equality of the two sets we need to show two things:

\begin{itemize}
\item for each $\rho' \in \Xi(h')$ such that $\rho \sqsubseteq \rho'$ there exists $\sigma' \in \Xi(k')$ such that $\sigma \sqsubseteq \sigma'$ and $\#(h', \theta, \rho') = \#(k', \theta, \sigma')$;
\item for each $\sigma' \in \Xi(k')$ such that $\sigma \sqsubseteq \sigma'$ there exists $\rho' \in \Xi(h')$ such that $\rho \sqsubseteq \rho'$ and $\#(k', \theta, \sigma') = \#(h', \theta, \rho')$.
\end{itemize}

Obviously once we have proved the first statement, we can prove the second in the exact same way, so proving the first statement is enough. Let then $\rho' \in \Xi(h')$ be such that $\rho \sqsubseteq \rho'$, we try to find $\sigma' \in \Xi(k')$ such that $\sigma \sqsubseteq \sigma'$ and $\#(h', \theta, \rho') = \#(k', \theta, \sigma')$.\\

We have that $\vartheta \in E_s(n,h)$, $u \in \mathcal{V} - var(h)$, $h' = h + <u, \vartheta>$, $h' \in K(n)$, $h' \ne \epsilon$, so $n \geqslant 2$. Using lemma~\ref{L:k-in-km+ext} we have
\[ \Xi(h') = \{ \delta + (u,s) | \, \delta \in \Xi(h), s \in \#(h,\vartheta,\delta) \}  . \]

So there exist $\delta \in \Xi(h)$, $s \in \#(h,\vartheta,\delta)$ such that $\rho' = \delta + (u,s)$. At this point we notice that 
\[ \delta = \rho'_{/dom(\delta)} = \rho'_{/dom(h)} = \rho'_{/dom(\rho)} = \rho \ . \]

\smallskip

So there exists $s \in \#(h,\vartheta,\rho)$ such that $\rho' = \rho + (u,s)$.\\

Let now $\sigma' = \sigma + (u,s)$ and we want to verify that $\sigma' \in \Xi(k')$ and that $\#(h', \theta, \rho') = \#(k', \theta, \sigma')$.\\

We have that $\vartheta \in E_s(n,k)$, $u \in \mathcal{V} - var(k)$, $k' = k + <u, \vartheta>$, $k' \in K(n)$, $k' \ne \epsilon$, so $n \geqslant 2$. Using lemma~\ref{L:k-in-km+ext} we have
\[ \Xi(k') = \{ \delta + (u,s) | \, \delta \in \Xi(k), s \in \#(k,\vartheta,\delta) \}  . \]

Since $\sigma \in \Xi(k)$, in order to show that $\sigma' \in \Xi(k')$ we just need to show that $s \in \#(k,\vartheta,\sigma)$. We know that $s \in \#(h,\vartheta,\rho)$.\\

We can consider that $k,h \in K(n)$, $\vartheta \in E(n,k) \cap E(n,h)$, $\sigma \in \Xi(k)$, $\rho \in \Xi(h)$, so we can apply the inductive hypothesis and obtain that $\#(k,\vartheta,\sigma) = \#(h,\vartheta,\rho)$, so $s \in \#(k,\vartheta,\sigma)$ and $\sigma' \in \Xi(k')$.\\

We now want to prove that $\#(h', \theta, \rho') = \#(k', \theta, \sigma')$.\\

We first notice that $k' = k + <u, \vartheta>$, $h' = h + <u, \vartheta>$ where $u \in \mathcal{V} - var(k)$, $u \in \mathcal{V} - var(h)$. Using lemma~\ref{L:kappa-acca-sigma-ro-aux1} we have that there exist $x'_1, \dots , x'_p \in \mathcal{V}$, $\varphi'_1, \dots , \varphi'_p \in \Sigma^*$ such that $k' =  <<x'_1, \varphi'_1> \dots <x'_p, \varphi'_p>>$, $y'_1, \dots , y'_r \in \mathcal{V}$, $\psi'_1, \dots , \psi'_r \in \Sigma^*$ such that $h' = <<y'_1, \psi'_1> \dots <y'_r, \psi'_r>>$ and for each $i \in dom(k')$, $j \in dom(h')$ $x'_i = y'_j \to \varphi'_i = \psi'_j$.\\

We also notice that $\rho' = \rho + (u,s)$, $\sigma' = \sigma + (u,s)$, so we can apply lemma~\ref{L:useful-slp-ext} and obtain that, if we redefine $\sigma' = (x', \varphi')$ and $\rho' = (y', \psi')$, then for each $i \in dom(\sigma')$, $j \in dom(\rho')$ $x'_i = y'_j \to \varphi'_i = \psi'_j$.\\

At this point we can consider that $k',h' \in K(n)$, $\theta \in S(n,k') \cap S(n,h')$, $\sigma' \in \Xi(k')$, $\rho' \in \Xi(h')$ and we can apply our inductive hypothesis to obtain that $\#(h', \theta, \rho') = \#(k', \theta, \sigma')$, and this completes our proof.

\begin{flushright}
$\diamond$
\end{flushright}

\bigskip

We now consider the subcase $t \in E_\exists(n+1,k)$ (which implies $k \in K(n)$). Clearly $t \notin E_a(n+1,h) \cup E_\forall(n+1,h) \cup \bigcup_{c \in \mathcal{C}'} E^c(n+1,h) \cup \bigcup_{f \in \mathcal{F}} E^f(n+1,h)$, so $t \in E_\exists(n+1,h)$ which implies $h \in K(n)$.\\

Then $t \in H_\exists(n+1,k)$, so there exist $\vartheta \in E_s(n,k)$, $u \in \mathcal{V} - var(k)$ such that if we define $k' = k + <u, \vartheta>$ then $k' \in K(n)$ and there also exists $\theta \in S(n,k')$ such that $t = \exists(u: \vartheta, \theta)$.\\

For each $\sigma \in \Xi(k)$
\[ \#(k, t, \sigma) = \text{exists} \ \sigma' \in \Xi(k'): \sigma \sqsubseteq \sigma' \ \#(k', \theta, \sigma') \ . \]

Moreover $t \in H_\exists(n+1,h)$, so there exist $\xi \in E_s(n,h)$, $v \in \mathcal{V} - var(h)$ such that if we define $h' = h + <v, \xi>$ then $h' \in K(n)$ and there also exists $\phi \in S(n,h')$ such that $t = \exists(v: \xi, \phi)$.\\

For each $\rho \in \Xi(h)$
\[ \#(h, t, \rho) = \text{exists} \ \rho' \in \Xi(h'): \rho \sqsubseteq \rho' \ \#(h', \phi, \rho') \ . \]

\smallskip

Given $t \in H_\exists(n+1,k)$, by lemma~\ref{IL:HQ-recursive-aux1} it follows that 

\begin{itemize}
\item $u \in \mathcal{V} - var(k)$,
\item the set of the positive integers $r$ such that $4 < r < \ell(t)$, $t[r]$ = `,' and $d(t, r) = 1$ has just one member $r_1$.
\end{itemize}

\smallskip

If we define the following:\\
if $r_1 = 5$ then $\eta = \epsilon$, else $\eta = t[5, r_1 - 1]$,\\
if $r_1 = \ell(t) - 1$ then $\zeta = \epsilon$ else $\zeta = t[r_1 + 1, \ell(t) - 1]$,\\

then we have also
\begin{itemize}
\item $\eta = \vartheta \in E_s(n,k)$.
\item if we define $\kappa = k + <u, \eta>$ then $\kappa = k' \in K(n)$ and $\zeta = \theta \in S(n,\kappa)$ .
\end{itemize}

\smallskip

Given $t \in H_\exists(n+1,h)$, by lemma~\ref{IL:HQ-recursive-aux1} it also follows that 

\begin{itemize}
\item $u = v \in \mathcal{V} - var(h)$,
\item $\eta = \xi \in E_s(n,h)$.
\item if we define $g = h + <v, \eta>$ then $g = h' \in K(n)$ and $\zeta = \phi \in S(n,g)$ .
\end{itemize}

\smallskip

Finally it also follows that $\xi = \vartheta$ and $\phi = \theta$.\\

We want now verify that $\#(k,t,\sigma)$ = $\#(h, t, \rho)$.\\

We have just seen that 
\[ \#(k, t, \sigma) = \text{exists} \ \sigma' \in \Xi(k'): \sigma \sqsubseteq \sigma' \ \#(k', \theta, \sigma') \ . \]

and 
\[ \#(h, t, \rho) = \text{exists} \ \rho' \in \Xi(h'): \rho \sqsubseteq \rho' \ \#(h', \phi, \rho') \ . \]

So $\#(k, t, \sigma)$ is obtained by applying an `exists' predicate to the set 
\[ \{ \#(k', \theta, \sigma') | \  \sigma' \in \Xi(k'), \sigma \sqsubseteq \sigma'  \} \ . \] 

Similarly $\#(h, t, \rho)$ is obtained by applying an `exists' predicate to the set 
\[ \{  \#(h', \phi, \rho') | \  \rho' \in \Xi(h'), \rho \sqsubseteq \rho'  \}  \ . \]

Clearly if we can show that 
\[ \{ \#(k', \theta, \sigma') | \  \sigma' \in \Xi(k'), \sigma \sqsubseteq \sigma'  \} = \{  \#(h', \phi, \rho') | \  \rho' \in \Xi(h'), \rho \sqsubseteq \rho'  \} \ . \]

then we have shown that $\#(k,t,\sigma)$ = $\#(h, t, \rho)$.\\

Actually we  proved this statement in the case of `forall' and it can be proved in the exact same way for `exists'.

\begin{flushright}
$\diamond$
\end{flushright}

\bigskip

We now consider the subcase $t \in \bigcup_{c \in \mathcal{C}'} E^c(n+1,k)$. This implies there exists $c_1 \in \mathcal{C}'$ such that $t \in E^{c_1}(n+1,k)$. This also implies $k \in K(n)$ and we have $t \in H_{c_1}(n+1,k)$.\\

Clearly $t \notin E_a(n+1,h) \cup E_\forall(n+1,h) \cup E_\exists(n+1,h) \cup \bigcup_{f \in \mathcal{F}} E^f(n+1,h)$, so $t \in \bigcup_{c \in \mathcal{C}'} E^c(n+1,h)$. This implies there exists $c_2 \in \mathcal{C}'$ such that $t \in E^{c_2}(n+1,h)$. Clearly we must have $h \in K(n)$ and we have also $t \in H_{c_2}(n+1,h)$.\\

As we have seen in lemmas~\ref{IL:Hc-recursive-aux1} and~\ref{IL:Hc-recursive-aux2} since $t = (c_1)(\psi) = (c_2)(\psi)$ then $c_2 = c_1$ and we can unequivocally determine $\psi_1, \dots , \psi_{u}$ such that $t$ can be written as $(c_1)( \psi_1, \dots , \psi_{u})$. Since $t \in H_{c_1}(n+1,k) \cap H_{c_2}(n+1,h)$ then for each $i = 1 \dots u$ $\psi_i \in E(n,k) \cap E(n,h)$. By the inductive hypothesis it follows immediately that for each $i = 1 \dots u$ $\#(k, \psi_i, \sigma) = \#(h, \psi_i, \rho)$.\\

Finally it follows that 

\begin{align*}
\#(k, t, \sigma) &= \#(c_1)( \#(k, \psi_1, \sigma), \dots , \#(k, \psi_u, \sigma) )\\
&= \#(c_2)( \#(h, \psi_1, \rho), \dots , \#(h, \psi_u, \rho) )  = \#(h, t, \rho) 
\end{align*}

\begin{flushright}
$\diamond$
\end{flushright}

\bigskip

We now consider the subcase $t \in \bigcup_{f \in \mathcal{F}} E^f(n+1,k)$. This implies there exists $f_1 \in \mathcal{F}$ such that $t \in E^{f_1}(n+1,k)$. This also implies $k \in K(n)$ and we have $t \in H_{f_1}(n+1,k)$.\\

Clearly $t \notin E_a(n+1,h) \cup E_\forall(n+1,h) \cup E_\exists(n+1,h) \cup \bigcup_{c \in \mathcal{C}'} E^c(n+1,k)$, so $t \in \bigcup_{f \in \mathcal{F}} E^f(n+1,h)$. This implies there exists $f_2 \in \mathcal{F}$ such that $t \in E^{f_2}(n+1,h)$. This also implies $h \in K(n)$ and we have $t \in H_{f_2}(n+1,h)$.\\

Since $t \in H_{f_1}(n+1,k)$ then $t = f_1(\psi)$ where $\psi \in \Sigma^*$. Since $t \in H_{f_2}(n+1,k)$ then $t = f_2(\varphi)$ where $\varphi \in \Sigma^*$. It follows that $f_2 = f_1$.\\

If $f_1$ has multiplicity 1 then there exists $\psi \in E(n,k)$ such that $t = f_1(\psi)$ and there exists $\varphi \in E(n,h)$ such that $t = f_2(\varphi)$. It follows that $\varphi = \psi$ and by the inductive hypothesis $\#(k, \psi, \sigma) = \#(h, \psi , \rho) = \#(h, \varphi , \rho)$. It also follows that 
\[ \#(k, t, \sigma) = P_{f_1}(\#(k, \psi, \sigma)) = P_{f_2}(\#(h, \varphi, \rho)) = \#(h, t, \rho) . \]

\smallskip

If $f_1$ has multiplicity 2 we can consider that $t = f_1(\psi)$ where $\psi \in \Sigma^*$ and that $t \in H_{f_1}(n+1,k)$, so by lemma \ref{IL:Hf-recursive-aux2} we can determine $\psi_1, \psi_2 \in E(n,k)$ such that $t = f_1(\psi_1, \psi_2)$.\\

We have also $t = f_2(\psi)$ where $\psi \in \Sigma^*$ and $t \in H_{f_2}(n+1,h)$, so using the same lemma we can determine that $\psi_1, \psi_2 \in E(n,h)$ and $t = f_2(\psi_1, \psi_2)$.\\

By the inductive hypotesis $\#(k, \psi_1, \sigma) = \#(h, \psi_1, \rho)$ and $\#(k, \psi_2, \sigma) = \#(h, \psi_2, \rho)$, so 
\[  \#(k, t, \sigma) = P_{f_1}(\#(k, \psi_1, \sigma), \#(k, \psi_2, \sigma) ) =  P_{f_2}(\#(h, \psi_1, \rho), \#(h, \psi_2, \rho) ) = \#(h, t, \rho) .  \]

\begin{flushright}
$\diamond$
\end{flushright}

\medskip

Let's consider the second case, which as we recall is the following:
\begin{itemize}
\item $t \in E_a(n+1,k) \cup E_\forall(n+1,k) \cup E_\exists(n+1,k) \cup \bigcup_{c \in \mathcal{C}'} E^c(n+1,k) \cup \bigcup_{f \in \mathcal{F}} E^f(n+1,k)$ and
\item $n+1 > 2$ and there exist $\nu$ positive integer such that $2 \leqslant \nu < n+1$, $g \in K(\nu)$ such that $g \sqsubseteq h$, $t \in E_a(\nu,g) \cup E_\forall(\nu,g) \cup E_\exists(\nu,g) \cup \bigcup_{c \in \mathcal{C}'} E^c(\nu,g) \cup \bigcup_{f \in \mathcal{F}} E^f(\nu,g)$ and for each $\delta \in \Xi(h)$  $\delta_{/dom(g)} \in \Xi(g)$ and $\#(h,t,\delta) = \#(g,t,\delta_{/dom(g)})$.
\end{itemize}

\smallskip

Initially we consider the same five different subcases of the first case:\\
$t \in E_a(n+1,k) \cup E_\forall(n+1,k) \cup E_\exists(n+1,k) \cup \bigcup_{c \in \mathcal{C}'} E^c(n+1,k) \cup \bigcup_{f \in \mathcal{F}} E^f(n+1,k)$.\\

We start with the subcase $t \in E_a(n+1,k)$. We must have $t \in E_a(\nu,g)$. 

\medskip

If $k \in K(n)$ then $E_a(n+1,k) = \emptyset$ so $k \in K(n)^+$ and there exist $\kappa \in 
K(n), \ \theta \in E_s(n,\kappa), \  u \in (\mathcal{V}-var(\kappa))$ such that $k = \kappa 
+ <u,\theta>$, $E_a(n+1,k) = \{ u \}$. Since $\sigma \in \Xi(k)$ there exist $\xi \in \Xi
(\kappa), \ s \in \#(\kappa, \theta, \xi)$ such that $\sigma = \xi + (u,s)$,\\
$\#(k,t,\sigma) = \#(k,t,\sigma)_{(n+1,k,a)} = s$.

\medskip

If $g \in K(\nu-1)$ then $E_a(\nu,g) = \emptyset$ so $g \in K(\nu-1)^+$ and there exist $\vartheta \in K(\nu-1), \ \mu \in E_s(\nu-1,\vartheta), \  v \in (\mathcal{V}-var(\vartheta))$ such that $g = \vartheta + <v,\mu>$, $E_a(\nu,g) = \{ v \}$. We have $\rho \in \Xi(h)$ and $\rho_{/dom(g)} \in \Xi(g)$. Let $\eta = \rho_{/dom(g)}$, then there exist $\zeta \in \Xi(\vartheta), \ q \in \#(\vartheta, \mu, \zeta)$ such that $\eta = \zeta + (v,q)$, $\#(g,t,\eta) = \#(g,t,\eta)_{(\nu,g,a)} = q$.

\medskip

We have to prove that $\#(k,t,\sigma) = \#(h,t,\rho)$, and since $\#(h,t,\rho) = \#(g,t,\eta)$ it is enough to prove that $\#(k,t,\sigma) = \#(g,t,\eta)$.

\medskip

Since $t \in E_a(n+1,k)$ we have $t = u$, since $t \in E_a(\nu,g)$ we have $t = v$, therefore $u = v$.

\medskip

Since $\eta \sqsubseteq \rho$ we can apply lemma~\ref{L:useful-sub-slp} to show that if $\eta = (w, \delta)$ then for each $i \in dom(\sigma), \ j \in dom(\eta)$ $x_i = w_j \to z_i = \delta_j$. 

\medskip

There exists $i \in dom(\sigma)$ such that $u = x_i, \ s = z_i$, there exists $j \in dom(\eta)$ such that $v = w_j, \ q = \delta_j$.

\medskip

Therefore $x_i = u = v = w_j$ and $\#(k,t,\sigma) = s = z_i = \delta_j = q = \#(g,t,\eta)$.

\begin{flushright}
$\diamond$
\end{flushright}

\bigskip

We now consider the subcase $t \in E_\forall(n+1,k)$ (which implies $k \in K(n)$). Clearly $t \notin E_a(\nu,g) \cup E_\exists(\nu,g) \cup \bigcup_{c \in \mathcal{C}'} E^c(\nu,g) \cup \bigcup_{f \in \mathcal{F}} E^f(\nu,g)$, so $t \in E_\forall(\nu,g)$ which implies $g \in K(\nu-1)$.\\

Then $t \in H_\forall(n+1,k)$, so there exist $\vartheta \in E_s(n,k)$, $u \in \mathcal{V} - var(k)$ such that if we define $k' = k + <u, \vartheta>$ then $k' \in K(n)$ and there also exists $\theta \in S(n,k')$ such that $t = \forall(u: \vartheta, \theta)$.\\

For each $\sigma \in \Xi(k)$
\[ \#(k, t, \sigma) = \text{for each} \ \sigma' \in \Xi(k'): \sigma \sqsubseteq \sigma' \ \#(k', \theta, \sigma') \ . \]

Moreover $t \in H_\forall(\nu,g)$, so there exist $\xi \in E_s(\nu-1,g)$, $v \in \mathcal{V} - var(g)$ such that if we define $g' = g + <v, \xi>$ then $g' \in K(\nu-1)$ and there also exists $\phi \in S(\nu-1,g')$ such that $t = \forall(v: \xi, \phi)$.\\

For each $\eta \in \Xi(g)$
\[ \#(g, t, \eta) = \text{for each} \ \eta' \in \Xi(g'): \eta \sqsubseteq \eta' \ \#(g', \phi, \eta') \ . \]

\smallskip

Given $t \in H_\forall(n+1,k)$, by lemma~\ref{IL:HQ-recursive-aux1} it follows that 

\begin{itemize}
\item $u \in \mathcal{V} - var(k)$,
\item the set of the positive integers $r$ such that $4 < r < \ell(t)$, $t[r]$ = `,' and $d(t, r) = 1$ has just one member $r_1$.
\end{itemize}

\smallskip

If we define the following:\\
if $r_1 = 5$ then $\lambda = \epsilon$, else $\lambda = t[5, r_1 - 1]$,\\
if $r_1 = \ell(t) - 1$ then $\zeta = \epsilon$ else $\zeta = t[r_1 + 1, \ell(t) - 1]$,\\

then we have also
\begin{itemize}
\item $\lambda = \vartheta \in E_s(n,k)$.
\item if we define $\kappa = k + <u, \lambda>$ then $\kappa = k' \in K(n)$ and $\zeta = \theta \in S(n,\kappa)$ .
\end{itemize}

\smallskip

Given $t \in H_\forall(\nu,g)$, by lemma~\ref{IL:HQ-recursive-aux1} it also follows that 

\begin{itemize}
\item $u = v \in \mathcal{V} - var(g)$,
\item $\lambda = \xi \in E_s(\nu-1,g)$.
\item if we define $\chi = g + <v, \lambda>$ then $\chi = g' \in K(\nu-1)$ and $\zeta = \phi \in S(\nu-1,\chi)$ .
\end{itemize}

\smallskip

Finally it also follows that $\xi = \vartheta$ and $\phi = \theta$.\\

We want now verify that $\#(k,t,\sigma)$ = $\#(h, t, \rho)$.\\

As an assumption we have that $\rho_{/dom(g)} \in \Xi(g)$ and $\#(h,t,\rho) = \#(g,t,\rho_{/dom(g)})$.\\

So, if we define $\eta = \rho_{/dom(g)}$, then $\eta \in \Xi(g)$ and $\#(h,t,\rho) = \#(g,t,\eta)$, and what we need to prove is $\#(k,t,\sigma)$ = $\#(g, t, \eta)$.\\

We have just seen that 
\[ \#(k, t, \sigma) = \text{for each} \ \sigma' \in \Xi(k'): \sigma \sqsubseteq \sigma' \ \#(k', \theta, \sigma') \ . \]

and 
\[ \#(g, t, \eta) = \text{for each} \ \eta' \in \Xi(g'): \eta \sqsubseteq \eta' \ \#(g', \phi, \eta') \ . \]

So $\#(k, t, \sigma)$ is obtained by applying a `for each' predicate to the set 
\[ \{ \#(k', \theta, \sigma') | \  \sigma' \in \Xi(k'), \sigma \sqsubseteq \sigma'  \} \ . \] 

Similarly $\#(g, t, \eta)$ is obtained by applying a `for each' predicate to the set 
\[ \{  \#(g', \phi, \eta') | \  \eta' \in \Xi(g'), \eta \sqsubseteq \eta'  \}  \ . \]

Clearly if we can show that 
\[ \{ \#(k', \theta, \sigma') | \  \sigma' \in \Xi(k'), \sigma \sqsubseteq \sigma'  \} = \{  \#(g', \phi, \eta') | \  \eta' \in \Xi(g'), \eta \sqsubseteq \eta'  \} \ . \]

then we have shown that $\#(k,t,\sigma)$ = $\#(g, t, \eta)$.\\

Given that $\phi = \theta$, in order to show the equality of the two sets we need to show two things:

\begin{itemize}
\item for each $\eta' \in \Xi(g')$ such that $\eta \sqsubseteq \eta'$ there exists $\sigma' \in \Xi(k')$ such that $\sigma \sqsubseteq \sigma'$ and $\#(g', \theta, \eta') = \#(k', \theta, \sigma')$;
\item for each $\sigma' \in \Xi(k')$ such that $\sigma \sqsubseteq \sigma'$ there exists $\eta' \in \Xi(g')$ such that $\eta \sqsubseteq \eta'$ and $\#(k', \theta, \sigma') = \#(g', \theta, \eta')$.
\end{itemize}

\smallskip

Let then $\eta' \in \Xi(g')$ be such that $\eta \sqsubseteq \eta'$, we try to find $\sigma' \in \Xi(k')$ such that $\sigma \sqsubseteq \sigma'$ and $\#(g', \theta, \eta') = \#(k', \theta, \sigma')$.\\

We have that $\xi \in E_s(\nu-1,g)$, $v \in \mathcal{V} - var(g)$, $g' = g + <v, \xi> \, \in K(\nu - 1)$, $g' \ne \epsilon$ so $\nu - 1 \geqslant 2$, so using lemma~\ref{L:k-in-km+ext} we have
\[ \Xi(g') = \{ \delta + (v,s) | \, \delta \in \Xi(g), s \in \#(g,\xi,\delta) \}  . \]

So, given $\eta' \in \Xi(g'): \eta \sqsubseteq \eta'$ there exist $\delta \in \Xi(g)$, $s \in \#(g,\xi,\delta)$ such that $\eta' = \delta + (v,s)$. At this point we notice that 
\[ \delta = \eta'_{/dom(\delta)} = \eta'_{/dom(g)} = \eta'_{/dom(\eta)} = \eta \ . \]

\smallskip

So there exists $s \in \#(g,\xi,\eta)$ such that $\eta' = \eta + (v,s)$, and in other words there exists $s \in \#(g,\vartheta,\eta)$ such that $\eta' = \eta + (u,s)$.\\

Let now $\sigma' = \sigma + (u,s)$ and we want to verify that $\sigma' \in \Xi(k')$.\\ 

We have that $\vartheta \in E_s(n,k)$, $u \in \mathcal{V} - var(k)$, $k' = k + <u, \vartheta>$, $k' \in K(n)$, $k' \ne \epsilon$, so $n \geqslant 2$. Using lemma~\ref{L:k-in-km+ext} we have
\[ \Xi(k') = \{ \delta + (u,s) | \, \delta \in \Xi(k), s \in \#(k,\vartheta,\delta) \}  . \]

Since $\sigma \in \Xi(k)$, in order to show that $\sigma' \in \Xi(k')$ we just need to show that $s \in \#(k,\vartheta,\sigma)$. We know that $s \in \#(g,\vartheta,\eta)$.\\

We have $\vartheta \in E(n,k)$, $\vartheta \in E(\nu - 1, g) \subseteq E(n,g)$.\\

We also notice that by lemma~\ref{L:kappa-acca-sigma-ro-aux2}, since $g \sqsubseteq h$, $k = \epsilon$ or $g = \epsilon$ or 
\begin{itemize}
\item $k, g \ne \epsilon$ and so $h \ne \epsilon$, $k = <<x_1, \varphi_1> \dots <x_m, \varphi_m>>, \ h = <<y_1, \psi_1> \dots <y_q, \psi_q>>$ and for each $i \in dom(k)$, $j \in dom(h)$ $x_i = y_j \to \varphi_i = \psi_j$;
\item there exists $v$ positive integer such that $v \leqslant q$, $g = <<y_1, \psi_1> \dots <y_v, \psi_v>>$ and for each $i \in dom(k)$, $j \in dom(g)$ $x_i = y_j \to \varphi_i = \psi_j$.
\end{itemize}

\medskip

Since $\eta \sqsubseteq \rho$ we can apply lemma~\ref{L:useful-sub-slp} to also show that if $\eta = (w, \mu)$ then for each $i \in dom(\sigma), \ j \in dom(\eta)$ $x_i = w_j \to z_i = \mu_j$.\\

With this we can apply the inductive hypothesis and obtain that $\#(k, \vartheta, \sigma) = \#(g,\vartheta,\eta)$, therefore $s \in \#(k,\vartheta,\sigma)$ and $\sigma' \in \Xi(k')$.\\

We now want to prove that $\#(g', \theta, \eta') = \#(k', \theta, \sigma')$.\\

We first notice that $k' = k + <u, \vartheta>$, $g' = g+ <u, \vartheta>$, with $u \in \mathcal{V} - var(k)$, $u \in \mathcal{V} - var(g)$.\\

By lemma~\ref{L:kappa-acca-sigma-ro-aux1} we have that there exist $x'_1, \dots , x'_p \in \mathcal{V}$, $\varphi'_1, \dots , \varphi'_p \in \Sigma^*$ such that $k' =  <<x'_1, \varphi'_1> \dots <x'_p, \varphi'_p>>$, $y'_1, \dots , y'_r \in \mathcal{V}$, $\psi'_1, \dots , \psi'_r \in \Sigma^*$ such that $g' = <<y'_1, \psi'_1> \dots <y'_r, \psi'_r>>$ and for each $i \in dom(k')$, $j \in dom(g')$ $x'_i = y'_j \to \varphi'_i = \psi'_j$.\\

We can also notice that $\sigma' = \sigma + (u,s)$, $\eta' = \eta + (u,s)$, $u \in \mathcal{V} - var(\sigma)$, $u \in \mathcal{V} - var(\eta)$ so by lemma~\ref{L:useful-slp-ext} if we set $\sigma' = (x',z')$, $\eta' = (w', \mu')$ then for each $i \in dom(\sigma'), \ j \in dom(\eta')$ $x'_i = w'_j \to z'_i = \mu'_j$.\\

At this point we can consider that $k',g' \in K(n)$, $\theta \in S(n,k') \cap S(n,g')$, $\sigma' \in \Xi(k')$, $\eta' \in \Xi(g')$ and we can apply our inductive hypothesis to obtain that $\#(g', \theta, \eta') = \#(k', \theta, \sigma')$, and this completes our proof.\\

For the other side of the proof, let $\sigma' \in \Xi(k')$ such that $\sigma \sqsubseteq \sigma'$, we try to find $\eta' \in \Xi(g')$ such that $\eta \sqsubseteq \eta'$ and $\#(k', \theta, \sigma') = \#(g', \theta, \eta')$.\\

We have that $\vartheta \in E_s(n,k)$, $u \in \mathcal{V} - var(k)$, $k' = k + <u, \vartheta>$, $k' \in K(n)$, $k' \ne \epsilon$, so $n \geqslant 2$. Using lemma~\ref{L:k-in-km+ext} we have
\[ \Xi(k') = \{ \delta + (u,s) | \, \delta \in \Xi(k), s \in \#(k,\vartheta,\delta) \}  . \]

So, given $\sigma' \in \Xi(k')$:  $\sigma \sqsubseteq \sigma'$ there exist $\delta \in \Xi(k)$, $s \in \#(k,\vartheta,\delta)$ such that $\sigma' = \delta + (u,s)$. At this point we notice that 
\[ \delta = \sigma'_{/dom(\delta)} = \sigma'_{/dom(k)} = \sigma'_{/dom(\sigma)} = \sigma \ . \]

So there exists $s \in \#(k,\vartheta,\sigma)$ such that $\sigma' = \sigma + (u,s)$.\\

Let now $\eta' = \eta + (u,s)$ and we want to verify that $\eta' \in \Xi(g')$.\\

We have that $\vartheta \in E_s(\nu-1,g)$, $u \in \mathcal{V} - var(g)$, $g' = g + <u, \vartheta> \, \in K(\nu - 1)$, $g' \ne \epsilon$ so $\nu - 1 \geqslant 2$, so using lemma~\ref{L:k-in-km+ext} we have
\[ \Xi(g') = \{ \delta + (u,s) | \, \delta \in \Xi(g), s \in \#(g,\vartheta,\delta) \}  . \]

Since $\eta \in \Xi(g)$, in order to show that $\eta' \in \Xi(g')$ we just need to show that $s \in \#(g,\vartheta,\eta)$. We know that $s \in \#(k,\vartheta,\sigma)$.\\

We have $\vartheta \in E(n,k)$, $\vartheta \in E(\nu - 1, g) \subseteq E(n,g)$.\\

We also notice that by lemma~\ref{L:kappa-acca-sigma-ro-aux2}, since $g \sqsubseteq h$, $k = \epsilon$ or $g = \epsilon$ or 
\begin{itemize}
\item $k, g \ne \epsilon$ and so $h \ne \epsilon$, $k = <<x_1, \varphi_1> \dots <x_m, \varphi_m>>, \ h = <<y_1, \psi_1> \dots <y_q, \psi_q>>$ and for each $i \in dom(k)$, $j \in dom(h)$ $x_i = y_j \to \varphi_i = \psi_j$;
\item there exists $v$ positive integer such that $v \leqslant q$, $g = <<y_1, \psi_1> \dots <y_v, \psi_v>>$ and for each $i \in dom(k)$, $j \in dom(g)$ $x_i = y_j \to \varphi_i = \psi_j$.
\end{itemize}

\medskip

Since $\eta \sqsubseteq \rho$ we can apply lemma~\ref{L:useful-sub-slp} to also show that if $\eta = (w, \mu)$ then for each $i \in dom(\sigma), \ j \in dom(\eta)$ $x_i = w_j \to z_i = \mu_j$.\\

With this we can apply the inductive hypothesis and obtain that $\#(k, \vartheta, \sigma) = \#(g,\vartheta,\eta)$, therefore $s \in \#(g,\vartheta,\eta)$ and $\eta' \in \Xi(g')$.\\

We now want to prove that $\#(k', \theta, \sigma') = \#(g', \theta, \eta')$.\\

We can notice that $k' = k + <u, \vartheta>$, $g' = g+ <u, \vartheta>$, with $u \in \mathcal{V} - var(k)$, $u \in \mathcal{V} - var(g)$.\\

By lemma~\ref{L:kappa-acca-sigma-ro-aux1} we have that there exist $x'_1, \dots , x'_p \in \mathcal{V}$, $\varphi'_1, \dots , \varphi'_p \in \Sigma^*$ such that $k' =  <<x'_1, \varphi'_1> \dots <x'_p, \varphi'_p>>$, $y'_1, \dots , y'_r \in \mathcal{V}$, $\psi'_1, \dots , \psi'_r \in \Sigma^*$ such that $g' = <<y'_1, \psi'_1> \dots <y'_r, \psi'_r>>$ and for each $i \in dom(k')$, $j \in dom(g')$ $x'_i = y'_j \to \varphi'_i = \psi'_j$.\\

We can also notice that $\sigma' = \sigma + (u,s)$, $\eta' = \eta + (u,s)$, $u \in \mathcal{V} - var(\sigma)$, $u \in \mathcal{V} - var(\eta)$ so by lemma~\ref{L:useful-slp-ext} if we set $\sigma' = (x',z')$, $\eta' = (w', \mu')$ then for each $i \in dom(\sigma'), \ j \in dom(\eta')$ $x'_i = w'_j \to z'_i = \mu'_j$.\\

At this point we can consider that $k',g' \in K(n)$, $\theta \in S(n,k') \cap S(n,g')$, $\sigma' \in \Xi(k')$, $\eta' \in \Xi(g')$ and we can apply our inductive hypothesis to obtain that $\#(k', \theta, \sigma') = \#(g', \theta, \eta')$, and this completes our proof.\\

\begin{flushright}
$\diamond$
\end{flushright}

\bigskip

We now consider the subcase $t \in E_\exists(n+1,k)$ (which implies $k \in K(n)$). Clearly $t \notin E_a(\nu,g) \cup E_\forall(\nu,g) \cup \bigcup_{c \in \mathcal{C}'} E^c(\nu,g) \cup \bigcup_{f \in \mathcal{F}} E^f(\nu,g)$, so $t \in E_\exists(\nu,g)$ which implies $g \in K(\nu-1)$.\\

Then $t \in H_\exists(n+1,k)$, so there exist $\vartheta \in E_s(n,k)$, $u \in \mathcal{V} - var(k)$ such that if we define $k' = k + <u, \vartheta>$ then $k' \in K(n)$ and there also exists $\theta \in S(n,k')$ such that $t = \exists(u: \vartheta, \theta)$.\\

For each $\sigma \in \Xi(k)$
\[ \#(k, t, \sigma) = \text{exists} \ \sigma' \in \Xi(k'): \sigma \sqsubseteq \sigma' \ \#(k', \theta, \sigma') \ . \]

Moreover $t \in H_\exists(\nu,g)$, so there exist $\xi \in E_s(\nu-1,g)$, $v \in \mathcal{V} - var(g)$ such that if we define $g' = g + <v, \xi>$ then $g' \in K(\nu-1)$ and there also exists $\phi \in S(\nu-1,g')$ such that $t = \exists(v: \xi, \phi)$.\\

For each $\eta \in \Xi(g)$
\[ \#(g, t, \eta) = \text{exists} \ \eta' \in \Xi(g'): \eta \sqsubseteq \eta' \ \#(g', \phi, \eta') \ . \]

\smallskip

Given $t \in H_\exists(n+1,k)$, by lemma~\ref{IL:HQ-recursive-aux1} it follows that 

\begin{itemize}
\item $u \in \mathcal{V} - var(k)$,
\item the set of the positive integers $r$ such that $4 < r < \ell(t)$, $t[r]$ = `,' and $d(t, r) = 1$ has just one member $r_1$.
\end{itemize}

\smallskip

If we define the following:\\
if $r_1 = 5$ then $\lambda = \epsilon$, else $\lambda = t[5, r_1 - 1]$,\\
if $r_1 = \ell(t) - 1$ then $\zeta = \epsilon$ else $\zeta = t[r_1 + 1, \ell(t) - 1]$,\\

then we have also
\begin{itemize}
\item $\lambda = \vartheta \in E_s(n,k)$.
\item if we define $\kappa = k + <u, \lambda>$ then $\kappa = k' \in K(n)$ and $\zeta = \theta \in S(n,\kappa)$ .
\end{itemize}

\smallskip

Given $t \in H_\exists(\nu,g)$, by lemma~\ref{IL:HQ-recursive-aux1} it also follows that 

\begin{itemize}
\item $u = v \in \mathcal{V} - var(g)$,
\item $\lambda = \xi \in E_s(\nu-1,g)$.
\item if we define $\chi = g + <v, \lambda>$ then $\chi = g' \in K(\nu-1)$ and $\zeta = \phi \in S(\nu-1,\chi)$ .
\end{itemize}

\smallskip

Finally it also follows that $\xi = \vartheta$ and $\phi = \theta$.\\

We want now verify that $\#(k,t,\sigma)$ = $\#(h, t, \rho)$.\\

As an assumption we have that $\rho_{/dom(g)} \in \Xi(g)$ and $\#(h,t,\rho) = \#(g,t,\rho_{/dom(g)})$.\\

So, if we define $\eta = \rho_{/dom(g)}$, then $\eta \in \Xi(g)$ and $\#(h,t,\rho) = \#(g,t,\eta)$, and what we need to prove is $\#(k,t,\sigma)$ = $\#(g, t, \eta)$.\\

We have just seen that 
\[ \#(k, t, \sigma) = \text{exists} \ \sigma' \in \Xi(k'): \sigma \sqsubseteq \sigma' \ \#(k', \theta, \sigma') \ . \]

and 
\[ \#(g, t, \eta) = \text{exists} \ \eta' \in \Xi(g'): \eta \sqsubseteq \eta' \ \#(g', \phi, \eta') \ . \]

So $\#(k, t, \sigma)$ is obtained by applying an `exists' predicate to the set 
\[ \{ \#(k', \theta, \sigma') | \  \sigma' \in \Xi(k'), \sigma \sqsubseteq \sigma'  \} \ . \] 

Similarly $\#(g, t, \eta)$ is obtained by applying an `exists' predicate to the set 
\[ \{  \#(g', \phi, \eta') | \  \eta' \in \Xi(g'), \eta \sqsubseteq \eta'  \}  \ . \]

Clearly if we can show that 
\[ \{ \#(k', \theta, \sigma') | \  \sigma' \in \Xi(k'), \sigma \sqsubseteq \sigma'  \} = \{  \#(g', \phi, \eta') | \  \eta' \in \Xi(g'), \eta \sqsubseteq \eta'  \} \ . \]

then we have shown that $\#(k,t,\sigma)$ = $\#(g, t, \eta)$.\\

Actually we  proved this statement in the case of `forall' and it can be proved in the exact same way for `exists'.

\begin{flushright}
$\diamond$
\end{flushright}

\bigskip

We now consider the subcase $t \in \bigcup_{c \in \mathcal{C}'} E^c(n+1,k)$. This implies there exists $c_1 \in \mathcal{C}'$ such that $t \in E^{c_1}(n+1,k)$. This also implies $k \in K(n)$ and we have $t \in H_{c_1}(n+1,k)$.\\

Clearly $t \notin E_a(\nu,g) \cup E_\forall(\nu,g) \cup E_\exists(\nu,g) \cup \bigcup_{f \in \mathcal{F}} E^f(\nu,g)$, so $t \in \bigcup_{c \in \mathcal{C}'} E^c(\nu,g)$. This implies there exists $c_2 \in \mathcal{C}'$ such that $t \in E^{c_2}(\nu, g)$. Clearly we must have $g \in K(\nu - 1)$ and we have also $t \in H_{c_2}(\nu,g)$.\\

As we have seen in lemmas~\ref{IL:Hc-recursive-aux1} and~\ref{IL:Hc-recursive-aux2} since $t = (c_1)(\psi) = (c_2)(\psi)$ then $c_2 = c_1$ and $t$ can be written as $(c_1)( \psi_1, \dots , \psi_{u})$. Since $t \in H_{c_1}(n+1,k) \cap H_{c_2}(\nu,g)$ then for each $i = 1 \dots u$ $\psi_i \in E(n,k) \cap E(\nu - 1, g)$. Clearly $g \in K(n)$ and $\psi_i \in E(n,g)$.\\

Let $\eta = \rho_{/dom(g)} \in \Xi(g)$. By lemma~\ref{L:kappa-acca-sigma-ro-aux2} we have that $k = \epsilon$ or $g = \epsilon$ or 
\begin{itemize}
\item $k, g \ne \epsilon$ and so $k, h \ne \epsilon$,
\item there exist $p$ positive integer such that $p \leqslant q$, $g = <<y_1, \psi_1> \dots <y_r, \psi_p>>$ and for each $i \in dom(k)$, $j \in dom(g)$ $x_i = y_j \to \varphi_i = \psi_j$.
\end{itemize}

\medskip

Since $\eta \sqsubseteq \rho$ we can apply lemma~\ref{L:useful-sub-slp} to also show that if $\eta = (w, \mu)$ then for each $i \in dom(\sigma), \ j \in dom(\eta)$ $x_i = w_j \to z_i = \mu_j$.\\

We can apply the inductive hypothesis and obtain that for each $i = 1 \dots u$ $\#(k, \psi_i, \sigma) = \#(g, \psi_i, \eta)$.\\

Finally it follows that 
\begin{align*}
\#(k, t, \sigma) &= \#(c_1)( \#(k, \psi_1, \sigma), \dots , \#(k, \psi_u, \sigma) )\\
&= \#(c_2)( \#(g, \psi_1, \eta), \dots , \#(g, \psi_u, \eta) ) = \#(g, t, \eta)  = \#(h, t, \rho) 
\end{align*}

\begin{flushright}
$\diamond$
\end{flushright}

\bigskip

We now consider the subcase $t \in \bigcup_{f \in \mathcal{F}} E^f(n+1,k)$. This implies there exists $f_1 \in \mathcal{F}$ such that $t \in E^{f_1}(n+1,k)$. This also implies $k \in K(n)$ and we have $t \in H_{f_1}(n+1,k)$.\\

Clearly $t \notin E_a(\nu,g) \cup E_\forall(\nu,g) \cup E_\exists(\nu,g) \cup \bigcup_{c \in \mathcal{C}'} E^c(\nu,g)$ so $t  \in \bigcup_{f \in \mathcal{F}} E^f(\nu,g)$. This implies there exists $f_2 \in \mathcal{F}$ such that $t \in E^{f_2}(\nu,g)$. This also implies $g \in K(\nu - 1)$ and we have $t \in H_{f_2}(\nu,g)$.\\

Since $t \in H_{f_1}(n+1,k)$ then $t = f_1(\psi)$ where $\psi \in \Sigma^*$. Since $t \in H_{f_2}(\nu,g)$ then $t = f_2(\varphi)$ where $\varphi \in \Sigma^*$. It follows that $f_2 = f_1$.\\

Let $\eta = \rho_{/dom(g)} \in \Xi(g)$. By lemma~\ref{L:kappa-acca-sigma-ro-aux2} we have that $k = \epsilon$ or $g = \epsilon$ or 
\begin{itemize}
\item $k, g \ne \epsilon$ and so $k, h \ne \epsilon$,
\item there exist $p$ positive integer such that $p \leqslant q$, $g = <<y_1, \psi_1> \dots <y_r, \psi_p>>$ and for each $i \in dom(k)$, $j \in dom(g)$ $x_i = y_j \to \varphi_i = \psi_j$.
\end{itemize}

\medskip

Since $\eta \sqsubseteq \rho$ we can apply lemma~\ref{L:useful-sub-slp} to also show that if $\eta = (w, \mu)$ then for each $i \in dom(\sigma), \ j \in dom(\eta)$ $x_i = w_j \to z_i = \mu_j$.\\

If $f_1$ has multiplicity 1 then there exists $\psi \in E(n,k)$ such that $t = f_1(\psi)$ and there exists $\varphi \in E(\nu - 1,g)$ such that $t = f_2(\varphi)$. It follows that $\varphi = \psi$. Clearly $g \in K(n)$ and $\varphi \in E(n,g)$.\\

So we can apply the inductive hypothesis and obtain that $\#(k, \psi, \sigma) = \#(g, \psi , \eta)$. It also follows that 
\[ \#(k, t, \sigma) = P_{f_1}(\#(k, \psi, \sigma)) = P_{f_2}(\#(g, \varphi, \eta)) = \#(g, t, \eta) = \#(h, t, \rho) . \]

\smallskip 

If $f_1$ has multiplicity 2 we can consider that $t = f_1(\psi)$ where $\psi \in \Sigma^*$ and that $t \in H_{f_1}(n+1,k)$, so by lemma \ref{IL:Hf-recursive-aux2} we can determine $\psi_1, \psi_2 \in E(n,k)$ such that $t = f_1(\psi_1, \psi_2)$.\\

We have also $t = f_2(\psi)$ where $\psi \in \Sigma^*$ and $t \in H_{f_2}(\nu,g)$, so using the same lemma we can determine that $\psi_1, \psi_2 \in E(\nu - 1,g)$ and $t = f_2(\psi_1, \psi_2)$. Clearly $g \in K(n)$ and $\psi_1, \psi_2 \in E(n,g)$.\\

By the inductive hypotesis $\#(k, \psi_1, \sigma) = \#(g, \psi_1, \eta)$ and $\#(k, \psi_2, \sigma) = \#(g, \psi_2, \eta)$, so 
\[  \#(k, t, \sigma) = P_{f_1}(\#(k, \psi_1, \sigma), \#(k, \psi_2, \sigma) ) =  P_{f_2}(\#(g, \psi_1, \eta), \#(g, \psi_2, \eta) ) = \#(g, t, \eta) , \]
\[ \#(k, t, \sigma) = \#(g, t, \eta) = \#(h, t, \rho) . \]

\begin{flushright}
$\diamond$
\end{flushright}

\bigskip

Finally, let's consider the third case, which, we recall, is the following. 
\begin{itemize}
\item $n+1 > 2$ and there exist $\mu$ positive integer such that $2 \leqslant \mu < n+1$, $\kappa \in K(\mu)$ such that $\kappa \sqsubseteq k$, $t \in E_a(\mu,\kappa) \cup E_\forall(\mu,\kappa) \cup E_\exists(\mu,\kappa) \cup \bigcup_{c \in \mathcal{C}'} E^c(\mu,\kappa) \cup \bigcup_{f \in \mathcal{F}} E^f(\mu,\kappa)$ and for each $\delta \in \Xi(k)$  $\delta_{/dom(\kappa)} \in \Xi(\kappa)$ and $\#(k,t,\delta) = \#(\kappa,t,\delta_{/dom(\kappa)})$ and
\item $n+1 > 2$ and there exist $\nu$ positive integer such that $2 \leqslant \nu < n+1$, $g \in K(\nu)$ such that $g \sqsubseteq h$, $t \in E_a(\nu,g) \cup E_\forall(\nu,g) \cup E_\exists(\nu,g) \cup \bigcup_{c \in \mathcal{C}'} E^c(\nu,g) \cup \bigcup_{f \in \mathcal{F}} E^f(\nu,g)$ and for each $\delta \in \Xi(h)$  $\delta_{/dom(g)} \in \Xi(g)$ and $\#(h,t,\delta) = \#(g,t,\delta_{/dom(g)})$.
\end{itemize}

\medskip

We have $t \in E(\mu, \kappa) \cap E(\nu, g)$, with $\mu, \nu < n + 1$.\\

We have also $\sigma = (x,z) \in \Xi(k)$, $\rho = (y,r) \in \Xi(h)$ such that for each $i \in dom(\sigma)$, $j \in dom(\rho)$ $x_i = y_j \to z_i = r_j$ and we want to show that $\#(k,t,\sigma) = \#(h,t,\rho)$. So we just need to show that $\#(\kappa,t,\sigma_{/dom(\kappa)}) = \#(g,t,\rho_{/dom(g)})$.\\

We can have $\kappa = \epsilon$ or $g = \epsilon$. Otherwise $\kappa, g \ne \epsilon$, $k, h \ne \epsilon$, $k = <<x_1, \varphi_1> \dots <x_m, \varphi_m>>, \ h = <<y_1, \psi_1> \dots <y_q, \psi_q>>$ and for each $i \in dom(k)$, $j \in dom(h)$ $x_i = y_j \to \varphi_i = \psi_j$. By lemma \ref{L:kappa-acca-sigma-ro-aux2} there exist $p, v$ positive integers such that $p \leqslant m$, $v \leqslant q$, $\kappa = <<x_1, \varphi_1> \dots <x_p, \varphi_p>>$, $g = <<y_1, \psi_1> \dots <y_v, \psi_v>>$ and for each $i \in dom(\kappa)$, $j \in dom(g)$ $x_i = y_j \to \varphi_i = \psi_j$.\\

If we define $u = max\{\mu,\nu\}$ then $\kappa, g \in K(u)$, $t \in E(u,\kappa) \cap E(u,g)$ and $u < n+1$.\\
Moreover let $\sigma' = \sigma_{/dom(\kappa)}$, $\sigma' = (x', z')$, $\rho' = \rho_{/dom(g)}$, $\rho' = (y', r')$. Since $\sigma' \sqsubseteq \sigma$ and $\rho' \sqsubseteq \rho$ by lemma~\ref{L:useful-sub-slp} we obtain that for each $i \in dom(\sigma')$, $j \in dom(\rho')$ $x'_i = y'_j \to z'_i = r'_j$.\\
By the inductive hypothesis we then obtain $\#(\kappa,t,\sigma') = \#(g,t,\rho')$, and so we have proved $\#(k,t,\sigma) = \#(h,t,\rho)$.
\end{proof}

\bigskip

\begin{lemma}\label{L:meaning-kept-f2}
Given
\begin{itemize}
\item a positive integer $n$;
\item $k \in K(n)$;
\item $f \in \mathcal{F}$ such that $f$ has multplicity $2$;
\item $\varphi_1, \varphi_2 \in E(n,k)$;
\end{itemize}

such that for each $\sigma \in \Xi(k)$ $A_f( \#(k,\varphi_1, \sigma), \#(k,\varphi_2, \sigma) )$ is true,\\
we have that $t = f( \varphi_1, \varphi_2) \in E(n+1,k)$.

\medskip

Given $\sigma \in \Xi(k)$ we have also
\[
\#(k,t,\sigma) = P_f( \#(k, \varphi_1, \sigma), \#(k, \varphi_2, \sigma) ) \ .
\]

\end{lemma}

\begin{proof}

If $t \in E(n,k) \cup E_b(n+1,k)$ then $t \in E(n+1,k)$, else $t \in E^f(n+1,k) \subseteq E(n+1,k)$.\\

Using lemma~\ref{L:exp-in-sets-abcd} we have that one of the following alternatives holds: 
\begin{itemize}
\item $t \in E_a(n+1,k) \cup E_\forall(n+1,k) \cup E_\exists(n+1,k) \cup \bigcup_{c \in \mathcal{C}'} E^c(n+1,k) \cup \bigcup_{g \in \mathcal{F}} E^g(n+1,k)$;
\item there exist $m$ positive integer such that $2 \leqslant m < n+1$, $h \in K(m)$ such that $h \sqsubseteq k$, $t \in E_a(m,h) \cup E_\forall(m,h) \cup E_\exists(m,h) \cup \bigcup_{c \in \mathcal{C}'} E^c(m,h) \cup \bigcup_{g \in \mathcal{F}} E^g(m,h)$ and for each $\sigma \in \Xi(k)$  $\sigma_{/dom(h)} \in \Xi(h)$ and $\#(k,t,\sigma) = \#(h,t,\sigma_{/dom(h)})$.
\end{itemize}

\medskip

If the first alternative holds, that is $t \in E_a(n+1,k) \cup E_\forall(n+1,k) \cup E_\exists(n+1,k) \cup \bigcup_{c \in \mathcal{C}'} E^c(n+1,k) \cup \bigcup_{g \in \mathcal{F}} E^g(n+1,k)$, then clearly $t \in E^f(n+1,k)$. This implies that $\#(k,t,\sigma) = \#(k,t,\sigma)_{(n+1,k,<f>)}$, so in this case our proof is finished.\\

Otherwise it must be $t \in E^f(m,h)$. This implies that there exist $\psi_1, \psi_2 \in E(m-1,h)$ such that $t = f(\psi_1, \psi_2)$, for each $\rho \in \Xi(h)$ 
\begin{itemize}
\item $A_f( \#(h,\psi_1, \rho), \#(h,\psi_2, \rho) )$ is true;
\item $\#(h,t,\rho) = P_f( \#(h, \psi_1, \rho), \#(h, \psi_2, \rho) )$.
\end{itemize}

\medskip

We now consider what we have seen in lemma~\ref{IL:Hf-recursive-aux2}. We have $t = f(\psi)$ with $\psi \in \Sigma^*$. Since $t \in H_f(n+1,k)$ the set of the positive integers $r$ such that $2 < r < \ell(t)$, $t[r]$ = `,' and $d(t, r) = 1$ has just one member $r_1$. We can define $\chi_1$ as follows: if $r_1 = 3$ then $\chi_1 = \epsilon$, else $\chi_1 = t[3, r_1 - 1]$. We also define $\chi_2$ as follows: if $r_1 = \ell(t) - 1$ then $\chi_2 = \epsilon$ else $\chi_2 = t[r_1 + 1, \ell(t) - 1]$. The lemma tells us that $\chi_1, \chi_2 \in E(n,k)$, and we can notice that $t = f(\chi_1, \chi_2)$.\\

Using lemma~\ref{IL:Hf-aux4} we obtain that $\varphi_1 = \chi_1$ and $\varphi_2 = \chi_2$.\\

We can apply again lemma~\ref{IL:Hf-recursive-aux2} using the fact that $t \in H_f(m,h)$, to obtain that $\chi_1, \chi_2 \in E(m-1,h)$ and $t = f(\chi_1, \chi_2)$ still holds. Using lemma~\ref{IL:Hf-aux4} we obtain that $\psi_1 = \chi_1 = \varphi_1$ and $\psi_2 = \chi_2 = \varphi_2$. Therefore $\varphi_1, \varphi_2 \in E(m-1,h) \subseteq E(n,h)$ and for each $\rho \in \Xi(h)$ $\#(h,t,\rho) = P_f( \#(h, \varphi_1, \rho), \#(h, \varphi_2, \rho) )$.\\

So given $\sigma \in \Xi(k)$ if we define $\rho = \sigma_{/dom(h)} \in \Xi(h)$ then
\[ \#(k,t,\sigma) = \#(h,t, \rho) = P_f( \#(h, \varphi_1, \rho), \#(h, \varphi_2, \rho ) ) \ . \]

\medskip

So we want to prove that 
\[
P_f( \#(h, \varphi_1, \rho), \#(h, \varphi_2, \rho) ) = P_f( \#(k, \varphi_1, \sigma), \#(k, \varphi_2, \sigma) ) \ ,
\]

\smallskip

and to prove this it is enough to prove that for each $\alpha \in \{1, 2 \}$
\[ \#(h, \varphi_\alpha, \rho) = \#(k, \varphi_\alpha, \sigma) \ . \] 

It is not difficult to prove this. In fact, by lemma~\ref{L:context-and-subcontext-satisfy}, if $k = <<x_1, \theta_1> \dots <x_u, \theta_u>>, \ h = <<y_1, \vartheta_1> \dots <y_q, \vartheta_q>> \in K(n) - \{ \epsilon \}$ since $h \sqsubseteq k$ then for each $i \in dom(k)$, $j \in dom(h)$ $x_i = y_j \to \theta_i = \vartheta_j$. If $\sigma = (x,z)$, $\rho = (y,r)$ then using lemma~\ref{L:state-and-substate-satisfy} we obtain that for each $i \in dom(\sigma)$, $j \in dom(\rho)$ $x_i = y_j \to z_i = r_j$. With this we can apply lemma~\ref{L:kappa-acca-sigma-ro} and obtain that $\#(h,\varphi_\alpha,\rho) = \#(k, \varphi_\alpha, \sigma)$.
\end{proof}

\bigskip

\begin{lemma}\label{L:meaning-kept-f1}
Given
\begin{itemize}
\item a positive integer $n$;
\item $k \in K(n)$;
\item $f \in \mathcal{F}$ such that $f$ has multplicity $1$;
\item $\varphi_1 \in E(n,k)$;
\end{itemize}

such that for each $\sigma \in \Xi(k)$ $A_f( \#(k,\varphi_1, \sigma) )$ is true,\\
we have that $t = f( \varphi_1) \in E(n+1,k)$.

\medskip

Given $\sigma \in \Xi(k)$ we have also
\[
\#(k,t,\sigma) = P_f( \#(k, \varphi_1, \sigma) ) \ .
\]

\end{lemma}

\begin{proof}

If $t \in E(n,k) \cup E_b(n+1,k)$ then $t \in E(n+1,k)$, else $t \in E^f(n+1,k) \subseteq E(n+1,k)$.\\

Using lemma~\ref{L:exp-in-sets-abcd} we have that one of the following alternatives holds: 
\begin{itemize}
\item $t \in E_a(n+1,k) \cup E_\forall(n+1,k) \cup E_\exists(n+1,k) \cup \bigcup_{c \in \mathcal{C}'} E^c(n+1,k) \cup \bigcup_{g \in \mathcal{F}} E^g(n+1,k)$;
\item there exist $m$ positive integer such that $2 \leqslant m < n+1$, $h \in K(m)$ such that $h \sqsubseteq k$, $t \in E_a(m,h) \cup E_\forall(m,h) \cup E_\exists(m,h) \cup \bigcup_{c \in \mathcal{C}'} E^c(m,h) \cup \bigcup_{g \in \mathcal{F}} E^g(m,h)$ and for each $\sigma \in \Xi(k)$  $\sigma_{/dom(h)} \in \Xi(h)$ and $\#(k,t,\sigma) = \#(h,t,\sigma_{/dom(h)})$.
\end{itemize}

\medskip

If the first alternative holds, that is $t \in E_a(n+1,k) \cup E_\forall(n+1,k) \cup E_\exists(n+1,k) \cup \bigcup_{c \in \mathcal{C}'} E^c(n+1,k) \cup \bigcup_{g \in \mathcal{F}} E^g(n+1,k)$, then clearly $t \in E^f(n+1,k)$. This implies that $\#(k,t,\sigma) = \#(k,t,\sigma)_{(n+1,k,<f>)}$, so in this case our proof is finished.\\

Otherwise it must be $t \in E^f(m,h)$. This implies that there exist $\psi_1 \in E(m-1,h)$ such that $t = f(\psi_1)$, for each $\rho \in \Xi(h)$ 
\begin{itemize}
\item $A_f( \#(h,\psi_1, \rho) )$ is true;
\item $\#(h,t,\rho) = \#(h,t,\rho)_{(m,h,<f>)} = P_f( \#(h, \psi_1, \rho) )$.\\
\end{itemize}

Clearly $\psi_1 = \varphi_1$ so for each $\rho \in \Xi(h)$ $\#(h,t,\rho) = P_f( \#(h, \varphi_1, \rho) )$.\\

Moreover given $\sigma \in \Xi(k)$ if we define $\rho = \sigma_{/dom(h)} \in \Xi(h)$ then
\[ \#(k,t,\sigma) = \#(h,t, \rho) = P_f( \#(h, \varphi_1, \rho) ) \ . \]

\medskip

So we want to prove that 
\[
P_f( \#(h, \varphi_1, \rho) ) = P_f( \#(k, \varphi_1, \sigma) ) \ ,
\]

\smallskip

and to prove this it is enough to prove that 
\[ \#(h, \varphi_1, \rho) = \#(k, \varphi_1, \sigma) \ . \] 

It is not difficult to prove this. In fact, by lemma~\ref{L:context-and-subcontext-satisfy}, if $k = <<x_1, \theta_1> \dots <x_u, \theta_u>>, \ h = <<y_1, \vartheta_1> \dots <y_q, \vartheta_q>> \in K(n) - \{ \epsilon \}$ since $h \sqsubseteq k$ then for each $i \in dom(k)$, $j \in dom(h)$ $x_i = y_j \to \theta_i = \vartheta_j$. If $\sigma = (x,z)$, $\rho = (y,r)$ then using lemma~\ref{L:state-and-substate-satisfy} we obtain that for each $i \in dom(\sigma)$, $j \in dom(\rho)$ $x_i = y_j \to z_i = r_j$. With this we can apply lemma~\ref{L:kappa-acca-sigma-ro} and obtain that $\#(h,\varphi_1,\rho) = \#(k, \varphi_1, \sigma)$.
\end{proof}

\bigskip

\begin{lemma}\label{L:meaning-kept-c1}
Given
\begin{itemize}
\item a positive integer $n$;
\item $k \in K(n)$;
\item $c \in \mathcal{C}'$ such that $\#(c)$ has a domain $M_1 \times \dots \times M_m$;
\item $\varphi_1, \dots , \varphi_m \in E(n,k)$;
\end{itemize}

such that for each $j = 1 \dots m$, $\sigma \in \Xi(k)$ $\#(k, \varphi_j, \sigma) \in M_j$,\\
we have that $t = (c)( \varphi_1, \dots , \varphi_m ) \in E(n+1,k)$.

\medskip

Given $\sigma \in \Xi(k)$ we have also
\[
\#(k,t,\sigma) = \#(c)( \#(k, \varphi_1, \sigma), \dots , \#(k, \varphi_m, \sigma) ) \ .
\]

\end{lemma}

\begin{proof}

If $t \in E(n,k) \cup E_b(n+1,k)$ then $t \in E(n+1,k)$, else $t \in E^c(n+1,k) \subseteq E(n+1,k)$.\\

Using lemma~\ref{L:exp-in-sets-abcd} we have that one of the following alternatives holds: 
\begin{itemize}
\item $t \in E_a(n+1,k) \cup E_\forall(n+1,k) \cup E_\exists(n+1,k) \cup \bigcup_{c \in \mathcal{C}'} E^c(n+1,k) \cup \bigcup_{g \in \mathcal{F}} E^g(n+1,k)$;
\item there exist $\nu$ positive integer such that $2 \leqslant \nu < n+1$, $h \in K(\nu)$ such that $h \sqsubseteq k$, $t \in E_a(\nu,h) \cup E_\forall(\nu,h) \cup E_\exists(\nu,h) \cup \bigcup_{c \in \mathcal{C}'} E^c(\nu,h) \cup \bigcup_{g \in \mathcal{F}} E^g(\nu,h)$ and for each $\sigma \in \Xi(k)$  $\sigma_{/dom(h)} \in \Xi(h)$ and $\#(k,t,\sigma) = \#(h,t,\sigma_{/dom(h)})$.
\end{itemize}

\medskip

If the first alternative holds, that is $t \in E_a(n+1,k) \cup E_\forall(n+1,k) \cup E_\exists(n+1,k) \cup \bigcup_{d \in \mathcal{C}'} E^d(n+1,k) \cup \bigcup_{g \in \mathcal{F}} E^g(n+1,k)$, then clearly $t \in E^c(n+1,k)$. This implies that $\#(k,t,\sigma) = \#(k,t,\sigma)_{(n+1,k,<c>)}$, so in this case our proof is finished.\\

Otherwise it must be $t \in E^c(\nu,h)$. This implies that there exist $\psi_1, \dots , \psi_m \in E(\nu-1,h)$ such that $t = (c)(\psi_1, \dots , \psi_m)$, for each $\rho \in \Xi(h)$ 
\begin{itemize}
\item for each $j = 1 \dots m$ $\#(h, \psi_j, \rho) \in M_j$;
\item $\#(h,t,\rho) = \#(c)( \#(h, \psi_1, \rho), \dots , \#(h, \psi_m, \rho) )$.\\
\end{itemize}

We have $t = (c)(\chi)$ with $\chi \in \Sigma^*$. Consider the set of the positive integers $r$ such that $4 < r < \ell(t)$, $t[r]$ = `,' and $d(t, r) = 1$. If this set is empty we can call $\chi_1 = \chi$ and using lemma~\ref{IL:Hc-recursive-aux1}, given that $t \in H_c(n+1,k)$ and $t \in H_c(\nu,h)$ we obtain that $\chi_1 \in E(n,k)$ and $\chi_1 \in E(\nu-1, h)$.\\

If the mentioned set is not empty let's name $r_1, \dots, r_u$ its members (in increasing order).\\ 
Let's also define $\chi_1$ as follows: if $r_1 = 5$ then let $\chi_1 = \epsilon$ else $r_1 > 5$ and let $\chi_1 = t[5, r_1 - 1]$.\\
If $u > 1$ then for each $i = 1 \dots u-1$: 
if $r_{i+1} = r_i + 1$ then let $\chi_{i+1} = \epsilon$ else $r_{i+1} > r_i + 1$ and let $\chi_{i+1} = t[r_i + 1, r_{i+1} - 1]$.\\
Finally if $r_u = \ell(t) - 1$ we define $\chi_{u+1} = \epsilon$, else $r_u < \ell(t) - 1$ and let $\chi_{u+1} = t[r_u + 1, \ell(t) - 1]$ .\\

Using lemma~\ref{IL:Hc-recursive-aux2} we obtain $t = (c)( \psi_1, \dots , \psi_{u+1})$ and for each $i = 1 \dots u+1$ $\chi_i \in E(n,k)$, $\chi_i \in E(\nu-1,h)$.\\

Using lemma~\ref{IL:Hc-aux4} we obtain that in the first case $m = 1$ and $\chi_1 = \varphi_1$, $\chi_1 = \psi_1$, in the second case $m = u+1$ and for each $j = 1 \dots u+1$ $\chi_j = \varphi_j$, $\chi_j = \psi_j$. In both cases for each $j = 1 \dots m$ $\varphi_j = \psi_j$.\\

Moreover given $\sigma \in \Xi(k)$ if we define $\rho = \sigma_{/dom(h)} \in \Xi(h)$ then
\[ \#(k,t,\sigma) = \#(h,t, \rho) = \#(c)( \#(h, \psi_1, \rho), \dots , \#(h, \psi_m, \rho) ) \ . \]

\medskip

In order to prove that $\#(k,t,\sigma) = \#(c)( \#(k, \varphi_1, \sigma), \dots , \#(k, \varphi_m, \sigma) )$ we just need to prove that for each $j = 1 \dots m$ $\#(k, \varphi_j, \sigma) = \#(h, \psi_j, \rho)$.\\

It is not difficult to prove this. In fact, by lemma~\ref{L:context-and-subcontext-satisfy}, if $k = <<x_1, \theta_1> \dots <x_u, \theta_u>>, \ h = <<y_1, \vartheta_1> \dots <y_q, \vartheta_q>> \in K(n) - \{ \epsilon \}$ since $h \sqsubseteq k$ then for each $i \in dom(k)$, $\alpha \in dom(h)$ $x_i = y_\alpha \to \theta_i = \vartheta_\alpha$. If $\sigma = (x,z)$, $\rho = (y,r)$ then using lemma~\ref{L:state-and-substate-satisfy} we obtain that for each $i \in dom(\sigma)$, $\alpha \in dom(\rho)$ $x_i = y_\alpha \to z_i = r_\alpha$. Moreover $\varphi_j = \psi_j \in E(\nu-1, h) \subseteq E(n,h)$. With this we can apply lemma~\ref{L:kappa-acca-sigma-ro} and obtain that $\#(h,\psi_j,\rho) = \#(h,\varphi_j,\rho) = \#(k, \varphi_j, \sigma)$.
\end{proof}

\bigskip

\begin{lemma}\label{L:meaning-kept-forall}
Given
\begin{itemize}
\item a positive integer $n$;
\item $k \in K(n)$;
\item $\varphi \in E_s(n,k)$;
\item $x \in \mathcal{V} - var(k)$;
\end{itemize}

such that if we define $k' = k + <x, \varphi>$ then $k' \in K(n)$, and given also $\phi \in S(n,k')$, we have that $t = \forall(x:\varphi, \phi) \in E(n+1,k)$.\\

Given $\sigma \in \Xi(k)$ we have also
\[ \#(k, t, \sigma) = \text{for each} \ \sigma' \in \Xi(k'): \sigma \sqsubseteq \sigma' \ \#(k', \phi, \sigma') \ . \]
\end{lemma}

\begin{proof}

If $t \in E(n,k) \cup E_b(n+1,k)$ then $t \in E(n+1,k)$, else $t \in E_\forall(n+1,k) \subseteq E(n+1,k)$.\\

Using lemma~\ref{L:exp-in-sets-abcd} we have that one of the following alternatives holds: 
\begin{itemize}
\item $t \in E_a(n+1,k) \cup E_\forall(n+1,k) \cup E_\exists(n+1,k) \cup \bigcup_{c \in \mathcal{C}'} E^c(n+1,k) \cup \bigcup_{g \in \mathcal{F}} E^g(n+1,k)$;
\item there exist $m$ positive integer such that $2 \leqslant m < n+1$, $h \in K(m)$ such that $h \sqsubseteq k$, $t \in E_a(m,h) \cup E_\forall(m,h) \cup E_\exists(m,h) \cup \bigcup_{c \in \mathcal{C}'} E^c(m,h) \cup \bigcup_{g \in \mathcal{F}} E^g(m,h)$ and for each $\sigma \in \Xi(k)$  $\sigma_{/dom(h)} \in \Xi(h)$ and $\#(k,t,\sigma) = \#(h,t,\sigma_{/dom(h)})$.
\end{itemize}

\medskip

If the first alternative holds, that is $t \in E_a(n+1,k) \cup E_\forall(n+1,k) \cup E_\exists(n+1,k) \cup \bigcup_{c \in \mathcal{C}'} E^c(n+1,k) \cup \bigcup_{g \in \mathcal{F}} E^g(n+1,k)$, then clearly $t \in E_\forall(n+1,k)$. This implies that $\#(k,t,\sigma) = \#(k,t,\sigma)_{(n+1,k,\forall)}$, so in this case our proof is finished.\\

Otherwise it must be $t \in E_\forall(m,h)$. This implies that $t \in H_\forall(m,h)$, so there exist $\xi \in E_s(m-1,h)$, $v \in \mathcal{V} - var(h)$ such that if we define $h' = h + <v, \xi>$ then $h' \in K(m-1)$ and there also exists $\theta \in S(m-1,h')$ such that $t = \forall(v: \xi, \theta)$.\\

Moreover for each $\rho \in \Xi(h)$
\[ \#(h, t, \rho) = \text{for each} \ \rho' \in \Xi(h'): \rho \sqsubseteq \rho' \ \#(h', \theta, \rho') \ . \]

By our hypothesis we have that $t = \forall(x:\varphi, \phi) \in H_\forall(n+1,k)$, and by lemma~\ref{IL:HQ-recursive-aux1} it follows that 

\begin{itemize}
\item $x \in \mathcal{V} - var(k)$,
\item the set of the positive integers $r$ such that $4 < r < \ell(t)$, $t[r]$ = `,' and $d(t, r) = 1$ has just one member $r_1$.
\end{itemize}

\smallskip

If we define the following:\\
if $r_1 = 5$ then $\lambda = \epsilon$, else $\lambda = t[5, r_1 - 1]$,\\
if $r_1 = \ell(t) - 1$ then $\zeta = \epsilon$ else $\zeta = t[r_1 + 1, \ell(t) - 1]$,\\

then we have also
\begin{itemize}
\item $\lambda = \varphi \in E_s(n,k)$.
\item if we define $\kappa = k + <x, \lambda>$ then $\kappa = k' \in K(n)$ and $\zeta = \phi \in S(n,\kappa)$ .
\end{itemize}

\medskip

Given $t \in H_\forall(m,h)$, by lemma~\ref{IL:HQ-recursive-aux1} it also follows that 

\begin{itemize}
\item $x = v \in \mathcal{V} - var(h)$,
\item $\lambda = \xi \in E_s(m-1,h)$.
\item if we define $\chi = h + <x, \lambda>$ then $\chi = h' \in K(m-1)$ and $\zeta = \theta \in S(m-1,\chi)$ .
\end{itemize}

\smallskip

Finally it also follows that $\xi = \varphi$ and $\theta = \phi$.\\

Given $\sigma \in \Xi(k)$ if we define $\rho = \sigma_{/dom(h)}$ then 
\[ \#(k,t,\sigma) = \#(h,t,\rho) = \text{for each} \ \rho' \in \Xi(h'): \rho \sqsubseteq \rho' \ \#(h', \theta, \rho') \ . \]

In order to show that 
\[ \#(k, t, \sigma) = \text{for each} \ \sigma' \in \Xi(k'): \sigma \sqsubseteq \sigma' \ \#(k', \phi, \sigma') \ , \]

we just need to show that
\[ \{ \#(k', \phi, \sigma') | \  \sigma' \in \Xi(k'), \sigma \sqsubseteq \sigma'  \} = \{  \#(h', \phi, \rho') | \  \rho' \in \Xi(h'), \rho \sqsubseteq \rho'  \} \ . \]

In order to show the equality of the two sets we need to show two things:

\begin{itemize}
\item for each $\rho' \in \Xi(h')$ such that $\rho \sqsubseteq \rho'$ there exists $\sigma' \in \Xi(k')$ such that $\sigma \sqsubseteq \sigma'$ and $\#(h', \phi, \rho') = \#(k', \phi, \sigma')$;
\item for each $\sigma' \in \Xi(k')$ such that $\sigma \sqsubseteq \sigma'$ there exists $\rho' \in \Xi(h')$ such that $\rho \sqsubseteq \rho'$ and $\#(k', \phi, \sigma') = \#(h', \phi, \rho')$.
\end{itemize}

Let then $\rho' \in \Xi(h')$ be such that $\rho \sqsubseteq \rho'$, we try to find $\sigma' \in \Xi(k')$ such that $\sigma \sqsubseteq \sigma'$ and $\#(h', \phi, \rho') = \#(k', \phi, \sigma')$.\\

We have that $\xi \in E_s(m-1,h)$, $v \in \mathcal{V} - var(h)$, $h' = h + <v, \xi> \, \in K(m - 1)$, $h' \ne \epsilon$ so $m - 1 \geqslant 2$, so using lemma~\ref{L:k-in-km+ext} we have
\[ \Xi(h') = \{ \delta + (v,s) | \, \delta \in \Xi(h), s \in \#(h,\xi,\delta) \}  . \]

So, given $\rho' \in \Xi(h'): \rho \sqsubseteq \rho'$ there exist $\delta \in \Xi(h)$, $s \in \#(h,\xi,\delta)$ such that $\rho' = \delta + (v,s)$. At this point we notice that 
\[ \delta = \rho'_{/dom(\delta)} = \rho'_{/dom(h)} = \rho'_{/dom(\rho)} = \rho \ . \]

\smallskip

So there exists $s \in \#(h,\xi,\rho)$ such that $\rho' = \rho + (v,s)$, and in other words there exists $s \in \#(h,\varphi,\rho)$ such that $\rho' = \rho + (x,s)$.\\

Let now $\sigma' = \sigma + (x,s)$ and we want to verify that $\sigma' \in \Xi(k')$.\\ 

We have that $\varphi \in E_s(n,k)$, $x \in \mathcal{V} - var(k)$, $k' = k + <x, \varphi>$, $k' \in K(n)$, $k' \ne \epsilon$, so $n \geqslant 2$. Using lemma~\ref{L:k-in-km+ext} we have
\[ \Xi(k') = \{ \delta + (x,s) | \, \delta \in \Xi(k), s \in \#(k,\varphi,\delta) \}  . \]

Since $\sigma \in \Xi(k)$, in order to show that $\sigma' \in \Xi(k')$ we just need to show that $s \in \#(k,\varphi,\sigma)$. We know that $s \in \#(h,\varphi,\rho)$.\\

We have $\varphi \in E(n,k)$, $\varphi \in E(m - 1, h) \subseteq E(n,h)$.\\

Since $h \sqsubseteq k$, if $h, k \ne \epsilon$ and $k = <<u_1, \psi_1> \dots <u_p, \psi_p>>$ then there exists $q = 1 \dots p$ such that $h = <<u_1, \psi_1> \dots <u_q, \psi_q>>$, and by lemma~\ref{L:context-and-subcontext-satisfy} for each $i \in dom(k)$, $j \in dom(h)$ $u_i = u_j \to \psi_i = \psi_j$.\\

\medskip

Since $\sigma \in \Xi(k)$, $\rho \in \Xi(h)$ and $\rho \sqsubseteq \sigma$ we can also apply lemma~\ref{L:state-and-substate-satisfy} and prove that if $\sigma=(y,z)$ and $\rho = (w, \mu)$ then for each $i \in dom(\sigma), \ j \in dom(\rho)$ $y_i = w_j \to z_i = \mu_j$.\\

With all this we can apply lemma~\ref{L:kappa-acca-sigma-ro} and obtain that $\#(k,\varphi,\sigma) = \#(h,\varphi,\rho)$, so we have proved that $s \in \#(k,\varphi,\sigma)$ and that $\sigma' \in \Xi(k')$.\\

We now want to prove that $\#(h', \phi, \rho') = \#(k', \phi, \sigma')$.\\

We first notice that $k' = k + <x, \varphi>$, $h' = h + <x, \varphi>$, with $x \in \mathcal{V} - var(k)$, $x \in \mathcal{V} - var(h)$.\\

We have seen that if $h, k \ne \epsilon$ then we can represent $k = <<u_1, \psi_1> \dots <u_p, \psi_p>>$, there exists $q = 1 \dots p$ such that $h = <<u_1, \psi_1> \dots <u_q, \psi_q>>$ and for each $i \in dom(k)$, $j \in dom(h)$ $u_i = u_j \to \psi_i = \psi_j$.\\

By lemma~\ref{L:kappa-acca-sigma-ro-aux1} we have that there exist $x'_1, \dots , x'_p \in \mathcal{V}$, $\varphi'_1, \dots , \varphi'_p \in \Sigma^*$ such that $k' =  <<x'_1, \varphi'_1> \dots <x'_p, \varphi'_p>>$, $y'_1, \dots , y'_r \in \mathcal{V}$, $\psi'_1, \dots , \psi'_r \in \Sigma^*$ such that $h' = <<y'_1, \psi'_1> \dots <y'_r, \psi'_r>>$ and for each $i \in dom(k')$, $j \in dom(h')$ $x'_i = y'_j \to \varphi'_i = \psi'_j$.\\

We can also notice that $\sigma' = \sigma + (x,s)$, $\rho' = \rho + (x,s)$, $x \in \mathcal{V} - var(\sigma)$, $x \in \mathcal{V} - var(\rho)$ so by lemma~\ref{L:useful-slp-ext} if we set $\sigma' = (x',z')$, $\rho' = (w', \mu')$ then for each $i \in dom(\sigma'), \ j \in dom(\rho')$ $x'_i = w'_j \to z'_i = \mu'_j$.\\

At this point we can consider that $k',h' \in K(n)$, $\phi \in S(n,k') \cap S(n,h')$, $\sigma' \in \Xi(k')$, $\rho' \in \Xi(h')$ and we can apply lemma~\ref{L:kappa-acca-sigma-ro} to obtain that $\#(h', \phi, \rho') = \#(k', \phi, \sigma')$, and this completes our proof.\\

For the other side of the proof, let $\sigma' \in \Xi(k')$ such that $\sigma \sqsubseteq \sigma'$, we try to find $\rho' \in \Xi(h')$ such that $\rho \sqsubseteq \rho'$ and $\#(k', \phi, \sigma') = \#(h', \phi, \rho')$.\\

We have that $\varphi \in E_s(n,k)$, $x \in \mathcal{V} - var(k)$, $k' = k + <x, \varphi>$, $k' \in K(n)$, $k' \ne \epsilon$, so $n \geqslant 2$. Using lemma~\ref{L:k-in-km+ext} we have
\[ \Xi(k') = \{ \delta + (x,s) | \, \delta \in \Xi(k), s \in \#(k,\varphi,\delta) \}  . \]

So, given $\sigma' \in \Xi(k')$:  $\sigma \sqsubseteq \sigma'$ there exist $\delta \in \Xi(k)$, $s \in \#(k,\varphi,\delta)$ such that $\sigma' = \delta + (x,s)$. At this point we notice that 
\[ \delta = \sigma'_{/dom(\delta)} = \sigma'_{/dom(k)} = \sigma'_{/dom(\sigma)} = \sigma \ . \]

So there exists $s \in \#(k,\varphi,\sigma)$ such that $\sigma' = \sigma + (x,s)$.\\

Let now $\rho' = \rho + (x,s)$ and we want to verify that $\rho' \in \Xi(h')$.\\

We have that $\varphi \in E_s(m-1,h)$, $x \in \mathcal{V} - var(h)$, $h' = h + <x, \varphi> \, \in K(m - 1)$, $h' \ne \epsilon$ so $m - 1 \geqslant 2$, so using lemma~\ref{L:k-in-km+ext} we have
\[ \Xi(h') = \{ \delta + (x,s) | \, \delta \in \Xi(h), s \in \#(h,\varphi,\delta) \}  . \]

Since $\rho \in \Xi(h)$, in order to show that $\rho' \in \Xi(h')$ we just need to show that $s \in \#(h,\varphi,\rho)$. We know that $s \in \#(k,\varphi,\sigma)$.\\

We have already proved in the first side of the proof that $\#(k,\varphi,\sigma) = \#(h,\varphi,\rho)$, so $s \in \#(h,\varphi,\rho)$ holds true and $\rho' \in \Xi(h')$.\\

We now want to prove that $\#(h', \phi, \rho') = \#(k', \phi, \sigma')$, this proof is the same of the one we saw in the first side of the main proof, but to be safe we are repeating it again.\\

We first notice that $k' = k + <x, \varphi>$, $h' = h + <x, \varphi>$, with $x \in \mathcal{V} - var(k)$, $x \in \mathcal{V} - var(h)$.\\

We have seen that if $h, k \ne \epsilon$ then we can represent $k = <<u_1, \psi_1> \dots <u_p, \psi_p>>$, there exists $q = 1 \dots p$ such that $h = <<u_1, \psi_1> \dots <u_q, \psi_q>>$ and for each $i \in dom(k)$, $j \in dom(h)$ $u_i = u_j \to \psi_i = \psi_j$.\\

By lemma~\ref{L:kappa-acca-sigma-ro-aux1} we have that there exist $x'_1, \dots , x'_p \in \mathcal{V}$, $\varphi'_1, \dots , \varphi'_p \in \Sigma^*$ such that $k' =  <<x'_1, \varphi'_1> \dots <x'_p, \varphi'_p>>$, $y'_1, \dots , y'_r \in \mathcal{V}$, $\psi'_1, \dots , \psi'_r \in \Sigma^*$ such that $h' = <<y'_1, \psi'_1> \dots <y'_r, \psi'_r>>$ and for each $i \in dom(k')$, $j \in dom(h')$ $x'_i = y'_j \to \varphi'_i = \psi'_j$.\\

We can also notice that $\sigma' = \sigma + (x,s)$, $\rho' = \rho + (x,s)$, $x \in \mathcal{V} - var(\sigma)$, $x \in \mathcal{V} - var(\rho)$ so by lemma~\ref{L:useful-slp-ext} if we set $\sigma' = (x',z')$, $\rho' = (w', \mu')$ then for each $i \in dom(\sigma'), \ j \in dom(\rho')$ $x'_i = w'_j \to z'_i = \mu'_j$.\\

At this point we can consider that $k',h' \in K(n)$, $\phi \in S(n,k') \cap S(n,h')$, $\sigma' \in \Xi(k')$, $\rho' \in \Xi(h')$ and we can apply lemma~\ref{L:kappa-acca-sigma-ro} to obtain that $\#(h', \phi, \rho') = \#(k', \phi, \sigma')$, and this completes our proof.
\end{proof}

\bigskip

\begin{lemma}\label{L:meaning-kept-exists}
Given
\begin{itemize}
\item a positive integer $n$;
\item $k \in K(n)$;
\item $\varphi \in E_s(n,k)$;
\item $x \in \mathcal{V} - var(k)$;
\end{itemize}

such that if we define $k' = k + <x, \varphi>$ then $k' \in K(n)$, and given also $\phi \in S(n,k')$, we have that $t = \exists(x:\varphi, \phi) \in E(n+1,k)$.\\

Given $\sigma \in \Xi(k)$ we have also
\[ \#(k, t, \sigma) = \text{exists} \ \sigma' \in \Xi(k'): \sigma \sqsubseteq \sigma' \ \#(k', \phi, \sigma') \ . \]
\end{lemma}

\begin{proof}

If $t \in E(n,k) \cup E_b(n+1,k)$ then $t \in E(n+1,k)$, else $t \in E_\exists(n+1,k) \subseteq E(n+1,k)$.\\

Using lemma~\ref{L:exp-in-sets-abcd} we have that one of the following alternatives holds: 
\begin{itemize}
\item $t \in E_a(n+1,k) \cup E_\forall(n+1,k) \cup E_\exists(n+1,k) \cup \bigcup_{c \in \mathcal{C}'} E^c(n+1,k) \cup \bigcup_{g \in \mathcal{F}} E^g(n+1,k)$;
\item there exist $m$ positive integer such that $2 \leqslant m < n+1$, $h \in K(m)$ such that $h \sqsubseteq k$, $t \in E_a(m,h) \cup E_\forall(m,h) \cup E_\exists(m,h) \cup \bigcup_{c \in \mathcal{C}'} E^c(m,h) \cup \bigcup_{g \in \mathcal{F}} E^g(m,h)$ and for each $\sigma \in \Xi(k)$  $\sigma_{/dom(h)} \in \Xi(h)$ and $\#(k,t,\sigma) = \#(h,t,\sigma_{/dom(h)})$.
\end{itemize}

\medskip

If the first alternative holds, that is $t \in E_a(n+1,k) \cup E_\forall(n+1,k) \cup E_\exists(n+1,k) \cup \bigcup_{c \in \mathcal{C}'} E^c(n+1,k) \cup \bigcup_{g \in \mathcal{F}} E^g(n+1,k)$, then clearly $t \in E_\exists(n+1,k)$. This implies that $\#(k,t,\sigma) = \#(k,t,\sigma)_{(n+1,k,\exists)}$, so in this case our proof is finished.\\

Otherwise it must be $t \in E_\exists(m,h)$. This implies that $t \in H_\exists(m,h)$, so there exist $\xi \in E_s(m-1,h)$, $v \in \mathcal{V} - var(h)$ such that if we define $h' = h + <v, \xi>$ then $h' \in K(m-1)$ and there also exists $\theta \in S(m-1,h')$ such that $t = \exists(v: \xi, \theta)$.\\

Moreover for each $\rho \in \Xi(h)$
\[ \#(h, t, \rho) = \text{exists} \ \rho' \in \Xi(h'): \rho \sqsubseteq \rho' \ \#(h', \theta, \rho') \ . \]

By our hypothesis we have that $t = \exists(x:\varphi, \phi) \in H_\exists(n+1,k)$, and by lemma~\ref{IL:HQ-recursive-aux1} it follows that 

\begin{itemize}
\item $x \in \mathcal{V} - var(k)$,
\item the set of the positive integers $r$ such that $4 < r < \ell(t)$, $t[r]$ = `,' and $d(t, r) = 1$ has just one member $r_1$.
\end{itemize}

\smallskip

If we define the following:\\
if $r_1 = 5$ then $\lambda = \epsilon$, else $\lambda = t[5, r_1 - 1]$,\\
if $r_1 = \ell(t) - 1$ then $\zeta = \epsilon$ else $\zeta = t[r_1 + 1, \ell(t) - 1]$,\\

then we have also
\begin{itemize}
\item $\lambda = \varphi \in E_s(n,k)$.
\item if we define $\kappa = k + <x, \lambda>$ then $\kappa = k' \in K(n)$ and $\zeta = \phi \in S(n,\kappa)$ .
\end{itemize}

\medskip

Given $t \in H_\exists(m,h)$, by lemma~\ref{IL:HQ-recursive-aux1} it also follows that 

\begin{itemize}
\item $x = v \in \mathcal{V} - var(h)$,
\item $\lambda = \xi \in E_s(m-1,h)$.
\item if we define $\chi = h + <x, \lambda>$ then $\chi = h' \in K(m-1)$ and $\zeta = \theta \in S(m-1,\chi)$ .
\end{itemize}

\smallskip

Finally it also follows that $\xi = \varphi$ and $\theta = \phi$.\\

Given $\sigma \in \Xi(k)$ if we define $\rho = \sigma_{/dom(h)}$ then 
\[ \#(k,t,\sigma) = \#(h,t,\rho) = \text{exists} \ \rho' \in \Xi(h'): \rho \sqsubseteq \rho' \ \#(h', \theta, \rho') \ . \]

In order to show that 
\[ \#(k, t, \sigma) = \text{exists} \ \sigma' \in \Xi(k'): \sigma \sqsubseteq \sigma' \ \#(k', \phi, \sigma') \ , \]

we just need to show that
\[ \{ \#(k', \phi, \sigma') | \  \sigma' \in \Xi(k'), \sigma \sqsubseteq \sigma'  \} = \{  \#(h', \phi, \rho') | \  \rho' \in \Xi(h'), \rho \sqsubseteq \rho'  \} \ . \]

The equality of these two sets has already been shown in lemma~\ref{L:meaning-kept-forall} and can be shown here in the exact same way.
\end{proof}

\bigskip

\begin{lemma}\label{L:quantifiers_S_h}
Let $h \in K, \ \phi \in E_s(h), \ y \in (\mathcal{V} - var(h)), \ k = h + <y, \phi>$. We have $k \in K$, and if $\vartheta \in S(k)$ then
\begin{itemize}
\item $\forall( y:\phi, \vartheta ) \in S(h)$, $\exists( y:\phi, \vartheta )  \in S(h)$;
\item for each $\rho \in \Xi(h)$ $\#(h, \forall ( y:\phi, \vartheta ), \rho)  = \text{for each} \ \sigma \in \Xi(k): \rho \sqsubseteq \sigma \ \#(k, \vartheta, \sigma)$;
\item for each $\rho \in \Xi(h)$ $\#(h, \exists ( y:\phi, \vartheta ), \rho)  = \text{exists} \ \sigma \in \Xi(k): \rho \sqsubseteq \sigma \ \#(k, \vartheta, \sigma)$.
\end{itemize}
\end{lemma}

\begin{proof}

\smallskip

Since $\phi \in E_s(h)$  there is a positive integer $n$ such that $h \in K(n), \ \phi \in E_s(n,h)$. This implies that $k \in K(n)^+ \cup K(n) = K(n+1) \subseteq K$.

\medskip

Let $\vartheta \in S(k)$. There is a positive integer $m$ such that $k \in K(m)$ and $\vartheta \in S(m,k)$. We define $p = \text{ max}\{n+1,m\}$, then we have

\begin{itemize}
\item $h \in K(p)$
\item $y \in (\mathcal{V} - var(h))$
\item $\phi \in E_s(p,h)$
\item $k \in K(p), \ \vartheta \in S(p,k)$.
\end{itemize}

\smallskip

Here we can apply lemma~\ref{L:meaning-kept-forall}, and obtain that 
\begin{itemize}
\item $\forall( y:\phi, \vartheta ) \in E(p+1,h)$;
\item for each $\rho \in \Xi(h)$ $\#(h, \forall ( y:\phi, \vartheta ), \rho)  = \text{for each} \ \sigma \in \Xi(k): \rho \sqsubseteq \sigma \ \#(k, \vartheta, \sigma)$;
\end{itemize}

As a consequence $\forall( y:\phi, \vartheta ) \in S(h)$ also holds.\\

Similarly we can apply lemma~\ref{L:meaning-kept-exists}, and obtain that 
\begin{itemize}
\item $\exists( y:\phi, \vartheta ) \in E(p+1,h)$;
\item for each $\rho \in \Xi(h)$ $\#(h, \exists ( y:\phi, \vartheta ), \rho)  = \text{exists} \ \sigma \in \Xi(k): \rho \sqsubseteq \sigma \ \#(k, \vartheta, \sigma)$;
\end{itemize}

As a consequence $\exists( y:\phi, \vartheta ) \in S(h)$ also holds.

\end{proof}

\bigskip

\begin{lemma}\label{L:cons-gamma-def}
Let $m$ be a positive integer. Let $x_1, \dots , x_m \in \mathcal{V}$, with $x_i \ne x_j$ for $i \ne j$. Let $\varphi_1, \dots , \varphi_m \in E$ and assume $H[x_1:\varphi_1, \dots , x_m:\varphi_m ]$. Let $k_0 = \epsilon$ and for each $i = 1 \dots m$ $k_i = k[x_1:\varphi_1, \dots , x_i:\varphi_i]$. Let $\varphi \in S(k_m)$. Then for each $i = 1 \dots m$ $\gamma[x_i: \varphi_i, \dots , x_m: \varphi_m, \varphi] \in S(k_{i-1})$.
\end{lemma}

\begin{proof}
By definition we have $\gamma[x_m: \varphi_m, \varphi] = \forall(x_m:\varphi_m, \varphi )$.\\

Moreover $k_{m-1} \in K$, $k_m = k_{m-1} + <x_m, \varphi_m>$, $\varphi_m \in E_s(k_{m-1})$, $x_m \in \mathcal{V} - var(k_{m-1})$. So we can apply lemma~\ref{L:quantifiers_S_h} and obtain that $\gamma[x_m: \varphi_m, \varphi] \in S(k_{m-1})$.\\

If $m > 1$ for each $i = 2 \dots m$ we have defined $\gamma[x_i: \varphi_i, \dots , x_m: \varphi_m, \varphi]$ and we can assume it is a member of $S(k_{i-1})$, by our definitions we have also
\[ \gamma[x_{i-1}: \varphi_{i-1}, \dots , x_m: \varphi_m, \varphi] = \forall( x_{i-1}: \varphi_{i-1}, \gamma[x_i: \varphi_i, \dots , x_m: \varphi_m, \varphi] ) \ . \]

\smallskip

We have also $k_{i-2} \in K$, $k_{i-1} = k_{i-2} + <x_{i-1}, \varphi_{i-1}>$, $\varphi_{i-1} \in E_s(k_{i-2})$, $x_{i-1} \in \mathcal{V} - var(k_{i-2})$. So we can apply again lemma~\ref{L:quantifiers_S_h} and obtain that $\gamma[x_{i-1}: \varphi_{i-1}, \dots , x_m: \varphi_m, \varphi] \in S(k_{i-2})$.
\end{proof}

\bigskip

\begin{theorem}\label{T:dedmet-fund}
Let $m$ be a positive integer. Let $x_1, \dots , x_m \in \mathcal{V}$, with $x_i \ne x_j$ for $i \ne j$. Let $\varphi_1, \dots , \varphi_m \in E$ and assume $H[x_1:\varphi_1, \dots , x_m:\varphi_m ]$. Let $\varphi \in S(k[x_1:\varphi_1, \dots , x_m:\varphi_m ])$. Then $\gamma[x_1: \varphi_1, \dots , x_m: \varphi_m, \varphi] \in S(\epsilon)$ and
\begin{multline*} 
\#(\gamma[x_1: \varphi_1, \dots , x_m: \varphi_m, \varphi]) \leftrightarrow \\
\leftrightarrow \text{for each } \sigma \in \Xi( k[x_1:\varphi_1, \dots , x_m:\varphi_m ] ) \ \#( k[x_1: \varphi_1, \dots , x_m: \varphi_m ] , \varphi, \sigma) 
\end{multline*}
\end{theorem}

\begin{proof}
By lemma~\ref{L:cons-gamma-def} $\gamma[x_1: \varphi_1, \dots , x_m: \varphi_m, \varphi] \in S(\epsilon)$.\\

Let $k_0 = \epsilon$ and for each $i = 1 \dots m$ $k_i = k[x_1:\varphi_1, \dots , x_i:\varphi_i]$ as in remark~\ref{R:H[]-k[]}. Then what we need to show is:
\[
\#(\gamma[x_1: \varphi_1, \dots , x_m: \varphi_m, \varphi]) \leftrightarrow \text{for each } \sigma \in \Xi( k_m) \ \#( k_m, \varphi, \sigma) \ .
\]

Let's consider that, by lemma~\ref{L:cons-gamma-def}, for each $i = 1 \dots m$ $\gamma[x_i: \varphi_i, \dots , x_m: \varphi_m, \varphi] \in S(k_{i-1})$.\\

In order to prove our result we try to show that for each $i = m \dots 1$ and for each $\rho \in \Xi(k_{i-1})$
\[
\# ( k_{i-1} , \gamma[x_i: \varphi_i, \dots , x_m: \varphi_m, \varphi], \rho ) \leftrightarrow \text{for each } \sigma \in \Xi( k_m ): \rho \sqsubseteq \sigma \ \#( k_m , \varphi, \sigma)  \ .
\]

\smallskip

We prove this by induction on $i$, starting with the case where $i=m$. Here we need to show that for each $\rho \in \Xi(k_{m-1})$
\[
\# ( k_{m-1} , \gamma[x_m: \varphi_m, \varphi], \rho ) \leftrightarrow \text{for each } \sigma \in \Xi( k_m ): \rho \sqsubseteq \sigma \ \#( k_m , \varphi, \sigma)  \ .
\]

\smallskip

Actually, by lemma~\ref{L:quantifiers_S_h},
\begin{align*}
\# ( k_{m-1} , \gamma[x_m: \varphi_m, \varphi], \rho ) &= \# ( k_{m-1} , \forall ( x_m : \varphi_m, \varphi ) , \rho ) =\\
&= \text{for each } \sigma \in \Xi( k_m ): \rho \sqsubseteq \sigma \ \#( k_m , \varphi, \sigma) \ .
\end{align*}

\smallskip

\medskip

Now suppose $m>1$, let $i = 2 \dots m$ and suppose the property holds for $i$, we show it also holds for $i - 1$. We need to prove that for each $\rho \in \Xi(k_{i-2})$
\[
\# ( k_{i-2} , \gamma[x_{i-1}: \varphi_{i-1}, \dots , x_m: \varphi_m, \varphi], \rho ) \leftrightarrow \text{for each } \sigma \in \Xi( k_m ): \rho \sqsubseteq \sigma \ \#( k_m , \varphi, \sigma)  \ .
\]

\smallskip

By our definitions we have
\begin{align*}
\# ( k_{i-2} , &\gamma[x_{i-1}: \varphi_{i-1}, \dots , x_m: \varphi_m, \varphi], \rho ) =\\
&= \# ( k_{i-2} , \forall ( x_{i-1}: \varphi_{i-1}, \gamma[x_i: \varphi_i, \dots , x_m: \varphi_m, \varphi] ) , \rho ) 
\end{align*}

By lemma~\ref{L:quantifiers_S_h} 
\begin{align*}
\# ( k_{i-2} , &\forall ( x_{i-1}: \varphi_{i-1}, \gamma[x_i: \varphi_i, \dots , x_m: \varphi_m, \varphi] ) , \rho ) =\\
&= \text{for each } \delta \in \Xi(k_{i-1}): \rho \sqsubseteq \delta \ \#(k_{i-1}, \gamma[x_i: \varphi_i, \dots , x_m: \varphi_m, \varphi], \delta)
\end{align*}

\medskip

By the inductive hypothesis given $\delta \in \Sigma(k_{i-1})$ we have 
\[
\#( k_{i-1} , \gamma[x_i: \varphi_i, \dots , x_m: \varphi_m, \varphi], \delta ) \leftrightarrow \text{for each } \sigma \in \Xi( k_m ): \delta \sqsubseteq \sigma \ \#( k_m , \varphi, \sigma)  \ .
\]

and so 
\begin{align*}
\#( k_{i-2} , &\forall ( x_{i-1}: \varphi_{i-1}, \gamma[x_i: \varphi_i, \dots , x_m: \varphi_m, \varphi] ) , \rho ) \leftrightarrow \\
& \leftrightarrow \text{for each } \delta \in \Xi(k_{i-1}): \rho \sqsubseteq \delta, \ \text{for each } \sigma \in \Xi( k_m ): \delta \sqsubseteq \sigma \ \#( k_m , \varphi, \sigma))
\end{align*}

In the end it's about proving that
\begin{align*}
\text{for each } &\delta \in \Xi(k_{i-1}): \rho \sqsubseteq \delta, \ \text{for each } \sigma \in \Xi( k_m ): \delta \sqsubseteq \sigma \ \#( k_m , \varphi, \sigma)) \leftrightarrow \\
& \leftrightarrow \text{for each } \sigma \in \Xi( k_m ): \rho \sqsubseteq \sigma \ \#( k_m , \varphi, \sigma)
\end{align*}

If we assume that for each $\sigma \in \Xi( k_m ): \rho \sqsubseteq \sigma \ \#( k_m , \varphi, \sigma)$ then clearly given $\delta \in \Xi(k_{i-1}): \rho \sqsubseteq \delta$ and given $\sigma \in \Xi( k_m ): \delta \sqsubseteq \sigma$ $\#( k_m , \varphi, \sigma)$ still holds, so this part of the proof is trivial.\\

Conversely we assume that for each $\delta \in \Xi(k_{i-1})$: $\rho \sqsubseteq \delta$ and for each $\sigma \in \Xi( k_m ): \delta \sqsubseteq \sigma$ $\#( k_m , \varphi, \sigma)$ holds. Let $\sigma \in \Xi( k_m ): \rho \sqsubseteq \sigma$, we need to prove $\#( k_m , \varphi, \sigma)$.\\

We define $\delta = \sigma_{/dom(k_{i-1})}$. By lemma~\ref{L:rest-of-state-is_state} $\delta \in \Xi(k_{i-1})$. Moreover $\delta, \rho \in \mathcal{R}(\sigma)$ and $dom(\rho) = dom(k_{i-2}) \subseteq dom(k_{i-1}) = dom(\delta)$. By lemma~\ref{L:useful3-slp} we obtain $\rho \sqsubseteq \delta$. Therefore $\#( k_m , \varphi, \sigma)$ holds.\\

This completes the proof that for each $\rho \in \Xi(k_{i-2})$
\[
\# ( k_{i-2} , \gamma[x_{i-1}: \varphi_{i-1}, \dots , x_m: \varphi_m, \varphi], \rho ) \leftrightarrow \text{for each } \sigma \in \Xi( k_m ): \rho \sqsubseteq \sigma \ \#( k_m , \varphi, \sigma)  \ .
\]

\smallskip 

And we have also completed the proof that for each $i = m \dots 1$ and for each $\rho \in \Xi(k_{i-1})$
\[
\# ( k_{i-1} , \gamma[x_i: \varphi_i, \dots , x_m: \varphi_m, \varphi], \rho ) \leftrightarrow \text{for each } \sigma \in \Xi( k_m ): \rho \sqsubseteq \sigma \ \#( k_m , \varphi, \sigma)  \ .
\]

\smallskip

It follows that for each $\rho \in \Xi(k_0)$
\[
\# ( k_0 , \gamma[x_1: \varphi_1, \dots , x_m: \varphi_m, \varphi], \rho ) \leftrightarrow \text{for each } \sigma \in \Xi( k_m ): \rho \sqsubseteq \sigma \ \#( k_m , \varphi, \sigma)  \ .
\]

\smallskip

and clearly this can be rewritten
\begin{align*}
\# ( \epsilon , \gamma[x_1: \varphi_1, \dots , x_m: \varphi_m, \varphi], \epsilon ) &\leftrightarrow \text{for each } \sigma \in \Xi( k_m ): \epsilon \sqsubseteq \sigma \ \#( k_m , \varphi, \sigma)  \ , \\
\# ( \gamma[x_1: \varphi_1, \dots , x_m: \varphi_m, \varphi] ) &\leftrightarrow \text{for each } \sigma \in \Xi( k_m ) \ \#( k_m , \varphi, \sigma)  \ .
\end{align*}
\end{proof}

\bigskip

We now need to prove the following result, which is in some way similar to~\ref{L:quantifiers_S_h} but involves the other logical connectives. After that we will be able to discuss the consistency and the completeness of our system.\\

\begin{lemma}\label{L:connectives_S_h_new}
Let $h \in K$, $\varphi_1, \varphi_2 \in S(h)$. Then

\begin{itemize}
\item $\wedge(\varphi_1, \varphi_2), \vee(\varphi_1, \varphi_2), \to(\varphi_1, \varphi_2), \leftrightarrow(\varphi_1, \varphi_2), \neg(\varphi_1) \in S(h)$;
\item for each $\rho \in \Xi(h)$ $\#(h, \wedge(\varphi_1, \varphi_2), \rho) = P_{\wedge}( \#(h, \varphi_1, \rho), \#(h, \varphi_2, \rho) )$ ;
\item for each $\rho \in \Xi(h)$ $\#(h, \vee(\varphi_1, \varphi_2), \rho) = P_{\vee}( \#(h, \varphi_1, \rho), \#(h, \varphi_2, \rho) )$ ;
\item for each $\rho \in \Xi(h)$ $\#(h, \to(\varphi_1, \varphi_2), \rho) = P_{\to}( \#(h, \varphi_1, \rho), \#(h, \varphi_2, \rho) )$ ;
\item for each $\rho \in \Xi(h)$ $\#(h, \leftrightarrow(\varphi_1, \varphi_2), \rho) = P_{\leftrightarrow}( \#(h, \varphi_1, \rho), \#(h, \varphi_2, \rho) )$ ;
\item for each $\rho \in \Xi(h)$ $\#(h, \neg(\varphi_1), \rho) = P_{\neg}( \#(h, \varphi_1, \rho) )$ .
\end{itemize}
\end{lemma}

\begin{proof}

\smallskip

For each $\rho \in \Xi(h)$ $\#(h, \varphi_1, \rho)$ is true or $\#(h, \varphi_1, \rho)$ is false; $\#(h, \varphi_2, \rho)$ is true or $\#(h, \varphi_2, \rho)$ is false.

\medskip

We recall that for each $\rho \in \Xi(h)$ $A_{\wedge}( \#(h, \varphi_1, \rho), \#(h, \varphi_2, \rho) )$, \\
$A_{\vee}( \#(h, \varphi_1, \rho), \#(h, \varphi_2, \rho) )$, $A_{\to}( \#(h, \varphi_1, \rho), \#(h, \varphi_2, \rho) )$, $A_{\leftrightarrow}( \#(h, \varphi_1, \rho), \#(h, \varphi_2, \rho) )$\\
are all defined as\\
\qquad ($\#(h, \varphi_1, \rho)$ is true or $\#(h, \varphi_1, \rho)$ is false) and ($\#(h, \varphi_2, \rho)$ is true or $\#(h, \varphi_2, \rho)$ is false).

\medskip

Therefore $A_{\wedge}( \#(h, \varphi_1, \rho), \#(h, \varphi_2, \rho) )$, $A_{\vee}( \#(h, \varphi_1, \rho), \#(h, \varphi_2, \rho) )$,\\
$A_{\to}( \#(h, \varphi_1, \rho), \#(h, \varphi_2, \rho) )$, $A_{\leftrightarrow}( \#(h, \varphi_1, \rho), \#(h, \varphi_2, \rho) )$ are all true.

\medskip

And for each $\rho \in \Xi(h)$ $A_{\neg}( \#(h, \varphi_1, \rho) )$ is true.

\medskip

Then by lemmas~\ref{L:meaning-kept-f2} and~\ref{L:meaning-kept-f1} 
\[ \wedge(\varphi_1, \varphi_2), \vee(\varphi_1, \varphi_2), \to(\varphi_1, \varphi_2), \leftrightarrow(\varphi_1, \varphi_2), \neg(\varphi_1) \in E(h) \ . \]

\smallskip

Moreover for each $\rho \in \Xi(h)$
\begin{align*}
\#(h, \wedge(\varphi_1, \varphi_2), \rho) &= P_{\wedge}( \#(h, \varphi_1, \rho), \#(h, \varphi_2, \rho) );\\
\#(h, \vee(\varphi_1, \varphi_2), \rho) &= P_{\vee}( \#(h, \varphi_1, \rho), \#(h, \varphi_2, \rho) );\\
\#(h, \to(\varphi_1, \varphi_2), \rho) &= P_{\to}( \#(h, \varphi_1, \rho), \#(h, \varphi_2, \rho) );\\
\#(h, \leftrightarrow(\varphi_1, \varphi_2), \rho) &= P_{\leftrightarrow}( \#(h, \varphi_1, \rho), \#(h, \varphi_2, \rho) );\\
\#(h, \neg(\varphi_1 ), \rho) &= P_{\neg}( \#(h, \varphi_1, \rho) ) \ .
\end{align*}

\smallskip

so
\begin{align*}
&\#(h, \wedge(\varphi_1, \varphi_2), \rho) \text{ is true or false;}\\
&\#(h, \vee(\varphi_1, \varphi_2), \rho) \text{ is true or false;}\\
&\#(h, \to(\varphi_1, \varphi_2), \rho) \text{ is true or false;}\\
&\#(h, \leftrightarrow(\varphi_1, \varphi_2), \rho) \text{ is true or false;}\\
&\#(h, \neg(\varphi_1 ), \rho) \text{ is true or false} \ .
\end{align*}

\smallskip

Therefore we get 
\[ \wedge(\varphi_1, \varphi_2), \vee(\varphi_1, \varphi_2), \to(\varphi_1, \varphi_2), \leftrightarrow(\varphi_1, \varphi_2), \neg(\varphi_1) \in S(h) \ . \]
\end{proof}

\medskip

\subsection{Consistency}\label{sec:consistency}

\smallskip

We have proved that a deductive system is sound, i.e. if we can derive a sentence $\varphi$ in our system then $\#(\varphi)$ holds. We now discuss the consistency of a deductive system.\\

Our definition of consistency implies that the symbol $\neg$ with the meaning we have associated to it in section~\ref{Ch:lang} is in the set $\mathcal{F}$ of our language. Actually in this section we have assumed that all of these symbols: $\neg, \wedge, \vee, \to, \leftrightarrow$ (with their meaning defined in section~\ref{Ch:lang}) are in our set $\mathcal{F}$, and we usually assume this since we expect in our deductions we'll frequently need these symbols.\\

A deductive system $\mathcal{D} = (\mathcal{A}, \mathcal{R})$ is said to be \emph{consistent} if and only if for each $\varphi$ sentence in $\mathcal{L}$ $(\vdash_{\mathcal{D}} \varphi)$ and $(\vdash_{\mathcal{D}} \neg (\varphi))$ aren't both true.\\

\begin{lemma}\label{L:ds-consistent}
Let $\mathcal{D} = (\mathcal{A}, \mathcal{R})$ be a deductive system in $\mathcal{L}$. Then $\mathcal{D}$ is consistent.
\end{lemma}

\begin{proof}

\medskip

Suppose there exists a sentence $\varphi$ such that $\vdash_{\mathcal{D}} \varphi$ and $\vdash_{\mathcal{D}} \neg (\varphi)$ both hold. By the soundness property we have $\#(\varphi)$ and $\#( \neg (\varphi) )$. Clearly by lemma~\ref{L:connectives_S_h_new}
\[ \#( \neg (\varphi) ) = \#(\epsilon, \neg (\varphi), \epsilon ) = P_\neg( \#(\varphi) ) = (\#(\varphi) \text{ is false} ) \ . \]

\medskip

So $\#(\varphi)$ would be true and false at the same time, a plain contradiction.
\end{proof}

\subsection{Completeness}\label{sec:completeness}

Let's now define the \emph{completeness} of a deductive system and talk a bit about this. Completeness is the converse property of soundness. A deductive system $\mathcal{D} = (\mathcal{A}, \mathcal{R})$ is said to be \emph{complete} if and only if for each $\varphi$ sentence in $\mathcal{L}$ if $\#(\varphi)$ holds then $\vdash_{\mathcal{D}} \varphi$. It was easy to prove the soundness of our system, unfortunately the topic of completeness is not that simple. Clearly, if we have defined a deductive system, there is no obvious reason to expect it is also complete.\\

Anyway, let's define a set $A$ as the set of all sentences $\varphi$ such that $\#(\varphi)$ holds. Assume $A$ is an axiom in $\mathcal{L}$ (this is a wrong assumption, but let's accept it for a moment). If we define $\mathcal{D} = ( \{ A \}, \emptyset)$ then $\mathcal{D}$ is a deductive system in $\mathcal{L}$. For each $\varphi$ sentence in $\mathcal{L}$ if $\#(\varphi)$ holds then $\varphi \in A$ and so $\vdash_{\mathcal{D}} \varphi$. In other words $\mathcal{D}$ is a complete deductive system. So, in the assumptions we made, a complete deductive system exists. Anyway as we said earlier, the assumption that $A$ is an axiom is clearly wrong, and it is wrong because there is no proof or evidence that $A$ is r.e..\\

Another trivial attempt we could make to define a complete system is the following. For each sentence $\varphi$ such that $\#(\varphi)$ holds let $\{ \varphi \}$ be an axiom in our deductive system $\mathcal{D}$. In this case for each $\varphi$ sentence in $\mathcal{L}$ if $\#(\varphi)$ holds then $\vdash_{\mathcal{D}} \varphi$ and so the system is complete. However, even in this case we have violated a requirement in the definition of a deductive system. In fact, there is no proof or evidence that our set of axioms is finite.\\

So we cannot trivially define a complete deductive system. It seems Cutland's book~\cite{Cutland} has interesting material with respect to the completeness or incompleteness of deductive systems, in chapter $8$. Actually Cutland introduces a notion of `recursively axiomatised formal system' and what he names a `simplified version of Goedel incompleteness theorem'. This theorem states that, given a recursively axiomatised formal system in which all provable statements are true, in this system there is a statement which is true but not provable (and so this system is not complete). The proof of this theorem is based on the fact that the set $\mathcal{P}$ of the provable statements of the system is recursively enumerable (r.e.) while the set $\mathcal{T}$ of the true statements of the system is not r.e.. Actually it seems to understand that Cutland refers to recursively axiomatised formal systems `of arithmetic' i.e. systems that are `adequate for making statements of ordinary arithmetic' and so include symbols like $0, 1, +, *, =$ and the logical connectives and quantifiers.\\

So, given a deductive system within our logic system, if we could describe it as a recursively axiomatised formal system of arithmetic, we would have proved that this same system is not complete. From another point of view, given a deductive system within our logic system, if one of the following conditions holds
\begin{itemize}
\item the system cannot be described as a recursively axiomatised formal system
\item the language does not include arithmetic
\end{itemize} 
we cannot state the incompleteness of the system.\\

This suggests two questions:
\begin{itemize}
\item can we describe a deductive system within our logic system as a recursively axiomatised formal system?
\item given a language that does not include arithmetic, under which conditions, if any, a deductive system within our logic system is complete?\\
\end{itemize} 

However these are non-trivial questions that I do not want to discuss in this manuscript, they are obviously of interest in further investigation of this approach.\\

In the next section we will build a deductive system and then use it to prove a given statement. This example system has many interesting and general features that can be applied also in other contexts in proving many statements. With our logic system we can certainly use many ideas to build powerful deductive systems and the example helps us to understand this. Anyway, looking at this single system, we just prove one single statement with it. We may want to prove other true statements in the same language, we may be able to do this with the axioms and rules we have provided or, to be able to do this, we may need to add other axioms or rules. However we will not make any statement about the completeness or incompleteness of the system.\\

We can also think to an alternative definition of completeness, let's call it \linebreak \emph{d-completeness}. Given a sentence $\varphi$ in $\mathcal{L}$ we say that $\varphi$ is \emph{derivable} in $\mathcal{L}$ if there exists a deductive system $\mathcal{D}$ in $\mathcal{L}$ such that $\vdash_{\mathcal{D}} \varphi$. We define the d-completeness of a deductive system $\mathcal{D}$ as follows: $\mathcal{D}$ is \emph{d-complete} if and only if for each $\varphi$ sentence in $\mathcal{L}$ if $\varphi$ is derivable in $\mathcal{L}$ then $\vdash_{\mathcal{D}} \varphi$.\\

Here we notice that if $\#(\varphi)$ holds then we can define $A = \{ \varphi \}$ and $A$ is clearly an axiom in $\mathcal{L}$. If we define $\mathcal{D} = ( \{ A \}, \emptyset)$ then $\mathcal{D}$ is a deductive system and $\vdash_{\mathcal{D}} \varphi$, so $\varphi$ is derivable in $\mathcal{L}$. Conversely is $\varphi$ is derivable in $\mathcal{L}$ then by soundess $\#(\varphi)$ holds. Therefore $\varphi$ is derivable in $\mathcal{L}$ if and only if $\#(\varphi)$ holds, so the notion of d-completeness is actually equivalent to the notion of completeness.

\section{Deductive methodology: further results}\label{Ch:dedMetFurther}

In this section we show some additional results, which can be referred to any language $\mathcal{L}=(\mathcal{V}, \mathcal{F}, \mathcal{C}, \#, \{D_1, \dots, D_n \})$ such that all of these symbols: $\neg, \wedge, \vee, \to, \leftrightarrow, \in, =$ are in our set $\mathcal{F}$. For each of these operators $f$ $A_f(x_1, \dots ,x_n)$ and $P_f(x_1, \dots ,x_n)$ are defined as specified at the beginning of section~\ref{Ch:lang}.\\

\begin{lemma}\label{L:xi_in_Ekj}
Let $m$ be a positive integer, $x_1, \dots , x_m \in \mathcal{V}$, with $x_i \ne x_j$ for $i \ne j$. Let $\varphi_1, \dots , \varphi_m \in E$, assume $H[x_1:\varphi_1, \dots , x_m:\varphi_m ]$, define $k = k[x_1:\varphi_1, \dots , x_m:\varphi_m ]$ and as usual $k_0 = \epsilon$ and for each $i = 1 \dots m$ $k_i = k[x_1:\varphi_1, \dots , x_i:\varphi_i ]$.\\
Then for each $i = 1 \dots m$, $j = i \dots m$ 
\begin{itemize}
\item $x_i \in E(k_j)$,
\item $\varphi_i \in E(k_j)$,
\item for each $\sigma \in \Xi(k_j)$ 
\begin{itemize}
\item $\sigma_{/dom(k_{i-1})} \in \Xi(k_{i-1})$,
\item $\#(k_j, x_i, \sigma) \in \#(k_{i-1}, \varphi_i, \sigma_{/dom(k_{i-1})})$,
\item $\#(k_j, \varphi_i, \sigma) = \#(k_{i-1}, \varphi_i, \sigma_{/dom(k_{i-1})})$,
\item $\#(k_j, x_i, \sigma) \in \#(k_j, \varphi_i, \sigma)$.
\end{itemize}
\end{itemize}
\end{lemma}

\begin{proof}
We prove our assertion by induction on $j$, so we begin by proving it at level $i$.\\

Since $k_i \in K$ there exists a positive integer $n$ such that $k_i \in K(n)$, and since $k_i \ne \epsilon$ we have $n \geqslant 2$. By lemma~\ref{L:k-in-km+} there exists a positive integer $q < n$ such that $k_i \in K(q)^+$. So there exist $h \in K(q), \phi \in E_s(q,h), y \in (\mathcal{V}-var(h))$ such that $k_i = h + <y,\phi>$. We have also $k_i = k_{i-1} + <x_i,\varphi_i>$ so
\[ x_i = y \in E_a(q+1,k_i) \subseteq E(q+1, k_i) \subseteq E(k_i) \ . \] 

\smallskip

For each $\sigma = \rho + (y,s) \in \Xi(k_i)$ 
\[
\#(k_i,x_i,\sigma)_{(q+1,k,a)} = s \in \#(h, \phi, \rho) = \#(k_{i-1}, \varphi_i, \rho).
\]

\smallskip

Clearly $\sigma_{/dom(\rho)} = \rho$, $dom(\rho) = dom(h) = dom(k_{i-1})$, therefore $\sigma_{/dom(k_{i-1})} = \rho$ and finally
\[
\#(k_i,x_i,\sigma) = \#(k_i,x_i,\sigma)_{(q+1,k,a)} = s \in \#(k_{i-1}, \varphi_i, \sigma_{/dom(k_{i-1})}).
\]

Since $\varphi_i \in E(k_{i-1})$ there exists a positive integer $q$ such that $k_{i-1} \in K(q)$ and $\varphi_i \in E(q, k_{i-1})$. 
Since $k_i \in K$ there also exists a positive integer $n$ such that $k_i \in K(n)$. Let $p = max \{ q,n \}$, then $\varphi_i \in E(p, k_{i-1})$ and $k_i \in K(p)$.\\

If $\varphi_i \in E(p, k_i)$ then clearly $\varphi_i \in E(k_i)$. Otherwise, since $k_i = k_{i-1} + <x_i, \varphi_i>$, $\varphi_i \in E_b(p+1, k_i) \subseteq E(p+1, k_i) \subseteq E(k_i)$.\\

At this point, given $\sigma \in \Xi(k_i)$, we observe that $k_{i-1}, k_i \in K(p+1)$, $k_{i-1} \sqsubseteq k_i$, $\varphi_i \in E(p+1,k_{i-1}) \cap E(p+1,k_i)$, $\sigma_{/dom(k_{i-1})} \in \Xi(k_{i-1})$, $\sigma_{/dom(k_{i-1})} \sqsubseteq \sigma$. Here we can apply lemma~\ref{L:kappa-acca-sigma-ro} and obtain that $\#(k_i, \varphi_i, \sigma) = \#(k_{i-1}, \varphi_i, \sigma_{/dom(k_{i-1})})$.\\

Now, in the case $i < m$, let $j = i \dots m-1$, we assume all of the following hold:
\begin{itemize}
\item $x_i \in E(k_j)$,
\item $\varphi_i \in E(k_j)$,
\item for each $\sigma \in \Xi(k_j)$ 
\begin{itemize}
\item $\sigma_{/dom(k_{i-1})} \in \Xi(k_{i-1})$,
\item $\#(k_j, x_i, \sigma) \in \#(k_{i-1}, \varphi_i, \sigma_{/dom(k_{i-1})})$,
\item $\#(k_j, \varphi_i, \sigma) = \#(k_{i-1}, \varphi_i, \sigma_{/dom(k_{i-1})})$,
\item $\#(k_j, x_i, \sigma) \in \#(k_j, \varphi_i, \sigma)$,
\end{itemize}
\end{itemize}

\smallskip

and we try to prove the same statements for $j+1$.\\

Since $k_{j+1} \in K$ there exists a positive integer $n$ such that $k_{j+1} \in K(n)$. There exists a positive integer $q$ such that $x_i \in E(q, k_j)$. Let $p = max \{ q,n \}$, then $k_{j+1} \in K(p)$ and $x_i \in E(p, k_j)$.\\

We can also observe that $k_{j+1} = k_j + (x_{j+1},\varphi_{j+1}) \in K(p) - \{ \epsilon \}$, so
\[
E_b(p+1,k_{j+1}) = \{ t | \ t \in E(p, k_j), \ t\notin E(p, k_{j+1}) \}.
\] 

\smallskip

Clearly if $x_i \in E(p, k_{j+1})$ then $x_i \in E(k_{j+1})$, otherwise $x_i \in E(p, k_j)$ and $x_i \notin E(p, k_{j+1})$, so $x_i \in E_b(p+1,k_{j+1}) \subseteq E(k_{j+1})$.\\

We now want to show that for each $\sigma \in \Xi(k_{j+1})$ $\sigma_{/dom(k_{i-1})} \in \Xi(k_{i-1})$ and
\[ \#(k_{j+1}, x_i, \sigma) \in \#(k_{i-1}, \varphi_i, \sigma_{/dom(k_{i-1})}) \ .\]

We define $\rho = \sigma_{/dom(k_j)}$, so (by lemma~\ref{L:rest-of-state-is_state}) $\rho \in \Xi(k_j)$ and by the inductive hypothesis $\rho_{/dom(k_{i-1})} \in \Xi(k_{i-1})$ $\#(k_j, x_i, \rho) \in \#(k_{i-1}, \varphi_i, \rho_{/dom(k_{i-1})})$.\\

It is also clear that $dom(k_{i-1}) \subseteq dom(k_j) \subseteq dom(k_{j+1}) = dom(\sigma)$ and therefore $\sigma_{/dom(k_{i-1})} = (\sigma_{/dom(k_j)})_{/dom(k_{i-1})} = \rho_{/dom(k_{i-1})} \in \Xi(k_{i-1})$.\\

Hence $\#(k_j, x_i, \rho) \in \#(k_{i-1}, \varphi_i, \sigma_{/dom(k_{i-1})})$ and to complete our proof that $\#(k_{j+1}, x_i, \sigma) \in \#(k_{i-1}, \varphi_i, \sigma_{/dom(k_{i-1})})$ we just need to show that $\#(k_{j+1}, x_i, \sigma) = \#(k_j, x_i, \rho)$.\\

In order to prove this we can use lemma~\ref{L:kappa-acca-sigma-ro}. In fact $k_j, k_{j+1} \in K(p+1)$, $k_j \sqsubseteq k_{j+1}$, $x_i \in E(p+1, k_j) \cap E(p+1, k_{j+1})$, $\sigma \in \Xi(k_{j+1})$, $\rho \in \Xi(k_j)$, $\rho \sqsubseteq \sigma$.\\

Since $\varphi_i \in E(k_j)$ there exists a positive integer $q$ such that $k_j \in K(q)$ and $\varphi_i \in E(q, k_j)$. 
Since $k_{j+1} \in K$ there also exists a positive integer $n$ such that $k_{j+1} \in K(n)$. Let $p = max \{ q,n \}$, then $\varphi_i \in E(p, k_j)$ and $k_{j+1} \in K(p)$.\\

If $\varphi_i \in E(p, k_{j+1})$ then clearly $\varphi_i \in E(k_{j+1})$. Otherwise, since $k_{j+1} = k_j + (x_{j+1}, \varphi_{j+1})$, $\varphi_i \in E_b(p+1, k_{j+1}) \subseteq E(p+1, k_{j+1}) \subseteq E(k_{j+1})$.\\

At this point we observe that $k_{i-1} \sqsubseteq k_{j+1}$, $\varphi_i \in E(k_{i-1}) \cap E(k_{j+1})$, $\sigma \in \Xi(k_{j+1})$, $\sigma_{/dom(k_{i-1})} \in \Xi(k_{i-1})$, $\sigma_{/dom(k_{i-1})} \sqsubseteq \sigma$. Here we can apply lemma~\ref{L:kappa-acca-sigma-ro} and obtain that $\#(k_{j+1}, \varphi_i, \sigma) = \#(k_{i-1}, \varphi_i, \sigma_{/dom(k_{i-1})})$.\\

\end{proof}

\bigskip

\begin{lemma}\label{L:in-t-varphi-in-Sk}
Suppose $k \in K$, $t, \varphi \in E(k)$ and for each $\sigma \in \Xi(k)$ $\#(k, \varphi, \sigma)$ is a set. Then 
\begin{itemize}
\item $\ \in(t, \varphi) \ \in S(k)$;
\item for each $\sigma \in \Xi(k)$ $\#(k, \, \in(t, \varphi), \sigma) =  P_{\in}(\#(k, t, \sigma), \#(k, \varphi, \sigma) )$.
\end{itemize}
\end{lemma}

\begin{proof}
This is a trivial consequence of lemma~\ref{L:meaning-kept-f2} .
\end{proof}

\bigskip

\begin{lemma}\label{L:sentence-propagation}
Let $m$ be a positive integer. Let $x_1, \dots , x_{m+1} \in \mathcal{V}$, 
with $x_i \ne x_j$ for $i \ne j$. Let $\varphi_1, \dots , \varphi_{m+1} 
\in E$ and assume $H[x_1:\varphi_1, \dots , x_{m+1}:\varphi_{m+1} ]$. 

\medskip

Define $k = k[x_1:\varphi_1, \dots , x_{m+1}:\varphi_{m+1} ]$. Of course 
$H[x_1:\varphi_1, \dots , x_m:\varphi_m ]$ also holds, we define $h = k
[x_1:\varphi_1, \dots , x_m:\varphi_m ]$. Let $\varphi \in E_s(h)$.

\medskip

Then $\varphi \in E_s(k)$ and for each $\sigma \in \Xi(k)$ $\sigma_{/dom(h)} \in \Xi(h)$, $\#(k, \varphi, \sigma) = \#(h, \varphi, \sigma_{/dom(h)})$.
\end{lemma}

\begin{proof}
Since $\varphi \in E(h)$ there exists a positive integer $q$ such that $h \in K(q)$ and $\varphi \in E(q, h)$. 
Since $k \in K$ there also exists a positive integer $n$ such that $k \in K(n)$. Let $p = max \{ q,n \}$, then $\varphi \in E(p, h)$ and $k \in K(p)$.\\

If $\varphi \in E(p, k)$ then clearly $\varphi \in E(k)$. Otherwise, since $k = k_{m+1} = k_m + <x_{m+1}, \varphi_{m+1}> = h + <x_{m+1}, \varphi_{m+1}>$, $\varphi \in E_b(p+1, k) \subseteq E(p+1, k) \subseteq E(k)$.\\

Let now $\sigma \in \Xi(k)$, by lemma~\ref{L:rest-of-state-is_state} we have $\sigma_{/dom(h)} \in \Xi(h)$, moreover $h \sqsubseteq k$, $\varphi \in E(h) \cap E(k)$, $\sigma_{/dom(h)} \sqsubseteq \sigma$. Here we can apply lemma~\ref{L:kappa-acca-sigma-ro} and obtain that $\#(k, \varphi, \sigma) = \#(h, \varphi, \sigma_{/dom(h)})$.
\end{proof}

\bigskip

\begin{lemma}\label{L:const}
Let $c \in \mathcal{C}$. For each positive integer $n$ and $k \in K(n)$ 
\begin{itemize}
\item $c \in E(n+1,k)$;
\item for each $\sigma \in \Xi(k)$ $\#(k,c,\sigma) = \#(c)$.
\end{itemize}

\end{lemma}

\begin{proof}

\medskip

The proof is by induction on $n$.

\medskip

For $n=1$ we have $k=\epsilon$ so $c \in E(1,\epsilon) = E(n,k) \subseteq E(n+1,k)$ and for each $\sigma \in \Xi(k)$ $\sigma = \epsilon$, so $\#(k,c,\sigma) = \#(\epsilon, c, \epsilon) = \#(c)$.

\medskip

Let $n$ be a positive integer and $k \in K(n+1) = K(n) \cup K(n)^+$.

\medskip

If $k \in K(n)$ then 
\begin{itemize}
\item $c \in E(n+1,k) \subseteq E(n+2,k)$;
\item for each $\sigma \in \Xi(k)$ $\#(k,c,\sigma) = \#(c)$.
\end{itemize}

\medskip

Otherwise $k \in K(n)^+$, so there exist $h \in K(n), \phi \in E_s(n,h), y \in (\mathcal{V}-var(h))$ such that $k = h + <y,\phi>$. By the inductive hypothesis 
\begin{itemize}
\item $c \in E(n+1,h)$;
\item for each $\rho \in \Xi(h)$ $\#(h,c,\rho) = \#(c)$.
\end{itemize}

\medskip

We have $c \notin E(n+1,k)$, $k = h + (y,\phi) \in K(n+1) - \{ \epsilon \}$, so $c \in E_b(n+2,k) \subseteq E(n+2,k)$ and for each $\sigma = \rho + (y,s) \in \Xi(k)$
\[ \#(k,c,\sigma) = \#(k,c,\sigma)_{(n+2,k,b)} = \#(h,c,\rho) = \#(c) \ . \]
\end{proof}

\section{Building a deductive system}\label{Ch:buildDedSystem}

In this section we will build a deductive sytem $\mathcal{D} = (\mathcal{A}, \mathcal{R})$, in order to be able to show an example of proof in the next section. The deductive system we are building can refer to any language $\mathcal{L}=(\mathcal{V}, \mathcal{F}, \mathcal{C}, \#, \{D_1, \dots, D_p \})$ such that all of these symbols: $\neg, \wedge, \vee, \to, \leftrightarrow, \in, =$ are in our set $\mathcal{F}$. For each of these operators $f$ $A_f(x_1, \dots ,x_n)$ and $P_f(x_1, \dots ,x_n)$ are defined as specified at the beginning of section~\ref{Ch:lang}.\\

We'll now list the set of axioms and rules of our deductive system. For every axiom/rule we first prove a result which ensures the soundness of the axiom/rule and then define properly the axiom/rule itself.\\

\begin{lemma}\label{L:axiom-cl-elim-new}
Let $m$ be a positive integer. Let $x_1, \dots , x_m \in \mathcal{V}$, with $x_i \ne x_j$ for $i \ne j$. Let $\varphi_1, \dots , \varphi_m \in E$ and assume $H[x_1:\varphi_1, \dots , x_m:\varphi_m ]$. Define $k = k[x_1:\varphi_1, \dots , x_m:\varphi_m ]$ and let $\varphi, \psi \in S(k)$.

\medskip

Under these assumptions we have 
\begin{itemize}
\item $\wedge(\varphi, \psi), \to \left( \wedge(\varphi, \psi), \varphi  \right), \to \left( \wedge(\varphi, \psi), \psi  \right) \in S(k)$,
\item $\gamma[ x_1: \varphi_1, \dots , x_m: \varphi_m, \to \left( \wedge(\varphi, \psi), \varphi  \right) ] \in S(\epsilon)$,
\item $\gamma[ x_1: \varphi_1, \dots , x_m: \varphi_m, \to \left( \wedge(\varphi, \psi), \psi  \right) ] \in S(\epsilon)$.
\end{itemize}

\smallskip

Moreover $\# \left( \gamma[ x_1: \varphi_1, \dots , x_m: \varphi_m, \to \left( \wedge(\varphi, \psi), \varphi  \right) ] \right)$ and\\ $\# \left( \gamma[ x_1: \varphi_1, \dots , x_m: \varphi_m, \to \left( \wedge(\varphi, \psi), \psi  \right) ] \right)$ are both true.
\end{lemma}

\begin{proof}

\medskip

Using theorem~\ref{T:dedmet-fund} and lemma~\ref{L:connectives_S_h_new} we can rewrite\\
$\# \left( \gamma[ x_1: \varphi_1, \dots , x_m: \varphi_m, \to \left( \wedge(\varphi, \psi), \varphi  \right) ] \right)$ as follows:
\begin{align} 
&\text{for each } \sigma \in \Xi(k) \ \#( k , \to \left( \wedge(\varphi, \psi), \varphi  \right), \sigma) \ , \notag  \\
&\text{for each } \sigma \in \Xi(k) \  P_{\to} \left( \# \left( k, \wedge(\varphi, \psi)  , \sigma   \right) ,  \# \left( k, \varphi, \sigma  \right)  \right) \ , \notag \\
&\text{for each } \sigma \in \Xi(k) \ P_{\to} \left( P_{\wedge} \left( \#(k, \varphi, \sigma), \#(k, \psi, \sigma)  \right) ,  \# \left( k, \varphi, \sigma  \right)  \right) \  \notag .
\end{align}

This can be expressed as

\begin{center}
for each $\sigma \in \Xi(k)$ if $\#(k, \varphi, \sigma)$ and $\#(k, \psi, \sigma)$ then $\#(k, \varphi, \sigma)$,
\end{center}

which is clearly true.\\

In the same way we can prove the truth of 
\[ \# \left( \gamma[ x_1: \varphi_1, \dots , x_m: \varphi_m, \to \left( \wedge(\varphi, \psi), \psi  \right) ] \right). \]
\end{proof}

\medskip

We can create a set~$A_{\ref{L:axiom-cl-elim-new}}$ which is the union of two sets of sentences.

\medskip

Let $G_1$ be the set of all the sentences $\gamma[ x_1: \varphi_1, \dots , x_m: \varphi_m, \to \left( \wedge(\varphi, \psi), \varphi  \right) ]$ such that
\begin{itemize}
\item $m$ is a positive integer, $x_1, \dots , x_m \in \mathcal{V}$, $x_i \ne x_j$ for $i \ne j$, $\varphi_1, \dots , \varphi_m \in E$, $H[x_1:\varphi_1, \dots , x_m:\varphi_m ]$,
\item $\varphi, \psi \in S(k[x_1:\varphi_1, \dots , x_m:\varphi_m ])$.
\end{itemize}

\smallskip

Let $G_2$ be the set of all the sentences $\gamma[ x_1: \varphi_1, \dots , x_m: \varphi_m, \to \left( \wedge(\varphi, \psi), \psi  \right) ]$ such that
\begin{itemize}
\item $m$ is a positive integer, $x_1, \dots , x_m \in \mathcal{V}$, $x_i \ne x_j$ for $i \ne j$, $\varphi_1, \dots , \varphi_m \in E$, $H[x_1:\varphi_1, \dots , x_m:\varphi_m ]$,
\item $\varphi, \psi \in S(k[x_1:\varphi_1, \dots , x_m:\varphi_m ])$.
\end{itemize}

\smallskip

Then $A_{\ref{L:axiom-cl-elim-new}}$ is the union of $G_1$ and $G_2$. Lemma~\ref{L:axiom-cl-elim-new} shows us that this set of sentences (which is a potential axiom) is `sound'. In order to use $A_{\ref{L:axiom-cl-elim-new}}$ as an axiom in our system we also need to show that $A_{\ref{L:axiom-cl-elim-new}}$ is r.e..\\

\begin{lemma}
$A_{\ref{L:axiom-cl-elim-new}}$ is r.e. .
\end{lemma}

\begin{proof}
Given a positive integer $m$ and $(x_1, \varphi_1, \dots , x_m, \varphi_m) \in R_m$ we can notice the following:

\begin{itemize}
\item $k[x_1:\varphi_1, \dots , x_m:\varphi_m ] \in K$;
\item $S(k[x_1:\varphi_1, \dots , x_m:\varphi_m ])$ is r.e.;
\item $\{(x_1, \varphi_1, \dots , x_m, \varphi_m)\} \times S(k[x_1:\varphi_1, \dots , x_m:\varphi_m ])^2$ is r.e..
\end{itemize}

So we can define the following

\[ Q_{m,2} = \bigcup_{(x_1, \varphi_1, \dots , x_m, \varphi_m) \in R_m} \{(x_1, \varphi_1, \dots , x_m, \varphi_m)\} \times S(k[x_1:\varphi_1, \dots , x_m:\varphi_m ])^2 \ . \]

Clearly $Q_{m,2} \subseteq (\Sigma^*)^{2m} \times (\Sigma^*)^2$ is r.e..\\

We can define a function $\chi$ over $(\Sigma^*)^{2m} \times (\Sigma^*)^2$ such that for each $((\psi_1, \varphi_1, \dots , \psi_m, \varphi_m), (\varphi, \psi)) \in (\Sigma^*)^{2m} \times (\Sigma^*)^2$ 
\[ \chi(((\psi_1, \varphi_1, \dots , \psi_m, \varphi_m), (\varphi, \psi))) = \gamma[\psi_1: \varphi_1, \dots , \psi_m: \varphi_m, \to \left( \wedge(\varphi, \psi), \varphi  \right)] \ . \]

\smallskip

Now $\chi$ clearly is a computable function and so the set $\{ \chi((x_1, \varphi_1, \dots , x_m, \varphi_m), (\varphi, \psi)) | \, ((x_1, \varphi_1, \dots , x_m, \varphi_m), (\varphi, \psi)) \in Q_{m,2} \}$  is a r.e. subset of $\Sigma^*$. And finally the set 

\[ \bigcup_{m \geqslant 1} \{ \chi(((x_1, \varphi_1, \dots , x_m, \varphi_m), (\varphi, \psi))) | \, (((x_1, \varphi_1, \dots , x_m, \varphi_m), (\varphi, \psi))) \in Q_{m,2} \}  \]

is itself a r.e. set. This set can obvioulsy be rewritten as follows

\[ \bigcup_{m \geqslant 1} \{ \gamma[x_1: \varphi_1, \dots , x_m: \varphi_m, \to \left( \wedge(\varphi, \psi), \varphi  \right)] | \, (((x_1, \varphi_1, \dots , x_m, \varphi_m), (\varphi, \psi))) \in Q_{m,2} \}  \]

and it should be clear at this point that this set is actually our axiom $G_1$, and so that $G_1$ is r.e..\\

In fact if $\xi \in G_1$ then there exist a positive integer $m$, $x_1, \dots , x_m \in \mathcal{V}$ such that $x_i \ne x_j$ for $i \ne j$, $\varphi_1, \dots , \varphi_m \in E$ such that $H[x_1:\varphi_1, \dots , x_m:\varphi_m ]$, $\varphi, \psi \in S(k[x_1:\varphi_1, \dots , x_m:\varphi_m ])$ such that $\xi = \gamma[ x_1: \varphi_1, \dots , x_m: \varphi_m, \to \left( \wedge(\varphi, \psi), \varphi  \right) ]$.\\

It follows that $(x_1, \varphi_1, \dots , x_m, \varphi_m) \in R_m$ and $((x_1, \varphi_1, \dots , x_m, \varphi_m), (\varphi, \psi)) \in Q_{m,2}$, so $\xi \in \{ \chi((y_1, \psi_1, \dots , y_m, \psi_m), (\phi, \theta)) | \, ((y_1, \psi_1, \dots , y_m, \psi_m), (\phi, \theta)) \in Q_{m,2} \}$, and so 

\[ \xi \in \bigcup_{p \geqslant 1} \{ \chi(((y_1, \psi_1, \dots , y_p, \psi_p), (\phi, \theta))) | \, ((y_1, \psi_1, \dots , y_p, \psi_p), (\phi, \theta)) \in Q_{p,2} \} \ . \]

Conversely if 
\[ \xi \in \bigcup_{p \geqslant 1} \{ \chi(((y_1, \psi_1, \dots , y_p, \psi_p), (\phi, \theta))) | \, ((y_1, \psi_1, \dots , y_p, \psi_p), (\phi, \theta)) \in Q_{p,2} \} \ , \]

there exist $p \geqslant 1$, $((y_1, \psi_1, \dots , y_p, \psi_p), (\phi, \theta)) \in Q_{p,2}$ such that\\
$\xi = \chi(((y_1, \psi_1, \dots , y_p, \psi_p), (\phi, \theta))) = \gamma[y_1: \psi_1, \dots , y_p: \psi_p, \to \left( \wedge(\phi, \theta), \phi  \right)]$, so $(y_1, \psi_1, \dots , y_p, \psi_p) \in R_p$, $\phi, \theta \in S(k[y_1:\psi_1, \dots , y_p:\psi_p ])$.\\

So $y_1, \dots , y_p \in \mathcal{V} \text{ with } y_i \ne y_j \text{ for } i \ne j, \ \psi_1, \dots , \psi_p \in E, H[y_1:\psi_1, \dots , y_p:\psi_p ].$\\

And this implies $\xi \in G_1$.\\

Similarly $G_2$ is r.e. and so $A_{\ref{L:axiom-cl-elim-new}}$ is r.e..
\end{proof}

\medskip

Then let $A_{\ref{L:axiom-cl-elim-new}} \in \mathcal{A}$.

\bigskip

\begin{lemma}\label{L:rule-il-trans-new}
Let $m$ be a positive integer. Let $x_1, \dots , x_m \in \mathcal{V}$, with $x_i \ne x_j$ for $i \ne j$. Let $\varphi_1, \dots , \varphi_m \in E$ and assume $H[x_1:\varphi_1, \dots , x_m:\varphi_m ]$. Define $k = k[x_1:\varphi_1, \dots , x_m:\varphi_m ]$ and let $\varphi, \psi, \chi \in S(k)$.

\medskip

Under these assumptions we have 
\begin{itemize}
\item $\to(\varphi, \psi), \to(\psi, \chi), \to(\varphi, \chi) \in S(k)$,
\item $\gamma[ x_1: \varphi_1, \dots , x_m: \varphi_m, \to(\varphi, \psi) ] \in S(\epsilon)$,
\item $\gamma[ x_1: \varphi_1, \dots , x_m: \varphi_m, \to(\psi, \chi) ] \in S(\epsilon)$,
\item $\gamma[ x_1: \varphi_1, \dots , x_m: \varphi_m, \to(\varphi, \chi) ] \in S(\epsilon)$.
\end{itemize}

\smallskip

Moreover if 
\begin{itemize}
\item $\#( \gamma[ x_1: \varphi_1, \dots , x_m: \varphi_m, \to(\varphi, \psi) ] )$,
\item $\#( \gamma[ x_1: \varphi_1, \dots , x_m: \varphi_m, \to(\psi, \chi) ] )$ 
\end{itemize}

then $\#( \gamma[ x_1: \varphi_1, \dots , x_m: \varphi_m, \to(\varphi, \chi) ] )$.
\end{lemma}

\begin{proof}

\medskip

We can rewrite $\#( \gamma[ x_1: \varphi_1, \dots , x_m: \varphi_m, \to(\varphi, \psi) ] )$ as follows:
\begin{align} 
&\text{for each } \sigma \in \Xi(k) \ \#( k , \to(\varphi, \psi), \sigma) \ , \notag  \\
&\text{for each } \sigma \in \Xi(k) \ P_{\to}\left( \#(k, \varphi , \sigma), \#(k, \psi , \sigma) \right)  \  \notag .
\end{align}

\medskip

And we can rewrite $\#( \gamma[ x_1: \varphi_1, \dots , x_m: \varphi_m, \to(\psi, \chi) ] )$ as follows:
\begin{align} 
&\text{for each } \sigma \in \Xi(k) \ \#( k , \to(\psi, \chi), \sigma) \ , \notag  \\
&\text{for each } \sigma \in \Xi(k) \ P_{\to}\left( \#(k, \psi , \sigma), \#(k, \chi , \sigma) \right)  \  \notag .
\end{align}

\medskip

In other words for each $\sigma \in \Xi(k)$ if $\#(k, \varphi, \sigma)$ then $\#(k, \psi, \sigma)$, and if $\#(k, \psi, \sigma)$ then $\#(k, \chi, \sigma)$. So, for each $\sigma \in \Xi(k)$, if $\#(k, \varphi, \sigma)$ then $\#(k, \chi, \sigma)$. This can be written as follows:
\begin{align} 
&\text{for each } \sigma \in \Xi(k) \ P_{\to}\left( \#(k, \varphi , \sigma), \#(k, \chi , \sigma) \right)  \ , \notag  \\
&\text{for each } \sigma \in \Xi(k) \ \#( k , \to(\varphi, \chi), \sigma) \ , \notag  \\
&\#( \gamma[ x_1: \varphi_1, \dots , x_m: \varphi_m, \to(\varphi, \chi) ] )  \  \notag .
\end{align}
\end{proof}

\medskip

We can create a set $R_{\ref{L:rule-il-trans-new}}$ as the set of all $3$-tuples
\[ \left(   
\begin{array} {l} 
{\gamma[ x_1: \varphi_1, \dots , x_m: \varphi_m, \to(\varphi, \psi) ], } \\
{\gamma[ x_1: \varphi_1, \dots , x_m: \varphi_m, \to(\psi, \chi) ], } \\
{\gamma[ x_1: \varphi_1, \dots , x_m: \varphi_m, \to(\varphi, \chi) ] } 
\end{array}
\right) \]

such that
\begin{itemize}
\item $m$ is a positive integer, $x_1, \dots , x_m \in \mathcal{V}$, $x_i \ne x_j$ for $i \ne j$, $\varphi_1, \dots , \varphi_m \in E$, $H[x_1:\varphi_1, \dots , x_m:\varphi_m ]$,
\item $\varphi, \psi, \chi \in S(k[x_1:\varphi_1, \dots , x_m:\varphi_m ])$.
\end{itemize}

\smallskip

Lemma~\ref{L:rule-il-trans-new} shows us that this set (which is a potential 2-ary rule) is `sound'. In order to use $R_{\ref{L:rule-il-trans-new}}$ as a rule in our system we also need to show that $R_{\ref{L:rule-il-trans-new}}$ is r.e..\\

\begin{lemma}
$R_{\ref{L:rule-il-trans-new}}$ is r.e. .
\end{lemma}

\begin{proof}
Given a positive integer $m$ and $(x_1, \varphi_1, \dots , x_m, \varphi_m) \in R_m$ we can notice the following:

\begin{itemize}
\item $k[x_1:\varphi_1, \dots , x_m:\varphi_m ] \in K$;
\item $S(k[x_1:\varphi_1, \dots , x_m:\varphi_m ])$ is r.e.;
\item $\{(x_1, \varphi_1, \dots , x_m, \varphi_m)\} \times S(k[x_1:\varphi_1, \dots , x_m:\varphi_m ])^3$ is r.e..
\end{itemize}

Let's define
\[ Q_{m,3} = \bigcup_{(x_1, \varphi_1, \dots , x_m, \varphi_m) \in R_m} \{(x_1, \varphi_1, \dots , x_m, \varphi_m)\} \times S(k[x_1:\varphi_1, \dots , x_m:\varphi_m ])^3 \ . \]

\smallskip

Clearly $Q_{m,3} \subseteq (\Sigma^*)^{2m} \times (\Sigma^*)^3$ is also r.e..\\

We now define three functions $\delta_{1,m}, \ \delta_{2,m}, \ \delta_{3,m}$ over $(\Sigma^*)^{2m} \times (\Sigma^*)^3$ as follows. Given $((\psi_1, \varphi_1, \dots , \psi_m, \varphi_m), (\varphi, \psi, \chi)) \in (\Sigma^*)^{2m} \times (\Sigma^*)^3$
\begin{align} 
&\delta_{1,m}((\psi_1, \varphi_1, \dots , \psi_{m}, \varphi_{m}), (\varphi, \psi, \chi)) = \gamma[ \psi_1: \varphi_1, \dots , \psi_m: \varphi_m, \to(\varphi, \psi) ] \ , \notag  \\
&\delta_{2,m}((\psi_1, \varphi_1, \dots , \psi_{m}, \varphi_{m}), (\varphi, \psi, \chi)) = \gamma[ \psi_1: \varphi_1, \dots , \psi_m: \varphi_m, \to(\psi, \chi) ] \ , \notag \\
&\delta_{3,m}((\psi_1, \varphi_1, \dots , \psi_{m}, \varphi_{m}), (\varphi, \psi, \chi)) = \gamma[ \psi_1: \varphi_1, \dots , \psi_m: \varphi_m, \to(\varphi, \chi) ] \  \notag .
\end{align}

\smallskip

All of the three functions we have defined are computable functions from $(\Sigma^*)^{2m} \times (\Sigma^*)^3$ to $\Sigma^*$. If we define a function $\delta_m$ over $(\Sigma^*)^{2m} \times (\Sigma^*)^3$ as follows: 
\[
\delta_m((\psi_1, \varphi_1, \dots , \psi_{m}, \varphi_{m}), (\varphi, \psi, \chi)) = 
\left(
\begin{array} {l}
{ \delta_{1,m}((\psi_1, \varphi_1, \dots , \psi_{m}, \varphi_{m}), (\varphi, \psi, \chi)) , } \\
{ \delta_{2,m}((\psi_1, \varphi_1, \dots , \psi_{m}, \varphi_{m}), (\varphi, \psi, \chi)) , } \\
{ \delta_{3,m}((\psi_1, \varphi_1, \dots , \psi_{m}, \varphi_{m}), (\varphi, \psi, \chi))  }
\end{array}
\right)
\]

then $\delta_m$ is a computable function from $(\Sigma^*)^{2m} \times (\Sigma^*)^3$ to $(\Sigma^*)^3$, therefore the set 
\[ D_m = \{ \delta_m((\psi_1, \varphi_1, \dots , \psi_{m}, \varphi_{m}), (\varphi, \psi, \chi)) | \ ((\psi_1, \varphi_1, \dots , \psi_{m}, \varphi_{m}), (\varphi, \psi, \chi)) \in Q_{m,3} \} \]
is a r.e. subset of $(\Sigma^*)^3$.\\

If we now consider the set $\bigcup_{m \geqslant 1} D_m$ then this is a r.e. subset of $(\Sigma^*)^3$ and actually this set is equal to our rule $R_{\ref{L:rule-il-trans-new}}$ which so is r.e. itself.\\

If $\xi = (\xi_1, \xi_2, \xi_3) \in R_{\ref{L:rule-il-trans-new}}$ then there exist a positive integer $m$, $x_1, \dots , x_m \in \mathcal{V}$, with $x_i \ne x_j$ for $i \ne j$, $\varphi_1, \dots , \varphi_m \in E$ such that $H[x_1:\varphi_1, \dots , x_m:\varphi_m ]$;
if we define $k = k[x_1:\varphi_1, \dots , x_m:\varphi_m ]$ there also exist $\varphi, \psi, \chi \in S(k)$ such that
\begin{itemize}
\item $\xi_1 = \gamma[ x_1: \varphi_1, \dots , x_m: \varphi_m, \to(\varphi, \psi) ]$,
\item $\xi_2 = \gamma[ x_1: \varphi_1, \dots , x_m: \varphi_m, \to(\psi, \chi) ]$,
\item $\xi_3 = \gamma[ x_1: \varphi_1, \dots , x_m: \varphi_m, \to(\varphi, \chi) ]$.
\end{itemize}

\smallskip

This means $(x_1, \varphi_1, \dots , x_m, \varphi_m) \in R_m$, $\varphi, \psi, \chi \in S(k[x_1:\varphi_1, \dots , x_m:\varphi_m])$, so $((x_1, \varphi_1, \dots , x_{m}, \varphi_{m}), (\varphi, \psi, \chi)) \in Q_{m,3}$ 
\begin{itemize}
\item $\xi_1 = \delta_{1,m}((x_1, \varphi_1, \dots , x_{m}, \varphi_{m}), (\varphi, \psi, \chi))$,
\item $\xi_2 = \delta_{2,m}((x_1, \varphi_1, \dots , x_{m}, \varphi_{m}), (\varphi, \psi, \chi))$,
\item $\xi_3 = \delta_{3,m}((x_1, \varphi_1, \dots , x_{m}, \varphi_{m}), (\varphi, \psi, \chi))$.
\end{itemize}

i.e. $\xi = \delta_m((x_1, \varphi_1, \dots , x_{m}, \varphi_{m}), (\varphi, \psi, \chi)) \in D_m$.\\

Conversely if there exists $p \geqslant 1$ such that $\xi \in D_p$ then there exists $((\psi_1, \varphi_1, \dots , \psi_{p}, \varphi_{p}), (\varphi, \psi, \chi)) \in Q_{p,3}$ such that $\xi = \delta_p((\psi_1, \varphi_1, \dots , \psi_{p}, \varphi_{p}), (\varphi, \psi, \chi))$.\\

It follows that $(\psi_1, \varphi_1, \dots , \psi_{p}, \varphi_{p}) \in R_p$, $\varphi, \psi, \chi \in S(k[\psi_1:\varphi_1, \dots , \psi_p:\varphi_p])$, so $\psi_1, \dots , \psi_{p} \in \mathcal{V}$, $\psi_i \ne \psi_j$ for $i \ne j$, $\varphi_1, \dots , \varphi_{p} \in E$, $H[\psi_1:\varphi_1, \dots , \psi_{p}:\varphi_{p} ]$.\\

Moreover 
\begin{multline*}
\xi = \delta_p((\psi_1, \varphi_1, \dots , \psi_{p}, \varphi_{p}), \varphi, \psi, \chi) = 
\left(
\begin{array} {l}
{ \delta_{1,p}((\psi_1, \varphi_1, \dots , \psi_{p}, \varphi_{p}), \varphi, \psi, \chi), } \\
{ \delta_{2,p}((\psi_1, \varphi_1, \dots , \psi_{p}, \varphi_{p}), \varphi, \psi, \chi), } \\
{ \delta_{3,p}((\psi_1, \varphi_1, \dots , \psi_{p}, \varphi_{p}), \varphi, \psi, \chi) }
\end{array}
\right) =\\ 
= \left(
\begin{array} {l}
{ \gamma[ \psi_1: \varphi_1, \dots , \psi_p: \varphi_p, \to(\varphi, \psi) ] , } \\
{ \gamma[ \psi_1: \varphi_1, \dots , \psi_p: \varphi_p, \to(\psi, \chi) ], } \\
{ \gamma[ \psi_1: \varphi_1, \dots , \psi_p: \varphi_p, \to(\varphi, \chi) ] }
\end{array}
\right) 
\end{multline*}

and so $\xi \in R_{\ref{L:rule-il-trans-new}}$.\\
\end{proof}

\medskip

Then let $R_{\ref{L:rule-il-trans-new}} \in \mathcal{R}$.

\bigskip

\begin{lemma}\label{L:axiom-xi-in-phii-new}
Let $m$ be a positive integer. Let $x_1, \dots , x_m \in \mathcal{V}$, with $x_i \ne x_j$ for $i \ne j$. Let $\varphi_1, \dots , \varphi_m \in E$ and assume $H[x_1:\varphi_1, \dots , x_m:\varphi_m ]$. Define $k = k[x_1:\varphi_1, \dots , x_m:\varphi_m ]$.

\medskip

Let $i = 1 \dots m$, then
\begin{itemize}
\item $\in(x_i, \varphi_i) \in S(k)$,
\item $\gamma[ x_1:\varphi_1, \dots , x_m:\varphi_m, \in(x_i, \varphi_i ) ] \in S(\epsilon)$,
\item $\#( \gamma[ x_1:\varphi_1, \dots , x_m:\varphi_m, \in(x_i, \varphi_i ) ] )$.
\end{itemize}
\end{lemma}

\begin{proof}

Using lemma~\ref{L:xi_in_Ekj} we obtain
\begin{itemize}
\item $x_i \in E(k)$,
\item $\varphi_i \in E(k)$,
\item for each $\sigma \in \Xi(k)$ 
\begin{itemize}
\item $\sigma_{/dom(k_{i-1})} \in \Xi(k_{i-1})$,
\item $\#(k, \varphi_i, \sigma) = \#(k_{i-1}, \varphi_i, \sigma_{/dom(k_{i-1})})$,
\item $\#(k, x_i, \sigma) \in \#(k, \varphi_i, \sigma)$.
\end{itemize}
\end{itemize}

\medskip

We have also that $\varphi_i \in E_s(k_{i-1})$, so for each $\sigma \in \Xi(k)$ $\#(k, \varphi_i, \sigma) = \#(k_{i-1}, \varphi_i, \sigma_{/dom(k_{i-1})})$ is a set. Therefore we can apply lemma~\ref{L:in-t-varphi-in-Sk} and obtain that $\in(x_i, \varphi_i) \in S(k)$. Consequently 
\[ \gamma[ x_1:\varphi_1, \dots , x_m:\varphi_m, (\in)(x_i, \varphi_i ) ] \in S(\epsilon) \ . \]

\medskip

Moreover we can rewrite $\#( \gamma[ x_1:\varphi_1, \dots , x_m:\varphi_m, (\in)(x_i, \varphi_i ) ] )$ as follows
\begin{align} 
&\text{for each } \sigma \in \Xi(k) \ \#( k , (\in)(x_i, \varphi_i), \sigma) \ , \notag  \\
&\text{for each } \sigma \in \Xi(k) \ P_{\in}( \#(k, x_i, \sigma), \#(k, \varphi_i, \sigma) ) \  \notag .
\end{align}

\smallskip

To show this we have to prove that for each $\sigma \in \Xi(k)$ $\#(k, x_i, \sigma)$ belongs to $\#(k, \varphi_i, \sigma)$. But we have just seen this is true.
\end{proof}

\medskip

We can create a set $A_{\ref{L:axiom-xi-in-phii-new}}$ which is the set of all sentences $\gamma[ x_1:\varphi_1, \dots , x_m:\varphi_m, \in(x_i, \varphi_i ) ]$ such that
\begin{itemize}
\item $m$ is a positive integer, $x_1, \dots , x_m \in \mathcal{V}$, $x_\alpha \ne x_\beta$ for $\alpha \ne \beta$, $\varphi_1, \dots , \varphi_m \in E$, $H[x_1:\varphi_1, \dots , x_m:\varphi_m ]$,
\item $i = 1 \dots m$.
\end{itemize}

\smallskip

Lemma~\ref{L:axiom-xi-in-phii-new} shows us that this set of sentences (which is a potential axiom) is `sound'. In order to use $A_{\ref{L:axiom-xi-in-phii-new}}$ as an axiom in our system we also need to show that $A_{\ref{L:axiom-xi-in-phii-new}}$ is r.e..\\

\begin{lemma}
$A_{\ref{L:axiom-xi-in-phii-new}}$ is r.e. .
\end{lemma}

\begin{proof}

Let $m$ be a positive integer and let $i = 1 \dots m$. We define a function $\chi_i$ over $(\Sigma^*)^{2m}$ such that for each $(\psi_1, \varphi_1, \dots , \psi_m, \varphi_m) \in (\Sigma^*)^{2m}$ 
\[ \chi_i(\psi_1, \varphi_1, \dots , \psi_m, \varphi_m) = \gamma[\psi_1: \varphi_1, \dots , \psi_m: \varphi_m, \in( \psi_i, \varphi_i) ] \ . \]

Now $\chi_i$ clearly is a computable function and so the set $\{ \chi_i(x_1, \varphi_1, \dots , x_m, \varphi_m) | \, (x_1, \varphi_1, \dots , x_m, \varphi_m) \in R_m \}$  is a r.e. subset of $\Sigma^*$. And moreover the set 

\[ \bigcup_{i = 1 \dots m} \{ \chi_i(x_1, \varphi_1, \dots , x_m, \varphi_m) | \, (x_1, \varphi_1, \dots , x_m, \varphi_m) \in R_m \}  \]

is itself a r.e. set. And finally the set 

\[ \bigcup_{m \geqslant 1} (\bigcup_{i = 1 \dots m} \{ \chi_i(x_1, \varphi_1, \dots , x_m, \varphi_m) | \, (x_1, \varphi_1, \dots , x_m, \varphi_m) \in R_m \} ) \]

is itself a r.e. set. This set can obviously be rewritten as follows:

\[ \bigcup_{m \geqslant 1} (\bigcup_{i = 1 \dots m} \{ \gamma[x_1: \varphi_1, \dots , x_m: \varphi_m, \in( x_i, \varphi_i) ] | \, (x_1, \varphi_1, \dots , x_m, \varphi_m) \in R_m \} ) \]

and it should be clear at this point that this set is actually our set $A_{\ref{L:axiom-xi-in-phii-new}}$.\\

In fact if $\xi \in A_{\ref{L:axiom-xi-in-phii-new}}$ then there exist a positive integer $m$, $x_1, \dots , x_m \in \mathcal{V}$ such that $x_\alpha \ne x_\beta$ for $\alpha \ne \beta$, $\varphi_1, \dots , \varphi_m \in E$ such that $H[x_1:\varphi_1, \dots , x_m:\varphi_m ]$,  $i = 1 \dots m$ such that $\xi = \gamma[ x_1:\varphi_1, \dots , x_m:\varphi_m, \in(x_i, \varphi_i ) ]$.\\

Of course this implies $(x_1, \varphi_1, \dots , x_m, \varphi_m) \in R_m$, so 
\[ \xi \in \{ \chi_i(x_1, \varphi_1, \dots , x_m, \varphi_m) | \, (x_1, \varphi_1, \dots , x_m, \varphi_m) \in R_m \} \ .  \]

And then
\begin{align} 
&\xi \in \bigcup_{j = 1 \dots m} \{ \chi_j(x_1, \varphi_1, \dots , x_m, \varphi_m) | \, (x_1, \varphi_1, \dots , x_m, \varphi_m) \in R_m \}  \ , \notag  \\
&\xi \in \bigcup_{p \geqslant 1} (\bigcup_{j = 1 \dots p} \{ \chi_j(x_1, \varphi_1, \dots , x_p, \varphi_p) | \, (x_1, \varphi_1, \dots , x_p, \varphi_p) \in R_p \}) \  \notag .
\end{align}

\medskip

Conversely if 
\[ \xi \in \bigcup_{p \geqslant 1} (\bigcup_{j = 1 \dots p} \{ \chi_j(x_1, \varphi_1, \dots , x_p, \varphi_p) | \, (x_1, \varphi_1, \dots , x_p, \varphi_p) \in R_p \})  \]

then there exists $p$ positive integer, $j = 1 \dots p$, $(x_1, \varphi_1, \dots , x_p, \varphi_p) \in R_p$ such that 
\[ \xi = \chi_j(x_1, \varphi_1, \dots , x_p, \varphi_p) = \gamma[ x_1:\varphi_1, \dots , x_p:\varphi_p, \in(x_j, \varphi_j ) ] \ . \]

Clearly we have $x_1, \dots , x_p \in \mathcal{V}$ such that $x_\alpha \ne x_\beta$ for $\alpha \ne \beta$, $\varphi_1, \dots , \varphi_p \in E$ such that $H[x_1:\varphi_1, \dots , x_p:\varphi_p ]$, so $\xi \in A_{\ref{L:axiom-xi-in-phii-new}}$.
\end{proof}

At this point let $A_{\ref{L:axiom-xi-in-phii-new}} \in \mathcal{A}$.

\bigskip

\begin{lemma}\label{L:rule-il-introd-new}
Let $m$ be a positive integer. Let $x_1, \dots , x_m \in \mathcal{V}$, with $x_i \ne x_j$ for $i \ne j$. Let $\varphi_1, \dots , \varphi_m \in E$ and assume $H[x_1:\varphi_1, \dots , x_m:\varphi_m ]$. Define $k = k[x_1:\varphi_1, \dots , x_m:\varphi_m ]$ and let $\varphi, \psi \in S(k)$.

\medskip

Under these assumptions we have 
\begin{itemize}
\item $\to(\psi, \varphi) \in S(k)$,
\item $\gamma[ x_1: \varphi_1, \dots , x_m: \varphi_m, \varphi ] \in S(\epsilon)$,
\item $\gamma[ x_1: \varphi_1, \dots , x_m: \varphi_m, \to(\psi, \varphi) ] \in S(\epsilon)$.
\end{itemize}

Moreover if $\#\left( \gamma[ x_1: \varphi_1, \dots , x_m: \varphi_m, \varphi ] \right)$ then $\#\left( \gamma[ x_1: \varphi_1, \dots , x_m: \varphi_m, \to(\psi, \varphi) ] \right)$ also holds.

\end{lemma}

\begin{proof}

\medskip

Suppose $\#\left( \gamma[ x_1: \varphi_1, \dots , x_m: \varphi_m, \varphi ] \right)$ holds. It can be rewritten as
\[ \text{for each } \sigma \in \Xi(k) \ \#( k , \varphi, \sigma) \ . \]

\smallskip

We can rewrite $\#\left( \gamma[ x_1: \varphi_1, \dots , x_m: \varphi_m, \to(\psi, \varphi) ] \right)$ as
\[ \text{for each } \sigma \in \Xi(k) \ \#( k, \to(\psi, \varphi), \sigma) \ , \]
\[ \text{for each } \sigma \in \Xi(k) \ P_{\to}( \#(k, \psi, \sigma), \#(k, \varphi, \sigma) ) \ . \]

\smallskip

For each $\sigma \in \Xi(k)$ $\#( k , \varphi, \sigma)$ holds, this implies that
\[ P_{\to}( \#(k, \psi, \sigma), \#(k, \varphi, \sigma) ) \]
holds too, therefore 
\[ \text{for each } \sigma \in \Xi(k) \ P_{\to}( \#(k, \psi, \sigma), \#(k, \varphi, \sigma) ) \ \]
also holds and this completes the proof.
\end{proof}

\medskip

We can create a set $R_{\ref{L:rule-il-introd-new}}$ as the set of all pairs
\[ \left( \gamma[ x_1: \varphi_1, \dots , x_m: \varphi_m, \varphi ], \gamma[ x_1: \varphi_1, \dots , x_m: \varphi_m, \to(\psi, \varphi) ] \right)  \]

such that
\begin{itemize}
\item $m$ is a positive integer, $x_1, \dots , x_m \in \mathcal{V}$, $x_i \ne x_j$ for $i \ne j$, $\varphi_1, \dots , \varphi_m \in E$, $H[x_1:\varphi_1, \dots , x_m:\varphi_m ]$,
\item $\varphi, \psi \in S(k[x_1:\varphi_1, \dots , x_m:\varphi_m ])$.
\end{itemize}

\smallskip

Lemma~\ref{L:rule-il-introd-new} shows us that this set (which is a potential 1-ary rule) is `sound'. In order to use $R_{\ref{L:rule-il-introd-new}}$ as a rule in our system we also need to show that $R_{\ref{L:rule-il-introd-new}}$ is r.e..\\

\begin{lemma}
$R_{\ref{L:rule-il-introd-new}}$ is r.e. .
\end{lemma}

\begin{proof}
Given a positive integer $m$ and $(x_1, \varphi_1, \dots , x_m, \varphi_m) \in R_m$ we can notice the following:

\begin{itemize}
\item $k[x_1:\varphi_1, \dots , x_m:\varphi_m ] \in K$;
\item $S(k[x_1:\varphi_1, \dots , x_m:\varphi_m ])$ is r.e.;
\item $\{(x_1, \varphi_1, \dots , x_m, \varphi_m)\} \times S(k[x_1:\varphi_1, \dots , x_m:\varphi_m ])^2$ is r.e..
\end{itemize}

So we can define the following

\[ Q_{m,2} = \bigcup_{(x_1, \varphi_1, \dots , x_m, \varphi_m) \in R_m} \{(x_1, \varphi_1, \dots , x_m, \varphi_m)\} \times S(k[x_1:\varphi_1, \dots , x_m:\varphi_m ])^2 \ . \]

Clearly $Q_{m,2} \subseteq (\Sigma^*)^{2m} \times (\Sigma^*)^2$ is r.e..\\

We now define two functions $\delta_{1,m}, \ \delta_{2,m}$ over $(\Sigma^*)^{2m} \times (\Sigma^*)^2$ as follows. Given $((\psi_1, \varphi_1, \dots , \psi_m, \varphi_m), (\varphi, \psi)) \in (\Sigma^*)^{2m} \times (\Sigma^*)^2$
\begin{align} 
&\delta_{1,m}((\psi_1, \varphi_1, \dots , \psi_{m}, \varphi_{m}), (\varphi, \psi)) = \gamma[ \psi_1: \varphi_1, \dots , \psi_m: \varphi_m, \varphi ] \ , \notag  \\
&\delta_{2,m}((\psi_1, \varphi_1, \dots , \psi_{m}, \varphi_{m}), (\varphi, \psi)) = \gamma[ \psi_1: \varphi_1, \dots , \psi_m: \varphi_m, \to(\psi, \varphi) ] \  \notag .
\end{align}

\smallskip

All of the two functions we have defined are computable functions from $(\Sigma^*)^{2m} \times (\Sigma^*)^2$ to $\Sigma^*$. If we define a function $\delta_m$ over $(\Sigma^*)^{2m} \times (\Sigma^*)^2$ as follows: 
\[
\delta_m((\psi_1, \varphi_1, \dots , \psi_{m}, \varphi_{m}), (\varphi, \psi)) = 
\left(
\begin{array} {l}
{ \delta_{1,m}((\psi_1, \varphi_1, \dots , \psi_{m}, \varphi_{m}), (\varphi, \psi)) , } \\
{ \delta_{2,m}((\psi_1, \varphi_1, \dots , \psi_{m}, \varphi_{m}), (\varphi, \psi)) , } 
\end{array}
\right)
\]

then $\delta_m$ is a computable function from $(\Sigma^*)^{2m} \times (\Sigma^*)^2$ to $(\Sigma^*)^2$, therefore the set 
\[ D_m = \{ \delta_m((\psi_1, \varphi_1, \dots , \psi_{m}, \varphi_{m}), (\varphi, \psi)) | \ ((\psi_1, \varphi_1, \dots , \psi_{m}, \varphi_{m}), (\varphi, \psi)) \in Q_{m,2} \} \]
is a r.e. subset of $(\Sigma^*)^2$.\\

If we now consider the set $\bigcup_{m \geqslant 1} D_m$ then this is a r.e. subset of $(\Sigma^*)^2$ and actually this set is equal to our rule $R_{\ref{L:rule-il-introd-new}}$ which so is r.e. itself.
\end{proof}

\medskip

Then let $R_{\ref{L:rule-il-introd-new}} \in \mathcal{R}$.

\bigskip

\begin{lemma}\label{L:rule-forall-elim-rest-new}
Let $m$ be a positive integer. Let $x_1, \dots , x_{m+1} \in \mathcal{V}$, 
with $x_i \ne x_j$ for $i \ne j$. Let $\varphi_1, \dots , \varphi_{m+1} 
\in E$ and assume $H[x_1:\varphi_1, \dots , x_{m+1}:\varphi_{m+1} ]$. 

\medskip

Define $k = k[x_1:\varphi_1, \dots , x_{m+1}:\varphi_{m+1} ]$. Of course 
$H[x_1:\varphi_1, \dots , x_m:\varphi_m ]$ also holds, we define $h = k
[x_1:\varphi_1, \dots , x_m:\varphi_m ]$. Let $\chi \in S(h)$, $t \in E
(h)$, $\varphi \in E_s(h)$.

\medskip

Under these assumptions
\begin{itemize}
\item $\in(x_{m+1}, \varphi) \in S(k)$,
\item $\forall( x_{m+1}: \varphi_{m+1}, \in(x_{m+1}, 
\varphi )) \in S(h)$,
\item $\gamma[ x_1: \varphi_1, \dots , x_m: \varphi_m, \to(\chi, 
\forall( x_{m+1}: \varphi_{m+1}, \in(x_{m+1}, \varphi) ) ) ]  \in S(\epsilon)$,
\item $\in(t, \varphi_{m+1}) \in S(h)$,
\item $\gamma[ x_1: \varphi_1, \dots , x_m: \varphi_m, \to(\chi, \in
(t, \varphi_{m+1}) ) ]  \in S(\epsilon)$,
\item $\in(t, \varphi) \in S(h)$,
\item $\gamma[ x_1: \varphi_1, \dots , x_m: \varphi_m, \to(\chi, \in
(t, \varphi) ) ]  \in S(\epsilon)$.
\end{itemize}

\smallskip

Moreover if 
\begin{itemize}
\item $\#( \gamma[ x_1: \varphi_1, \dots , x_m: \varphi_m, \to(\chi, 
\forall( x_{m+1}: \varphi_{m+1}, \in(x_{m+1}, \varphi) ) ) ] )$ and
\item $\#( \gamma[ x_1: \varphi_1, \dots , x_m: \varphi_m, \to(\chi, 
\in(t, \varphi_{m+1}) ) ] )$
\end{itemize}

\smallskip

then $\#( \gamma[ x_1: \varphi_1, \dots , x_m: \varphi_m, \to(\chi, 
\in(t, \varphi) ) ] )$.
\end{lemma}

\begin{proof}

\medskip

By lemma~\ref{L:xi_in_Ekj} we obtain that $x_{m+1} \in E(k)$.

\medskip

By lemma~\ref{L:sentence-propagation}, since $\varphi \in E_s(h)$, we obtain that $\varphi \in E_s(k)$ and for each $\sigma \in \Xi(k)$ $\sigma_{/dom(h)} \in \Xi(h)$, $\#(k, \varphi, \sigma) = \#(h, \varphi, \sigma_{/dom(h)})$.

\medskip

By lemma~\ref{L:in-t-varphi-in-Sk} we obtain that $\in(x_{m+1}, \varphi) \in S(k)$.

\medskip

By lemma~\ref{L:quantifiers_S_h} we obtain $\forall( x_{m+1}: 
\varphi_{m+1}, \in(x_{m+1}, \varphi) ) \in S(h)$.

\medskip

Clearly this implies that 
\[ \gamma[ x_1: \varphi_1, \dots , x_m: \varphi_m, \to(\chi, \forall( x_{m+1}: \varphi_{m+1}, \in(x_{m+1}, \varphi) ) ) ]  \in S
(\epsilon) . \]

\medskip

Furthermore we have $t \in E(h)$, $\varphi_{m+1} \in E_s(h)$, so $\in(t, 
\varphi_{m+1}) \in S(h)$. It clearly follows that $\gamma[ x_1: \varphi_1, \dots , 
x_m: \varphi_m, \to(\chi, \in(t, \varphi_{m+1}) ) ]  \in S(\epsilon)$.

\medskip

We have also $\varphi \in E_s(h)$, so $\in(t, \varphi) \in S(h)$. It follows 
that 
\[ \gamma[ x_1: \varphi_1, \dots , x_m: \varphi_m, \to(\chi, \in(t, \varphi) ) 
]  \in S(\epsilon) . \]

\medskip

We now assume
\begin{itemize}
\item $\#( \gamma[ x_1: \varphi_1, \dots , x_m: \varphi_m, \to(\chi, \forall
( x_{m+1}: \varphi_{m+1}, \in(x_{m+1}, \varphi) ) ) ] )$ and
\item $\#( \gamma[ x_1: \varphi_1, \dots , x_m: \varphi_m, \to(\chi, \in(t, 
\varphi_{m+1}) ) ] )$
\end{itemize}

\smallskip

both hold and we try to prove
$\#( \gamma[ x_1: \varphi_1, \dots , x_m: \varphi_m, \to(\chi, \in(t, \varphi) 
) ] )$.

\medskip

We can rewrite 
\[ \#( \gamma[ x_1: \varphi_1, \dots , x_m: \varphi_m, \to(\chi, \forall( x_{m+1}: \varphi_{m+1}, \in(x_{m+1}, \varphi) ) ) ] ) \]

as
\begin{align} 
&\text{for each } \rho \in \Xi(h) \ \#\left( h , \to \left(\chi, \forall \left( x_
{m+1}: \varphi_{m+1}, \in(x_{m+1}, \varphi) \right) \right), \rho 
\right) \ , \notag  \\
&\text{for each } \rho \in \Xi(h) \ P_{\to} \left( \#\left(h, \chi , \rho \right)  , \#\left(h, 
\forall \left( x_{m+1}: \varphi_{m+1}, \in(x_{m+1}, \varphi) \right) , \rho \right) \right) \ , \notag \\
&\text{for each } \rho \in \Xi(h) \ P_{\to} \left( \#\left(h, \chi , \rho \right)  , \text{ for each } \sigma \in \Xi(k): \rho \sqsubseteq \sigma \    \# \left(k, \in(x_{m+1}, \varphi), \sigma \right) \right) \ , \notag \\
&\text{for each } \rho \in \Xi(h) \ P_{\to} \left( \#\left(h, \chi , \rho \right)  , \text{ for each } \sigma \in \Xi(k): \rho \sqsubseteq \sigma \    P_\in \left( \#(k, x_{m+1}, \sigma) , \#(k, \varphi, \sigma) \right) \right) \  \notag .
\end{align}

\medskip

We can rewrite
\[ \#( \gamma[ x_1: \varphi_1, \dots , x_m: \varphi_m, \to(\chi, \in(t, 
\varphi_{m+1}) ) ] ) \]

as
\begin{align} 
&\text{for each } \rho \in \Xi( h ) \  \#( h , \to(\chi, \in(t, \varphi_{m+1}) ), \rho) \ , \notag  \\
&\text{for each } \rho \in \Xi( h ) \  P_{\to} ( \#(h, \chi, \rho) , \#(h, \in(t, \varphi_{m+1}), 
\rho) ) \ , \notag \\
&\text{for each } \rho \in \Xi( h ) \  P_{\to} ( \#(h, \chi, \rho), P_{\in}( \#(h, t, \rho), \#(h, 
\varphi_{m+1}, \rho) ) ) \  \notag .
\end{align}

\bigskip

We can rewrite
\[ \#( \gamma[ x_1: \varphi_1, \dots , x_m: \varphi_m, \to(\chi, \in(t, 
\varphi) ) ] ) \]

as
\begin{align} 
&\text{for each } \rho \in \Xi( h ) \  \#( h , \to(\chi, \in(t, \varphi) ), \rho) \ , \notag  \\
&\text{for each } \rho \in \Xi( h ) \  P_{\to} ( \#(h, \chi, \rho) , \#(h, \in(t, \varphi), 
\rho) ) \ , \notag \\
&\text{for each } \rho \in \Xi( h ) \  P_{\to} ( \#(h, \chi, \rho), P_{\in}( \#(h, t, \rho), \#(h, 
\varphi, \rho) ) ) \  \notag .
\end{align}

\bigskip

Let $\rho \in \Xi(h)$ and let $\#(h, \chi, \rho)$. We need to show that $\#(h, 
t, \rho)$ belongs to $\#(h, \varphi, \rho)$.

\medskip

There exists a positive integer $q$ such that $k \in K(q)^+$. So there exist $g \in K(q), \phi \in E_s(q,g), y \in (\mathcal{V}-var(g))$ such that $k = g + <y,\phi>$. At the same time 
\[ k = k_{m+1} = k_m + <x_{m+1}, \varphi_{m+1}> = h + <x_{m+1}, \varphi_{m+1}> \ . \]

\smallskip

Therefore 
\begin{align*}
\Xi(k) &= \{ \delta + (y,s) | \, \delta \in \Xi(g), s \in \#(g,\phi,\delta) \} = \\
&= \{ \delta + (x_{m+1},s) | \, \delta \in \Xi(h), s \in \#(h,\varphi_{m+1},\delta) \} \, .
\end{align*}

\medskip

We have $\rho \in \Xi(h)$, $\#(h,t,\rho) \in \#(h,\varphi_{m+1},\rho)$, so $\rho + (x_{m+1},\#(h,t,\rho)) \in \Xi(k)$.

\medskip

Let $\sigma = \rho + (x_{m+1},\#(h,t,\rho)) \in \Xi(k)$, clearly $\rho \sqsubseteq \sigma$, so $\#(k, x_{m+1}, \sigma)$ belongs to $\#(k, \varphi, \sigma)$. And we have also
\begin{align} 
&x_{m+1} = y \in E_a(q+1, k) \subseteq E(q+1, k) \ , \notag  \\
&\#(k, x_{m+1}, \sigma) = \#(k, x_{m+1}, \sigma)_{(q+1,k,a)} = \#(h,t,\rho) \ , \notag \\
&\#(k, \varphi, \sigma) = \#(h, \varphi, \sigma_{/dom(h)}) = \#(h, \varphi, \sigma_{/dom(\rho)}) = \#(h, \varphi, \rho) \  \notag .
\end{align}

\smallskip

Finally we obtain $\#(h,t,\rho) = \#(k, x_{m+1}, \sigma)$ belongs to $\#(k, \varphi, \sigma) = \#(h, \varphi, \rho)$.
\end{proof}

\medskip

We can create a set~$R_{\ref{L:rule-forall-elim-rest-new}}$ which is the set of all 3-tuples
\[ \left(   
\begin{array} {l} 
{ \gamma[ x_1: \varphi_1, \dots , x_m: \varphi_m, \to(\chi, \forall( x_{m+1}: \varphi_{m+1}, \in(x_{m+1}, \varphi) ) ) ], } \\
{ \gamma[ x_1: \varphi_1, \dots , x_m: \varphi_m, \to(\chi,\in(t, \varphi_{m+1}) ) ], } \\
{ \gamma[ x_1: \varphi_1, \dots , x_m: \varphi_m, \to(\chi, \in(t, \varphi) ) ] } 
\end{array}
\right) \]

such that 
\begin{itemize}
\item $m$ is a positive integer, $x_1, \dots , x_{m+1} \in \mathcal{V}$, with $x_i \ne x_j$ for $i \ne j$, $\varphi_1, \dots , \varphi_{m+1} \in E$, $H[x_1:\varphi_1, \dots , x_{m+1}:\varphi_{m+1} ]$;
\item if we define $k = k[x_1:\varphi_1, \dots , x_{m+1}:\varphi_{m+1} ]$ and $h = k[x_1:\varphi_1, \dots , x_m:\varphi_m ]$ then
\begin{itemize}
\item $\chi \in S(h)$,
\item $t \in E(h)$,
\item $\varphi \in E_s(h)$.
\end{itemize}
\end{itemize}

\smallskip 

Lemma~\ref{L:rule-forall-elim-rest-new} shows us that this set (which is a potential 2-ary rule) is `sound'. In order to use $R_{\ref{L:rule-forall-elim-rest-new}}$ as a rule in our system we also need to show that $R_{\ref{L:rule-forall-elim-rest-new}}$ is r.e..\\

\begin{lemma}
$R_{\ref{L:rule-forall-elim-rest-new}}$ is r.e. .
\end{lemma}

\begin{proof}
Given a positive integer $m$ and $(x_1, \varphi_1, \dots , x_{m+1}, \varphi_{m+1}) \in R_{m+1}$ all of the following sets are r.e.:
\begin{itemize}
\item $E(k[x_1:\varphi_1, \dots , x_m:\varphi_m ])$,
\item $S(k[x_1:\varphi_1, \dots , x_m:\varphi_m ])$,
\item $E_s(k[x_1:\varphi_1, \dots , x_m:\varphi_m ])$.
\end{itemize}

\smallskip

Therefore the following set is also r.e.:
\begin{multline*} 
\{ (x_1, \varphi_1, \dots , x_{m+1}, \varphi_{m+1}) \} \times S(k[x_1:\varphi_1, \dots , x_m:\varphi_m ]) \times E(k[x_1:\varphi_1, \dots , x_m:\varphi_m ])\\
 \times E_s(k[x_1:\varphi_1, \dots , x_m:\varphi_m ]) .
\end{multline*}

\smallskip

Let's use this temporary definition 

\begin{multline*} 
Q'_{m+1,3} = \bigcup_{(x_1, \varphi_1, \dots , x_{m+1}, \varphi_{m+1}) \in R_{m+1}} \{(x_1, \varphi_1, \dots , x_{m+1}, \varphi_{m+1})\} \times S(k[x_1:\varphi_1, \dots , x_m:\varphi_m ])\\ 
\times E(k[x_1:\varphi_1, \dots , x_m:\varphi_m ]) \times E_s(k[x_1:\varphi_1, \dots , x_m:\varphi_m ]).
\end{multline*}

With this $Q'_{m+1,3}$ is a r.e. subset of $(\Sigma^*)^{2(m + 1)} \times \Sigma^* \times \Sigma^* \times \Sigma^*$.\\

We now define three functions $\delta_{1,m}, \ \delta_{2,m}, \ \delta_{3,m}$ over $(\Sigma^*)^{2(m + 1)} \times \Sigma^* \times \Sigma^* \times \Sigma^*$ as follows. Given $((\psi_1, \varphi_1, \dots , \psi_{m+1}, \varphi_{m+1}), \chi, t, \varphi) \in (\Sigma^*)^{2(m + 1)} \times \Sigma^* \times \Sigma^* \times \Sigma^*$
\begin{multline*}
\delta_{1,m}((\psi_1, \varphi_1, \dots , \psi_{m+1}, \varphi_{m+1}), \chi, t, \varphi) =\\
\gamma[ \psi_1: \varphi_1, \dots , \psi_m: \varphi_m, \to(\chi, \forall( \psi_{m+1}: \varphi_{m+1}, \in(\psi_{m+1}, \varphi) ) ) ] \ ,
\end{multline*} 
\begin{align} 
&\delta_{2,m}((\psi_1, \varphi_1, \dots , \psi_{m+1}, \varphi_{m+1}), \chi, t, \varphi) = \gamma[ \psi_1: \varphi_1, \dots , \psi_m: \varphi_m, \to(\chi,\in(t, \varphi_{m+1}) ) ] \ , \notag  \\
&\delta_{3,m}((\psi_1, \varphi_1, \dots , \psi_{m+1}, \varphi_{m+1}), \chi, t, \varphi) = \gamma[ \psi_1: \varphi_1, \dots , \psi_m: \varphi_m, \to(\chi, \in(t, \varphi) ) ] \  \notag .
\end{align}

\smallskip

All of the three functions we have defined are computable functions from $(\Sigma^*)^{2(m + 1)} \times \Sigma^* \times \Sigma^* \times \Sigma^*$ to $\Sigma^*$. If we define a function $\delta_m$ over $(\Sigma^*)^{2(m + 1)} \times \Sigma^* \times \Sigma^* \times \Sigma^*$ as follows: 
\[
\delta_m((\psi_1, \varphi_1, \dots , \psi_{m+1}, \varphi_{m+1}), \chi, t, \varphi) = 
\left(
\begin{array} {l}
{ \delta_{1,m}((\psi_1, \varphi_1, \dots , \psi_{m+1}, \varphi_{m+1}), \chi, t, \varphi), } \\
{ \delta_{2,m}((\psi_1, \varphi_1, \dots , \psi_{m+1}, \varphi_{m+1}), \chi, t, \varphi), } \\
{ \delta_{3,m}((\psi_1, \varphi_1, \dots , \psi_{m+1}, \varphi_{m+1}), \chi, t, \varphi) }
\end{array}
\right)
\]

then $\delta_m$ is a computable function from $(\Sigma^*)^{2(m + 1)} \times \Sigma^* \times \Sigma^* \times \Sigma^*$ to $(\Sigma^*)^3$, therefore the set 
\[ D_m = \{ \delta_m((\psi_1, \varphi_1, \dots , \psi_{m+1}, \varphi_{m+1}), \chi, t, \varphi) | ((\psi_1, \varphi_1, \dots , \psi_{m+1}, \varphi_{m+1}), \chi, t, \varphi) \in Q'_{m+1,3} \} \]
is a r.e. subset of $(\Sigma^*)^3$.\\

If we now consider the set $\bigcup_{m \geqslant 1} D_m$ then this is a r.e. subset of $(\Sigma^*)^3$ and actually this set is equal to our rule $R_{\ref{L:rule-forall-elim-rest-new}}$ which so is r.e. itself.\\

If $\xi \in R_{\ref{L:rule-forall-elim-rest-new}}$ then there exist a positive integer $m$, $x_1, \dots , x_{m+1} \in \mathcal{V}$, with $x_i \ne x_j$ for $i \ne j$, $\varphi_1, \dots , \varphi_{m+1} \in E$ such that $H[x_1:\varphi_1, \dots , x_{m+1}:\varphi_{m+1} ]$;
if we define $k = k[x_1:\varphi_1, \dots , x_{m+1}:\varphi_{m+1} ]$ and $h = k[x_1:\varphi_1, \dots , x_m:\varphi_m ]$ 
there also exist $\chi \in S(h)$, $t \in E(h)$, $\varphi \in E_s(h)$, $\xi_1, \xi_2, \xi_3 \in \Sigma^*$ such that
\begin{itemize}
\item $\xi = (\xi_1, \xi_2, \xi_3)$
\item $\xi_1 = \gamma[ x_1: \varphi_1, \dots , x_m: \varphi_m, \to(\chi, \forall( x_{m+1}: \varphi_{m+1}, \in(x_{m+1}, \varphi) ) ) ]$,
\item $\xi_2 = \gamma[ x_1: \varphi_1, \dots , x_m: \varphi_m, \to(\chi,\in(t, \varphi_{m+1}) ) ]$,
\item $\xi_3 = \gamma[ x_1: \varphi_1, \dots , x_m: \varphi_m, \to(\chi, \in(t, \varphi) ) ]$.
\end{itemize}

\smallskip

This means that $(x_1, \varphi_1, \dots , x_{m+1}, \varphi_{m+1}) \in R_{m+1}$, $\chi \in S(k[x_1:\varphi_1, \dots , x_m:\varphi_m ])$, $t \in E(k[x_1:\varphi_1, \dots , x_m:\varphi_m ])$, $\varphi \in E_s(k[x_1:\varphi_1, \dots , x_m:\varphi_m ])$, so 
$((x_1, \varphi_1, \dots , x_{m+1}, \varphi_{m+1}), \chi, t, \varphi) \in Q'_{m+1,3}$.\\

Moreover 
\begin{itemize}
\item $\xi_1 = \delta_{1,m}((x_1, \varphi_1, \dots , x_{m+1}, \varphi_{m+1}), \chi, t, \varphi)$,
\item $\xi_2 = \delta_{2,m}((x_1, \varphi_1, \dots , x_{m+1}, \varphi_{m+1}), \chi, t, \varphi)$,
\item $\xi_3 = \delta_{3,m}((x_1, \varphi_1, \dots , x_{m+1}, \varphi_{m+1}), \chi, t, \varphi)$.
\end{itemize}

i.e. $\xi = \delta_m((x_1, \varphi_1, \dots , x_{m+1}, \varphi_{m+1}), \chi, t, \varphi) \in D_m$.\\

Conversely if there exists $p \geqslant 1$ such that $\xi \in D_p$ then there exists $((\psi_1, \varphi_1, \dots , \psi_{p+1}, \varphi_{p+1}), \chi, t, \varphi) \in Q'_{p+1,3}$ such that\\
$\xi = \delta_p((\psi_1, \varphi_1, \dots , \psi_{p+1}, \varphi_{p+1}), \chi, t, \varphi)$.\\

Since $((\psi_1, \varphi_1, \dots , \psi_{p+1}, \varphi_{p+1}), \chi, t, \varphi) \in Q'_{p+1,3}$ we have $(\psi_1, \varphi_1, \dots , \psi_{p+1}, \varphi_{p+1}) \in R_{p+1}$, $\chi \in S(k[\psi_1:\varphi_1, \dots , \psi_p:\varphi_p ])$, $t \in E(k[\psi_1:\varphi_1, \dots , \psi_p:\varphi_p ])$, $\varphi \in E_s(k[\psi_1:\varphi_1, \dots , \psi_p:\varphi_p ])$.\\

It follows that $\psi_1, \dots , \psi_{p+1} \in \mathcal{V}$, $\psi_i \ne \psi_j$ for $i \ne j$, $\varphi_1, \dots , \varphi_{p+1} \in E$, $H[\psi_1:\varphi_1, \dots , \psi_{p+1}:\varphi_{p+1} ]$.\\

Moreover 
\begin{multline*}
\xi = \delta_p((\psi_1, \varphi_1, \dots , \psi_{p+1}, \varphi_{p+1}), \chi, t, \varphi) = 
\left(
\begin{array} {l}
{ \delta_{1,p}((\psi_1, \varphi_1, \dots , \psi_{p+1}, \varphi_{p+1}), \chi, t, \varphi), } \\
{ \delta_{2,p}((\psi_1, \varphi_1, \dots , \psi_{p+1}, \varphi_{p+1}), \chi, t, \varphi), } \\
{ \delta_{3,p}((\psi_1, \varphi_1, \dots , \psi_{p+1}, \varphi_{p+1}), \chi, t, \varphi) }
\end{array}
\right) =\\ 
= \left(
\begin{array} {l}
{ \gamma[ \psi_1: \varphi_1, \dots , \psi_p: \varphi_p, \to(\chi, \forall( \psi_{p+1}: \varphi_{p+1}, \in(\psi_{p+1}, \varphi) ) ) ] , } \\
{ \gamma[ \psi_1: \varphi_1, \dots , \psi_p: \varphi_p, \to(\chi,\in(t, \varphi_{p+1}) ) ], } \\
{ \gamma[ \psi_1: \varphi_1, \dots , \psi_p: \varphi_p, \to(\chi, \in(t, \varphi) ) ] }
\end{array}
\right) 
\end{multline*}

and so $\xi \in R_{\ref{L:rule-forall-elim-rest-new}}$.\\
\end{proof}

\medskip

Then let $R_{\ref{L:rule-forall-elim-rest-new}} \in \mathcal{R}$.

\bigskip

\begin{lemma}\label{L:leverage-exists-simple-new}

Let $x_1 \in \mathcal{V}$, $\varphi_1 \in E$ and assume $H[x_1: \varphi_1]$. Define $k = k[x_1: \varphi_1]$. Let $\psi \in S(k)$ and $\varphi \in S(k) \cap S(\epsilon)$. Under these assumptions we have
\begin{itemize}
\item $\to(\psi, \varphi) \in S(k)$,
\item $\gamma [ x_1: \varphi_1, \to(\psi, \varphi) ] \in S(\epsilon)$,
\item $\exists \left( x_1: \varphi_1, \psi \right) \in S(\epsilon)$,
\item $\to \left( \exists \left( x_1: \varphi_1, \psi \right), \varphi \right) \in S(\epsilon)$.
\end{itemize}

\medskip

Moreover if $\#( \gamma [ x_1: \varphi_1, \to(\psi, \varphi) ] )$ then $\#( \to \left( \exists \left( x_1: \varphi_1, \psi \right), \varphi \right)  )$.
\end{lemma}

\begin{proof}

\medskip

Suppose $\#( \gamma [ x_1: \varphi_1, (\to)(\psi, \varphi) ] )$. We have
\begin{align} 
&\text{for each } \sigma \in \Xi(k) \ \#(k, \to(\psi, \varphi), \sigma) \ , \notag  \\
&\text{for each } \sigma \in \Xi(k) \ P_\to ( \#(k, \psi, \sigma), \#(k, \varphi, \sigma) ) \  \notag .
\end{align}

\medskip

In turn $\#( \to \left( \exists \left( \{ \} ( x_1: \varphi_1, \psi )  \right), \varphi \right)  )$ can be rewritten as
\begin{align} 
&\#(\epsilon, \to \left( \exists \left( x_1: \varphi_1, \psi \right), \varphi \right), \epsilon ) \ , \notag  \\
&P_\to ( \#(\epsilon, \exists \left( x_1: \varphi_1, \psi \right), \epsilon), \#(\epsilon, \varphi, \epsilon)  ) \ , \notag \\
&P_\to ( \#(\exists \left( x_1: \varphi_1, \psi \right)), \#(\varphi)  ) \ , \notag \\
&P_\to ( \text{exists } \sigma \in \Xi(k): \ \#(k, \psi, \sigma) , \#(\varphi)  ) \  \notag .
\end{align}

\medskip

In order to prove the last statement, we suppose there exists $\sigma \in \Xi(k)$ such that $\#(k, \psi, \sigma)$. This implies $\#(k, \varphi, \sigma)$, but we need to show that $\#(\varphi)$ holds.

\medskip

To perform this step we can use lemma~\ref{L:kappa-acca-sigma-ro}. In fact there exists a positive integer $q$ such that $\epsilon, k \in K(q)$, $\varphi \in E(q, \epsilon) \cap E(q,k)$. Moreover $\epsilon \sqsubseteq k$, $\epsilon \in \Xi(\epsilon)$, $\sigma \in \Xi(k)$, $\epsilon \sqsubseteq \sigma$ so by lemma~\ref{L:kappa-acca-sigma-ro} $\#(k, \varphi, \sigma) = \#(\epsilon, \varphi, \epsilon) = \#(\varphi)$.
\end{proof}

\medskip

We can create a set~$R_{\ref{L:leverage-exists-simple-new}}$ as the set of all pairs
\[ \left(   
\begin{array} {l} 
{ \gamma [ x_1: \varphi_1, \to(\psi, \varphi) ], } \\
{ \to \left( \exists \left( x_1: \varphi_1, \psi \right), \varphi \right) } 
\end{array}
\right) \]

such that $x_1 \in \mathcal{V}$, $\varphi_1 \in E$, $H[x_1: \varphi_1]$, $\psi \in S(k[x_1: \varphi_1])$ and $\varphi \in S(k[x_1: \varphi_1]) \cap S(\epsilon)$.\\

Lemma~\ref{L:leverage-exists-simple-new} shows us that this set (which is a potential 1-ary rule) is `sound'. In order to use $R_{\ref{L:leverage-exists-simple-new}}$ as a rule in our system we also need to show that $R_{\ref{L:leverage-exists-simple-new}}$ is r.e..\\

\begin{lemma}
$R_{\ref{L:leverage-exists-simple-new}}$ is r.e. .
\end{lemma}

\begin{proof}
Given $(x_1, \varphi_1) \in R_1$ all of the following sets are r.e.: 
\begin{itemize}
\item $S(k[x_1: \varphi_1])$,
\item $S(k[x_1: \varphi_1]) \cap S(\epsilon)$,
\item $\{ (x_1, \varphi_1) \} \times S(k[x_1: \varphi_1]) \times (S(k[x_1: \varphi_1]) \cap S(\epsilon))$.
\end{itemize}

\medskip

Let's use this temporary definition 
\[ Q'_{1,2} = \bigcup_{(x_1, \varphi_1) \in R_1} \{ (x_1, \varphi_1) \} \times S(k[x_1: \varphi_1]) \times (S(k[x_1: \varphi_1]) \cap S(\epsilon)) \ . \] 

With this $Q'_{1,2}$ is a r.e. subset of $(\Sigma^*)^2 \times \Sigma^* \times \Sigma^*$.\\

We now define two functions $\delta_{1,1}, \ \delta_{2,1}$ over $(\Sigma^*)^2 \times \Sigma^* \times \Sigma^*$ as follows:
\begin{align} 
&\delta_{1,1}((\psi_1, \varphi_1), \psi , \varphi) = \gamma [ \psi_1: \varphi_1, \to(\psi, \varphi) ] \ , \notag  \\
&\delta_{2,1}((\psi_1, \varphi_1), \psi , \varphi) = \ \to \left( \exists \left( \psi_1: \varphi_1, \psi \right), \varphi \right) \  \notag .
\end{align}

\smallskip

The two functions we have defined are both computable functions from $(\Sigma^*)^2 \times \Sigma^* \times \Sigma^*$ to $\Sigma^*$. If we define a function $\delta_1$ over $(\Sigma^*)^2 \times \Sigma^* \times \Sigma^*$ as follows
\[ \delta_1((\psi_1, \varphi_1), \psi , \varphi) = 
\left(
\begin{array} {l}
{ \delta_{1,1}((\psi_1, \varphi_1), \psi , \varphi), } \\
{ \delta_{2,1}((\psi_1, \varphi_1), \psi , \varphi) } \\
\end{array}
\right) \ , \]

then $\delta_1$ is a computable function from $(\Sigma^*)^2 \times \Sigma^* \times \Sigma^*$ to $(\Sigma^*)^2$, therefore the set 
\[ D_1 = \{ \delta_1((\psi_1, \varphi_1), \psi , \varphi) | ((\psi_1, \varphi_1), \psi , \varphi) \in Q'_{1,2} \} \]

is a r.e. subset of $(\Sigma^*)^2$, and $D_1$ is equal to our set $R_{\ref{L:leverage-exists-simple-new}}$ which so is r.e. itself.
\end{proof}

\medskip

Then let $R_{\ref{L:leverage-exists-simple-new}} \in \mathcal{R}$.

\bigskip

\begin{lemma}\label{L:AND-I-new}
Let $m$ be a positive integer. Let $x_1, \dots , x_m \in \mathcal{V}$, with $x_i 
\ne x_j$ for $i \ne j$. Let $\varphi_1, \dots , \varphi_m \in E$ and assume $H
[x_1:\varphi_1, \dots , x_m:\varphi_m ]$. Define $k = k[x_1:\varphi_1, \dots , 
x_m:\varphi_m ]$ and let $\varphi, \psi_1, \psi_2 \in S(k)$.

\medskip

Under these assumptions we have $\to(\varphi, \psi_1), \to(\varphi, \psi_2), 
\to(\varphi, \wedge(\psi_1, \psi_2) ) \in S(k)$.

\medskip

Moreover, if 
\[ \#(\gamma[x_1: \varphi_1, \dots , x_m: \varphi_m, \to(\varphi, \psi_1) ]), \ 
\#(\gamma[x_1: \varphi_1, \dots , x_m: \varphi_m, \to(\varphi, \psi_2) ]) \]

then
\[  \#(\gamma[x_1: \varphi_1, \dots , x_m: \varphi_m, \to(\varphi, \wedge
(\psi_1, \psi_2) ) ]) \ . \] 

\end{lemma}

\begin{proof}

\smallskip

We need to show 
\[  \#(\gamma[x_1: \varphi_1, \dots , x_m: \varphi_m, \to(\varphi, \wedge
(\psi_1, \psi_2) ) ]) \ , \] 

\smallskip

that is
\begin{align} 
&\text{for each } \sigma \in \Xi(k) \ \#( k , \to(\varphi, \wedge(\psi_1, \psi_2) ) , \sigma) \ , \notag \\
&\text{for each } \sigma \in \Xi(k) \ P_{\to}( \#(k, \varphi, \sigma)  , \#(k, \wedge(\psi_1, 
\psi_2) , \sigma)  ) \ \notag , \\
&\text{for each } \sigma \in \Xi(k) \ P_{\to}( \#(k, \varphi, \sigma)  , P_{\wedge}( \#(k, \psi_1, 
\sigma), \#(k, \psi_2, \sigma)  )  ) \ .
\label{E:AND-I-1-new}
\end{align}

\smallskip

But we have
\begin{align*} 
&\#(\gamma[x_1: \varphi_1, \dots , x_m: \varphi_m, \to(\varphi, \psi_1) ]) \ , 
\\
&\text{for each } \sigma \in \Xi(k) \ \#( k , \to(\varphi, \psi_1 ) , \sigma) \ , \notag \\
&\text{for each } \sigma \in \Xi(k) \ P_{\to}( \#(k, \varphi, \sigma)  , \#(k, \psi_1 , \sigma) ) \ . 
\end{align*}

\smallskip

And we have
\begin{align*} 
&\#(\gamma[x_1: \varphi_1, \dots , x_m: \varphi_m, \to(\varphi, \psi_2) ]) \ , 
\\
&\text{for each } \sigma \in \Xi(k) \ \#( k , \to(\varphi, \psi_2 ) , \sigma) \ , \notag \\
&\text{for each } \sigma \in \Xi(k) \ P_{\to}( \#(k, \varphi, \sigma)  , \#(k, \psi_2 , \sigma) ) \ . 
\end{align*}

\smallskip

So for each $\sigma \in \Xi(k)$ if $\#(k, \varphi, \sigma)$ holds true then both 
$\#(k, \psi_1 , \sigma)$ and $\#(k, \psi_2 , \sigma)$ hold. This implies \ref
{E:AND-I-1-new} holds true in turn.
\end{proof}

\medskip

We can create a set~$R_{\ref{L:AND-I-new}}$ which is the set of all $3$-tuples
\[ \left(   
\begin{array} {l} 
{\gamma[ x_1: \varphi_1, \dots , x_m: \varphi_m, \to(\varphi, \psi_1) ], } \\
{\gamma[ x_1: \varphi_1, \dots , x_m: \varphi_m, \to(\varphi, \psi_2) ], } \\
{\gamma[ x_1: \varphi_1, \dots , x_m: \varphi_m, \to(\varphi, \wedge (\psi_1, 
\psi_2) ) ] } 
\end{array}
\right) \]

such that
\begin{itemize}
\item $m$ is a positive integer, $x_1, \dots , x_m \in \mathcal{V}$, $x_i \ne x_j$ 
for $i \ne j$, $\varphi_1, \dots , \varphi_m \in E$, $H[x_1:\varphi_1, \dots , 
x_m:\varphi_m ]$,
\item $\varphi, \psi_1, \psi_2 \in S(k[x_1:\varphi_1, \dots , x_m:\varphi_m ])$.\\
\end{itemize}

\medskip

Lemma~\ref{L:AND-I-new} shows us that this set (which is a potential 2-ary rule) is `sound'. In order to use $R_{\ref{L:AND-I-new}}$ as a rule in our system we also need to show that $R_{\ref{L:AND-I-new}}$ is r.e..\\

\begin{lemma}
$R_{\ref{L:AND-I-new}}$ is r.e. .
\end{lemma}

\begin{proof}
Given a positive integer $m$ and $(x_1, \varphi_1, \dots , x_m, \varphi_m) \in R_m$ we can notice the following:

\begin{itemize}
\item $k[x_1:\varphi_1, \dots , x_m:\varphi_m ] \in K$;
\item $S(k[x_1:\varphi_1, \dots , x_m:\varphi_m ])$ is r.e.;
\item $\{(x_1, \varphi_1, \dots , x_m, \varphi_m)\} \times S(k[x_1:\varphi_1, \dots , x_m:\varphi_m ])^3$ is r.e..
\end{itemize}

Let's define
\[ Q_{m,3} = \bigcup_{(x_1, \varphi_1, \dots , x_m, \varphi_m) \in R_m} \{(x_1, \varphi_1, \dots , x_m, \varphi_m)\} \times S(k[x_1:\varphi_1, \dots , x_m:\varphi_m ])^3 \ . \]

\smallskip

Clearly $Q_{m,3} \subseteq (\Sigma^*)^{2m} \times (\Sigma^*)^3$ is also r.e..\\

We now define three functions $\delta_{1,m}, \ \delta_{2,m}, \ \delta_{3,m}$ over $(\Sigma^*)^{2m} \times (\Sigma^*)^3$ as follows. Given $((\theta_1, \varphi_1, \dots , \theta_m, \varphi_m), (\varphi, \psi_1, \psi_2)) \in (\Sigma^*)^{2m} \times (\Sigma^*)^3$
\begin{align} 
&\delta_{1,m}((\theta_1, \varphi_1, \dots , \theta_{m}, \varphi_{m}), (\varphi, \psi_1, \psi_2)) = \gamma[ \theta_1: \varphi_1, \dots , \theta_m: \varphi_m, \to(\varphi, \psi_1) ] \ , \notag  \\
&\delta_{2,m}((\theta_1, \varphi_1, \dots , \theta_{m}, \varphi_{m}), (\varphi, \psi_1, \psi_2)) = \gamma[ \theta_1: \varphi_1, \dots , \theta_m: \varphi_m, \to(\varphi, \psi_2) ] \ , \notag \\
&\delta_{3,m}((\theta_1, \varphi_1, \dots , \theta_{m}, \varphi_{m}), (\varphi, \psi_1, \psi_2)) = \gamma[ \theta_1: \varphi_1, \dots , \theta_m: \varphi_m, \to(\varphi, \wedge (\psi_1, \psi_2) ) ]  \  \notag .
\end{align}

\smallskip

All of the three functions we have defined are computable functions from $(\Sigma^*)^{2m} \times (\Sigma^*)^3$ to $\Sigma^*$. If we define a function $\delta_m$ over $(\Sigma^*)^{2m} \times (\Sigma^*)^3$ as follows: 
\[
\delta_m((\theta_1, \varphi_1, \dots , \theta_{m}, \varphi_{m}), (\varphi, \psi_1, \psi_2)) = 
\left(
\begin{array} {l}
{ \delta_{1,m}((\theta_1, \varphi_1, \dots , \theta_{m}, \varphi_{m}), (\varphi, \psi_1, \psi_2)) , } \\
{ \delta_{2,m}((\theta_1, \varphi_1, \dots , \theta_{m}, \varphi_{m}), (\varphi, \psi_1, \psi_2)) , } \\
{ \delta_{3,m}((\theta_1, \varphi_1, \dots , \theta_{m}, \varphi_{m}), (\varphi, \psi_1, \psi_2))  }
\end{array}
\right)
\]

then $\delta_m$ is a computable function from $(\Sigma^*)^{2m} \times (\Sigma^*)^3$ to $(\Sigma^*)^3$, therefore the set 
\[ D_m = \{ \delta_m((\theta_1, \varphi_1, \dots , \theta_{m}, \varphi_{m}), (\varphi, \psi_1, \psi_2)) | \ ((\theta_1, \varphi_1, \dots , \theta_{m}, \varphi_{m}), (\varphi, \psi_1, \psi_2)) \in Q_{m,3} \} \]
is a r.e. subset of $(\Sigma^*)^3$.\\

If we now consider the set $\bigcup_{m \geqslant 1} D_m$ then this is a r.e. subset of $(\Sigma^*)^3$ and actually this set is equal to our set $R_{\ref{L:AND-I-new}}$ which so is r.e. itself.
\end{proof}

\medskip

Then let $R_{\ref{L:AND-I-new}} \in \mathcal{R}$.

\bigskip

\begin{lemma}\label{L:contradict-elim-new}
Let $m$ be a positive integer. Let $x_1, \dots , x_m \in \mathcal{V}$, with $x_i \ne x_j$ for $i \ne j$. Let $\varphi_1, \dots , \varphi_m \in E$ and assume $H[x_1:\varphi_1, \dots , x_m:\varphi_m ]$. Define $k = k[x_1:\varphi_1, \dots , x_m:\varphi_m ]$ and let $\varphi, \psi \in S(k)$.

\medskip

Under these assumptions we have 
\begin{itemize}
\item $\to \left( \varphi, \wedge \left( \psi, \neg(\psi) \right) \right), \neg (\varphi) \in S(k)$,
\item $\gamma[ x_1: \varphi_1, \dots , x_m: \varphi_m, \to \left( \varphi, \wedge \left( \psi, \neg (\psi) \right) \right) ] \in S(\epsilon)$,
\item $\gamma[ x_1: \varphi_1, \dots , x_m: \varphi_m, \neg(\varphi) ] \in S(\epsilon)$.
\end{itemize}

\smallskip

Moreover if $\#( \gamma[ x_1: \varphi_1, \dots , x_m: \varphi_m, \to \left( \varphi, \wedge \left( \psi, \neg (\psi) \right) \right) ] )$ then\\
$\#( \gamma[ x_1: \varphi_1, \dots , x_m: \varphi_m, \neg(\varphi) ] )$.
\end{lemma}

\begin{proof}

\medskip

We can rewrite $\#( \gamma[ x_1: \varphi_1, \dots , x_m: \varphi_m, \to \left( \varphi, \wedge \left( \psi, \neg (\psi) \right) \right) ] )$ as
\begin{align} 
&\text{for each } \sigma \in \Xi(k) \ \#( k , \to \left( \varphi, \wedge \left( \psi, \neg (\psi) \right) \right), \sigma) , \sigma) \ , \notag  \\
&\text{for each } \sigma \in \Xi(k) \ P_{\to} \left( \#(k, \varphi , \sigma), \#(k, \wedge \left( \psi, \neg (\psi) \right), \sigma)   \right) \ , \notag \\
&\text{for each } \sigma \in \Xi(k) \ P_{\to} \left( \#(k, \varphi , \sigma), P_{\wedge} \left( \#(k, \psi, \sigma), \#(k, \neg (\psi) , \sigma)  \right) \right) \ , \notag \\
&\text{for each } \sigma \in \Xi(k) \ P_{\to} \left( \#(k, \varphi , \sigma), P_{\wedge} \left( \#(k, \psi, \sigma), P_{\neg} (\#(k, \psi, \sigma))  \right) \right) \ \notag . 
\end{align}

\medskip

This can be expressed as `for each $\sigma \in \Xi(k)$ either $\#(k, \varphi, \sigma)$ is false or both $\#(k, \psi, \sigma)$ and ($\#(k, \psi, \sigma)$ is false) are true'.

\medskip

Since $\#(k, \psi, \sigma)$ cannot be both true and false at the same time we have that `for each $\sigma \in \Xi(k)$ $\#(k, \varphi, \sigma)$  is false'. This is formally expressed as
\begin{align} 
&\text{for each } \sigma \in \Xi(k) \ P_{\neg} ( \#(k, \varphi, \sigma) ) \ , \notag  \\
&\text{for each } \sigma \in \Xi(k) \ \#(k, \neg (\varphi), \sigma) \ \notag .
\end{align}

\medskip

which we can finally rewrite as $\#( \gamma[ x_1: \varphi_1, \dots , x_m: \varphi_m, \neg (\varphi) ] )$.
\end{proof}

\medskip

We can create a set $R_{\ref{L:contradict-elim-new}}$ which is the set of all pairs
\[ \left( \gamma[ x_1: \varphi_1, \dots , x_m: \varphi_m, \to \left( \varphi, \wedge \left( \psi, \neg (\psi) \right) \right) ], \gamma[ x_1: \varphi_1, \dots , x_m: \varphi_m, \neg (\varphi) ] \right)  \]

such that
\begin{itemize}
\item $m$ is a positive integer, $x_1, \dots , x_m \in \mathcal{V}$, $x_i \ne x_j$ for $i \ne j$, $\varphi_1, \dots , \varphi_m \in E$, $H[x_1:\varphi_1, \dots , x_m:\varphi_m ]$,
\item $\varphi, \psi \in S(k[x_1:\varphi_1, \dots , x_m:\varphi_m ])$.
\end{itemize}

\medskip

Lemma~\ref{L:contradict-elim-new} shows us that this set (which is a potential 1-ary rule) is `sound'. In order to use $R_{\ref{L:contradict-elim-new}}$ as a rule in our system we also need to show that $R_{\ref{L:contradict-elim-new}}$ is r.e..\\

\begin{lemma}
$R_{\ref{L:contradict-elim-new}}$ is r.e..
\end{lemma}

\begin{proof}
Given a positive integer $m$ and $(x_1, \varphi_1, \dots , x_m, \varphi_m) \in R_m$ we can notice the following:

\begin{itemize}
\item $k[x_1:\varphi_1, \dots , x_m:\varphi_m ] \in K$;
\item $S(k[x_1:\varphi_1, \dots , x_m:\varphi_m ])$ is r.e.;
\item $\{(x_1, \varphi_1, \dots , x_m, \varphi_m)\} \times S(k[x_1:\varphi_1, \dots , x_m:\varphi_m ])^2$ is r.e..
\end{itemize}

Let's define
\[ Q_{m,2} = \bigcup_{(x_1, \varphi_1, \dots , x_m, \varphi_m) \in R_m} \{(x_1, \varphi_1, \dots , x_m, \varphi_m)\} \times S(k[x_1:\varphi_1, \dots , x_m:\varphi_m ])^2 \ . \]

\smallskip

Clearly $Q_{m,2} \subseteq (\Sigma^*)^{2m} \times (\Sigma^*)^2$ is also r.e..\\

We now define two functions $\delta_{1,m}, \ \delta_{2,m}$ over $(\Sigma^*)^{2m} \times (\Sigma^*)^2$ as follows. Given $((\psi_1, \varphi_1, \dots , \psi_m, \varphi_m), (\varphi, \psi)) \in (\Sigma^*)^{2m} \times (\Sigma^*)^2$
\begin{align}
&\delta_{1,m}((\psi_1, \varphi_1, \dots , \psi_{m}, \varphi_{m}), (\varphi, \psi)) = \gamma[ \psi_1: \varphi_1, \dots , \psi_m: \varphi_m, \to \left( \varphi, \wedge \left( \psi, \neg (\psi) \right) \right) ] \ , \notag \\
&\delta_{2,m}((\psi_1, \varphi_1, \dots , \psi_{m}, \varphi_{m}), (\varphi, \psi)) = \gamma[ \psi_1: \varphi_1, \dots , \psi_m: \varphi_m, \neg( \varphi ) ] \  \notag . 
\end{align}

\smallskip

All of the two functions we have defined are computable functions from $(\Sigma^*)^{2m} \times (\Sigma^*)^2$ to $\Sigma^*$. If we define a function $\delta_m$ over $(\Sigma^*)^{2m} \times (\Sigma^*)^2$ as follows: 
\[
\delta_m((\psi_1, \varphi_1, \dots , \psi_{m}, \varphi_{m}), (\varphi, \psi)) = 
\left(
\begin{array} {l}
{ \delta_{1,m}((\psi_1, \varphi_1, \dots , \psi_{m}, \varphi_{m}), (\varphi, \psi)) , } \\
{ \delta_{2,m}((\psi_1, \varphi_1, \dots , \psi_{m}, \varphi_{m}), (\varphi, \psi)) , }
\end{array}
\right)
\]

then $\delta_m$ is a computable function from $(\Sigma^*)^{2m} \times (\Sigma^*)^2$ to $(\Sigma^*)^2$, therefore the set 
\[ D_m = \{ \delta_m((\psi_1, \varphi_1, \dots , \psi_{m}, \varphi_{m}), (\varphi, \psi)) | \ ((\psi_1, \varphi_1, \dots , \psi_{m}, \varphi_{m}), (\varphi, \psi)) \in Q_{m,2} \} \]
is a r.e. subset of $(\Sigma^*)^2$.\\

If we now consider the set $\bigcup_{m \geqslant 1} D_m$ then this is a r.e. subset of $(\Sigma^*)^2$ and actually this set is equal to our set $R_{\ref{L:contradict-elim-new}}$ which so is r.e. itself.\\
\end{proof}

\medskip

Then let $R_{\ref{L:contradict-elim-new}} \in \mathcal{R}$.

\bigskip

\begin{lemma}\label{L:notand-to-impl-new}
Let $m$ be a positive integer. Let $x_1, \dots , x_m \in \mathcal{V}$, with $x_i \ne x_j$ for $i \ne 
j$. Let $\varphi_1, \dots , \varphi_m \in E$ and assume $H[x_1:\varphi_1, \dots , x_m:\varphi_m ]$. 
Define $k = k[x_1:\varphi_1, \dots , x_m:\varphi_m ]$ and let $\varphi, \psi \in S(k)$.

\medskip

Under these assumptions we have 
\begin{itemize}
\item $\neg \left( \wedge ( \varphi, \psi ) \right), \to ( \varphi, \neg (\psi) ) \in S(k)$,
\item $\gamma[ x_1: \varphi_1, \dots , x_m: \varphi_m, \neg \left( \wedge ( \varphi, \psi ) 
\right) ] \in S(\epsilon)$,
\item $\gamma[ x_1: \varphi_1, \dots , x_m: \varphi_m, \to ( \varphi, \neg (\psi) ) ] \in S
(\epsilon)$.
\end{itemize}

\smallskip

Moreover if $\#( \gamma[ x_1: \varphi_1, \dots , x_m: \varphi_m, \neg \left( \wedge ( \varphi, 
\psi ) \right) ] )$ then\\
$\#( \gamma[ x_1: \varphi_1, \dots , x_m: \varphi_m, \to ( \varphi, \neg (\psi) ) ] )$.
\end{lemma}

\begin{proof}

\medskip

We can rewrite $\#( \gamma[ x_1: \varphi_1, \dots , x_m: \varphi_m, \neg \left( \wedge ( \varphi, 
\psi ) \right) ] )$ as
\begin{align} 
&\text{for each } \sigma \in \Xi(k) \ \#( k , \neg \left( \wedge ( \varphi, \psi ) \right) , \sigma) \ , \notag  \\
&\text{for each } \sigma \in \Xi(k) \ P_{\neg} ( \#(k, \wedge ( \varphi, \psi ) , \sigma) )  \ , \notag  \\
&\text{for each } \sigma \in \Xi(k) \ P_{\neg} ( P_{\wedge} ( \#(k, \varphi, \sigma), \#(k, \psi, \sigma)  ) ) \  \notag .
\end{align}

\medskip

We can rewrite $\#( \gamma[ x_1: \varphi_1, \dots , x_m: \varphi_m, \to ( \varphi, \neg (\psi) ) 
] )$ as
\begin{align} 
&\text{for each } \sigma \in \Xi(k) \ \#( k , \to ( \varphi, \neg (\psi) ), \sigma) \ , \notag  \\
&\text{for each } \sigma \in \Xi(k) \ P_{\to} ( \#(k, \varphi, \sigma), \#(k, \neg (\psi), \sigma) ) \ , \notag  \\
&\text{for each } \sigma \in \Xi(k) \ P_{\to} ( \#(k, \varphi, \sigma), P_{\neg} ( \#(k, \psi, \sigma) ) ) \  \notag .
\end{align}

\medskip

Thus if $\#( \gamma[ x_1: \varphi_1, \dots , x_m: \varphi_m, \neg \left( \wedge ( \varphi, \psi ) 
\right) ] )$ we have that `for each $\sigma \in \Xi(k)$ it is false that $\#(k, \varphi, \sigma)$ and 
$\#(k, \psi, \sigma)$ are both true'.

\medskip

In other words for each $\sigma \in \Xi(k)$ ($\#(k, \varphi, \sigma)$ is false) or ($\#(k, \psi, 
\sigma)$ is false).

\medskip

In other words for each $\sigma \in \Xi(k)$ $P_{\to} ( \#(k, \varphi, \sigma), P_{\neg} ( \#(k, \psi, 
\sigma) ) )$.

\medskip

The last condition clearly implies $\#( \gamma[ x_1: \varphi_1, \dots , x_m: \varphi_m, \to ( 
\varphi, \neg (\psi) ) ] )$.
\end{proof}

\medskip

We can create a set $R_{\ref{L:notand-to-impl-new}}$ which is the set 
of all pairs
\[ \left( \gamma[ x_1: \varphi_1, \dots , x_m: \varphi_m, \neg \left( \wedge ( \varphi, \psi ) 
\right) ], \gamma[ x_1: \varphi_1, \dots , x_m: \varphi_m, \to ( \varphi, \neg (\psi) ) ] \right) 
 \]

such that
\begin{itemize}
\item $m$ is a positive integer, $x_1, \dots , x_m \in \mathcal{V}$, $x_i \ne x_j$ for $i \ne j$, $
\varphi_1, \dots , \varphi_m \in E$, $H[x_1:\varphi_1, \dots , x_m:\varphi_m ]$,
\item $\varphi, \psi \in S(k[x_1:\varphi_1, \dots , x_m:\varphi_m ])$.
\end{itemize}

\medskip

Lemma~\ref{L:notand-to-impl-new} shows us that this set (which is a potential 1-ary rule) is `sound'. In order to use $R_{\ref{L:notand-to-impl-new}}$ as a rule in our system we also need to show that $R_{\ref{L:notand-to-impl-new}}$ is r.e..\\

\begin{lemma}
$R_{\ref{L:notand-to-impl-new}}$ is r.e..
\end{lemma}

\begin{proof}
Given a positive integer $m$ and $(x_1, \varphi_1, \dots , x_m, \varphi_m) \in R_m$ we can notice the following:

\begin{itemize}
\item $k[x_1:\varphi_1, \dots , x_m:\varphi_m ] \in K$;
\item $S(k[x_1:\varphi_1, \dots , x_m:\varphi_m ])$ is r.e.;
\item $\{(x_1, \varphi_1, \dots , x_m, \varphi_m)\} \times S(k[x_1:\varphi_1, \dots , x_m:\varphi_m ])^2$ is r.e..
\end{itemize}

Let's define
\[ Q_{m,2} = \bigcup_{(x_1, \varphi_1, \dots , x_m, \varphi_m) \in R_m} \{(x_1, \varphi_1, \dots , x_m, \varphi_m)\} \times S(k[x_1:\varphi_1, \dots , x_m:\varphi_m ])^2 \ . \]

\smallskip

Clearly $Q_{m,2} \subseteq (\Sigma^*)^{2m} \times (\Sigma^*)^2$ is also r.e..\\

We now define two functions $\delta_{1,m}, \ \delta_{2,m}$ over $(\Sigma^*)^{2m} \times (\Sigma^*)^2$ as follows. Given $((\psi_1, \varphi_1, \dots , \psi_m, \varphi_m), (\varphi, \psi)) \in (\Sigma^*)^{2m} \times (\Sigma^*)^2$
\begin{align} 
& \delta_{1,m}((\psi_1, \varphi_1, \dots , \psi_{m}, \varphi_{m}), (\varphi, \psi)) = \gamma[ \psi_1: \varphi_1, \dots , \psi_m: \varphi_m, \neg \left( \wedge ( \varphi, \psi ) \right) ] \ , \notag  \\
& \delta_{2,m}((\psi_1, \varphi_1, \dots , \psi_{m}, \varphi_{m}), (\varphi, \psi)) = \gamma[ \psi_1: \varphi_1, \dots , \psi_m: \varphi_m, \to ( \varphi, \neg (\psi) ) ] \  \notag .
\end{align}

\smallskip

All of the two functions we have defined are computable functions from $(\Sigma^*)^{2m} \times (\Sigma^*)^2$ to $\Sigma^*$. If we define a function $\delta_m$ over $(\Sigma^*)^{2m} \times (\Sigma^*)^2$ as follows: 
\[
\delta_m((\psi_1, \varphi_1, \dots , \psi_{m}, \varphi_{m}), (\varphi, \psi)) = 
\left(
\begin{array} {l}
{ \delta_{1,m}((\psi_1, \varphi_1, \dots , \psi_{m}, \varphi_{m}), (\varphi, \psi)) , } \\
{ \delta_{2,m}((\psi_1, \varphi_1, \dots , \psi_{m}, \varphi_{m}), (\varphi, \psi)) , }
\end{array}
\right)
\]

then $\delta_m$ is a computable function from $(\Sigma^*)^{2m} \times (\Sigma^*)^2$ to $(\Sigma^*)^2$, therefore the set 
\[ D_m = \{ \delta_m((\psi_1, \varphi_1, \dots , \psi_{m}, \varphi_{m}), (\varphi, \psi)) | \ ((\psi_1, \varphi_1, \dots , \psi_{m}, \varphi_{m}), (\varphi, \psi)) \in Q_{m,2} \} \]
is a r.e. subset of $(\Sigma^*)^2$.\\

If we now consider the set $\bigcup_{m \geqslant 1} D_m$ then this is a r.e. subset of $(\Sigma^*)^2$ and actually this set is equal to our set $R_{\ref{L:notand-to-impl-new}}$ which so is r.e. itself.\\
\end{proof}

\medskip

Then let $R_{\ref{L:notand-to-impl-new}} \in \mathcal{R}$.

\bigskip

\begin{lemma}\label{L:notforall-to-existsnot-new}
Let $m$ be a positive integer. Let $x_1, \dots , x_{m+1} \in \mathcal{V}$, with $x_i \ne x_j$ for $i 
\ne j$. Let $\varphi_1, \dots , \varphi_{m+1} \in E$ and assume $H[x_1:\varphi_1, \dots , x_{m+1}:
\varphi_{m+1} ]$. 

\medskip

Define $k = k[x_1:\varphi_1, \dots , x_{m+1}:\varphi_{m+1} ]$. Of course $H[x_1:\varphi_1, \dots , 
x_m:\varphi_m ]$ also holds, we define $h = k[x_1:\varphi_1, \dots , x_m:\varphi_m ]$. Let $\chi \in 
S(h)$, $\varphi \in S(k)$.

\medskip

Under these assumptions we have 
\begin{itemize}
\item $\forall ( x_{m+1} : \varphi_{m+1}, \varphi ) \in S(h)$,
\item $\neg \left( \forall ( x_{m+1} : \varphi_{m+1}, \varphi ) \right) \in S(h)$,
\item $\to \left( \chi, \neg \left( \forall ( x_{m+1} : \varphi_{m+1}, \varphi ) 
\right) \right) \in S(h)$,
\item $\gamma[ x_1: \varphi_1, \dots , x_m: \varphi_m, \to \left( \chi, \neg \left( \forall ( x_{m+1} : \varphi_{m+1}, \varphi ) \right) \right) ] \in S(\epsilon)$,
\item $\neg(\varphi) \in S(k)$,
\item $\exists ( x_{m+1} : \varphi_{m+1}, \neg(\varphi) ) \in S(h)$,
\item $\to \left( \chi, \exists ( x_{m+1} : \varphi_{m+1}, \neg(\varphi) ) \right) \in 
S(h)$,
\item $\gamma[ x_1: \varphi_1, \dots , x_m: \varphi_m, \to \left( \chi, \exists ( x_{m+1} 
: \varphi_{m+1}, \neg(\varphi) ) \right) ] \in S(\epsilon)$.
\end{itemize}

\medskip

Moreover if $\#(\gamma[ x_1: \varphi_1, \dots , x_m: \varphi_m, \to \left( \chi, \neg \left( 
\forall ( x_{m+1} : \varphi_{m+1}, \varphi ) \right) \right) ])$ then 
\[ \#( \gamma[ x_1: \varphi_1, \dots , x_m: \varphi_m, \to \left( \chi, \exists ( x_{m+1} 
: \varphi_{m+1}, \neg(\varphi) ) \right) ] ) \ . \]

\end{lemma}

\begin{proof}

\medskip

We can rewrite $\#(\gamma[ x_1: \varphi_1, \dots , x_m: \varphi_m, \to \left( \chi, \neg \left( 
\forall ( x_{m+1} : \varphi_{m+1}, \varphi ) \right) \right) ])$ as

\medskip

for each $\rho \in \Xi(h)$ $\#( h , \to \left( \chi, \neg \left( \forall ( x_{m+1} : \varphi_{m
+1}, \varphi ) \right) \right) , \rho)$,

for each $\rho \in \Xi(h)$ $P_\to \left( \#(h, \chi, \rho), \#(h, \neg \left( \forall ( x_{m+1} : 
\varphi_{m+1}, \varphi ) \right), \rho) \right)$,

for each $\rho \in \Xi(h)$ $P_\to \left( \#(h, \chi, \rho), P_{\neg} \left( \#(h, \forall \left( x_{m+1} : \varphi_{m+1}, \varphi \right), \rho)  \right) \right)$,

for each $\rho \in \Xi(h)$ $P_\to \left( \#(h, \chi, \rho), P_{\neg} \left(  
\text{for each } \sigma \in \Xi(k): \rho \sqsubseteq \sigma \ \#(k, \varphi, \sigma)
  \right) \right)$,

\medskip

We can furtherly express this as

\medskip

for each $\rho \in \Xi(h)$ if $\#(h, \chi, \rho)$ then it is false that (for each $\sigma \in \Xi
(k)$ such that $\rho \sqsubseteq \sigma$ $\#(k, \varphi, \sigma)$ holds),

for each $\rho \in \Xi(h)$ if $\#(h, \chi, \rho)$ then (there exists $\sigma \in \Xi(k)$ such that 
$\rho \sqsubseteq \sigma$ and $\#(k, \varphi, \sigma)$ is false).

\bigskip

We can rewrite $\#( \gamma[ x_1: \varphi_1, \dots , x_m: \varphi_m, \to \left( \chi, \exists ( x_{m+1} : \varphi_{m+1}, \neg(\varphi) ) \right) ] )$ as

\medskip

for each $\rho \in \Xi(h)$ $\#( h , \to \left( \chi, \exists ( x_{m+1} : \varphi_{m+1}, \neg
(\varphi) ) \right) , \rho)$,

for each $\rho \in \Xi(h)$ $P_\to \left( \#(h, \chi, \rho), \ \#(h, \exists ( x_{m+1} : \varphi_{m+1}, 
\neg(\varphi) ), \rho) \right)$,

for each $\rho \in \Xi(h)$ $P_\to \left( \#(h, \chi, \rho), \ \text{there exists } \sigma \in \Xi(k): \rho \sqsubseteq \sigma, \ \#(k, \neg(\varphi), \sigma) \right)$,

for each $\rho \in \Xi(h)$ if $\#(h, \chi, \rho)$ then (there exists $\sigma \in \Xi(k)$ such that 
$\rho \sqsubseteq \sigma$ and $\#(k, \varphi, \sigma)$ is false).

\medskip

The last condition is clearly ensured by our hypothesis.
\end{proof}

\medskip

We can create a set~$R_{\ref{L:notforall-to-existsnot-new}}$ 
which is the set of all pairs
\[ \left(   
\begin{array} {l} 
{ \gamma[ x_1: \varphi_1, \dots , x_m: \varphi_m, \to \left( \chi, \neg \left( \forall ( x_{m+1} : \varphi_{m+1}, \varphi ) \right) \right) ], } \\
{ \gamma[ x_1: \varphi_1, \dots , x_m: \varphi_m, \to \left( \chi, \exists ( x_{m+1} : 
\varphi_{m+1}, \neg(\varphi) ) \right) ] } 
\end{array}
\right) \]

such that 
\begin{itemize}
\item $m$ is a positive integer, $x_1, \dots , x_{m+1} \in \mathcal{V}$, with $x_i \ne x_j$ for $i 
\ne j$, $\varphi_1, \dots , \varphi_{m+1} \in E$, $H[x_1:\varphi_1, \dots , x_{m+1}:\varphi_{m+1} ]$;
\item if we define $k = k[x_1:\varphi_1, \dots , x_{m+1}:\varphi_{m+1} ]$ and $h = k[x_1:\varphi_1, 
\dots , x_m:\varphi_m ]$ then $\chi \in S(h)$, $\varphi \in S(k)$.
\end{itemize}

\medskip

Lemma~\ref{L:notforall-to-existsnot-new} shows us that this set (which is a potential 1-ary rule) is `sound'. In order to use $R_{\ref{L:notforall-to-existsnot-new}}$ as a rule in our system we also need to show that $R_{\ref{L:notforall-to-existsnot-new}}$ is r.e..\\

\begin{lemma}
$R_{\ref{L:notforall-to-existsnot-new}}$ is r.e..
\end{lemma}

\begin{proof}

Given a positive integer $m$ and $(x_1, \varphi_1, \dots , x_{m+1}, \varphi_{m+1}) \in R_{m+1}$ all of the following sets are r.e.:
\begin{itemize}
\item $S(k[x_1:\varphi_1, \dots , x_m:\varphi_m ])$,
\item $S(k[x_1:\varphi_1, \dots , x_{m+1}:\varphi_{m+1} ])$.
\end{itemize}

\medskip

Therefore the following set is also r.e.:
\[ \{ (x_1, \varphi_1, \dots , x_{m+1}, \varphi_{m+1}) \} \times S(k[x_1:\varphi_1, \dots , x_m:\varphi_m ]) \times S(k[x_1:\varphi_1, \dots , x_{m+1}:\varphi_{m+1} ]) \ . \]

\medskip

Let's use this temporary definition 

\begin{multline*} 
Q'_{m+1,2} = \bigcup_{(x_1, \varphi_1, \dots , x_{m+1}, \varphi_{m+1}) \in R_{m+1}} \{(x_1, \varphi_1, \dots , x_{m+1}, \varphi_{m+1})\} \times S(k[x_1:\varphi_1, \dots , x_m:\varphi_m ])\\ 
\times S(k[x_1:\varphi_1, \dots , x_{m+1}:\varphi_{m+1} ]) ).
\end{multline*}

With this $Q'_{m+1,2}$ is a r.e. subset of $(\Sigma^*)^{2(m + 1)} \times \Sigma^* \times \Sigma^*$.\\

We now define two functions $\delta_{1,m}, \ \delta_{2,m}$ over $(\Sigma^*)^{2(m + 1)} \times \Sigma^* \times \Sigma^*$ as follows. Given $((\psi_1, \varphi_1, \dots , \psi_{m+1}, \varphi_{m+1}), \chi, \varphi) \in (\Sigma^*)^{2(m + 1)} \times \Sigma^* \times \Sigma^*$
\begin{multline*}
\delta_{1,m}((\psi_1, \varphi_1, \dots , \psi_{m+1}, \varphi_{m+1}), \chi, \varphi) =\\
\gamma[ \psi_1: \varphi_1, \dots , \psi_m: \varphi_m, \to \left( \chi, \neg \left( \forall ( \psi_{m+1} : \varphi_{m+1}, \varphi ) \right) \right) ] \ .
\end{multline*}
\begin{multline*}
\delta_{2,m}((\psi_1, \varphi_1, \dots , \psi_{m+1}, \varphi_{m+1}), \chi, \varphi) =\\
\gamma[ \psi_1: \varphi_1, \dots , \psi_m: \varphi_m, \to \left( \chi, \exists ( \psi_{m+1} : 
\varphi_{m+1}, \neg(\varphi) ) \right) ] \ .
\end{multline*}

\medskip

All of the two functions we have defined are computable functions from $(\Sigma^*)^{2(m + 1)} \times \Sigma^* \times \Sigma^*$ to $\Sigma^*$. If we define a function $\delta_m$ over $(\Sigma^*)^{2(m + 1)} \times \Sigma^* \times \Sigma^*$ as follows: 
\[
\delta_m((\psi_1, \varphi_1, \dots , \psi_{m+1}, \varphi_{m+1}), \chi, \varphi) = 
\left(
\begin{array} {l}
{ \delta_{1,m}((\psi_1, \varphi_1, \dots , \psi_{m+1}, \varphi_{m+1}), \chi, \varphi), } \\
{ \delta_{2,m}((\psi_1, \varphi_1, \dots , \psi_{m+1}, \varphi_{m+1}), \chi, \varphi) } 
\end{array}
\right)
\]

then $\delta_m$ is a computable function from $(\Sigma^*)^{2(m + 1)} \times \Sigma^* \times \Sigma^*$ to $(\Sigma^*)^2$, therefore the set 
\[ D_m = \{ \delta_m((\psi_1, \varphi_1, \dots , \psi_{m+1}, \varphi_{m+1}), \chi, \varphi) | ((\psi_1, \varphi_1, \dots , \psi_{m+1}, \varphi_{m+1}), \chi, \varphi) \in Q'_{m+1,2} \} \]
is a r.e. subset of $(\Sigma^*)^2$.\\

If we now consider the set $\bigcup_{m \geqslant 1} D_m$ then this is a r.e. subset of $(\Sigma^*)^2$ and actually this set is equal to our rule $R_{\ref{L:notforall-to-existsnot-new}}$ which so is r.e. itself.\\

If $\xi \in R_{\ref{L:notforall-to-existsnot-new}}$ then there exist a positive integer $m$, $x_1, \dots , x_{m+1} \in \mathcal{V}$, with $x_i \ne x_j$ for $i \ne j$, $\varphi_1, \dots , \varphi_{m+1} \in E$ such that $H[x_1:\varphi_1, \dots , x_{m+1}:\varphi_{m+1} ]$;
if we define $k = k[x_1:\varphi_1, \dots , x_{m+1}:\varphi_{m+1} ]$ and $h = k[x_1:\varphi_1, \dots , x_m:\varphi_m ]$ 
there also exist $\chi \in S(h)$, $\varphi \in S(k)$, $\xi_1, \xi_2 \in \Sigma^*$ such that
\begin{itemize}
\item $\xi = (\xi_1, \xi_2)$,
\item $\xi_1 = \gamma[ x_1: \varphi_1, \dots , x_m: \varphi_m, \to \left( \chi, \neg \left( \forall ( x_{m+1} : \varphi_{m+1}, \varphi ) \right) \right) ]$,
\item $\xi_2 = \gamma[ x_1: \varphi_1, \dots , x_m: \varphi_m, \to \left( \chi, \exists ( x_{m+1} : 
\varphi_{m+1}, \neg(\varphi) ) \right) ]$.
\end{itemize}

\smallskip

This means that $(x_1, \varphi_1, \dots , x_{m+1}, \varphi_{m+1}) \in R_{m+1}$, $\chi \in S(k[x_1:\varphi_1, \dots , x_m:\varphi_m ])$, $\varphi \in S(k[x_1:\varphi_1, \dots , x_{m+1}:\varphi_{m+1} ])$, so 
$((x_1, \varphi_1, \dots , x_{m+1}, \varphi_{m+1}), \chi, \varphi) \in Q'_{m+1,2}$.\\

Moreover 
\begin{itemize}
\item $\xi_1 = \delta_{1,m}((x_1, \varphi_1, \dots , x_{m+1}, \varphi_{m+1}), \chi, \varphi)$,
\item $\xi_2 = \delta_{2,m}((x_1, \varphi_1, \dots , x_{m+1}, \varphi_{m+1}), \chi, \varphi)$.
\end{itemize}

i.e. $\xi = \delta_m((x_1, \varphi_1, \dots , x_{m+1}, \varphi_{m+1}), \chi, \varphi) \in D_m$.\\

Conversely if there exists $p \geqslant 1$ such that $\xi \in D_p$ then there exists $((\psi_1, \varphi_1, \dots , \psi_{p+1}, \varphi_{p+1}), \chi, \varphi) \in Q'_{p+1,2}$ such that\\
$\xi = \delta_p((\psi_1, \varphi_1, \dots , \psi_{p+1}, \varphi_{p+1}), \chi, \varphi)$.\\

Since $((\psi_1, \varphi_1, \dots , \psi_{p+1}, \varphi_{p+1}), \chi, \varphi) \in Q'_{p+1,2}$ we have $(\psi_1, \varphi_1, \dots , \psi_{p+1}, \varphi_{p+1}) \in R_{p+1}$, $\chi \in S(k[\psi_1:\varphi_1, \dots , \psi_p:\varphi_p ])$, $\varphi \in S(k[\psi_1:\varphi_1, \dots , \psi_{p+1}:\varphi_{p+1} ])$.\\

It follows that $\psi_1, \dots , \psi_{p+1} \in \mathcal{V}$, $\psi_i \ne \psi_j$ for $i \ne j$, $\varphi_1, \dots , \varphi_{p+1} \in E$, $H[\psi_1:\varphi_1, \dots , \psi_{p+1}:\varphi_{p+1} ]$.\\

Moreover 
\begin{multline*}
\xi = \delta_p((\psi_1, \varphi_1, \dots , \psi_{p+1}, \varphi_{p+1}), \chi, \varphi) = 
\left(
\begin{array} {l}
{ \delta_{1,p}((\psi_1, \varphi_1, \dots , \psi_{p+1}, \varphi_{p+1}), \chi, \varphi), } \\
{ \delta_{2,p}((\psi_1, \varphi_1, \dots , \psi_{p+1}, \varphi_{p+1}), \chi, \varphi) } 
\end{array}
\right) =\\ 
= \left(
\begin{array} {l}
{ \gamma[ \psi_1: \varphi_1, \dots , \psi_p: \varphi_p, \to \left( \chi, \neg \left( \forall ( \psi_{p+1} : \varphi_{p+1}, \varphi ) \right) \right) ], } \\
{ \gamma[ \psi_1: \varphi_1, \dots , \psi_p: \varphi_p, \to \left( \chi, \exists ( \psi_{p+1} : 
\varphi_{p+1}, \neg(\varphi) ) \right) ] } 
\end{array}
\right) 
\end{multline*}

\medskip

and so $\xi \in R_{\ref{L:notforall-to-existsnot-new}}$.\\
\end{proof}

\medskip

Then let $R_{\ref{L:notforall-to-existsnot-new}} \in \mathcal{R}$.

\bigskip

\begin{lemma}\label{L:rule-cl-to-il-new}
Let $m$ be a positive integer. Let $x_1, \dots , x_m \in \mathcal{V}$, with $x_i \ne x_j$ for $i \ne 
j$. Let $\varphi_1, \dots , \varphi_m \in E$ and assume $H[x_1:\varphi_1, \dots , x_m:\varphi_m ]$. 
Define $k = k[x_1:\varphi_1, \dots , x_m:\varphi_m ]$ and let $\varphi, \psi, \chi \in S(k)$.

\medskip

Under these assumptions we have 
\begin{itemize}
\item $\to( \wedge(\varphi, \psi) , \chi), \to(\varphi, \to(\psi, \chi) ) \in S(k)$,
\item $\gamma[ x_1: \varphi_1, \dots , x_m: \varphi_m, \to( \wedge(\varphi, \psi) , \chi) ] \in 
S(\epsilon)$,
\item $\gamma[ x_1: \varphi_1, \dots , x_m: \varphi_m, \to(\varphi, \to(\psi, \chi) ) ] \in S
(\epsilon)$.
\end{itemize}

\smallskip

Moreover if $\#( \gamma[ x_1: \varphi_1, \dots , x_m: \varphi_m, \to( \wedge(\varphi, \psi) , 
\chi) ] )$ then\\
$\#( \gamma[ x_1: \varphi_1, \dots , x_m: \varphi_m, \to(\varphi, \to(\psi, \chi) ) ] )$.

\end{lemma}

\begin{proof}

\medskip

We assume $\#( \gamma[ x_1: \varphi_1, \dots , x_m: \varphi_m, \to( \wedge(\varphi, \psi) , \chi) 
] )$ which can be rewritten
\begin{align} 
&\text{for each } \sigma \in \Xi(k) \ \#( k , \to( \wedge(\varphi, \psi) , \chi), \sigma) \ , \notag  \\
&\text{for each } \sigma \in \Xi(k) \ P_{\to}( \#(k, \wedge(\varphi, \psi), \sigma) , \#(k, \chi , \sigma) ) \ , \notag  \\
&\text{for each } \sigma \in \Xi(k) \ P_{\to}( P_{\wedge}( \#(k, \varphi , \sigma), \#(k, \psi , \sigma) ) , \#(k, \chi 
, \sigma) ) \  \notag .
\end{align}

\medskip

Of course we now try to show $\#( \gamma[ x_1: \varphi_1, \dots , x_m: \varphi_m, \to(\varphi, 
\to(\psi, \chi) ) ] )$ which in turn can be rewritten
\begin{align} 
&\text{for each } \sigma \in \Xi(k) \ \#( k , \to(\varphi, \to(\psi, \chi) ), \sigma) \ , \notag  \\
&\text{for each } \sigma \in \Xi(k) \ P_{\to}( \#(k, \varphi , \sigma), \#(k, \to(\psi, \chi) , \sigma)  ) \ , \notag  \\
&\text{for each } \sigma \in \Xi(k) \ P_{\to}( \#(k, \varphi , \sigma), P_{\to}( \#(k, \psi, \sigma), \#(k, \chi, 
\sigma) )  ) \  \notag .
\end{align}

\medskip

Let $\sigma \in \Xi(k)$, suppose $\#(k, \varphi, \sigma)$ and $\#(k, \psi, \sigma)$, then we have $
\#(k, \chi, \sigma)$ and this completes the proof.
\end{proof}

\medskip

We can create a set $R_{\ref{L:rule-cl-to-il-new}}$ which is the set of all pairs
\[ \left( \gamma[ x_1: \varphi_1, \dots , x_m: \varphi_m, \to( \wedge(\varphi, \psi) , \chi) ], 
\gamma[ x_1: \varphi_1, \dots , x_m: \varphi_m, \to(\varphi, \to(\psi, \chi) ) ] \right)  \]

such that
\begin{itemize}
\item $m$ is a positive integer, $x_1, \dots , x_m \in \mathcal{V}$, $x_i \ne x_j$ for $i \ne j$, $
\varphi_1, \dots , \varphi_m \in E$, $H[x_1:\varphi_1, \dots , x_m:\varphi_m ]$,
\item $\varphi, \psi, \chi \in S(k[x_1:\varphi_1, \dots , x_m:\varphi_m ])$.
\end{itemize}

\medskip

Lemma~\ref{L:rule-cl-to-il-new} shows us that this set (which is a potential 1-ary rule) is `sound'. In order to use $R_{\ref{L:rule-cl-to-il-new}}$ as a rule in our system we also need to show that $R_{\ref{L:rule-cl-to-il-new}}$ is r.e..\\

\begin{lemma}
$R_{\ref{L:rule-cl-to-il-new}}$ is r.e..
\end{lemma}

\begin{proof}
Given a positive integer $m$ and $(x_1, \varphi_1, \dots , x_m, \varphi_m) \in R_m$ we can notice the following:

\begin{itemize}
\item $k[x_1:\varphi_1, \dots , x_m:\varphi_m ] \in K$;
\item $S(k[x_1:\varphi_1, \dots , x_m:\varphi_m ])$ is r.e.;
\item $\{(x_1, \varphi_1, \dots , x_m, \varphi_m)\} \times S(k[x_1:\varphi_1, \dots , x_m:\varphi_m ])^3$ is r.e..
\end{itemize}

Let's define
\[ Q_{m,3} = \bigcup_{(x_1, \varphi_1, \dots , x_m, \varphi_m) \in R_m} \{(x_1, \varphi_1, \dots , x_m, \varphi_m)\} \times S(k[x_1:\varphi_1, \dots , x_m:\varphi_m ])^3 \ . \]

\smallskip

Clearly $Q_{m,3} \subseteq (\Sigma^*)^{2m} \times (\Sigma^*)^3$ is also r.e..\\

We now define two functions $\delta_{1,m}, \ \delta_{2,m}$ over $(\Sigma^*)^{2m} \times (\Sigma^*)^3$ as follows. Given $((\psi_1, \varphi_1, \dots , \psi_m, \varphi_m), (\varphi, \psi, \chi)) \in (\Sigma^*)^{2m} \times (\Sigma^*)^3$
\begin{align} 
& \delta_{1,m}((\psi_1, \varphi_1, \dots , \psi_{m}, \varphi_{m}), (\varphi, \psi, \chi)) = \gamma[ \psi_1: \varphi_1, \dots , \psi_m: \varphi_m, \to( \wedge(\varphi, \psi) , \chi) ] \ , \notag  \\
& \delta_{2,m}((\psi_1, \varphi_1, \dots , \psi_{m}, \varphi_{m}), (\varphi, \psi, \chi)) = \gamma[ \psi_1: \varphi_1, \dots , \psi_m: \varphi_m, \to(\varphi, \to(\psi, \chi) ) ] \  \notag .
\end{align}

\smallskip

All of the two functions we have defined are computable functions from $(\Sigma^*)^{2m} \times (\Sigma^*)^3$ to $\Sigma^*$. If we define a function $\delta_m$ over $(\Sigma^*)^{2m} \times (\Sigma^*)^3$ as follows: 
\[
\delta_m((\psi_1, \varphi_1, \dots , \psi_{m}, \varphi_{m}), (\varphi, \psi, \chi)) = 
\left(
\begin{array} {l}
{ \delta_{1,m}((\psi_1, \varphi_1, \dots , \psi_{m}, \varphi_{m}), (\varphi, \psi, \chi)) , } \\
{ \delta_{2,m}((\psi_1, \varphi_1, \dots , \psi_{m}, \varphi_{m}), (\varphi, \psi, \chi)) , }
\end{array}
\right)
\]

then $\delta_m$ is a computable function from $(\Sigma^*)^{2m} \times (\Sigma^*)^3$ to $(\Sigma^*)^2$, therefore the set 
\[ D_m = \{ \delta_m((\psi_1, \varphi_1, \dots , \psi_{m}, \varphi_{m}), (\varphi, \psi, \chi)) | \ ((\psi_1, \varphi_1, \dots , \psi_{m}, \varphi_{m}), (\varphi, \psi, \chi)) \in Q_{m,3} \} \]
is a r.e. subset of $(\Sigma^*)^2$.\\

If we now consider the set $\bigcup_{m \geqslant 1} D_m$ then this is a r.e. subset of $(\Sigma^*)^2$ and actually this set is equal to our set $R_{\ref{L:rule-cl-to-il-new}}$ which so is r.e. itself.
\end{proof}

\medskip

Then let $R_{\ref{L:rule-cl-to-il-new}} \in \mathcal{R}$.

\bigskip

\begin{lemma}\label{L:rule-il-to-cl-simple-new}
Let $\varphi, \psi, \chi \in S(\epsilon)$. We have 
\begin{itemize}
\item $\to (\varphi, \to (\psi, \chi ) ) \in S(\epsilon)$,
\item $\to ( \wedge(\varphi, \psi) , \chi  ) \in S(\epsilon)$.
\end{itemize}

\smallskip

Moreover if $\#( \to (\varphi, \to (\psi, \chi ) ) )$ then $\#( \to ( \wedge(\varphi, \psi) , \chi  ) )$.

\end{lemma}

\begin{proof}

\medskip

Suppose $\#( \to (\varphi, \to (\psi, \chi ) ) )$ holds. It can be rewritten
\begin{align} 
& P_{\to} ( \#(\varphi), \#( \to (\psi, \chi ) )  ) \ , \notag  \\
& P_{\to} ( \#(\varphi), P_{\to} ( \#(\psi),  \#(\chi) )  ) \  \notag .
\end{align}

\medskip

In turn, $\#( \to ( \wedge(\varphi, \psi) , \chi  ) )$ can be rewritten
\begin{align} 
& P_{\to} ( \#( \wedge(\varphi, \psi) ), \#( \chi )  ) \ , \notag  \\
& P_{\to} ( P_{\wedge} ( \#(\varphi), \#(\psi) ) , \#( \chi )  ) \  \notag .
\end{align}

\medskip

Suppose $\#(\varphi)$ and $\#(\psi)$ both hold, we need to show that $\#(\chi)$ holds. This is granted by 
\[ P_{\to} ( \#(\varphi), P_{\to} ( \#(\psi),  \#(\chi) )  ) \ . \]
\end{proof}

\medskip

We can create a set~$R_{\ref{L:rule-il-to-cl-simple-new}}$ which is the set of all pairs
\[ \left(   
\begin{array} {l} 
{ \to (\varphi, \to (\psi, \chi ) ), } \\
{ \to ( \wedge(\varphi, \psi) , \chi  ) } 
\end{array}
\right) \]

such that $\varphi, \psi, \chi \in S(\epsilon)$. 

\medskip

Lemma~\ref{L:rule-il-to-cl-simple-new} shows us that this set (which is a potential 1-ary rule) is `sound'. In order to use $R_{\ref{L:rule-il-to-cl-simple-new}}$ as a rule in our system we also need to show that $R_{\ref{L:rule-il-to-cl-simple-new}}$ is r.e..\\

\begin{lemma}
$R_{\ref{L:rule-il-to-cl-simple-new}}$ is r.e.
\end{lemma}

\begin{proof}
Clearly $S(\epsilon)$ is .r.e. and so is $S(\epsilon)^3$.\\

Let's define two functions $\delta_{1,1}$, $\delta_{2,1}$ over $(\Sigma^*)^3$ as follows: 
\begin{align} 
& \delta_{1,1}(\varphi, \psi, \chi) = \ \to (\varphi, \to (\psi, \chi ) ) \ , \notag  \\
& \delta_{2,1}(\varphi, \psi, \chi) = \ \to ( \wedge(\varphi, \psi) , \chi  ) \  \notag .
\end{align}

\smallskip

The two functions we have defined are both computable functions from $(\Sigma^*)^3$ to $\Sigma^*$. If we define a function $\delta_1$ over $(\Sigma^*)^3$ as follows
\[ \delta_1(\varphi, \psi, \chi) = 
\left(
\begin{array} {l}
{ \delta_{1,1}( \varphi, \psi, \chi ), } \\
{ \delta_{2,1}( \varphi, \psi, \chi ) } \\
\end{array}
\right) \ , \]

then $\delta_1$ is a computable function from $(\Sigma^*)^3$ to $(\Sigma^*)^2$, therefore the set 
\[ D_1 = \{ \delta_1(\varphi, \psi, \chi) | (\varphi, \psi, \chi) \in S(\epsilon)^3 \} \]

is a r.e. subset of $(\Sigma^*)^2$, and $D_1$ is equal to our set $R_{\ref{L:rule-il-to-cl-simple-new}}$ which so is r.e. itself.
\end{proof}

\medskip

Then let $R_{\ref{L:rule-il-to-cl-simple-new}} \in \mathcal{R}$.

\section{Example of a proof}\label{Ch:proofExample}

As an example of proof, we want to prove a form of the Bocardo syllogism. In Ferreir\'{o}s' referenced paper (\cite{Ferreiros}), on paragraph 3.1, the syllogism is expressed as follows:\\

Some $A$ are not $B$. All $C$ are $B$. Therefore, some $A$ are not $C$.\\

Suppose $A$, $B$ and $C$ represent sets, the statement we actually want to prove is the following:\\

If ( (there exists $x \in A$ such that $x \notin B$) and (for each $y \in C$ $y \in B$) ) then\\
(there exists $z \in A$ such that $z \notin C$).\\

In order to formalize this, we will use a language $(\mathcal{V}, \mathcal{F}, \mathcal{C}, \#, \{D_1, \dots, D_p \})$ which must be as follows
\begin{align} 
&\mathcal{V} = \{ x, y, z \} \ , \notag  \\
&\mathcal{F} = \{ \neg, \wedge, \vee, \to, \leftrightarrow, \in, =  \} \ , \notag  \\
&\mathcal{C} = \{ A, B, C \} \  \notag .
\end{align}

\smallskip

where $A, B, C$ are constants each representing a set.\\

Moreover, we do not need the additional sets $\{D_1, \dots, D_p \}$ so we can set $p = 0$.\\

Here we notice that we set as a constraint that for each $c \in \mathcal{C}$ and for any positive integer $q$ we must be able to decide all of the following conditions

\begin{itemize}
\item $Set_q( \#(c))$;
\item $Event_q( \#(c))$;
\item $\#(c) \in D_i$;
\item $\#(c) \in \mathcal{P}^q(D_i)$;
\item if ($Set_q( \#(c))$) then ($NotEmpty_q( \#(c))$).
\end{itemize}

And moreover, the last condition must be decided as true.\\

With respect to our specific set of constants $\mathcal{C}$, each of its members represents a set and has nothing else as a specific constraint. For instance with respect to $A$ we can take the following decisions. 

\begin{itemize}
\item $Set_1( \#(A))$: true;
\item for $q > 1$: $Set_q( \#(A))$: false;
\item $Event_q( \#(A))$: false;
\item if ($Set_q( \#(A))$) then ($NotEmpty_q( \#(A))$): true.
\end{itemize}

We can omit the decisions related to sets $D_i$ because we don't have such sets in our language. The exact same decisions are taken for $B$ and $C$.\\

At this point we suppose we can formalize the statement as
\begin{equation}\label{E:proof-ex}
\to \left( \wedge \left( 
\begin{array} {l}
{  \exists \left( x:A, \neg \left( \in(x, B) \right) \right) , } \\
{ \forall \left( y:C, \in(y,B) \right) }
\end{array}
\right) , 
\exists \left( z:A, \neg \left(  \in (z, C)  \right) \right) \right)     \tag{$Th_1$}  \ .
\end{equation}

\medskip

We'll soon see a proof of this statement and of course if we can show a proof of a statement then we have also proved the statement is a sentence in our language.\\

First of all we need the following lemma, that can be applied to any language which includes all the symbols $\neg, \wedge, \vee, \to, \leftrightarrow, \in, =$ in the set $\mathcal{F}$, and therefore it can also be applied to our current language.\\

\begin{lemma}\label{L:xi-in-D-new}
Let $m$ be a positive integer, $x_1, \dots , x_m \in \mathcal{V}$, with $x_i \ne x_j$ for $i \ne j$. Let $A_1, \dots , A_m \in \mathcal{C}$ such that for each $i = 1 \dots m$ $\#(A_i)$ is a set. Let $D \in \mathcal{C}$ such that $\#(D)$ is a set. We have $H[x_1:A_1, \dots , x_m:A_m ]$. If we define $k = k[x_1:A_1, \dots , x_m:A_m ]$ then for each $i = 1 \dots m$ 
\begin{itemize}
\item $\ \in(x_i, D) \ \in S(k)$,
\item for each $\sigma \in \Xi(k)$ $\#(k, \, \in(x_i, D), \sigma) =  P_{\in}(\#(k, x_i, \sigma), \#(D) )$.
\end{itemize}
\end{lemma}

\begin{proof}

\medskip

We first consider that $A_1 \in E(\epsilon)$ and $\#(A_1)$ is a set, so $A_1 \in E_s(\epsilon)$ and $H[x_1:A_1]$. Let $k_1 = k[x_1:A_1]$.

\medskip

If $m>1$ then for each $i = 1 \dots m-1$ we suppose $H[x_1:A_1, \dots , x_i:A_i ]$ holds and we define $k_i = k[x_1:A_1, \dots , x_i:A_i ]$.\\
Clearly by lemma~\ref{L:const} $A_{i+1} \in E(k_i)$ and for each $\rho \in \Xi(k_i)$ $\#(k_i, A_{i+1}, \rho) = \#(A_{i+1})$ is a set.\\
So $A_{i+1} \in E_s(k_i)$, which implies $H[x_1:A_1, \dots , x_{i+1}:A_{i+1} ]$ (and we can define \linebreak
$k_{i+1} = k[x_1:A_1, \dots , x_{i+1}:A_{i+1} ]$).

\medskip

This proves that $H[x_1:A_1, \dots , x_m:A_m ]$ holds.

\medskip

Let $i = 1 \dots m$. Using lemma~\ref{L:xi_in_Ekj} we obtain that $x_i \in E(k)$.

\medskip

Moreover $D \in E(k)$ and for each $\sigma \in \Xi(k)$ $\#(k, D, \sigma) = \#(D)$ is a set. By lemma~\ref{L:in-t-varphi-in-Sk} we have
\begin{itemize}
\item $\ \in(x_i, D) \ \in S(k)$,
\item for each $\sigma \in \Xi(k)$ $\#(k, \, \in(x_i, D), \sigma) =  P_{\in}(\#(k, x_i, \sigma), \#(D) )$.
\end{itemize}
\end{proof}

\bigskip

In order to provide a proof of statement~\ref{E:proof-ex} we'll make use of a deductive system 
which includes all the axioms and rules listed in section~\ref{Ch:buildDedSystem}.

\medskip

Using the former lemma we can derive $H[x:A]$ and we can define $h = k[x:A]$. Moreover $\in(x,B) \in S(h)$, so $\neg(\in(x,B)) \in S(h)$.

\medskip

We also have $H[x:A, y:C]$ and we define $k_y = k[x:A, y:C]$.\\
We have $\in(y,B) \in S(k_y)$ and by lemma~\ref{L:quantifiers_S_h} $\forall( y:C, \in(y,B) ) \in S(h)$.

\medskip

Thus $\wedge \left( \neg( \in(x,B) ),  \forall ( y:C, \in(y,B) ) \right)$ also belongs to $S(h)$.

\medskip

Moreover $H[x:A, z:A]$ and we define $k_z = k[x:A, z:A]$.\\
We have $\in(z,C) \in S(k_z)$ and by lemma~\ref{L:quantifiers_S_h} $\forall( z:A, \in(z,C) ) \in S(h)$.

\medskip

The first sentence in our proof is an instance of axiom~$A_{\ref{L:axiom-cl-elim-new}}$.
\scriptsize
\begin{equation}\label{E:ex2-a-new}
 \gamma \left[x:A, \to \left( 
\wedge \left(  
\begin{array} {l}
{
\wedge \left(
\begin{array} {l}
{  \neg(\in(x,B) ) , } \\
{ \forall ( y:C, \in(y,B) ) }
\end{array}
\right), }\\
{ \forall ( z:A, \in(z,C) ) }
\end{array}
\right) , 
\wedge \left(
\begin{array} {l}
{ \neg( \in(x,B) ) , } \\
{ \forall ( y:C, \in(y,B) ) }
\end{array}
\right)
\right) \right] .
\end{equation}
\normalsize

\medskip

By~$A_{\ref{L:axiom-cl-elim-new}}$ we also obtain
\begin{equation}\label{E:ex2-b-new}
 \gamma \left[x:A, \to \left( 
\wedge \left(
\begin{array} {l}
{  \neg( \in(x,B) ) , } \\
{ \forall ( y:C, \in(y,B) ) }
\end{array}
\right) , 
\neg( \in(x,B) )
\right) \right] .
\end{equation}

\medskip

By~\ref{E:ex2-a-new}, \ref{E:ex2-b-new} and rule~$R_{\ref{L:rule-il-trans-new}}$
\begin{equation}\label{E:ex2-c-new}
 \gamma \left[x:A, \to \left( 
\wedge \left(  
\begin{array} {l}
{
\wedge \left(
\begin{array} {l}
{  \neg( \in(x,B) ) , } \\
{ \forall ( y:C, \in(y,B) ) }
\end{array}
\right), }\\
{ \forall ( z:A, \in(z,C) ) }
\end{array}
\right) , 
\neg( \in(x,B) )
\right) \right] .
\end{equation}

\medskip

Another instance of $A_{\ref{L:axiom-cl-elim-new}}$ is the following
\small
\begin{equation}\label{E:ex2-d-new}
 \gamma \left[x:A, \to \left( 
\wedge \left(  
\begin{array} {l}
{
\wedge \left(
\begin{array} {l}
{  \neg( \in(x,B) ) , } \\
{ \forall ( y:C, \in(y,B) ) }
\end{array}
\right), }\\
{ \forall ( z:A, \in(z,C) ) }
\end{array}
\right) , 
\forall ( z:A, \in(z,C) )
\right) \right] .
\end{equation}
\normalsize

\medskip

By axiom~$A_{\ref{L:axiom-xi-in-phii-new}}$ we obtain 
\begin{equation}\label{E:ex2-e-new}
\gamma[ x:A, \in(x, A ) ] .
\end{equation}

\medskip

By \ref{E:ex2-e-new} and rule~$R_{\ref{L:rule-il-introd-new}}$ we also get
\begin{equation}\label{E:ex2-f-new}
 \gamma \left[x:A, \to \left( 
\wedge \left(  
\begin{array} {l}
{
\wedge \left(
\begin{array} {l}
{  \neg( \in(x,B) ) , } \\
{ \forall ( y:C, \in(y,B) ) }
\end{array}
\right), }\\
{ \forall ( z:A, \in(z,C) ) }
\end{array}
\right) , 
\in(x, A )
\right) \right] .
\end{equation}

\medskip

Since $x \in E(h)$, $C \in E_s(h)$ etc. we can apply rule~$R_{\ref{L:rule-forall-elim-rest-new}}$ to~\ref{E:ex2-d-new} and~\ref{E:ex2-f-new} and obtain
\begin{equation}\label{E:ex2-g-new}
 \gamma \left[x:A, \to \left( 
\wedge \left(  
\begin{array} {l}
{
\wedge \left(
\begin{array} {l}
{  \neg( \in(x,B) ) , } \\
{ \forall ( y:C, \in(y,B) ) }
\end{array}
\right), }\\
{ \forall ( z:A, \in(z,C) ) }
\end{array}
\right) , 
\in(x, C )
\right) \right] .
\end{equation}

\medskip

By axiom~$A_{\ref{L:axiom-cl-elim-new}}$ 
\begin{equation}\label{E:ex2-h-new}
 \gamma \left[x:A, \to \left( 
\wedge \left(
\begin{array} {l}
{  \neg( \in(x,B) ) , } \\
{ \forall ( y:C, \in(y,B) ) }
\end{array}
\right) , 
\forall( y:C, \in(y,B) )
\right) \right] .
\end{equation}

\medskip

By~\ref{E:ex2-a-new}, \ref{E:ex2-h-new} and rule~$R_{\ref{L:rule-il-trans-new}}$
\small
\begin{equation}\label{E:ex2-i-new}
 \gamma \left[x:A, \to \left( 
\wedge \left(  
\begin{array} {l}
{
\wedge \left(
\begin{array} {l}
{  \neg( \in(x,B) ) , } \\
{ \forall( y:C, \in(y,B) ) }
\end{array}
\right), }\\
{ \forall ( z:A, \in(z,C) ) }
\end{array}
\right) , 
\forall ( y:C, \in(y,B) )
\right) \right] .
\end{equation}
\normalsize

\medskip

Since $x \in E(h)$, $B \in E_s(h)$ etc. we can apply rule~$R_{\ref{L:rule-forall-elim-rest-new}}$ to~\ref{E:ex2-g-new} and~\ref{E:ex2-i-new} and obtain
\begin{equation}\label{E:ex2-l-new}
 \gamma \left[x:A, \to \left( 
\wedge \left(  
\begin{array} {l}
{
\wedge \left(
\begin{array} {l}
{  \neg( \in(x,B) ) , } \\
{ \forall ( y:C, \in(y,B) ) }
\end{array}
\right), }\\
{ \forall ( z:A, \in(z,C) ) }
\end{array}
\right) , 
\in(x, B )
\right) \right] .
\end{equation}

\medskip

By \ref{E:ex2-l-new}, \ref{E:ex2-c-new} and $R_{\ref{L:AND-I-new}}$
\small
\begin{equation}\label{E:ex2-m-new}
 \gamma \left[x:A, \to \left( 
\wedge \left(  
\begin{array} {l}
{
\wedge \left(
\begin{array} {l}
{  \neg( \in(x,B) ) , } \\
{ \forall ( y:C, \in(y,B) ) }
\end{array}
\right), }\\
{ \forall ( z:A, \in(z,C) ) }
\end{array}
\right) , 
\wedge \left(
\begin{array} {l}
{  \in(x,B), } \\
{  \neg( \in(x,B) ) }
\end{array}
\right)
\right) \right] .
\end{equation}
\normalsize

\medskip

By $R_{\ref{L:contradict-elim-new}}$ 
\begin{equation}\label{E:ex2-n-new}
 \gamma \left[x:A, \neg \left( 
\wedge \left(  
\begin{array} {l}
{
\wedge \left(
\begin{array} {l}
{  \neg( \in(x,B) ) , } \\
{ \forall ( y:C, \in(y,B) ) }
\end{array}
\right), }\\
{ \forall ( z:A, \in(z,C) ) }
\end{array}
\right) \right) \right] .
\end{equation}

\medskip

By $R_{\ref{L:notand-to-impl-new}}$
\begin{equation}\label{E:ex2-o-new}
 \gamma \left[x:A, \to \left( 
\wedge 
\left(
\begin{array} {l}
{  \neg( \in(x,B) ) , } \\
{  \forall ( y:C, \in(y,B) ) }
\end{array}
\right),
\neg \left( \forall ( z:A, \in(z,C) ) \right)
\right) \right] .
\end{equation}

\medskip

By $R_{\ref{L:notforall-to-existsnot-new}}$
\begin{equation}\label{E:ex2-p-new}
 \gamma \left[x:A, \to \left( 
\wedge 
\left(
\begin{array} {l}
{  \neg( \in(x,B) ) , } \\
{  \forall ( y:C, \in(y,B) ) }
\end{array}
\right),
\exists ( z:A, \neg (\in(z,C)) ) 
\right) \right] .
\end{equation}

\medskip

Since $\exists ( z:A, \neg (\in(z,C)) ) \in S(h)$ we can apply $R_{\ref{L:rule-cl-to-il-new}}$ and obtain
\begin{equation}\label{E:ex2-q-new}
 \gamma \left[x:A, \to \left( 
\neg( \in(x,B) ) ,
\to
\left(
\begin{array} {l}
{ \forall ( y:C, \in(y,B) ) , } \\
{ \exists ( z:A, \neg (\in(z,C)) ) }
\end{array}
\right)
\right) \right] .
\end{equation}

\medskip

Using lemma~\ref{L:xi-in-D-new} we obtain that $\in(y,B) \in S(k[y:C])$ and $\in(z,C) \in S(k[z:A])$.

\medskip

By lemma~\ref{L:quantifiers_S_h} we obtain that $\forall ( y:C, \in(y,B) ) \in S(\epsilon)$ and similarly \newline $\exists ( z:A, \neg (\in(z,C)) ) \in S(\epsilon)$.

\medskip

We can apply rule~$R_{\ref{L:leverage-exists-simple-new}}$ to~\ref{E:ex2-q-new} and obtain
\begin{equation}\label{E:ex2-r-new}
\to \left( 
\exists \left( x: A, \neg( \in(x,B) )  \right),
\to
\left(
\begin{array} {l}
{ \forall ( y:C, \in(y,B) ) , } \\
{ \exists ( z:A, \neg (\in(z,C)) ) }
\end{array}
\right)
\right)
\end{equation}

\medskip

Finally, by~$R_{\ref{L:rule-il-to-cl-simple-new}}$, we obtain
\begin{equation}\label{E:ex2-s-new}
\to \left( 
\wedge
\left(
\begin{array} {l}
{ \exists \left( x: A, \neg( \in(x,B) )  \right), } \\
{ \forall ( y:C, \in(y,B) ) }
\end{array}
\right)
,
\exists ( z:A, \neg (\in(z,C)) ) 
\right)
\end{equation}

\begin{flushright}
$\qed$\\
\end{flushright}

\bigskip

We have proved statement~\ref{E:proof-ex}, this also means that \ref{E:proof-ex} is a sentence in our language.\\

\section{Extending our deductive system}\label{Ch:extendDedSystem}

In this section we are going to extend our deductive systems, in other words we are going to add axioms and rules to a subset of the deductive sytems $\mathcal{D} = (\mathcal{A}, \mathcal{R})$ which we have built in section~\ref{Ch:buildDedSystem}. We are going to do this in order to be able to show another example of proof in the next section. Our new deductive systems can refer to any language $\mathcal{L}=(\mathcal{V}, \mathcal{F}, \mathcal{C}, \#, \{D_1, \dots, D_p \})$ such that all of these symbols: $N, *$ are in our set $\mathcal{C}$, all of these symbols $x, y, z, u, v, w$ are in our set $\mathcal{V}$, all of these symbols: $\neg, \wedge, \vee, \to, \leftrightarrow, \in, =$ are in our set $\mathcal{F}$. For each of these operators $f$ $A_f(x_1, \dots ,x_n)$ and $P_f(x_1, \dots ,x_n)$ are defined as specified at the beginning of section~\ref{Ch:lang}. Moreover we require $p \geqslant 1$ and $D_1 = \mathbb{N}$.\\

The constant symbol $N$ represents the set of natural numbers $\mathbb{N}$, so that we have $\#(N) = \mathbb{N}$.\\

The symbol~$*$ stands for the product (or multiplication) operation in the domain $\mathbb{N}$ of natural numbers. Therefore $\#(*)$ is a function defined on $\mathbb{N} \times \mathbb{N}$ and for each $\alpha, \beta \in \mathbb{N}$ $\#(*)(\alpha, \beta)$ is the product of $\alpha$ and $\beta$, in other words $\#(*)(\alpha, \beta) = \alpha \cdot \beta$.\\ 

Given a language $\mathcal{L}=(\mathcal{V}, \mathcal{F}, \mathcal{C}, \#, \{D_1, \dots, D_p \})$ as above, in section~\ref{Ch:buildDedSystem} we have defined a deductive system for this language, and we assume that all the axioms and rules we have defined for that deductive system apply to our new deductive system. We are now going to add new axioms and rules to our new deductive system.\\

\begin{lemma}\label{L:preliminary-results-ch6-1}
$H[x:N, y:N, z:N, u:N, v:N]$ holds.

\end{lemma}

\begin{proof}

Follows from lemma~\ref{L:xi-in-D-new}.

\end{proof}

\medskip

\begin{lemma}\label{L:preliminary-results-ch6-2}
Let $k \in K$ and $\varphi, \psi \in E(k)$. Then

\begin{itemize}
\item $=(\varphi, \psi) \in S(k)$
\item for each $\sigma \in \Xi(k)$ $\#(k, =(\varphi, \psi), \sigma) = (\#(k,\varphi, \sigma) = \#(k,\psi, \sigma)) \ .$ 
\end{itemize}
\end{lemma}

\begin{proof}
It's a simply a case of lemma~\ref{L:meaning-kept-f2}. 
\end{proof}

\medskip

\begin{lemma}\label{L:preliminary-results-ch6-3}
Let $k \in K$ and let $\varphi, \psi \in E(k)$. Assume for each $\sigma \in \Xi(k)$ $\#(k, \varphi, \sigma) \in \mathbb{N}$ and $\#(k, \psi, \sigma) \in \mathbb{N}$, then

\begin{itemize}
\item $(*)(\varphi, \psi) \in E(k)$
\item for each $\sigma \in \Xi(k)$ $\#(k, (*)(\varphi, \psi), \sigma) = (\#(k,\varphi, \sigma) \cdot \#(k,\psi, \sigma)) \ .$ 
\end{itemize}
\end{lemma}

\begin{proof}
It's simply a case of lemma~\ref{L:meaning-kept-c1}.
\end{proof}

\medskip

\begin{lemma}\label{L:preliminary-results-ch6-4}
Let $k \in K$ and let $\varphi \in E(k)$, then

\begin{itemize}
\item $\in(\varphi, N) \in S(k)$,
\item for each $\sigma \in \Xi(k)$ $\#(k, \in(\varphi, N), \sigma) = P_{\in}(\#(k,\varphi, \sigma), \mathbb{N})$. 
\end{itemize}
\end{lemma}

\begin{proof}
It's a simply a case of lemma~\ref{L:meaning-kept-f2}. 
\end{proof}

\medskip

\begin{lemma}\label{L:axiom-substitute-exp-in-mult}
Let $k = k[x:N, y:N, z:N, u:N, v:N]$, $\varphi \in E(k)$ such that for each $\sigma \in \Xi(k)$ $\#(k, \varphi, \sigma) \in \mathbb{N}$ then 

\begin{itemize}
\item $=(y,\varphi) \in S(k)$
\item $=(z,*(y,v)) \in S(k)$
\item $=(z,*(\varphi,v)) \in S(k)$
\end{itemize}

\smallskip

Moreover
\[ \# \left( \gamma[ x:N, y:N, z:N, u:N, v:N, \to \left( \wedge(=(y,\varphi), =(z,yv)), =(z,\varphi v)  \right) ] \right) \]
is true.
\end{lemma}

\begin{proof}

\smallskip

By lemma~\ref{L:xi_in_Ekj} $y \in E(k)$ and by~\ref{L:preliminary-results-ch6-2} $=(y,\varphi) \in S(k)$.

\medskip

Moreover by~\ref{L:xi_in_Ekj} $v \in E(k)$. If we define $k_v = k[x:N, y:N, z:N, u:N]$ then for each $\sigma \in \Xi(k)$ $\sigma_{/dom(k_v)} \in \Xi(k_v)$ and $\#(k,v,\sigma) \in \#(k_v, N, \sigma_{/dom(k_v)}) = \#(N) = \mathbb{N}$.

\medskip

Similarly by lemma~\ref{L:xi_in_Ekj} if we define $k_y = k[x:N]$ then for each $\sigma \in \Xi(k)$ $\sigma_{/dom(k_y)} \in \Xi(k_y)$ and $\#(k,y,\sigma) \in \#(k_y, N, \sigma_{/dom(k_y)}) = \#(N) = \mathbb{N}$.

\medskip 

Similarly by the same lemma $z \in E(k)$ and if we define $k_z = k[x:N, y:N]$ then for each $\sigma \in \Xi(k)$ $\sigma_{/dom(k_z)} \in \Xi(k_z)$ and $\#(k,z,\sigma) \in \#(k_z, N, \sigma_{/dom(k_z)}) = \#(N) = \mathbb{N}$.

\medskip

By lemma~\ref{L:preliminary-results-ch6-3} it follows that $(*)(y, v) \in E(k)$ and for each $\sigma \in \Xi(k)$ $\#(k, (*)(y, v), \sigma) = (\#(k,y, \sigma) \cdot \#(k,v, \sigma)) \in \mathbb{N} \ .$

\medskip

Similarly $(*)(\varphi, v) \in E(k)$ and for each $\sigma \in \Xi(k)$ $\#(k, (*)(\varphi, v), \sigma) = (\#(k,\varphi, \sigma) \cdot \#(k,v, \sigma)) \in \mathbb{N} \ .$

\medskip

By lemma~\ref{L:preliminary-results-ch6-2} $=(z,*(y,v)) \in S(k)$ and $=(z,*(\varphi,v)) \in S(k)$. 

\medskip 

Clearly it follows that $\to \left( \wedge(=(y,\varphi), =(z,yv)), =(z,\varphi v)  \right) \in S(k)$ and we can rewrite 
\[ \# \left( \gamma[ x:N, y:N, z:N, u:N, v:N, \to \left( \wedge(=(y,\varphi), =(z,yv)), =(z,\varphi v)  \right) ] \right) \]

as follows
\begin{align} 
&\text{for each } \sigma \in \Xi(k) \ \#( k , \to \left( \wedge(=(y,\varphi), =(z,yv)), =(z,\varphi v) \right), \sigma) \ , \notag  \\
&\text{for each } \sigma \in \Xi(k) \ P_{\to} (\#(k, \wedge(=(y,\varphi), =(z,yv)), \sigma) , \#(k, =(z,\varphi v) ,\sigma) ) \ , \notag \\
&\text{for each } \sigma \in \Xi(k) \ \#(k, \wedge(=(y,\varphi), =(z,yv)), \sigma) \text{ is false or } \#(k, =(z,\varphi v) ,\sigma) \notag .
\end{align}

\medskip

Given $\sigma \in \Xi(k)$ we assume $\#(k, \wedge(=(y,\varphi), =(z,yv)), \sigma)$ holds and want to show that $\#(k, =(z,\varphi v) ,\sigma)$ then holds. 

\medskip

We have 
\begin{align} 
&P_{\wedge}( \#(k, =(y,\varphi), \sigma), \#(k, =(z,yv), \sigma) ) \ , \notag  \\
&P_{\wedge}( \#(k, y, \sigma) = \#(k, \varphi, \sigma), \#(k, z, \sigma) = \#(k, yv, \sigma) ) \ , \notag  \\
&P_{\wedge}( \#(k, y, \sigma) = \#(k, \varphi, \sigma), \#(k, z, \sigma) = \#(k, y, \sigma) \cdot \#(k, v, \sigma) ) \  \notag .
\end{align}

\smallskip

From there it follows that $\#(k, z, \sigma) = \#(k, \varphi, \sigma) \cdot \#(k, v, \sigma)$.

\medskip 

We have shown that $\#(k, =(z,\varphi v) ,\sigma)$ holds, in fact it can be rewritten 
\begin{align} 
&\#(k, z ,\sigma) = \#(k, \varphi v ,\sigma) \ , \notag  \\
&\#(k, z, \sigma) = \#(k, \varphi, \sigma) \cdot \#(k, v, \sigma) \  \notag .
\end{align}
\end{proof}

\medskip

We can create a set~$A_{\ref{L:axiom-substitute-exp-in-mult}}$ which is the set of all sentences 
\[ \gamma[ x:N, y:N, z:N, u:N, v:N, \to \left( \wedge(=(y,\varphi), =(z,yv)), =(z,\varphi v)  \right) ] \] 
such that
\begin{itemize}
\item $\varphi \in E(k[x:N, y:N, z:N, u:N, v:N])$ ,
\item for each $\sigma \in \Xi(k)$ $\#(k, \varphi, \sigma) \in \mathbb{N}$ .
\end{itemize}

\smallskip

Lemma~\ref{L:axiom-substitute-exp-in-mult} shows us that this set of sentences (which is a potential axiom) is `sound'. In order to use $A_{\ref{L:axiom-substitute-exp-in-mult}}$ as an axiom in our system we also need to show that $A_{\ref{L:axiom-substitute-exp-in-mult}}$ is r.e..\\

\begin{lemma}
$A_{\ref{L:axiom-substitute-exp-in-mult}}$ is r.e..
\end{lemma}

\begin{proof}
$A_{\ref{L:axiom-substitute-exp-in-mult}}$ is the set of all sentences 
\[ \gamma[ x:N, y:N, z:N, u:N, v:N, \to \left( \wedge(=(y,\varphi), =(z,yv)), =(z,\varphi v)  \right) ] \] 
such that $\varphi \in E_{\mathbb{N}}(k[x:N, y:N, z:N, u:N, v:N])$.\\

Let's define a function $\eta$ over $\Sigma^*$ with $\eta(\varphi) =  \gamma[ x:N, y:N, z:N, u:N, v:N, \to \left( \wedge(=(y,\varphi), =(z,yv)), =(z,\varphi v)  \right) ]$.\\

Then $A_{\ref{L:axiom-substitute-exp-in-mult}}$ is simply the set  $\{  \eta(\varphi) | \ \varphi \in E_{\mathbb{N}}(k[x:N, y:N, z:N, u:N, v:N]) \}$.\\

Since $E_{\mathbb{N}}(k[x:N, y:N, z:N, u:N, v:N])$ is r.e. then $A_{\ref{L:axiom-substitute-exp-in-mult}}$ is also r.e..

\end{proof}

\medskip

Then let $A_{\ref{L:axiom-substitute-exp-in-mult}} \in \mathcal{A}$.

\bigskip

\begin{lemma}\label{L:axiom-commutative-mult}
Let $k = k[x:N, y:N, z:N, u:N, v:N]$, $\chi \in S(k)$ then 

\begin{itemize}
\item $=((xu)v,x(uv)) \in S(k)$
\end{itemize}

\smallskip

Moreover
\[ \# \left( \gamma[ x:N, y:N, z:N, u:N, v:N, \to \left( \chi, =((xu)v,x(uv)) \right) ] \right) \]
is true.
\end{lemma}

\begin{proof}

\smallskip

By lemma~\ref{L:xi_in_Ekj} $u \in E(k)$. If we define $k_u = k[x:N, y:N, z:N]$ then for each $\sigma \in \Xi(k)$ $\sigma_{/dom(k_u)} \in \Xi(k_u)$ and $\#(k,u,\sigma) \in \#(k_u, N, \sigma_{/dom(k_u)}) = \#(N) = \mathbb{N}$.\\

Similarly by~\ref{L:xi_in_Ekj} $x \in E(k)$. If we define $k_x = \epsilon$ then for each $\sigma \in \Xi(k)$ $\sigma_{/dom(k_x)} \in \Xi(k_x)$ and $\#(k,x,\sigma) \in \#(k_x, N, \sigma_{/dom(k_x)}) = \#(N) = \mathbb{N}$.\\

Similarly by~\ref{L:xi_in_Ekj} $v \in E(k)$. If we define $k_v = k[x:N, y:N, z:N, u:N]$ then for each $\sigma \in \Xi(k)$ $\sigma_{/dom(k_v)} \in \Xi(k_v)$ and $\#(k,v,\sigma) \in \#(k_v, N, \sigma_{/dom(k_v)}) = \#(N) = \mathbb{N}$.\\

By lemma~\ref{L:preliminary-results-ch6-3} $(*)(x, u) \in E(k)$ and for each $\sigma \in \Xi(k)$ $\#(k, (*)(x, u), \sigma) = (\#(k,x, \sigma) \cdot \#(k,u, \sigma)) \in \mathbb{N}$.\\

Also by lemma~\ref{L:preliminary-results-ch6-3} $(*)(xu, v) \in E(k)$ and for each $\sigma \in \Xi(k)$ 
\small
\[ \#(k, (*)(xu, v), \sigma) = \#(k,xu, \sigma) \cdot \#(k,v, \sigma) = (\#(k,x, \sigma) \cdot \#(k,u, \sigma)) \cdot \#(k,v, \sigma) \ . \]
\normalsize

\medskip

By lemma~\ref{L:preliminary-results-ch6-3} $(*)(u, v) \in E(k)$ and for each $\sigma \in \Xi(k)$ $\#(k, (*)(u, v), \sigma) = (\#(k,u, \sigma) \cdot \#(k,v, \sigma)) \in \mathbb{N}$.\\

Also by lemma~\ref{L:preliminary-results-ch6-3} $(*)(x, uv) \in E(k)$ and for each $\sigma \in \Xi(k)$ 
\small
\[ \#(k, (*)(x, uv), \sigma) = \#(k,x, \sigma) \cdot \#(k,uv, \sigma) = \#(k,x, \sigma) \cdot (\#(k,u, \sigma) \cdot \#(k,v, \sigma)) \ . \]
\normalsize

\medskip

Clearly it follows that for each $\sigma \in \Xi(k)$ 
\[ \#(k, (*)(x, uv), \sigma) = \#(k, (*)(xu, v), \sigma) \ . \]

\medskip

By lemma~\ref{L:preliminary-results-ch6-2} it also follows that $=((xu)v,x(uv)) \in S(k)$ and that for each $\sigma \in \Xi(k)$ $\#(k, =((xu)v,x(uv)) , \sigma )$ is true.\\

Finally we observe that 
\[ \# \left( \gamma[ x:N, y:N, z:N, u:N, v:N, \to \left( \chi, =((xu)v,x(uv)) \right) ] \right) \]

\medskip

can be rewritten as 
\begin{align} 
&\text{for each } \sigma \in \Xi(k) \ \#( k , \to \left( \chi, =((xu)v,x(uv)) \right), \sigma) \ , \notag  \\
&\text{for each } \sigma \in \Xi(k) \ P_{\to} ( \#(k, \chi, \sigma), \#(k, =((xu)v,x(uv)) , \sigma) ) \ \notag .
\end{align}

\smallskip

for each $\sigma \in \Xi(k)$ $\#(k, \chi, \sigma)$ is false or $\#(k, =((xu)v,x(uv)) , \sigma )$.\\

So we have proved it.
\end{proof}

\medskip

We can create a set~$A_{\ref{L:axiom-commutative-mult}}$ which is the set of all sentences 
\[ \gamma[ x:N, y:N, z:N, u:N, v:N, \to \left( \chi, =((xu)v,x(uv)) \right) ] \]
such that
\begin{itemize}
\item $\chi \in S(k[x:N, y:N, z:N, u:N, v:N])$ .
\end{itemize}

\smallskip

Lemma~\ref{L:axiom-commutative-mult} shows us that this set of sentences (which is a potential axiom) is `sound'. In order to use $A_{\ref{L:axiom-commutative-mult}}$ as an axiom in our system we also need to show that $A_{\ref{L:axiom-commutative-mult}}$ is r.e..\\

\begin{lemma}
$A_{\ref{L:axiom-commutative-mult}}$ is r.e..
\end{lemma}

\begin{proof}

$A_{\ref{L:axiom-commutative-mult}}$ is the set of all sentences 
\[ \gamma[ x:N, y:N, z:N, u:N, v:N, \to \left( \chi, =((xu)v,x(uv)) \right) ] \]
such that
\begin{itemize}
\item $\chi \in S(k[x:N, y:N, z:N, u:N, v:N])$ .\\
\end{itemize}

Let's define a function $\eta$ over $\Sigma^*$ with 
\[ \eta(\chi) = \gamma[ x:N, y:N, z:N, u:N, v:N, \to \left( \chi, =((xu)v,x(uv)) \right) ] \ .  \] 

Then $A_{\ref{L:axiom-commutative-mult}}$ is simply the set  $\{  \eta(\chi) | \ \chi \in S(k[x:N, y:N, z:N, u:N, v:N]) \}$.\\

Since $S(k[x:N, y:N, z:N, u:N, v:N])$ is r.e. then $A_{\ref{L:axiom-commutative-mult}}$ is also r.e..
\end{proof}

\medskip

Then let $A_{\ref{L:axiom-commutative-mult}} \in \mathcal{A}$.

\bigskip

\begin{lemma}\label{L:ax-eq-trans}
Let $m$ be a positive integer. Let $x_1, \dots , x_m \in \mathcal{V}$, with $x_i \ne x_j$ for $i \ne j$. Let $\varphi_1, \dots , \varphi_m \in E$ and assume $H[x_1:\varphi_1, \dots , x_m:\varphi_m ]$. Define $k = k[x_1:\varphi_1, \dots , x_m:\varphi_m ]$ and let $\chi \in S(k)$, $\varphi, \psi, \theta \in E(k)$.

\medskip

Then\\ 
$\#( \gamma[ x_1: \varphi_1, \dots , x_m: \varphi_m, \to \left( \chi, \to \left( \wedge \left(=(\varphi,\psi),=(\psi,\theta) \right), =(\varphi, \theta) \right) \right) ] )$.
\end{lemma}

\begin{proof}

\smallskip

We can rewrite 
\[ \#( \gamma[ x_1: \varphi_1, \dots , x_m: \varphi_m, \to \left( \chi, \to \left( \wedge \left(=(\varphi,\psi),=(\psi,\theta) \right), =(\varphi, \theta) \right) \right) ] ) \] 

\smallskip

as
\begin{align} 
&\text{for each } \sigma \in \Xi(k) \ \#( k , \to \left( \chi, \to \left( \wedge \left(=(\varphi,\psi),=(\psi,\theta) \right), =(\varphi, \theta) \right) \right), \sigma) \ , \notag  \\
&\text{for each } \sigma \in \Xi(k) \ P_{\to} ( \#(k, \chi , \sigma), \#(k, \to \left( \wedge \left(=(\varphi,\psi),=(\psi,\theta) \right), =(\varphi, \theta) \right), \sigma)  ) \ , \notag  \\
&\text{for each } \sigma \in \Xi(k) \ \#(k, \chi , \sigma) \text{ is false or } \#(k, \to \left( \wedge \left(=(\varphi,\psi),=(\psi,\theta) \right), =(\varphi, \theta) \right), \sigma) \notag .
\end{align}

We can rewrite $\#(k, \to \left( \wedge \left(=(\varphi,\psi),=(\psi,\theta) \right), =(\varphi, \theta) \right), \sigma)$ as
\begin{align} 
&P_{\to}( \#(k, \wedge \left(=(\varphi,\psi),=(\psi,\theta) \right), \sigma), \#(k, =(\varphi, \theta), \sigma) ) \ , \notag  \\
&P_{\to}( P_{\wedge} ( \#(k, =(\varphi,\psi), \sigma), \#(k, =(\psi,\theta), \sigma) ), \#(k, =(\varphi, \theta), \sigma) ) \ , \notag \\
&P_{\to}( P_{\wedge} ( \#(k, \varphi, \sigma) = \#(k, \psi, \sigma), \#(k, \psi, \sigma) = \#(k, \theta, \sigma) ), \#(k, \varphi, \sigma) = \#(k, \theta, \sigma) ) \  \notag .
\end{align}

\smallskip

($\#(k, \varphi, \sigma) = \#(k, \psi, \sigma)$ and $\#(k, \psi, \sigma) = \#(k, \theta, \sigma)$) is false or\\
$\#(k, \varphi, \sigma) = \#(k, \theta, \sigma)$.\\

If ($\#(k, \varphi, \sigma) = \#(k, \psi, \sigma)$ and $\#(k, \psi, \sigma) = \#(k, \theta, \sigma)$) is false then $\#(k, \to \left( \wedge \left(=(\varphi,\psi),=(\psi,\theta) \right), =(\varphi, \theta) \right), \sigma)$ is true.\\

Otherwise clearly $\#(k, \varphi, \sigma) = \#(k, \theta, \sigma)$ holds and so\\
$\#(k, \to \left( \wedge \left(=(\varphi,\psi),=(\psi,\theta) \right), =(\varphi, \theta) \right), \sigma)$ is true all the same.

\end{proof}

\medskip

We can create a set~$A_{\ref{L:ax-eq-trans}}$ which is the set of all sentences 
\[ \gamma[ x_1: \varphi_1, \dots , x_m: \varphi_m, \to \left( \chi, \to \left( \wedge \left(=(\varphi,\psi),=(\psi,\theta) \right), =(\varphi, \theta) \right) \right) ] \]

such that
\begin{itemize}
\item $m$ is a positive integer, $x_1, \dots , x_m \in \mathcal{V}$, $x_i \ne x_j$ for $i \ne j$, $\varphi_1, \dots , \varphi_m \in E$, $H[x_1:\varphi_1, \dots , x_m:\varphi_m ]$,
\item $\varphi, \psi, \theta \in E(k[x_1:\varphi_1, \dots , x_m:\varphi_m ])$,
\item $\chi \in S(k[x_1:\varphi_1, \dots , x_m:\varphi_m ])$.
\end{itemize}

\smallskip

Lemma~\ref{L:ax-eq-trans} shows us that this set of sentences (which is a potential axiom) is `sound'. In order to use $A_{\ref{L:ax-eq-trans}}$ as an axiom in our system we also need to show that $A_{\ref{L:ax-eq-trans}}$ is r.e..\\

\begin{lemma}
$A_{\ref{L:ax-eq-trans}}$ is r.e..
\end{lemma}

\begin{proof}
Given a positive integer $m$ and $(x_1, \varphi_1, \dots , x_m, \varphi_m) \in R_m$ we can notice the following:

\begin{itemize}
\item $k[x_1:\varphi_1, \dots , x_m:\varphi_m ] \in K$;
\item $E(k[x_1:\varphi_1, \dots , x_m:\varphi_m ])$ is r.e.;
\item $S(k[x_1:\varphi_1, \dots , x_m:\varphi_m ])$ is r.e.;
\item $\{(x_1, \varphi_1, \dots , x_m, \varphi_m)\} \times E(k[x_1:\varphi_1, \dots , x_m:\varphi_m ])^3 \times S(k[x_1:\varphi_1, \dots , x_m:\varphi_m ])$ is r.e..
\end{itemize}

So we can define the following
\small
\[ Q_{m,4} = \bigcup_{(x_1, \varphi_1, \dots , x_m, \varphi_m) \in R_m} \{(x_1, \varphi_1, \dots , x_m, \varphi_m)\} \times E(k[x_1:\varphi_1, \dots , x_m:\varphi_m ])^3 \times S(k[x_1:\varphi_1, \dots , x_m:\varphi_m ]) \ . \]
\normalsize

Clearly $Q_{m,4} \subseteq (\Sigma^*)^{2m} \times (\Sigma^*)^3 \times \Sigma^*$ is r.e..\\

We can define a function $\eta$ over $(\Sigma^*)^{2m} \times (\Sigma^*)^3 \times \Sigma^*$ such that for each $((\psi_1, \varphi_1, \dots , \psi_m, \varphi_m), (\varphi, \psi, \theta), \chi) \in (\Sigma^*)^{2m} \times (\Sigma^*)^3 \times \Sigma^*$ 
\small
\[ \eta(((\psi_1, \varphi_1, \dots , \psi_m, \varphi_m), (\varphi, \psi, \theta), \chi)) = \gamma[ x_1: \varphi_1, \dots , x_m: \varphi_m, \to \left( \chi, \to \left( \wedge \left(=(\varphi,\psi),=(\psi,\theta) \right), =(\varphi, \theta) \right) \right) ] \ . \]
\normalsize

\smallskip

Now $\eta$ clearly is a computable function and so the set $\{ \eta(((x_1, \varphi_1, \dots , x_m, \varphi_m), (\varphi, \psi, \theta), \chi)) | \, ((x_1, \varphi_1, \dots , x_m, \varphi_m), (\varphi, \psi, \theta), \chi) \in Q_{m,4} \}$  is a r.e. subset of $\Sigma^*$. And finally the set 
\[ \bigcup_{m \geqslant 1} \{ \eta(((x_1, \varphi_1, \dots , x_m, \varphi_m), (\varphi, \psi, \theta), \chi)) | \, ((x_1, \varphi_1, \dots , x_m, \varphi_m), (\varphi, \psi, \theta), \chi) \in Q_{m,4} \}  \]

is itself a r.e. set. It should be clear at this point that this set is actually our set $A_{\ref{L:ax-eq-trans}}$, and so that $A_{\ref{L:ax-eq-trans}}$ is r.e..
\end{proof}

\medskip

Then let $A_{\ref{L:ax-eq-trans}} \in \mathcal{A}$.

\bigskip

\begin{lemma}\label{L:rule-mp}
Let $m$ be a positive integer. Let $x_1, \dots , x_m \in \mathcal{V}$, with $x_i \ne x_j$ for $i \ne j$. Let $\varphi_1, \dots , \varphi_m \in E$ and assume $H[x_1:\varphi_1, \dots , x_m:\varphi_m ]$. Define $k = k[x_1:\varphi_1, \dots , x_m:\varphi_m ]$ and let $\varphi, \psi, \chi \in S(k)$. Under these assumptions, if
\begin{align} 
&\#( \gamma[ x_1: \varphi_1, \dots , x_m: \varphi_m, \to \left( \chi, \varphi \right) ]) \ , \notag  \\
&\#( \gamma[ x_1: \varphi_1, \dots , x_m: \varphi_m, \to \left( \chi, \to \left( \varphi, \psi \right) \right) ]) \  \notag .
\end{align}

\smallskip

then
\[ \#( \gamma[ x_1: \varphi_1, \dots , x_m: \varphi_m, \to \left( \chi, \psi \right) ]). \]

\end{lemma}

\begin{proof}

\smallskip

We can rewrite 
\[ \#( \gamma[ x_1: \varphi_1, \dots , x_m: \varphi_m, \to \left( \chi, \varphi \right) ]) \]

\smallskip 

as 
\begin{align} 
&\text{for each } \sigma \in \Xi(k) \ \#( k , \to \left( \chi, \varphi \right), \sigma) \ , \notag  \\
&\text{for each } \sigma \in \Xi(k) \ P_{\to} ( \#(k, \chi , \sigma), \#(k, \varphi , \sigma)  ) \notag .
\end{align}

\medskip

We can rewrite 
\[ \#( \gamma[ x_1: \varphi_1, \dots , x_m: \varphi_m, \to \left( \chi, \to \left( \varphi, \psi \right) \right) ]) \]

\smallskip 

as 
\begin{align} 
&\text{for each } \sigma \in \Xi(k) \ \#( k , \to \left( \chi, \to \left( \varphi, \psi \right) \right), \sigma) \ , \notag  \\
&\text{for each } \sigma \in \Xi(k) \ P_{\to} \left( \#(k, \chi , \sigma), \#(k, \to \left( \varphi, \psi \right) , \sigma) \right) \ , \notag \\
&\text{for each } \sigma \in \Xi(k) \ P_{\to} \left( \#(k, \chi , \sigma), P_{\to} \left( \#(k, \varphi, \sigma ), \#(k, \psi, \sigma ) \right) \right) \notag .
\end{align}

\medskip

Finally we can rewrite 
\[ \#( \gamma[ x_1: \varphi_1, \dots , x_m: \varphi_m, \to \left( \chi, \psi \right) ]) \]

\smallskip 

as 
\begin{align} 
&\text{for each } \sigma \in \Xi(k) \ \#( k , \to \left( \chi, \psi \right), \sigma) \ , \notag  \\
&\text{for each } \sigma \in \Xi(k) \ P_{\to} ( \#(k, \chi , \sigma), \#(k, \psi , \sigma)  ) \notag .
\end{align}

\bigskip

If we assume both
\begin{align} 
&\#( \gamma[ x_1: \varphi_1, \dots , x_m: \varphi_m, \to \left( \chi, \varphi \right) ]) \ , \notag  \\
&\#( \gamma[ x_1: \varphi_1, \dots , x_m: \varphi_m, \to \left( \chi, \to \left( \varphi, \psi \right) \right) ]) \  \notag .
\end{align}

\smallskip

then for each $\sigma \in \Xi(k)$ 
\begin{itemize}
\item $\#(k, \chi , \sigma)$ is false or $\#(k, \varphi , \sigma)$ is true,
\item $\#(k, \chi , \sigma)$ is false or $\#(k, \varphi , \sigma)$ is false or $\#(k, \psi , \sigma)$ is true.
\end{itemize}

\smallskip

Clearly this implies $\#(k, \chi , \sigma)$ is false or $\#(k, \psi , \sigma)$ is true.

\medskip

Therefore in our assumptions $\#( \gamma[ x_1: \varphi_1, \dots , x_m: \varphi_m, \to \left( \chi, \psi \right) ])$ holds.
\end{proof}

\medskip

We can create a set~$R_{\ref{L:rule-mp}}$ which is the set of all $3$-tuples
\[ \left(   
\begin{array} {l} 
{\gamma[ x_1: \varphi_1, \dots , x_m: \varphi_m, \to \left( \chi, \varphi \right) ], } \\
{\gamma[ x_1: \varphi_1, \dots , x_m: \varphi_m, \to \left( \chi, \to \left( \varphi, \psi \right) \right) ], } \\
{\gamma[ x_1: \varphi_1, \dots , x_m: \varphi_m, \to \left( \chi, \psi \right) ] } 
\end{array}
\right) \]

such that
\begin{itemize}
\item $m$ is a positive integer, $x_1, \dots , x_m \in \mathcal{V}$, $x_i \ne x_j$ 
for $i \ne j$, $\varphi_1, \dots , \varphi_m \in E$, $H[x_1:\varphi_1, \dots , 
x_m:\varphi_m ]$,
\item $\varphi, \psi, \chi \in S(k[x_1:\varphi_1, \dots , x_m:\varphi_m ])$.\\
\end{itemize}

\medskip

Lemma~\ref{L:rule-mp} shows us that this set (which is a potential 2-ary rule) is `sound'. In order to use $R_{\ref{L:rule-mp}}$ as a rule in our system we also need to show that $R_{\ref{L:rule-mp}}$ is r.e..\\

\begin{lemma}
$R_{\ref{L:rule-mp}}$ is r.e.
\end{lemma}

\begin{proof}
Given a positive integer $m$ and $(x_1, \varphi_1, \dots , x_m, \varphi_m) \in R_m$ we can notice the following:

\begin{itemize}
\item $k[x_1:\varphi_1, \dots , x_m:\varphi_m ] \in K$;
\item $S(k[x_1:\varphi_1, \dots , x_m:\varphi_m ])$ is r.e.;
\item $\{(x_1, \varphi_1, \dots , x_m, \varphi_m)\} \times S(k[x_1:\varphi_1, \dots , x_m:\varphi_m ])^3$ is r.e..
\end{itemize}

Let's define
\[ Q_{m,3} = \bigcup_{(x_1, \varphi_1, \dots , x_m, \varphi_m) \in R_m} \{(x_1, \varphi_1, \dots , x_m, \varphi_m)\} \times S(k[x_1:\varphi_1, \dots , x_m:\varphi_m ])^3 \ . \]

\smallskip

Clearly $Q_{m,3} \subseteq (\Sigma^*)^{2m} \times (\Sigma^*)^3$ is also r.e..\\

We now define three functions $\delta_{1,m}, \ \delta_{2,m}, \ \delta_{3,m}$ over $(\Sigma^*)^{2m} \times (\Sigma^*)^3$ as follows.
Given $((\psi_1, \varphi_1, \dots , \psi_m, \varphi_m), (\varphi, \psi, \chi)) \in (\Sigma^*)^{2m} \times (\Sigma^*)^3$
\begin{align} 
&\delta_{1,m}((\psi_1, \varphi_1, \dots , \psi_{m}, \varphi_{m}), (\varphi, \psi, \chi)) = \gamma[ \psi_1: \varphi_1, \dots , \psi_m: \varphi_m, \to \left( \chi, \varphi \right) ] \ , \notag  \\
&\delta_{2,m}((\psi_1, \varphi_1, \dots , \psi_{m}, \varphi_{m}), (\varphi, \psi, \chi)) = \gamma[ \psi_1: \varphi_1, \dots , \psi_m: \varphi_m, \to \left( \chi, \to \left( \varphi, \psi \right) \right) ] \ , \notag \\
&\delta_{3,m}((\psi_1, \varphi_1, \dots , \psi_{m}, \varphi_{m}), (\varphi, \psi, \chi)) = \gamma[ \psi_1: \varphi_1, \dots , \psi_m: \varphi_m, \to \left( \chi, \psi \right) ] \  \notag .
\end{align}

\smallskip

All of the three functions we have defined are computable functions from $(\Sigma^*)^{2m} \times (\Sigma^*)^3$ to $\Sigma^*$. If we define a function $\delta_m$ over $(\Sigma^*)^{2m} \times (\Sigma^*)^3$ as follows: 
\[
\delta_m((\psi_1, \varphi_1, \dots , \psi_{m}, \varphi_{m}), (\varphi, \psi, \chi)) = 
\left(
\begin{array} {l}
{ \delta_{1,m}((\psi_1, \varphi_1, \dots , \psi_{m}, \varphi_{m}), (\varphi, \psi, \chi)) , } \\
{ \delta_{2,m}((\psi_1, \varphi_1, \dots , \psi_{m}, \varphi_{m}), (\varphi, \psi, \chi)) , } \\
{ \delta_{3,m}((\psi_1, \varphi_1, \dots , \psi_{m}, \varphi_{m}), (\varphi, \psi, \chi)) , }
\end{array}
\right)
\]

then $\delta_m$ is a computable function from $(\Sigma^*)^{2m} \times (\Sigma^*)^3$ to $(\Sigma^*)^3$, therefore the set 
\[ D_m = \{ \delta_m((\psi_1, \varphi_1, \dots , \psi_{m}, \varphi_{m}), (\varphi, \psi, \chi)) | \ ((\psi_1, \varphi_1, \dots , \psi_{m}, \varphi_{m}), (\varphi, \psi, \chi)) \in Q_{m,3} \} \]
is a r.e. subset of $(\Sigma^*)^3$.\\

If we now consider the set $\bigcup_{m \geqslant 1} D_m$ then this is a r.e. subset of $(\Sigma^*)^3$ and actually this set is equal to our set $R_{\ref{L:rule-mp}}$ which so is r.e. itself.
\end{proof}
 
\medskip

Then let $R_{\ref{L:rule-mp}} \in \mathcal{R}$.

\bigskip

\begin{lemma}\label{L:rule-introd-exists-mult}
Let $k = k[x:N, y:N, z:N, u:N, v:N]$, $\chi \in S(k)$, $\varphi \in E(k)$ such that for each $\sigma \in \Xi(k)$ $\#(k, \varphi, \sigma) \in \mathbb{N}$ then 

\begin{itemize}
\item $=(z, x \varphi) \in S(k)$,
\item $\exists(w:N, =(z, xw)) \in S(k)$.
\end{itemize}

\smallskip

Under these assumptions if
\[ \# \left( \gamma[ x:N, y:N, z:N, u:N, v:N, \to \left( \chi,  =(z, x \varphi) \right) ] \right) \]

then
\[ \# \left( \gamma[ x:N, y:N, z:N, u:N, v:N, \to \left( \chi, \exists( w:N, =(z, xw) ) \right) ] \right) \]
is true.
\end{lemma}

\begin{proof}

\smallskip

By ~\ref{L:xi_in_Ekj} $x \in E(k)$. If we define $k_x = \epsilon$ then for each $\sigma \in \Xi(k)$ $\sigma_{/dom(k_x)} \in \Xi(k_x)$ and $\#(k,x,\sigma) \in \#(k_x, N, \sigma_{/dom(k_x)}) = \#(N) = \mathbb{N}$.

\medskip

By lemma~\ref{L:preliminary-results-ch6-3} it follows that $(*)(x, \varphi) \in E(k)$ and for each $\sigma \in \Xi(k)$ $\#(k, (*)(x, \varphi), \sigma) = (\#(k,x, \sigma) \cdot \#(k,\varphi, \sigma)) \in \mathbb{N}$.

\medskip

By ~\ref{L:xi_in_Ekj} $z \in E(k)$, so we can apply lemma~\ref{L:preliminary-results-ch6-2} and obtain that $=(z, x \varphi)$ belongs to $S(k)$.

\medskip

Let $h = k + <w,N>$, we have $N \in E(k)$ and for each $\sigma \in \Xi(k)$ $\#(k,N,\sigma) = \#(N) = \mathbb{N}$. So $N \in E_s(k)$. Moreover $w \in (\mathcal{V} - var(k))$ so by lemma~\ref{L:quantifiers_S_h} $h \in K$. 

\medskip

We now want to show that $=(z, xw)$ belongs to $S(h)$. Since $N \in E_s(k)$ we have $H[x:N, y:N, z:N, u:N, v:N, w:N]$. We have \[ k[x:N, y:N, z:N, u:N, v:N, w:N] = k + <w,N> = h \ . \]

\medskip

Using lemma~\ref{L:xi_in_Ekj} we obtain that $z, x, w \in E(h)$. If we define $h_x = \epsilon$ then for each $\rho \in \Xi(h)$ $\rho_{/dom(h_x)} \in \Xi(h_x)$ and $\#(h,x,\rho) \in \#(h_x, N, \rho_{/dom(h_x)}) = \#(N) = \mathbb{N}$.

\medskip

Moreover for each $\rho \in \Xi(h)$ $\rho_{/dom(k)} \in \Xi(k)$ $\#(h,w,\rho) \in \#(k, N, \rho_{/dom(k)}) = \#(N) = \mathbb{N}$.

\medskip

By lemma~\ref{L:preliminary-results-ch6-3} it follows that $(*)(x, w) \in E(h)$ and for each $\rho \in \Xi(h)$ $\#(h, (*)(x, w), \rho) = (\#(h,x, \rho) \cdot \#(h,w, \rho)) \in \mathbb{N}$.

\medskip

By lemma~\ref{L:preliminary-results-ch6-2} $=(z, xw)$ belongs to $S(h)$. We can now apply lemma~\ref{L:quantifiers_S_h} and obtain that $\exists( w:N, =(z, xw) )  \in S(k)$ and for each $\sigma \in \Xi(k)$
\[ \#(k, \exists ( w:N, =(z, xw) ), \sigma)  = \text{ exists } \rho \in \Xi(h):  \sigma \sqsubseteq \rho, \  \#(h, =(z, xw), \rho)  \ . \]

\smallskip

We can rewrite 
\[ \# \left( \gamma[ x:N, y:N, z:N, u:N, v:N, \to \left( \chi, =(z, x \varphi) \right) ] \right) \]
as:
\begin{align} 
&\text{for each } \sigma \in \Xi(k) \ \#( k , \to \left( \chi,  =(z, x \varphi) \right), \sigma) \ , \notag  \\
&\text{for each } \sigma \in \Xi(k) \ P_{\to}( \#(k, \chi, \sigma),  \#(k, =(z, x \varphi) , \sigma) ) \ , \notag \\
&\text{for each } \sigma \in \Xi(k) \ \#(k, \chi, \sigma) \text{ is false or } \#(k, =(z, x \varphi) , \sigma) \notag .
\end{align}

\medskip

We can rewrite
\[ \# \left( \gamma[ x:N, y:N, z:N, u:N, v:N, \to \left( \chi, \exists( w:N, =(z, xw) ) \right) ] \right) \]
as:
\begin{align} 
&\text{for each } \sigma \in \Xi(k) \ \#( k , \to \left( \chi, \exists( w:N, =(z, xw) ) \right), \sigma) \ , \notag  \\
&\text{for each } \sigma \in \Xi(k) \ P_{\to}( \#(k, \chi, \sigma), \#(k, \exists( w:N, =(z, xw) ), \sigma)  ) \ , \notag \\
&\text{for each } \sigma \in \Xi(k) \ \#(k, \chi, \sigma) \text{ is false or } \#(k, \exists( w:N, =(z, xw) ), \sigma) \ , \notag \\ 
&\text{for each } \sigma \in \Xi(k) \ \#(k, \chi, \sigma) \text{ is false or exists } \rho \in \Xi(h):  \sigma \sqsubseteq \rho, \  \#(h, =(z, xw), \rho) \notag .
\end{align}

\medskip

We now assume 
\[ \# \left( \gamma[ x:N, y:N, z:N, u:N, v:N, \to \left( \chi, =(z, x \varphi) \right) ] \right) \]
and try to prove
\[ \# \left( \gamma[ x:N, y:N, z:N, u:N, v:N, \to \left( \chi, \exists( w:N, =(z, xw) ) \right) ] \right) \ . \]

\medskip

Let $\sigma \in \Xi(k)$, if $\#(k, \chi, \sigma)$ is false then our proof is already finished. So we assume $\#(k, \chi, \sigma)$ is true. In this case $\#(k, =(z, x \varphi), \sigma)$ holds.

\medskip

It follows that  
\[ \#(k, z, \sigma) = \#(k, x \varphi, \sigma) = \#(k,x, \sigma) \cdot \#(k,\varphi, \sigma) \ . \]

\smallskip

We have to show there exists $\rho \in \Xi(h)$ such that $\sigma \sqsubseteq \rho$ and $\#(h, =(z, xw), \rho)$. We can rewrite $\#(h, =(z, xw), \rho)$ as 
\[ \#(h, z, \rho) = \#(h, xw, \rho) = \#(h, x, \rho) \cdot \#(h, w, \rho) \ . \]

\smallskip

Let's define $\rho = \sigma + (w, \#(k, \varphi, \sigma))$. There exists a positive integer $n$ such that $h \in K(n)$. Since $h \ne \epsilon$ we have $n \geqslant 2$ and by lemma~\ref{L:k-in-km+} there exists $q < n$ such that $h \in K(q)^+$. Then there exist $g \in K(q)$, $\phi \in E_s(q,g)$, $\alpha \in (\mathcal{V}-var(g))$ such that $h = g + <\alpha, \phi>$ and 
\[
\Xi(h) = \{ \delta + (\alpha,s) | \, \delta \in \Xi(g), s \in \#(g,\phi,\delta) \} \, .
\]

\smallskip

Now we have also $h = k + <w,N>$ therefore 
\[
\Xi(h) = \{ \delta + (w,s) | \, \delta \in \Xi(k), s \in \mathbb{N} \} \, .
\]

\smallskip

It follows that $\rho \in \Xi(h)$ and moreover

\[ \#(h, w, \rho) = \#(h, w, \rho)_{(q+1,h,a)} = \#(k, \varphi, \sigma) \ . \]

\smallskip

We have also $\#(h, z, \rho) = \#(k, z, \sigma)$. In fact $z \in E(h) \cap E(k)$, $k \sqsubseteq h$, $\sigma \sqsubseteq \rho$ and we can use lemma~\ref{L:kappa-acca-sigma-ro}. Similarly we obtain $\#(h, x, \rho) = \#(k, x, \sigma)$. Since 
\[ \#(k, z, \sigma) = \#(k,x, \sigma) \cdot \#(k,\varphi, \sigma) \ . \]

we have 
\[ \#(h, z, \rho) = \#(h, x, \rho) \cdot \#(h, w, \rho) \ . \]

\smallskip

and then of course $\#(h, =(z, xw), \rho)$.
\end{proof}

\medskip

We can create a set~$R_{\ref{L:rule-introd-exists-mult}}$ which is the set of all pairs
\[ \left(   
\begin{array} {l} 
{\gamma[ x:N, y:N, z:N, u:N, v:N, \to \left( \chi, =(z, x \varphi) \right) ], } \\
{\gamma[ x:N, y:N, z:N, u:N, v:N, \to \left( \chi, \exists( w:N, =(z, xw) ) \right) ] } \\
\end{array}
\right) \]
such that $\chi \in S(k)$, $\varphi \in E(k)$ such that for each $\sigma \in \Xi(k)$ $\#(k, \varphi, \sigma) \in \mathbb{N}$.\\

Lemma~\ref{L:rule-introd-exists-mult} shows us that this set (which is a potential 1-ary rule) is `sound'. In order to use $R_{\ref{L:rule-introd-exists-mult}}$ as a rule in our system we also need to show that $R_{\ref{L:rule-introd-exists-mult}}$ is r.e..\\

\begin{lemma}
$R_{\ref{L:rule-introd-exists-mult}}$ is r.e..
\end{lemma}

\begin{proof}

Our set~$R_{\ref{L:rule-introd-exists-mult}}$ is the set of all pairs
\[ \left(   
\begin{array} {l} 
{\gamma[ x:N, y:N, z:N, u:N, v:N, \to \left( \chi, =(z, x \varphi) \right) ], } \\
{\gamma[ x:N, y:N, z:N, u:N, v:N, \to \left( \chi, \exists(w:N, =(z, xw) ) \right) ] } \\
\end{array}
\right) \]
such that $\chi \in S(k[x:N, y:N, z:N, u:N, v:N])$, $\varphi \in E_{\mathbb{N}}(k[x:N, y:N, z:N, u:N, v:N])$.\\

We now define two functions $\delta_1, \delta_2$ over $(\Sigma^*)^2$ as follows. Given $\chi, \varphi \in \Sigma^*$
\begin{align} 
&\delta_1(\chi, \varphi) = \gamma[ x:N, y:N, z:N, u:N, v:N, \to \left( \chi, =(z, x \varphi) \right) ] \ , \notag  \\
&\delta_2(\chi, \varphi) = \gamma[ x:N, y:N, z:N, u:N, v:N, \to \left( \chi, \exists(w:N, =(z, xw) ) \right) ] \notag .
\end{align}
 
 All of the two functions we have defined are computable functions from $(\Sigma^*)^2$ to $\Sigma^*$. If we define a function $\delta$ from $(\Sigma^*)^2$ to $(\Sigma^*)^2$ as follows: 
\[
\delta(\chi, \varphi) = 
\left(
\begin{array} {l}
{ \delta_1(\chi, \varphi) , } \\
{ \delta_2(\chi, \varphi) }
\end{array}
\right)
\]

then $\delta$ is a computable function from $(\Sigma^*)^2$ to $(\Sigma^*)^2$. We can actually rewrite $R_{\ref{L:rule-introd-exists-mult}}$ as
\[ \{ \delta(\chi, \varphi) | \ (\chi, \varphi) \in S(k[x:N, y:N, z:N, u:N, v:N]) \times E_{\mathbb{N}}(k[x:N, y:N, z:N, u:N, v:N] \} \ . \]

Since $S(k[x:N, y:N, z:N, u:N, v:N])$ and $E_{\mathbb{N}}(k[x:N, y:N, z:N, u:N, v:N]$ are r.e. then $R_{\ref{L:rule-introd-exists-mult}}$ is r.e. itself.
\end{proof}

\medskip

Then let $R_{\ref{L:rule-introd-exists-mult}} \in \mathcal{R}$.\\

\begin{lemma}\label{L:move-inside-universal-quantifier}
Let $m$ be a positive integer. Let $x_1, \dots , x_{m+1} \in \mathcal{V}$, with $x_i \ne x_j$ for $i 
\ne j$. Let $\varphi_1, \dots , \varphi_{m+1} \in E$ and assume $H[x_1:\varphi_1, \dots , x_{m+1}:
\varphi_{m+1} ]$. 

\medskip

Define $k = k[x_1:\varphi_1, \dots , x_{m+1}:\varphi_{m+1} ]$. Of course $H[x_1:\varphi_1, \dots , 
x_m:\varphi_m ]$ also holds, we define $h = k[x_1:\varphi_1, \dots , x_m:\varphi_m ]$. Let $\chi \in 
S(h) \cap S(k)$, $\varphi \in S(k)$.

\medskip

Under these assumptions we have 
\begin{itemize}
\item $\forall ( x_{m+1} : \varphi_{m+1}, \varphi ) \in S(h)$,
\item $\to \left( \chi, \forall ( x_{m+1} : \varphi_{m+1}, \varphi ) \right) \in S(h)$,
\item $\gamma[ x_1: \varphi_1, \dots , x_m: \varphi_m, \to \left( \chi, \forall ( x_{m+1} : \varphi_{m+1}, \varphi ) \right) ] \in S(\epsilon)$,
\item $\gamma[ x_1: \varphi_1, \dots , x_{m+1}: \varphi_{m+1}, \to \left( \chi, \varphi \right) ] \in S(\epsilon)$.
\end{itemize}

\medskip

Moreover if $\#(\gamma[ x_1: \varphi_1, \dots , x_{m+1}: \varphi_{m+1}, \to \left( \chi, \varphi \right) ])$ then 
\[ \#( \gamma[ x_1: \varphi_1, \dots , x_m: \varphi_m, \to \left( \chi, \forall ( x_{m+1} : \varphi_{m+1}, \varphi ) \right) ] ) \ . \]

\end{lemma}

\begin{proof}

By lemma~\ref{L:quantifiers_S_h} $\forall ( x_{m+1} : \varphi_{m+1}, \varphi ) \in S(h)$, and clearly all the other `preliminary' results hold.

\medskip

We can rewrite 
\[ \#(\gamma[ x_1: \varphi_1, \dots , x_{m+1}: \varphi_{m+1}, \to \left( \chi, \varphi \right) ]) \] 

as
\begin{align} 
&\text{for each } \sigma \in \Xi(k) \ \#\left( k , \to \left(\chi, \varphi \right), \sigma 
\right) \ , \notag  \\
&\text{for each } \sigma \in \Xi(k) \ P_{\to} \left( \#\left(k, \chi , \sigma \right)  , \#\left(k, 
\varphi , \sigma \right) \right) \ , \notag \\
&\text{for each } \sigma \in \Xi(k) \ \#\left(k, \chi , \sigma \right) \text{ is false or } \#\left(k, 
\varphi , \sigma \right) \notag .
\end{align}

\medskip

We can rewrite 
\[ \#(\gamma[ x_1: \varphi_1, \dots , x_m: \varphi_m, \to \left( \chi, \forall ( x_{m+1} : \varphi_{m+1}, \varphi ) \right) ]) \] 

as
\begin{align} 
&\text{for each } \rho \in \Xi(h) \ \#\left( h , \to \left(\chi, \forall \left( x_
{m+1}: \varphi_{m+1}, \varphi \right) \right), \rho 
\right) \ , \notag  \\
&\text{for each } \rho \in \Xi(h) \ P_{\to} \left( \#\left(h, \chi , \rho \right)  , \#\left(h, 
\forall \left( x_{m+1}: \varphi_{m+1}, \varphi 
\right) , \rho \right) \right) \ , \notag \\
&\text{for each } \rho \in \Xi(h) \ \#\left(h, \chi , \rho \right) \text{ is false or }
\text{for each } \sigma \in \Xi(k) \text{ such that } \rho \sqsubseteq \sigma \ \#(k, \varphi, \sigma) \notag .
\end{align}

\medskip

Let $\rho \in \Xi( h )$ and $\#\left(h, \chi , \rho \right)$, let $\sigma \in \Xi(k)$ such that $\rho \sqsubseteq \sigma$, we want to show that $\# \left(k, \varphi, \sigma \right)$ holds. To show this it is clearly enough to show that $\#\left(k, \chi , \sigma \right)$ holds. To do this we can use lemma~\ref{L:kappa-acca-sigma-ro}. In fact there exists a positive integer $n$ such that $h \in K(n)$, $\chi \in E(n, h)$, $k \in K(n)$, $\chi \in E(n, k)$. Given that $\rho \in \Xi( h )$, $\sigma \in \Xi(k)$, $\rho \sqsubseteq \sigma$ we can apply that lemma and get $\#(h, \chi, \rho) = \#(k, \chi, \sigma)$, so $\#(k, \chi, \sigma)$ is proved.
\end{proof}

\medskip

We can create a set~$R_{\ref{L:move-inside-universal-quantifier}}$ 
which is the set of all pairs
\[ \left(   
\begin{array} {l} 
{ \gamma[ x_1: \varphi_1, \dots , x_{m+1}: \varphi_{m+1}, \to \left( \chi, \varphi \right) ], } \\
{ \gamma[ x_1: \varphi_1, \dots , x_m: \varphi_m, \to \left( \chi, \forall ( x_{m+1} : \varphi_{m+1}, \varphi ) \right) ] } 
\end{array}
\right) \]

such that 
\begin{itemize}
\item $m$ is a positive integer, $x_1, \dots , x_{m+1} \in \mathcal{V}$, with $x_i \ne x_j$ for $i 
\ne j$, $\varphi_1, \dots , \varphi_{m+1} \in E$, $H[x_1:\varphi_1, \dots , x_{m+1}:\varphi_{m+1} ]$;
\item if we define $k = k[x_1:\varphi_1, \dots , x_{m+1}:\varphi_{m+1} ]$ and $h = k[x_1:\varphi_1, 
\dots , x_m:\varphi_m ]$ then $\chi \in S(h) \cap S(k)$, $\varphi \in S(k)$.
\end{itemize}

\medskip

Lemma~\ref{L:move-inside-universal-quantifier} shows us that this set (which is a potential 1-ary rule) is `sound'. In order to use $R_{\ref{L:move-inside-universal-quantifier}}$ as a rule in our system we also need to show that $R_{\ref{L:move-inside-universal-quantifier}}$ is r.e..\\

\begin{lemma}
$R_{\ref{L:move-inside-universal-quantifier}}$ is r.e..
\end{lemma}

\begin{proof}
Given a positive integer $m$ and $(x_1, \varphi_1, \dots , x_{m+1}, \varphi_{m+1}) \in R_{m+1}$ all of the following sets are r.e.:
\begin{itemize}
\item $S(k[x_1:\varphi_1, \dots , x_m:\varphi_m ])$,
\item $S(k[x_1:\varphi_1, \dots , x_{m+1}:\varphi_{m+1} ])$,
\item $S(k[x_1:\varphi_1, \dots , x_m:\varphi_m ]) \cap S(k[x_1:\varphi_1, \dots , x_{m+1}:\varphi_{m+1} ])$.
\end{itemize}

\medskip

Therefore the following set is also r.e.:
\begin{multline*}
\{ (x_1, \varphi_1, \dots , x_{m+1}, \varphi_{m+1}) \} \times (S(k[x_1:\varphi_1, \dots , x_m:\varphi_m ]) \cap S(k[x_1:\varphi_1, \dots , x_{m+1}:\varphi_{m+1} ])) \\ 
 \times S(k[x_1:\varphi_1, \dots , x_{m+1}:\varphi_{m+1} ]) \ .
\end{multline*}

\medskip

Let's use this temporary definition 
\begin{multline*} 
Q'_{m+1,2} = \bigcup_{(x_1, \varphi_1, \dots , x_{m+1}, \varphi_{m+1}) \in R_{m+1}} \{(x_1, \varphi_1, \dots , x_{m+1}, \varphi_{m+1})\} \\ \times (S(k[x_1:\varphi_1, \dots , x_m:\varphi_m ]) \cap S(k[x_1:\varphi_1, \dots , x_{m+1}:\varphi_{m+1} ]))\\ 
\times S(k[x_1:\varphi_1, \dots , x_{m+1}:\varphi_{m+1} ]) ).
\end{multline*}

With this $Q'_{m+1,2}$ is a r.e. subset of $(\Sigma^*)^{2(m + 1)} \times \Sigma^* \times \Sigma^*$.\\

We now define two functions $\delta_{1,m}, \ \delta_{2,m}$ over $(\Sigma^*)^{2(m + 1)} \times \Sigma^* \times \Sigma^*$ as follows. Given $((\psi_1, \varphi_1, \dots , \psi_{m+1}, \varphi_{m+1}), \chi, \varphi) \in (\Sigma^*)^{2(m + 1)} \times \Sigma^* \times \Sigma^*$
\begin{multline*}
\delta_{1,m}((\psi_1, \varphi_1, \dots , \psi_{m+1}, \varphi_{m+1}), \chi, \varphi) =
\gamma[ \psi_1: \varphi_1, \dots , \psi_{m+1}: \varphi_{m+1}, \to \left( \chi, \varphi \right) ] \ .
\end{multline*}
\begin{multline*}
\delta_{2,m}((\psi_1, \varphi_1, \dots , \psi_{m+1}, \varphi_{m+1}), \chi, \varphi) =\\
\gamma[ \psi_1: \varphi_1, \dots , \psi_m: \varphi_m, \to \left( \chi, \forall ( \psi_{m+1} : \varphi_{m+1}, \varphi ) \right) ] \ .
\end{multline*}

\medskip

All of the two functions we have defined are computable functions from $(\Sigma^*)^{2(m + 1)} \times \Sigma^* \times \Sigma^*$ to $\Sigma^*$. If we define a function $\delta_m$ over $(\Sigma^*)^{2(m + 1)} \times \Sigma^* \times \Sigma^*$ as follows: 
\[
\delta_m((\psi_1, \varphi_1, \dots , \psi_{m+1}, \varphi_{m+1}), \chi, \varphi) = 
\left(
\begin{array} {l}
{ \delta_{1,m}((\psi_1, \varphi_1, \dots , \psi_{m+1}, \varphi_{m+1}), \chi, \varphi), } \\
{ \delta_{2,m}((\psi_1, \varphi_1, \dots , \psi_{m+1}, \varphi_{m+1}), \chi, \varphi) } 
\end{array}
\right)
\]

then $\delta_m$ is a computable function from $(\Sigma^*)^{2(m + 1)} \times \Sigma^* \times \Sigma^*$ to $(\Sigma^*)^2$, therefore the set 
\[ D_m = \{ \delta_m((\psi_1, \varphi_1, \dots , \psi_{m+1}, \varphi_{m+1}), \chi, \varphi) | ((\psi_1, \varphi_1, \dots , \psi_{m+1}, \varphi_{m+1}), \chi, \varphi) \in Q'_{m+1,2} \} \]
is a r.e. subset of $(\Sigma^*)^2$.\\

If we now consider the set $\bigcup_{m \geqslant 1} D_m$ then this is a r.e. subset of $(\Sigma^*)^2$ and actually this set is equal to our set $R_{\ref{L:move-inside-universal-quantifier}}$ which so is r.e. itself.
\end{proof}

\medskip

Then let $R_{\ref{L:move-inside-universal-quantifier}} \in \mathcal{R}$.\\

\begin{lemma}\label{L:universal-to-existential-first}
Let $m$ be a positive integer. Let $x_1, \dots , x_{m+1} \in \mathcal{V}$, with $x_i \ne x_j$ for $i 
\ne j$. Let $\varphi_1, \dots , \varphi_{m+1} \in E$ and assume $H[x_1:\varphi_1, \dots , x_{m+1}:
\varphi_{m+1} ]$. 

\medskip

Define $k = k[x_1:\varphi_1, \dots , x_{m+1}:\varphi_{m+1} ]$. Of course $H[x_1:\varphi_1, \dots , 
x_m:\varphi_m ]$ also holds, we define $h = k[x_1:\varphi_1, \dots , x_m:\varphi_m ]$. Let $\chi \in 
S(h)$, $\varphi \in S(k)$, $\psi \in S(h) \cap S(k)$.

\medskip

Under these assumptions we have 
\begin{itemize}
\item $\forall ( x_{m+1} : \varphi_{m+1}, \to (\varphi, \psi) ) \in S(h)$,
\item $\to \left( \chi, \forall ( x_{m+1} : \varphi_{m+1}, \to (\varphi, \psi) ) \right) \in S(h)$,
\item $\exists ( x_{m+1} : \varphi_{m+1}, \varphi ) \in S(h)$,
\item $\to (\chi,  \to( \exists ( x_{m+1} : \varphi_{m+1}, \varphi ), \psi)) \in S(h)$

\item $\gamma[ x_1: \varphi_1, \dots , x_m: \varphi_m, \to \left( \chi, \forall ( x_{m+1} : \varphi_{m+1}, \to (\varphi, \psi) ) \right) ] \in S(\epsilon)$,
\item $\gamma[ x_1: \varphi_1, \dots , x_{m}: \varphi_{m}, \to (\chi,  \to( \exists ( x_{m+1} : \varphi_{m+1}, \varphi ), \psi)) ] \in S(\epsilon)$.
\end{itemize}

\medskip

Moreover if $\#(\gamma[ x_1: \varphi_1, \dots , x_m: \varphi_m, \to \left( \chi, \forall ( x_{m+1} : \varphi_{m+1}, \to (\varphi, \psi) ) \right) ])$ then 
\[ \#( \gamma[ x_1: \varphi_1, \dots , x_{m}: \varphi_{m}, \to (\chi,  \to( \exists ( x_{m+1} : \varphi_{m+1}, \varphi ), \psi)) ] ) \ . \]

\end{lemma}

\begin{proof}

\medskip

Clearly $\to (\varphi, \psi) \in S(k)$ and by lemma~\ref{L:quantifiers_S_h}
\[ \forall ( x_{m+1} : \varphi_{m+1}, \to (\varphi, \psi) ) \in S(h) . \]

Similarly $\exists ( x_{m+1} : \varphi_{m+1}, \varphi ) \in S(h)$ and all the other `preliminary' results hold.

\medskip

We can rewrite 
\[ \#(\gamma[ x_1: \varphi_1, \dots , x_m: \varphi_m, \to \left( \chi, \forall ( x_{m+1} : \varphi_{m+1}, \to (\varphi, \psi) ) \right) ]) \] 
as
\small
\begin{align} 
&\text{for each } \rho \in \Xi(h) \ \#\left( h , \to \left(\chi, \forall \left( x_
{m+1}: \varphi_{m+1}, \to(\varphi, \psi) \right) \right), \rho 
\right) \ , \notag  \\
&\text{for each } \rho \in \Xi(h) \ P_{\to} \left( \#\left(h, \chi , \rho \right)  , \#\left(h, 
\forall \left( x_{m+1}: \varphi_{m+1}, \to(\varphi, \psi) \right) , \rho \right) \right) \ , \notag \\
&\text{for each } \rho \in \Xi(h) \ \#\left(h, \chi , \rho \right) \text{ is false or }  \#\left(h, 
\forall \left( x_{m+1}: \varphi_{m+1}, \to(\varphi, \psi) \right) , \rho \right) \ , \notag \\
&\text{for each } \rho \in \Xi(h) \ \#\left(h, \chi , \rho \right) \text{ is false or for each } \sigma \in \Xi(k): \rho \sqsubseteq \sigma \ \#(k, \to(\varphi, \psi), \sigma ) \ , \notag \\
&\text{for each } \rho \in \Xi(h) \ \#\left(h, \chi , \rho \right) \text{ is false or for each } \sigma \in \Xi(k): \rho \sqsubseteq \sigma \ \#(k, \varphi, \sigma) \text{ is false or } \#(k, \psi, \sigma) \ \notag .
\end{align}
\normalsize

\medskip

We can rewrite 
\[ \#( \gamma[ x_1: \varphi_1, \dots , x_{m}: \varphi_{m}, \to (\chi,  \to( \exists ( x_{m+1} : \varphi_{m+1}, \varphi ), \psi)) ] ) \]

as
\begin{align} 
&\text{for each } \rho \in \Xi(h) \ \#\left( h , \to \left(\chi, \to( \exists ( x_{m+1} : \varphi_{m+1}, \varphi ), \psi) \right), \rho 
\right) \ , \notag  \\
&\text{for each } \rho \in \Xi(h) \ P_{\to} \left( \#\left(h, \chi , \rho \right)  , \#\left(h, 
\to( \exists ( x_{m+1} : \varphi_{m+1}, \varphi ), \psi) , \rho \right) \right) \ , \notag \\
&\text{for each } \rho \in \Xi(h) \ P_{\to} \left( \#\left(h, \chi , \rho \right)  , P_{\to} \left( \#(h, \exists ( x_{m+1} : \varphi_{m+1}, \varphi ), \rho), \#(h, \psi, \rho)  \right) \right) \ , \notag \\
&\text{for each } \rho \in \Xi(h) \ \#\left(h, \chi , \rho \right) \text{ is false or } (  \#(h, \exists ( x_{m+1} : \varphi_{m+1}, \varphi ), \rho) \text{ is false or } \#(h, \psi, \rho) ) \ \notag .
\end{align}

\medskip

We can furtherly express this as

\smallskip

`for each $\rho \in \Xi( h )$ $\#\left(h, \chi , \rho \right)$ is false or\\
((there exists $\sigma \in \Xi(k)$ such that $\rho \sqsubseteq \sigma$ and $\# \left(k, \varphi, \sigma \right)$) is false or $\#(h, \psi, \rho)$)'.\\

We now assume 
\[ \#(\gamma[ x_1: \varphi_1, \dots , x_m: \varphi_m, \to \left( \chi, \forall ( x_{m+1} : \varphi_{m+1}, \to (\varphi, \psi) ) \right) ]) \] and try to prove 
\[ \#( \gamma[ x_1: \varphi_1, \dots , x_{m}: \varphi_{m}, \to (\chi,  \to( \exists ( x_{m+1} : \varphi_{m+1}, \varphi ), \psi)) ] ) \ . \]

\medskip

Let $\rho \in \Xi( h )$ and $\#\left(h, \chi , \rho \right)$, suppose there exists $\sigma \in \Xi(k)$ such that $\rho \sqsubseteq \sigma$ and $\# \left(k, \varphi, \sigma \right)$. Clearly under our assumptions $\#(k, \psi, \sigma)$ holds. We need to prove  $\#(h, \psi, \rho)$, and to do this we can use lemma~\ref{L:kappa-acca-sigma-ro}. In fact there exists a positive integer $n$ such that $h \in K(n)$, $\psi \in E(n, h)$, $k \in K(n)$, $\psi \in E(n, k)$. Given that $\rho \in \Xi( h )$, $\sigma \in \Xi(k)$, $\rho \sqsubseteq \sigma$ we can apply that lemma and get $\#(h, \psi, \rho) = \#(k, \psi, \sigma)$, so $\#(h, \psi, \rho)$ is proved.
\end{proof}

\medskip

We can create a set~$R_{\ref{L:universal-to-existential-first}}$ 
which is the set of all pairs
\[ \left(   
\begin{array} {l} 
{ \gamma[ x_1: \varphi_1, \dots , x_m: \varphi_m, \to \left( \chi, \forall ( x_{m+1} : \varphi_{m+1}, \to (\varphi, \psi) ) \right) ], } \\
{ \gamma[ x_1: \varphi_1, \dots , x_{m}: \varphi_{m}, \to (\chi,  \to( \exists ( x_{m+1} : \varphi_{m+1}, \varphi ), \psi)) ] } 
\end{array}
\right) \]

such that 
\begin{itemize}
\item $m$ is a positive integer, $x_1, \dots , x_{m+1} \in \mathcal{V}$, with $x_i \ne x_j$ for $i 
\ne j$, $\varphi_1, \dots , \varphi_{m+1} \in E$, $H[x_1:\varphi_1, \dots , x_{m+1}:\varphi_{m+1} ]$;
\item if we define $k = k[x_1:\varphi_1, \dots , x_{m+1}:\varphi_{m+1} ]$ and $h = k[x_1:\varphi_1, 
\dots , x_m:\varphi_m ]$ then $\chi \in S(h)$, $\varphi \in S(k)$, $\psi \in S(h) \cap S(k)$.
\end{itemize}
 
 \medskip
 
Lemma~\ref{L:universal-to-existential-first} shows us that this set (which is a potential 1-ary rule) is `sound'. In order to use $R_{\ref{L:universal-to-existential-first}}$ as a rule in our system we also need to show that $R_{\ref{L:universal-to-existential-first}}$ is r.e..\\

\begin{lemma}
$R_{\ref{L:universal-to-existential-first}}$ is r.e..
\end{lemma}

\begin{proof}
Given a positive integer $m$ and $(x_1, \varphi_1, \dots , x_{m+1}, \varphi_{m+1}) \in R_{m+1}$ all of the following sets are r.e.:
\begin{itemize}
\item $S(k[x_1:\varphi_1, \dots , x_m:\varphi_m ])$,
\item $S(k[x_1:\varphi_1, \dots , x_{m+1}:\varphi_{m+1} ])$,
\item $S(k[x_1:\varphi_1, \dots , x_m:\varphi_m ]) \cap S(k[x_1:\varphi_1, \dots , x_{m+1}:\varphi_{m+1} ])$.
\end{itemize}

\medskip

Therefore the following set is also r.e.:
\begin{multline*}
\{ (x_1, \varphi_1, \dots , x_{m+1}, \varphi_{m+1}) \} \times S(k[x_1:\varphi_1, \dots , x_m:\varphi_m ]) \\
\times S(k[x_1:\varphi_1, \dots , x_{m+1}:\varphi_{m+1} ]) \\
\times (S(k[x_1:\varphi_1, \dots , x_m:\varphi_m ]) \cap S(k[x_1:\varphi_1, \dots , x_{m+1}:\varphi_{m+1} ])) \ .
\end{multline*}

\medskip

Let's use this temporary definition 

\begin{multline*} 
Q'_{m+1,3} = \bigcup_{(x_1, \varphi_1, \dots , x_{m+1}, \varphi_{m+1}) \in R_{m+1}} \{ (x_1, \varphi_1, \dots , x_{m+1}, \varphi_{m+1}) \} \times S(k[x_1:\varphi_1, \dots , x_m:\varphi_m ]) \\
\times S(k[x_1:\varphi_1, \dots , x_{m+1}:\varphi_{m+1} ]) \\
\times (S(k[x_1:\varphi_1, \dots , x_m:\varphi_m ]) \cap S(k[x_1:\varphi_1, \dots , x_{m+1}:\varphi_{m+1} ])) \ .
\end{multline*}

With this $Q'_{m+1,3}$ is a r.e. subset of $(\Sigma^*)^{2(m + 1)} \times \Sigma^* \times \Sigma^* \times \Sigma^*$.\\

We now define two functions $\delta_{1,m}, \ \delta_{2,m}$ over $(\Sigma^*)^{2(m + 1)} \times \Sigma^* \times \Sigma^* \times \Sigma^*$ as follows. Given $((\psi_1, \varphi_1, \dots , \psi_{m+1}, \varphi_{m+1}), \chi, \varphi, \psi) \in (\Sigma^*)^{2(m + 1)} \times \Sigma^* \times \Sigma^* \times \Sigma^*$
\begin{multline*}
\delta_{1,m}((\psi_1, \varphi_1, \dots , \psi_{m+1}, \varphi_{m+1}), \chi, \varphi, \psi) =\\
\gamma[ \psi_1: \varphi_1, \dots , \psi_m: \varphi_m, \to \left( \chi, \forall ( \psi_{m+1} : \varphi_{m+1}, \to (\varphi, \psi) ) \right) ] \ .
\end{multline*}
\begin{multline*}
\delta_{2,m}((\psi_1, \varphi_1, \dots , \psi_{m+1}, \varphi_{m+1}), \chi, \varphi, \psi) =\\
\gamma[ \psi_1: \varphi_1, \dots , \psi_m: \varphi_m, \to (\chi,  \to( \exists ( \psi_{m+1} : \varphi_{m+1}, \varphi ), \psi)) ] \ .
\end{multline*}

\medskip

All of the two functions we have defined are computable functions from $(\Sigma^*)^{2(m + 1)} \times \Sigma^* \times \Sigma^* \times \Sigma^*$ to $\Sigma^*$. If we define a function $\delta_m$ over $(\Sigma^*)^{2(m + 1)} \times \Sigma^* \times \Sigma^* \times \Sigma^*$ as follows: 
\[
\delta_m((\psi_1, \varphi_1, \dots , \psi_{m+1}, \varphi_{m+1}), \chi, \varphi, \psi) = 
\left(
\begin{array} {l}
{ \delta_{1,m}((\psi_1, \varphi_1, \dots , \psi_{m+1}, \varphi_{m+1}), \chi, \varphi, \psi), } \\
{ \delta_{2,m}((\psi_1, \varphi_1, \dots , \psi_{m+1}, \varphi_{m+1}), \chi, \varphi, \psi) } 
\end{array}
\right)
\]

then $\delta_m$ is a computable function from $(\Sigma^*)^{2(m + 1)} \times \Sigma^* \times \Sigma^* \times \Sigma^*$ to $(\Sigma^*)^2$, therefore the set 
\[ D_m = \{ \delta_m((\psi_1, \varphi_1, \dots , \psi_{m+1}, \varphi_{m+1}), \chi, \varphi, \psi) | ((\psi_1, \varphi_1, \dots , \psi_{m+1}, \varphi_{m+1}), \chi, \varphi, \psi) \in Q'_{m+1,3} \} \]
is a r.e. subset of $(\Sigma^*)^2$.\\

If we now consider the set $\bigcup_{m \geqslant 1} D_m$ then this is a r.e. subset of $(\Sigma^*)^2$ and actually this set is equal to our set $R_{\ref{L:universal-to-existential-first}}$ which so is r.e. itself.
\end{proof}

\medskip

Then let $R_{\ref{L:universal-to-existential-first}} \in \mathcal{R}$.\\

\begin{lemma}\label{L:universal-to-existential-second}
Let $m$ be a positive integer. Let $x_1, \dots , x_{m+1} \in \mathcal{V}$, with $x_i \ne x_j$ for $i 
\ne j$. Let $\varphi_1, \dots , \varphi_{m+1} \in E$ and assume $H[x_1:\varphi_1, \dots , x_{m+1}:
\varphi_{m+1} ]$. 

\medskip

Define $k = k[x_1:\varphi_1, \dots , x_{m+1}:\varphi_{m+1} ]$. Of course $H[x_1:\varphi_1, \dots , 
x_m:\varphi_m ]$ also holds, we define $h = k[x_1:\varphi_1, \dots , x_m:\varphi_m ]$. Let $\varphi \in S(k)$, $\psi \in S(h) \cap S(k)$.

\medskip

Under these assumptions we have 
\begin{itemize}
\item $\forall ( x_{m+1} : \varphi_{m+1}, \to (\varphi, \psi) ) \in S(h)$,
\item $\exists ( x_{m+1} : \varphi_{m+1}, \varphi ) \in S(h)$,

\item $\gamma[ x_1: \varphi_1, \dots , x_m: \varphi_m, \forall ( x_{m+1} : \varphi_{m+1}, \to (\varphi, \psi) ) ] \in S(\epsilon)$,
\item $\gamma[ x_1: \varphi_1, \dots , x_{m}: \varphi_{m}, \to( \exists ( x_{m+1} : \varphi_{m+1}, \varphi ), \psi) ] \in S(\epsilon)$.
\end{itemize}

\medskip

Moreover if $\#(\gamma[ x_1: \varphi_1, \dots , x_m: \varphi_m, \forall ( x_{m+1} : \varphi_{m+1}, \to (\varphi, \psi) ) ])$ then 
\[ \#( \gamma[ x_1: \varphi_1, \dots , x_{m}: \varphi_{m}, \to( \exists ( x_{m+1} : \varphi_{m+1}, \varphi ), \psi) ] ) \ . \]

\end{lemma}

\begin{proof}
 
\medskip

Clearly $\to (\varphi, \psi) \in S(k)$ and by lemma~\ref{L:quantifiers_S_h}
\[ \forall ( x_{m+1} : \varphi_{m+1}, \to (\varphi, \psi) ) \in S(h) . \]

Similarly $\exists ( x_{m+1} : \varphi_{m+1}, \varphi ) \in S(h)$ and all the other `preliminary' results hold.

\medskip

We can rewrite 
\[ \#(\gamma[ x_1: \varphi_1, \dots , x_m: \varphi_m, \forall ( x_{m+1} : \varphi_{m+1}, \to (\varphi, \psi) )  ]) \] 

as
\begin{align} 
&\text{for each } \rho \in \Xi(h) \ \#\left( h , \forall \left( x_
{m+1}: \varphi_{m+1}, \to(\varphi, \psi) \right) , \rho 
\right) \ , \notag  \\
&\text{for each } \rho \in \Xi(h) \text{ for each } \sigma \in \Xi(k) \text{ such that } \rho \sqsubseteq \sigma \ \#(k, \to( \varphi, \psi), \sigma) \ , \notag \\
&\text{for each } \rho \in \Xi(h) \text{ for each } \sigma \in \Xi(k) \text{ such that } \rho \sqsubseteq \sigma \ P_{\to}( \#(k, \varphi, \sigma), \#(k, \psi, \sigma) )  \ , \notag \\
&\text{for each } \rho \in \Xi(h) \text{ for each } \sigma \in \Xi(k) \text{ such that } \rho \sqsubseteq \sigma \  \#(k, \varphi, \sigma) \text{ is false or } \#(k, \psi, \sigma) \ \notag .
\end{align}

\medskip

We can rewrite 
\[ \#( \gamma[ x_1: \varphi_1, \dots , x_{m}: \varphi_{m}, \to( \exists ( x_{m+1} : \varphi_{m+1}, \varphi ), \psi) ] )  \]

as
\begin{align} 
&\text{for each } \rho \in \Xi(h) \ \#\left( h , \to( \exists ( x_{m+1} : \varphi_{m+1}, \varphi ), \psi), \rho 
\right) \ , \notag  \\
&\text{for each } \rho \in \Xi(h) \ P_{\to} \left( \#(h, \exists ( x_{m+1} : \varphi_{m+1}, \varphi ), \rho), \#(h, \psi, \rho)  \right) \ , \notag \\
&\text{for each } \rho \in \Xi(h) \  \#(h, \exists ( x_{m+1} : \varphi_{m+1}, \varphi ), \rho) \text{ is false or } \#(h, \psi, \rho)  \ \notag .
\end{align}

\medskip

We can furtherly express this as

\smallskip

`for each $\rho \in \Xi( h )$\\
((there exists $\sigma \in \Xi(k)$ such that $\rho \sqsubseteq \sigma$ and $\# \left(k, \varphi, \sigma \right)$) is false or $\#(h, \psi, \rho)$)'.\\

We now assume 
\[ \#(\gamma[ x_1: \varphi_1, \dots , x_m: \varphi_m, \forall ( x_{m+1} : \varphi_{m+1}, \to (\varphi, \psi) ) ]) \] and try to prove 
\[ \#( \gamma[ x_1: \varphi_1, \dots , x_{m}: \varphi_{m}, \to( \exists ( x_{m+1} : \varphi_{m+1}, \varphi ), \psi) ] ) \ . \]

\medskip

Let $\rho \in \Xi( h )$, suppose there exists $\sigma \in \Xi(k)$ such that $\rho \sqsubseteq \sigma$ and $\# \left(k, \varphi, \sigma \right)$. Clearly under our assumptions $\#(k, \psi, \sigma)$ holds. We need to prove  $\#(h, \psi, \rho)$, and to do this we can use lemma~\ref{L:kappa-acca-sigma-ro}. In fact there exists a positive integer $n$ such that $h \in K(n)$, $\psi \in E(n, h)$, $k \in K(n)$, $\psi \in E(n, k)$. Given that $\rho \in \Xi( h )$, $\sigma \in \Xi(k)$, $\rho \sqsubseteq \sigma$ we can apply that lemma and get $\#(h, \psi, \rho) = \#(k, \psi, \sigma)$, so $\#(h, \psi, \rho)$ is proved.
 \end{proof}

\medskip

We can create a set~$R_{\ref{L:universal-to-existential-second}}$ 
which is the set of all pairs
\[ \left(   
\begin{array} {l} 
{ \gamma[ x_1: \varphi_1, \dots , x_m: \varphi_m, \forall ( x_{m+1} : \varphi_{m+1}, \to (\varphi, \psi) ) ], } \\
{ \gamma[ x_1: \varphi_1, \dots , x_{m}: \varphi_{m}, \to( \exists ( x_{m+1} : \varphi_{m+1}, \varphi ), \psi) ] } 
\end{array}
\right) \]

such that 
\begin{itemize}
\item $m$ is a positive integer, $x_1, \dots , x_{m+1} \in \mathcal{V}$, with $x_i \ne x_j$ for $i 
\ne j$, $\varphi_1, \dots , \varphi_{m+1} \in E$, $H[x_1:\varphi_1, \dots , x_{m+1}:\varphi_{m+1} ]$;
\item if we define $k = k[x_1:\varphi_1, \dots , x_{m+1}:\varphi_{m+1} ]$ and $h = k[x_1:\varphi_1, 
\dots , x_m:\varphi_m ]$ then $\varphi \in S(k)$, $\psi \in S(h) \cap S(k)$.
\end{itemize}

\medskip

Lemma~\ref{L:universal-to-existential-second} shows us that this set (which is a potential 1-ary rule) is `sound'. In order to use $R_{\ref{L:universal-to-existential-second}}$ as a rule in our system we also need to show that $R_{\ref{L:universal-to-existential-second}}$ is r.e..\\

\begin{lemma}
$R_{\ref{L:universal-to-existential-second}}$ is r.e..
\end{lemma}

\begin{proof}
Given a positive integer $m$ and $(x_1, \varphi_1, \dots , x_{m+1}, \varphi_{m+1}) \in R_{m+1}$ all of the following sets are r.e.:
\begin{itemize}
\item $S(k[x_1:\varphi_1, \dots , x_m:\varphi_m ])$,
\item $S(k[x_1:\varphi_1, \dots , x_{m+1}:\varphi_{m+1} ])$,
\item $S(k[x_1:\varphi_1, \dots , x_m:\varphi_m ]) \cap S(k[x_1:\varphi_1, \dots , x_{m+1}:\varphi_{m+1} ])$.
\end{itemize}

\medskip

Therefore the following set is also r.e.:
\begin{multline*}
\{ (x_1, \varphi_1, \dots , x_{m+1}, \varphi_{m+1}) \} \times S(k[x_1:\varphi_1, \dots , x_{m+1}:\varphi_{m+1} ]) \\
\times (S(k[x_1:\varphi_1, \dots , x_m:\varphi_m ]) \cap S(k[x_1:\varphi_1, \dots , x_{m+1}:\varphi_{m+1} ])) \ .
\end{multline*}

\medskip

Let's use this temporary definition 

\begin{multline*} 
Q'_{m+1,2} = \bigcup_{(x_1, \varphi_1, \dots , x_{m+1}, \varphi_{m+1}) \in R_{m+1}} \{ (x_1, \varphi_1, \dots , x_{m+1}, \varphi_{m+1}) \} \times S(k[x_1:\varphi_1, \dots , x_{m+1}:\varphi_{m+1} ]) \\
\times (S(k[x_1:\varphi_1, \dots , x_m:\varphi_m ]) \cap S(k[x_1:\varphi_1, \dots , x_{m+1}:\varphi_{m+1} ])) \ .
\end{multline*}

With this $Q'_{m+1,2}$ is a r.e. subset of $(\Sigma^*)^{2(m + 1)} \times \Sigma^* \times \Sigma^*$.\\

We now define two functions $\delta_{1,m}, \ \delta_{2,m}$ over $(\Sigma^*)^{2(m + 1)} \times \Sigma^* \times \Sigma^*$ as follows. Given $((\psi_1, \varphi_1, \dots , \psi_{m+1}, \varphi_{m+1}), \varphi, \psi) \in (\Sigma^*)^{2(m + 1)} \times \Sigma^* \times \Sigma^*$
\begin{multline*}
\delta_{1,m}((\psi_1, \varphi_1, \dots , \psi_{m+1}, \varphi_{m+1}), \varphi, \psi) =\\
\gamma[ \psi_1: \varphi_1, \dots , \psi_m: \varphi_m, \forall ( \psi_{m+1} : \varphi_{m+1}, \to (\varphi, \psi) ) ] \ .
\end{multline*}
\begin{multline*}
\delta_{2,m}((\psi_1, \varphi_1, \dots , \psi_{m+1}, \varphi_{m+1}), \varphi, \psi) =\\
\gamma[ \psi_1: \varphi_1, \dots , \psi_m: \varphi_m, \to( \exists ( \psi_{m+1} : \varphi_{m+1}, \varphi ), \psi) ] \ .
\end{multline*}

\medskip

All of the two functions we have defined are computable functions from $(\Sigma^*)^{2(m + 1)} \times \Sigma^* \times \Sigma^*$ to $\Sigma^*$. If we define a function $\delta_m$ over $(\Sigma^*)^{2(m + 1)} \times \Sigma^* \times \Sigma^*$ as follows: 
\[
\delta_m((\psi_1, \varphi_1, \dots , \psi_{m+1}, \varphi_{m+1}), \varphi, \psi) = 
\left(
\begin{array} {l}
{ \delta_{1,m}((\psi_1, \varphi_1, \dots , \psi_{m+1}, \varphi_{m+1}), \varphi, \psi), } \\
{ \delta_{2,m}((\psi_1, \varphi_1, \dots , \psi_{m+1}, \varphi_{m+1}), \varphi, \psi) } 
\end{array}
\right)
\]

then $\delta_m$ is a computable function from $(\Sigma^*)^{2(m + 1)} \times \Sigma^* \times \Sigma^*$ to $(\Sigma^*)^2$, therefore the set 
\[ D_m = \{ \delta_m((\psi_1, \varphi_1, \dots , \psi_{m+1}, \varphi_{m+1}), \varphi, \psi) | ((\psi_1, \varphi_1, \dots , \psi_{m+1}, \varphi_{m+1}), \varphi, \psi) \in Q'_{m+1,2} \} \]
is a r.e. subset of $(\Sigma^*)^2$.\\

If we now consider the set $\bigcup_{m \geqslant 1} D_m$ then this is a r.e. subset of $(\Sigma^*)^2$ and actually this set is equal to our set $R_{\ref{L:universal-to-existential-second}}$ which so is r.e. itself.
\end{proof}

\medskip

Then let $R_{\ref{L:universal-to-existential-second}} \in \mathcal{R}$.\\

\begin{lemma}\label{L:rule-il-to-cl-new}
Let $m$ be a positive integer. Let $x_1, \dots , x_m \in \mathcal{V}$, with $x_i \ne x_j$ for $i \ne 
j$. Let $\varphi_1, \dots , \varphi_m \in E$ and assume $H[x_1:\varphi_1, \dots , x_m:\varphi_m ]$. 
Define $k = k[x_1:\varphi_1, \dots , x_m:\varphi_m ]$ and let $\varphi, \psi, \chi \in S(k)$.

\medskip

Under these assumptions we have 
\begin{itemize}
\item $\to( \wedge(\varphi, \psi) , \chi), \to(\varphi, \to(\psi, \chi) ) \in S(k)$,
\item $\gamma[ x_1: \varphi_1, \dots , x_m: \varphi_m, \to(\varphi, \to(\psi, \chi) ) ] \in S
(\epsilon)$,
\item $\gamma[ x_1: \varphi_1, \dots , x_m: \varphi_m, \to( \wedge(\varphi, \psi) , \chi) ] \in 
S(\epsilon)$.
\end{itemize}

\smallskip

Moreover if $\#( \gamma[ x_1: \varphi_1, \dots , x_m: \varphi_m, \to(\varphi, \to(\psi, \chi) ) ] )$ then\\
$\#( \gamma[ x_1: \varphi_1, \dots , x_m: \varphi_m, \to( \wedge(\varphi, \psi) , 
\chi) ] )$

\end{lemma}

\begin{proof}

\medskip

We assume $\#( \gamma[ x_1: \varphi_1, \dots , x_m: \varphi_m, \to(\varphi, 
\to(\psi, \chi) ) ] )$ which can be rewritten
\begin{align} 
&\text{for each } \sigma \in \Xi(k) \ \#( k , \to(\varphi, \to(\psi, \chi) ), \sigma) \ , \notag  \\
&\text{for each } \sigma \in \Xi(k) \ P_{\to}( \#(k, \varphi , \sigma), \#(k, \to(\psi, \chi) , \sigma)  ) \ , \notag \\
&\text{for each } \sigma \in \Xi(k) \ \#(k, \varphi , \sigma) \text{ is false or (} \#(k, \psi , \sigma) \text{ is false or } \#(k, \chi, \sigma) \text{)}  \notag .
\end{align}

\medskip

We now try to show $\#( \gamma[ x_1: \varphi_1, \dots , x_m: \varphi_m, \to( \wedge(\varphi, \psi) , \chi) 
] )$ which in turn can be rewritten
\begin{align} 
&\text{for each } \sigma \in \Xi(k) \ \#( k , \to( \wedge(\varphi, \psi) , \chi), \sigma) \ , \notag  \\
&\text{for each } \sigma \in \Xi(k) \ P_{\to}( \#(k, \wedge(\varphi, \psi), \sigma) , \#(k, \chi , \sigma) ) \ , \notag \\
&\text{for each } \sigma \in \Xi(k) \ P_{\to}( P_{\wedge}( \#(k, \varphi , \sigma), \#(k, \psi , \sigma) ) , \#(k, \chi 
, \sigma) ) \ , \notag \\
&\text{for each } \sigma \in \Xi(k) \ \text{ it is false that (} \#(k, \varphi , \sigma) \text{ and } \#(k, \psi, \sigma) \text{) or } \#(k, \chi, \sigma) \notag .
\end{align}

\medskip

Let $\sigma \in \Xi(k)$, let's also keep in mind that $\#(k, \varphi , \sigma)$ is false or $\#(k, \psi, \sigma)$ is false or $\#(k, \chi, \sigma)$. If $\#(k, \varphi , \sigma)$ is false then it is false that ($\#(k, \varphi , \sigma)$ and $\#(k, \psi, \sigma)$). Similarly if $\#(k, \psi, \sigma)$ is false then it is false that ($\#(k, \varphi , \sigma)$ and $\#(k, \psi, \sigma)$). Finally if $\#(k, \chi, \sigma)$ holds then it holds itself and what we wanted to show is true.
\end{proof}

\medskip

We can create a set $R_{\ref{L:rule-il-to-cl-new}}$ which is the set of all pairs
\[ \left( 
\gamma[ x_1: \varphi_1, \dots , x_m: \varphi_m, \to(\varphi, \to(\psi, \chi) ) ] ,
\gamma[ x_1: \varphi_1, \dots , x_m: \varphi_m, \to( \wedge(\varphi, \psi) , \chi) ]
\right)  \]

such that
\begin{itemize}
\item $m$ is a positive integer, $x_1, \dots , x_m \in \mathcal{V}$, $x_i \ne x_j$ for $i \ne j$, $
\varphi_1, \dots , \varphi_m \in E$, $H[x_1:\varphi_1, \dots , x_m:\varphi_m ]$,
\item $\varphi, \psi, \chi \in S(k[x_1:\varphi_1, \dots , x_m:\varphi_m ])$.
\end{itemize}

\medskip

Lemma~\ref{L:rule-il-to-cl-new} shows us that this set (which is a potential 1-ary rule) is `sound'. In order to use $R_{\ref{L:rule-il-to-cl-new}}$ as a rule in our system we also need to show that $R_{\ref{L:rule-il-to-cl-new}}$ is r.e..\\

\begin{lemma}
$R_{\ref{L:rule-il-to-cl-new}}$ is r.e..
\end{lemma}

\begin{proof}
Given a positive integer $m$ and $(x_1, \varphi_1, \dots , x_m, \varphi_m) \in R_m$ we can notice the following:

\begin{itemize}
\item $k[x_1:\varphi_1, \dots , x_m:\varphi_m ] \in K$;
\item $S(k[x_1:\varphi_1, \dots , x_m:\varphi_m ])$ is r.e.;
\item $\{(x_1, \varphi_1, \dots , x_m, \varphi_m)\} \times S(k[x_1:\varphi_1, \dots , x_m:\varphi_m ])^3$ is r.e..
\end{itemize}

Let's define
\[ Q_{m,3} = \bigcup_{(x_1, \varphi_1, \dots , x_m, \varphi_m) \in R_m} \{(x_1, \varphi_1, \dots , x_m, \varphi_m)\} \times S(k[x_1:\varphi_1, \dots , x_m:\varphi_m ])^3 \ . \]

\smallskip

Clearly $Q_{m,3} \subseteq (\Sigma^*)^{2m} \times (\Sigma^*)^3$ is also r.e..\\

We now define two functions $\delta_{1,m}, \ \delta_{2,m}$ over $(\Sigma^*)^{2m} \times (\Sigma^*)^3$ as follows. Given $((\psi_1, \varphi_1, \dots , \psi_m, \varphi_m), (\varphi, \psi, \chi)) \in (\Sigma^*)^{2m} \times (\Sigma^*)^3$
\begin{align} 
&\delta_{1,m}((\psi_1, \varphi_1, \dots , \psi_{m}, \varphi_{m}), (\varphi, \psi, \chi)) = \gamma[ \psi_1: \varphi_1, \dots , \psi_m: \varphi_m, \to(\varphi, \to(\psi, \chi) ) ]  \ , \notag  \\
&\delta_{2,m}((\psi_1, \varphi_1, \dots , \psi_{m}, \varphi_{m}), (\varphi, \psi, \chi)) = \gamma[ \psi_1: \varphi_1, \dots , \psi_m: \varphi_m, \to( \wedge(\varphi, \psi) , \chi) ]  \  \notag .
\end{align}

\smallskip

All of the two functions we have defined are computable functions from $(\Sigma^*)^{2m} \times (\Sigma^*)^3$ to $\Sigma^*$. If we define a function $\delta_m$ over $(\Sigma^*)^{2m} \times (\Sigma^*)^3$ as follows: 
\[
\delta_m((\psi_1, \varphi_1, \dots , \psi_{m}, \varphi_{m}), (\varphi, \psi, \chi)) = 
\left(
\begin{array} {l}
{ \delta_{1,m}((\psi_1, \varphi_1, \dots , \psi_{m}, \varphi_{m}), (\varphi, \psi, \chi)) , } \\
{ \delta_{2,m}((\psi_1, \varphi_1, \dots , \psi_{m}, \varphi_{m}), (\varphi, \psi, \chi)) , }
\end{array}
\right)
\]

then $\delta_m$ is a computable function from $(\Sigma^*)^{2m} \times (\Sigma^*)^3$ to $(\Sigma^*)^2$, therefore the set 
\[ D_m = \{ \delta_m((\psi_1, \varphi_1, \dots , \psi_{m}, \varphi_{m}), (\varphi, \psi, \chi)) | \ ((\psi_1, \varphi_1, \dots , \psi_{m}, \varphi_{m}), (\varphi, \psi, \chi)) \in Q_{m,3} \} \]
is a r.e. subset of $(\Sigma^*)^2$.\\

If we now consider the set $\bigcup_{m \geqslant 1} D_m$ then this is a r.e. subset of $(\Sigma^*)^2$ and actually this set is equal to our set $R_{\ref{L:rule-il-to-cl-new}}$ which so is r.e. itself.
\end{proof}

\medskip

Then let $R_{\ref{L:rule-il-to-cl-new}} \in \mathcal{R}$.

\section{Another proof}\label{Ch:anotherProof}

For each $x, y$ natural numbers we say that $x$ divides $y$ if there exists a natural number $\alpha$ such that $y = x \alpha$.\\

In our example we want to show that for each $x, y, z$ natural numbers if $x$ divides $y$ and $y$ divides $z$ then $x$ divides $z$.\\

Of course, we first need to build an expression in our language to express this. To build that expression we must add to our language a constant symbol $N$ to represent the set of natural numbers $\mathbb{N}$, so that we have $\#(N) = \mathbb{N}$.\\

And we need to add another constant symbol in our language. This is the symbol~$*$ that stands for the product (or multiplication) operation in the domain $\mathbb{N}$ of natural numbers. Therefore $\#(*)$ is a function defined on $\mathbb{N} \times \mathbb{N}$ and for each $\alpha, \beta \in \mathbb{N}$ $\#(*)(\alpha, \beta)$ is the product of $\alpha$ and $\beta$, in other words $\#(*)(\alpha, \beta) = \alpha \cdot \beta$.\\ 

The set $\mathcal{F}$ of operators is the same we have assumed in our former example, so it must contain all of these symbols: $\neg, \wedge, \vee, \to, \leftrightarrow, \in, =$.\\

So, in order to formalize our statement and a proof of it, we will use a language $(\mathcal{V}, \mathcal{F}, \mathcal{C}, \#, \{D_1, \dots, D_p \})$ which must be as follows
\begin{align} 
&\mathcal{V} = \{x, y, z, u, v, w \} \ , \notag  \\
&\mathcal{F} = \{ \neg, \wedge, \vee, \to, \leftrightarrow, \in, =  \} \ , \notag  \\
&\mathcal{C} = \{ N, *  \} \  \notag .
\end{align}

\smallskip

Moreover, we need to include the set $\mathbb{N}$ of natural numbers in our additional sets, so let $p = 1$ and $D_1 = \mathbb{N}$.\\

Here we notice that we set as a constraint that for each $c \in \mathcal{C}$ and for any positive integer $q$ we must be able to decide all of the following conditions

\begin{itemize}
\item $Set_q( \#(c))$;
\item $Event_q( \#(c))$;
\item $\#(c) \in D_i$;
\item $\#(c) \in \mathcal{P}^q(D_i)$;
\item if ($Set_q( \#(c))$) then ($NotEmpty_q( \#(c))$).
\end{itemize}

And moreover, the last condition must be decided as true.\\

For the constant $N$ we can take the following decisions: 

\begin{itemize}
\item $Set_1( \#(N))$: true;
\item for $q > 1$, $Set_q( \#(N))$: false;
\item $Event_q( \#(N))$: false;
\item $\#(N) \in \mathbb{N}$: false;
\item $\#(N) \in \mathcal{P}^1(\mathbb{N})$: true;
\item for $q > 1$ $\#(N) \in \mathcal{P}^q(\mathbb{N})$: false;
\item if ($Set_q( \#(N))$) then ($NotEmpty_q( \#(N))$): true.
\end{itemize}

\medskip

For the constant $*$ we can take the following decisions: 

\begin{itemize}
\item $Set_1( \#(*))$: true;
\item for $q > 1$, $Set_q( \#(*))$: false;
\item $Event_q( \#(*))$: false;
\item $\#(*) \in \mathbb{N}$: false;
\item $\#(*) \in \mathcal{P}^q(\mathbb{N})$: false;
\item if ($Set_q( \#(*))$) then ($NotEmpty_q( \#(*))$): true.
\end{itemize}

\medskip

At this point, the statement we wish to prove is the following:

\scriptsize
\begin{equation}\label{E:th-second-ex}
 \gamma \left[x:N, y:N , z:N, \to \left( 
\wedge 
\left(
\begin{array} {l}
{  \exists( u:N, =(y, *(x, u) ) ), } \\
{  \exists( v:N, =(z, *(y, v) ) )}
\end{array}
\right),
\exists( w:N, =(z, *(x, w) ) )
\right) \right] \tag{$Th_1$} \, .
\end{equation}
\normalsize

\medskip

Let $k = k[x:N, y:N, z:N, u:N, v:N]$. By lemma~\ref{L:xi_in_Ekj} $u \in E(k)$. If we define $k_u = k[x:N, y:N, z:N]$ then for each $\sigma \in \Xi(k)$ $\sigma_{/dom(k_u)} \in \Xi(k_u)$ and $\#(k,u,\sigma) \in \#(k_u, N, \sigma_{/dom(k_u)}) = \#(N) = \mathbb{N}$.\\

Similarly by ~\ref{L:xi_in_Ekj} $x \in E(k)$. If we define $k_x = \epsilon$ then for each $\sigma \in \Xi(k)$ $\sigma_{/dom(k_x)} \in \Xi(k_x)$ and $\#(k,x,\sigma) \in \#(k_x, N, \sigma_{/dom(k_x)}) = \#(N) = \mathbb{N}$.\\

By lemma~\ref{L:preliminary-results-ch6-3} it follows that $(*)(x, u) \in E(k)$ and for each $\sigma \in \Xi(k)$ $\#(k, (*)(x, u), \sigma) = (\#(k,x, \sigma) \cdot \#(k,u, \sigma)) \in \mathbb{N}$.\\

The first sentence in our proof is an instance of axiom~$A_{\ref{L:axiom-substitute-exp-in-mult}}$.

\begin{equation}\label{E:ex3-a}
\gamma[ x:N, y:N, z:N, u:N, v:N, \to \left( \wedge(=(y,xu), =(z,yv)), =(z,(xu)v)  \right) ]
\end{equation}

\medskip

The following also hold: 
\begin{itemize}
\item $\wedge(=(y,xu), =(z,yv)) \in S(k)$.
\end{itemize}

\smallskip

By~$A_{\ref{L:axiom-commutative-mult}}$ we obtain
\begin{equation}\label{E:ex3-b}
\gamma[ x:N, y:N, z:N, u:N, v:N, \to \left( \wedge(=(y,xu), =(z,yv)), =((xu)v,x(uv)) \right) ]
\end{equation}

\medskip

The following also hold: 
\begin{itemize}
\item $=(z,(xu)v) \in S(k)$,
\item $ =((xu)v,x(uv)) \in S(k)$.
\end{itemize}

\smallskip

By \ref{E:ex3-a}, \ref{E:ex3-b} and rule~$R_{\ref{L:AND-I-new}}$
\small
\begin{equation}\label{E:ex3-c}
 \gamma \left[x:N, y:N, z:N, u:N, v:N, \to \left( 
\wedge \left(  
\begin{array} {l}
{  =(y,xu), } \\
{  =(z,yv)  }
\end{array}
\right) , 
\wedge \left(
\begin{array} {l}
{ =(z,(xu)v), } \\
{  =((xu)v,x(uv)) }
\end{array}
\right)
\right) \right] .
\end{equation}
\normalsize

\medskip

The following also hold: 
\begin{itemize}
\item $z \in E(k)$,
\item $(xu)v \in E(k)$,
\item $x(uv) \in E(k)$.
\end{itemize}

\smallskip

By axiom~$A_{\ref{L:ax-eq-trans}}$ 

\small
\begin{equation}\label{E:ex3-d}
 \gamma \left[x:N, \dots , v:N, \to \left( 
\wedge \left(  
\begin{array} {l}
{  =(y,xu), } \\
{  =(z,yv)  }
\end{array}
\right) , \to \left(
\wedge \left(
\begin{array} {l}
{ =(z,(xu)v), } \\
{  =((xu)v,x(uv)) }
\end{array}
\right), =(z,x(uv)) \right)
\right) \right] .
\end{equation}
\normalsize

\medskip

The following also hold: 
\begin{itemize}
\item $=(z, x(uv)) \in S(k)$.
\end{itemize}

\smallskip

By \ref{E:ex3-c}, \ref{E:ex3-d} and rule~$R_{\ref{L:rule-mp}}$

\small
\begin{equation}\label{E:ex3-e}
 \gamma \left[x:N, y:N, z:N, u:N, v:N, \to \left( 
\wedge \left(  
\begin{array} {l}
{  =(y,xu), } \\
{  =(z,yv)  }
\end{array}
\right) , 
=(z,x(uv))
\right) \right] .
\end{equation}
\normalsize

\medskip

The following also hold: $(*)(u, v) \in E(k)$ and for each $\sigma \in \Xi(k)$\\
$\#(k, (*)(u, v), \sigma) = (\#(k,u, \sigma) \cdot \#(k,v, \sigma)) \in \mathbb{N}$ (cfr lemma~\ref{L:axiom-commutative-mult}). 

\medskip

By \ref{E:ex3-e} and rule~$R_{\ref{L:rule-introd-exists-mult}}$

\small
\begin{equation}\label{E:ex3-h}
 \gamma \left[x:N, y:N, z:N, u:N, v:N, \to \left( 
\wedge \left(  
\begin{array} {l}
{  =(y,xu), } \\
{  =(z,yv)  }
\end{array}
\right) , 
\exists( w:N, =(z, xw) )
\right)
\right] .
\end{equation}
\normalsize

\medskip

The following also holds: $\exists( w:N, =(z, xw) ) \in S(k)$ (cfr lemma~\ref{L:rule-introd-exists-mult}).

\medskip

By \ref{E:ex3-h} and rule~$R_{\ref{L:rule-cl-to-il-new}}$

\small
\begin{equation}\label{E:ex3-i}
 \gamma \left[x:N, y:N, z:N, u:N, v:N, \to \left( 
 =(y,xu) , 
\to \left(
\begin{array} {l}
{  =(z,yv), } \\
{ \exists( w:N, =(z, xw) )   }
\end{array}
\right)
\right)
\right] .
\end{equation}
\normalsize

\medskip

Let $h = k[x:N, y:N, z:N, u:N]$. By lemma~\ref{L:xi_in_Ekj} $u \in E(h)$. If we define $h_u = k[x:N, y:N, z:N]$ then for each $\rho \in \Xi(h)$ $\rho_{/dom(h_u)} \in \Xi(h_u)$ and $\#(h,u,\rho) \in \#(h_u, N, \rho_{/dom(h_u)}) = \#(N) = \mathbb{N}$.\\

Similarly by ~\ref{L:xi_in_Ekj} $x \in E(h)$. If we define $h_x = \epsilon$ then for each $\rho \in \Xi(h)$ $\rho_{/dom(h_x)} \in \Xi(h_x)$ and $\#(h,x,\rho) \in \#(h_x, N, \rho_{/dom(h_x)}) = \#(N) = \mathbb{N}$.\\

By lemma~\ref{L:preliminary-results-ch6-3} it follows that $(*)(x, u) \in E(h)$ and for each $\rho \in \Xi(h)$ $\#(h, (*)(x, u), \rho) = (\#(h,x, \rho) \cdot \#(h,u, \rho)) \in \mathbb{N}$.\\

Still by ~\ref{L:xi_in_Ekj} $y \in E(h)$ and by lemma~\ref{L:preliminary-results-ch6-2} $=(y, xu) \in S(h)$.\\

By \ref{E:ex3-i} and rule~$R_{\ref{L:move-inside-universal-quantifier}}$

\scriptsize
\begin{equation}\label{E:ex3-l}
 \gamma \left[x:N, y:N, z:N, u:N, \to \left( 
 =(y,xu) , 
\forall \left( v:N, 
\to \left(
\begin{array} {l}
{  =(z,yv), } \\
{ \exists( w:N, =(z, xw) )   }
\end{array}
\right)
\right)
\right)
\right] .
\end{equation}
\normalsize

\medskip

We now want to prove that $\exists( w:N, =(z, xw) ) \in S(h)$. We start by defining $g = h + <w,N>$.

\medskip

We have $N \in E(h)$ and for each $\rho \in \Xi(h)$ $\#(h,N,\rho) = \#(N) = \mathbb{N}$. So $N \in E_s(h)$. Moreover $w \in (\mathcal{V} - var(h))$ so by lemma~\ref{L:quantifiers_S_h} $g \in K$. 

\medskip

We now want to show that $=(z, xw)$ belongs to $S(g)$. Since $N \in E_s(h)$ we have $H[x:N, y:N, z:N, u:N, w:N]$. We have \[ k[x:N, y:N, z:N, u:N, w:N] = h + (w,N) = g \ . \]

\medskip

Using lemma~\ref{L:xi_in_Ekj} we obtain that $z, x, w \in E(g)$. If we define $g_x = \epsilon$ then for each $\sigma \in \Xi(g)$ $\sigma_{/dom(g_x)} \in \Xi(g_x)$ and $\#(g,x,\sigma) \in \#(g_x, N, \sigma_{/dom(g_x)}) = \#(N) = \mathbb{N}$.

\medskip

Moreover for each $\sigma \in \Xi(g)$ $\sigma_{/dom(h)} \in \Xi(h)$ $\#(g,w,\sigma) \in \#(h, N, \sigma_{/dom(h)}) = \#(N) = \mathbb{N}$.

\medskip

By lemma~\ref{L:preliminary-results-ch6-3} it follows that $(*)(x, w) \in E(g)$ and for each $\sigma \in \Xi(g)$ $\#(g, (*)(x, w), \sigma) = (\#(g,x, \sigma) \cdot \#(g,w, \sigma)) \in \mathbb{N}$.

\medskip

By lemma~\ref{L:preliminary-results-ch6-2} $=(z, xw)$ belongs to $S(g)$. We can now apply lemma~\ref{L:quantifiers_S_h} and obtain that $\exists( w:N, =(z, xw) ) \in S(h)$.

\medskip

To sum up we have $=(y,xu) \in S(h)$, $=(z,yv) \in S(k)$,\\
$\exists( w:N, =(z, xw) ) \in S(h) \cap S(k)$. 

\medskip

By \ref{E:ex3-l} and rule~$R_{\ref{L:universal-to-existential-first}}$

\small
\begin{equation}\label{E:ex3-m}
 \gamma \left[x:N, y:N, z:N, u:N, \to \left( 
 =(y,xu) , 
\to \left(
\begin{array} {l}
{ \exists( v:N, =(z,yv) ), } \\
{ \exists( w:N, =(z, xw) )   }
\end{array}
\right)
\right)
\right] .
\end{equation}
\normalsize

Using lemma \ref{L:gamma-expand-last}, we can rewrite \ref{E:ex3-m} as 

\scriptsize
\begin{equation}\label{E:ex3-o}
 \gamma \left[x:N, y:N, z:N, \forall \left( u:N, \to \left( 
 =(y,xu) , 
\to \left(
\begin{array} {l}
{ \exists( v:N, =(z,yv) ), } \\
{ \exists( w:N, =(z, xw) )   }
\end{array}
\right)
\right)
\right)
\right] .
\end{equation}
\normalsize

Let $\kappa = k[x:N, y:N, z:N]$. We have proved that\\
$\exists( w:N, =(z, xw) ) \in S(h)$ and $\exists( v:N, =(z,yv) ) \in S(h)$.

\medskip 

We also need to prove that $\exists( w:N, =(z, xw) ) \in S(\kappa)$\\ and $\exists( v:N, =(z,yv) ) \in S(\kappa)$.

\medskip

In order to prove $\exists( w:N, =(z, xw) ) \in S(\kappa)$ we redefine $g$ as $\kappa + <w,N>$.

\medskip 

We have $N \in E(\kappa)$ and for each $\rho \in \Xi(\kappa)$ $\#(\kappa,N,\rho) = \#(N) = \mathbb{N}$. So $N \in E_s(\kappa)$. Moreover $w \in (\mathcal{V} - var(\kappa))$ so by lemma~\ref{L:quantifiers_S_h} $g \in K$. 

\medskip

We now want to show that $=(z, xw)$ belongs to $S(g)$. It follows from lemma~\ref{L:xi-in-D-new} that $H[x:N, y:N, z:N, w:N]$. We have 
\[ k[x:N, y:N, z:N, w:N] = \kappa + <w,N> = g \ . \]

\medskip

Using lemma~\ref{L:xi_in_Ekj} we obtain that $z, x, w \in E(g)$. If we define $g_x = \epsilon$ then for each $\sigma \in \Xi(g)$ $\sigma_{/dom(g_x)} \in \Xi(g_x)$ and $\#(g,x,\sigma) \in \#(g_x, N, \sigma_{/dom(g_x)}) = \#(N) = \mathbb{N}$.

\medskip

Moreover for each $\sigma \in \Xi(g)$ $\sigma_{/dom(\kappa)} \in \Xi(\kappa)$ $\#(g,w,\sigma) \in \#(\kappa, N, \sigma_{/dom(\kappa)}) = \#(N) = \mathbb{N}$.

\medskip

By lemma~\ref{L:preliminary-results-ch6-3} it follows that $(*)(x, w) \in E(g)$ and for each $\sigma \in \Xi(g)$ $\#(g, (*)(x, w), \sigma) = (\#(g,x, \sigma) \cdot \#(g,w, \sigma)) \in \mathbb{N}$.

\medskip

By lemma~\ref{L:preliminary-results-ch6-2} $=(z, xw)$ belongs to $S(g)$. We can now apply lemma~\ref{L:quantifiers_S_h} and obtain that $\exists( w:N, =(z, xw) ) \in S(\kappa)$.

\medskip

In order to prove $\exists( v:N, =(z,yv) ) \in S(\kappa)$ we redefine $g$ as $\kappa + <v,N>$.

\medskip 

We have $N \in E(\kappa)$ and for each $\rho \in \Xi(\kappa)$ $\#(\kappa,N,\rho) = \#(N) = \mathbb{N}$. So $N \in E_s(\kappa)$. Moreover $v \in (\mathcal{V} - var(\kappa))$ so by lemma~\ref{L:quantifiers_S_h} $g \in K$. 

\medskip

We now want to show that $=(z, yv)$ belongs to $S(g)$. It follows from lemma~\ref{L:xi-in-D-new} that $H[x:N, y:N, z:N, v:N]$. We have \[ k[x:N, y:N, z:N, v:N] = \kappa + <v,N> = g \ . \]

\medskip

Using lemma~\ref{L:xi_in_Ekj} we obtain that $z, y, v \in E(g)$. If we define $g_y = k[x:N]$ then for each $\sigma \in \Xi(g)$ $\sigma_{/dom(g_y)} \in \Xi(g_y)$ and $\#(g,y,\sigma) \in \#(g_y, N, \sigma_{/dom(g_y)}) = \#(N) = \mathbb{N}$.

\medskip 

Moreover for each $\sigma \in \Xi(g)$ $\sigma_{/dom(\kappa)} \in \Xi(\kappa)$ $\#(g,v,\sigma) \in \#(\kappa, N, \sigma_{/dom(\kappa)}) = \#(N) = \mathbb{N}$.

\medskip

By lemma~\ref{L:preliminary-results-ch6-3} it follows that $(*)(y, v) \in E(g)$ and for each $\sigma \in \Xi(g)$ $\#(g, (*)(y, v), \sigma) = (\#(g,y, \sigma) \cdot \#(g,v, \sigma)) \in \mathbb{N}$.

\medskip

By lemma~\ref{L:preliminary-results-ch6-2} $=(z, yv)$ belongs to $S(g)$. We can now apply lemma~\ref{L:quantifiers_S_h} and obtain that $\exists( v:N, =(z, yv) ) \in S(\kappa)$.

\medskip

Then if we apply rule~$R_{\ref{L:universal-to-existential-second}}$ to \ref{E:ex3-o} we obtain

\scriptsize
\begin{equation}\label{E:ex3-p}
 \gamma \left[x:N, y:N, z:N, \to \left( 
\exists( u:N, =(y,xu) ) , 
\to \left(
\begin{array} {l}
{ \exists( v:N, =(z,yv) ), } \\
{ \exists( w:N, =(z, xw) )   }
\end{array}
\right)
\right)
\right] .
\end{equation}
\normalsize

We have also $\exists( u:N, =(y,xu) ) \in S(\kappa)$, so if we apply rule~$R_{\ref{L:rule-il-to-cl-new}}$ we finally obtain

\scriptsize
\begin{equation}\label{E:ex3-q}
 \gamma \left[x:N, y:N, z:N, \to \left( 
 \wedge \left( 
\begin{array} {l}
{ \exists( u:N, =(y,xu) ) } \\
{ \exists( v:N, =(z,yv) )   }
\end{array}
 \right)
 , 
\exists( w:N, =(z, xw) )  
\right)
\right] .
\end{equation}
\normalsize

\section{Thing we can express in our system}\label{Ch:expMixed}

In this section we will wee some interesting things we can express with our approach. First of all we mentioned in the introduction that in our system we can express statements in which both quantifiers over individuals and quantifiers over sets of individuals occur. We made the simple example of the following statement:\\

for each subset X of $\mathbb{N}$ and for each x $\in \mathbb{N}$ we have x $\in$ X or x $\notin$ X .\\

Let's see how we can map the statement within our system. In our language we need two constants: $N$ whose meaning is the set $\mathbb{N}$ of natural numbers, $P$ which represents the set of the subsets of $\mathbb{N}$ i.e. $\mathcal{P}(\mathbb{N})$. In our system we have a constraint that constants cannot represent the empty set or a set that contains the empty set, so in our case the empty set is not a member of $\mathcal{P}(\mathbb{N})$.\\

The set $\mathcal{F}$ of operators is the same we have assumed in our other examples, so it must contain all of these symbols: $\neg, \wedge, \vee, \to, \leftrightarrow, \in, =$.\\

So, in order to formalize our statement and a proof of it, we will use a language $(\mathcal{V}, \mathcal{F}, \mathcal{C}, \#, \{D_1, \dots, D_p \})$ which must be as follows
\begin{align} 
&\mathcal{V} = \{x, X \} \ , \notag  \\
&\mathcal{F} = \{ \neg, \wedge, \vee, \to, \leftrightarrow, \in, =  \} \ , \notag  \\
&\mathcal{C} = \{ N, P  \} \  \notag .
\end{align}

\smallskip

Moreover, in order to write down this statament we do not need to include the set $\mathbb{N}$ or another set in our additional sets, so let $p = 0$.\\

Here we notice that we also set as a constraint that for each $c \in \mathcal{C}$ and for any positive integer $q$ we must be able to decide all of the following conditions

\begin{itemize}
\item $Set_q( \#(c))$;
\item $Event_q( \#(c))$;
\item $\#(c) \in D_i$;
\item $\#(c) \in \mathcal{P}^q(D_i)$;
\item if ($Set_q( \#(c))$) then ($NotEmpty_q( \#(c))$).
\end{itemize}

And moreover, the last condition must be decided as true.\\

For the constant $N$ we can take the following decisions: 

\begin{itemize}
\item $Set_1( \#(N))$: true;
\item for $q > 1$, $Set_q( \#(N))$: false;
\item $Event_q( \#(N))$: false;
\item if ($Set_q( \#(N))$) then ($NotEmpty_q( \#(N))$): true.
\end{itemize}

\medskip

For the constant $P$ we can take the following decisions: 

\begin{itemize}
\item $Set_1( \#(P))$: true;
\item $Set_2( \#(P))$: true;
\item for $q > 2$, $Set_q( \#(P))$: false;
\item $Event_q( \#(P))$: false;
\item if ($Set_q( \#(P))$) then ($NotEmpty_q( \#(P))$): true.
\end{itemize}

\medskip

With this setup, we can express the statement as follows
\[ \gamma[ x: N, X: P, \vee ( \in(x, X), \neg( \in(x, X)) ) ] \ . \]

\medskip

Let's now verify this is an expression of our language.\\

First of all, by lemma \ref{L:xi-in-D-new}, $H[x: N, X: P]$ holds.\\

Let now $k = k[x: N, X: P]$, we try to show that $\in(x, X) \in S(k)$.\\

Using lemma~\ref{L:xi_in_Ekj} we obtain that $x,X \in E(k)$.\\

Moreover, let $h = k[x: N]$, then for each $\sigma \in \Xi(k)$, $\sigma_{/dom(h)} \in \Xi(h)$ and $\#(k, X, \sigma) \in \#(h, P, \sigma_{/dom(h)}) = \mathcal{P} (\mathbb{N})$. Therefore $\#(k, X, \sigma)$ is a set, we can apply lemma~\ref{L:in-t-varphi-in-Sk} and obtain that $\in(x, X) \in S(k)$.\\

As a consequence of this $\vee ( \in(x, X), \neg( \in(x, X)) ) \in S(k)$ and finally
\[ \gamma[ x: N, X: P, \vee ( \in(x, X), \neg( \in(x, X)) ) ] \in S(\epsilon) \ . \]

We will now see other interesting things we can express in our approach, which are related to the liar paradox that we will discuss later in the manuscript.\\

First of all, given a set $A$, we need to be able to express the condition $\delta$ = `for each $x$ in $A$ $x$ is false'. And we also wat to express the condition `$\delta$ belongs to $A$'.\\

Let's now see how we can express the mentioned conditions in our system. Suppose we have a language $(\mathcal{V}, \mathcal{F}, \mathcal{C}, \#, \{D_1, \dots, D_p \})$ which must be as follows
\begin{align} 
&\mathcal{V} = \{ x \} \ , \notag  \\
&\mathcal{F} = \{ \neg, \wedge, \vee, \to, \leftrightarrow, \in, =  \} \  \notag .
\end{align}

We can have in our language a constant $A$ which represents a set, but we refer to the more general case of an expression $\psi \in E_s(\epsilon)$.\\

At this point we have $H[x:\psi]$ and if we define $k = k[x:\psi]$ then by lemma~\ref{L:xi_in_Ekj} $x \in E(k)$. Clearly by lemma~\ref{L:meaning-kept-f1} $\neg(x) \in E(k)$.\\

Given $\sigma \in \Xi(k)$ $\#(k, \neg(x), \sigma) = (\#(k, x, \sigma) \text{ is false})$, so $\#(k, \neg(x), \sigma)$ itself is true or false, and this implies that $\neg(x) \in S(k)$.\\

Given that $\epsilon \in K$, $\psi \in E_s(\epsilon)$, $x \in \mathcal{V} - var(\epsilon)$, $k = \epsilon + <x, \psi>$ we can apply lemma~\ref{L:quantifiers_S_h} and obtain that $\forall(x: \psi, \neg(x)) \in S(\epsilon)$.\\

Let $\delta = \forall(x: \psi, \neg(x))$, we also want to show that $\in( \delta , \psi ) \in S(\epsilon)$.\\

We have that $\#( \epsilon, \psi, \epsilon)$ is a set, so $A_{\in}(\#( \epsilon, \delta, \epsilon), \#( \epsilon, \psi, \epsilon))$, and by lemma~\ref{L:meaning-kept-f2} $\in( \delta , \psi ) \in E(\epsilon)$.\\

Moreover $\#( \epsilon, \in( \delta , \psi ), \epsilon) = \#( \epsilon, \delta, \epsilon) \text{ belongs to } \#( \epsilon, \psi, \epsilon)$, so $\#( \epsilon, \in( \delta , \psi ), \epsilon)$ is true or false, so $\in( \delta , \psi ) \in S(\epsilon)$.

\section{Further study}\label{Ch:conspar}

Of course, further investigations about our approach to logic can be performed. We have mentioned in section~\ref{sec:completeness} the topic on the completeness or incompleteness of our deductive systems. Then we have introduced some example of a deductive system. Some questions that I have not investigated in depth are the following:
\begin{itemize}
\item can we describe a deductive system within our logic system as a recursively axiomatised formal system?
\item given a language that does not include arithmetic, under which conditions, if any, a deductive system within our logic system is complete?
\end{itemize} 

\medskip

Another interesting (and not extremely easy) topic is about comparing the expressive power of our system with the one of standard logic systems.\\

Another topic to consider is substitution. First-order logic features the notion of `substitution' (see e.g. Enderton's book~\cite{Enderton}). Under appropriate assumptions, we can apply substitution to a formula $\varphi$ and obtain a new formula $\varphi^x_t$, by replacing the free occurrences of the variable $x$ by the term $t$. In our approach we could be able to define a similar notion, with the difference that for us $t$ could be a generic expression. I have somehow studied how the topic of substitution could be applied to this type of system, but with respect to a former version of my system. I am rather confident that general substitution mechanisms can be introduced for this type of logic, but I'm not sure how much work this would require. After all I suppose the introduction of general substitution mechanisms could be considered as not being properly a core topic about this approach, since for instance we can use simplified substitution mechanisms.\\

Finally, let's also briefly talk about paradoxes. A paradox is usually a situation in which a contradiction or inconsistency occurs, in other words a paradox arises when we can build a sentence $\varphi$ such that both $\varphi$ and $\neg (\varphi)$ can be derived. Since our system is consistent it shouldn't be possible to have true paradoxes in it. If we have proved the consistency of our system, what can we do more than this to exclude that the system is vulnerable to paradoxes?\\

It could anyway not be wrong to discuss some of the most known paradoxical arguments to ask ourselves if our system could be vulnerable to one of them.\\

We begin with Russell's paradox. Assume we can build the set $A$ of all those sets $X$ such that $X$ is not a member of $X$. Clearly, if $A \in A$ then $A \notin A$ and conversely if $A \notin A$ then $A \in A$. We have proved both $A \in A$ and its negation, and this is the Russell's paradox.\\

Postulating that a certain set contains all sets that don't belong to themselves would mean stating something like this:\\

For every set $X$, $X \in A$ if and only if $X \notin X$.\\

And here we're using a quantifier on the class of `all sets', the existence of which is not assumed in our system.

\vspace{18pt}

We also want to examine the liar paradox and related topics. Let's consider how the paradox is stated in Mendelson's book.\\

A man says, `I am lying'. If he is lying, then what he says is true, so he is not lying. If he is not lying, then what he says is false, so he is lying. In any case, he is lying and he is not lying.\\

Mendelson classifies this paradox as a `semantic paradox' because it makes use of concepts which
need not occur within our standard mathematical language. I agree that, in his formulation, the
paradox has some step which seems not mathematically rigorous.\\

We'll try to provide a more rigorous wording of the paradox.\\

Let $A$ be a set, and let $\delta$ be the condition `for each $x$ in $A$ $x$ is false'. Suppose $\delta$ is the only member of $A$. In this case if $\delta$ is true then it is false; if on the contrary $\delta$ is false then it is true.\\

The explanation of the paradox is the following: simply $\delta$ cannot be the only item in set $A$. In fact, suppose $A$ has only one element, and let's call it $\varphi$. This implies $\delta$ is equivalent to `$\varphi$ is false' so it seems
acceptable that $\delta$ is not $\varphi$.\\

Another approach to the explanation is the following.\\

If $\delta$ is true then for each $x$ in $A$ $x$ is false, so $\delta$ is not in $A$. By contraposition if $\delta$ is in $A$ then $\delta$ is false.\\

Moreover if the uniqueness condition `for each $x$ in $A$ $x=\delta$' is true, then $\delta$ can be true or false, but if it is false then it is true, so it is true in both cases.\\

Therefore if $\delta$ is the only element in $A$ then $\delta$ is true and false at the same time. This implies $\delta$
cannot be the only item in $A$.\\

On the basis of this argument I consider the liar paradox as an apparent paradox that actually has an explanation. What is the relation between our approach to logic and the liar paradox?\\

Standard logic isn't very suitable to express this paradox. In fact first-order logic is not designed to construct a condition like our condition $\delta$ (= `for each $x$ in $A$ $x$ is false'), and moreover, it is clearly not designed to say `$\delta$ belongs to set $A$'. These conditions aren't plainly leading to inconsistency, so it is desirable they can be expressed in a general approach to logic. And our system permits to express them, as we have seen in chapter~\ref{Ch:expMixed}. The paradox isn't ought to simply using these conditions, it is due to an assumption that is clearly false, and the so-called paradox is simply the proof of its falseness.

\bigskip

Related to the liar paradox is the Cretan `paradox', which is actually not a proper paradox, but is perhaps even more `unsettling' and we quote again Mendelson in this regard: (\cite{Mendelson}).\\

\begin{quotation}
The Cretan “paradox”, known in antiquity, is similar to the Liar Paradox. The Cretan philosopher Epimenides said, “All Cretans are liars”. If what he said is true, then, since Epimenides is a Cretan, it must be false. Hence, what he said is false. Thus, there must be some Cretan who is not a liar. This is not logically impossible, so we do not have a genuine paradox. However, the fact that the utterance by Epimenides of that false sentence could imply the existence of some Cretan who is not a liar is rather unsettling.
\end{quotation}

\medskip

If we try to put this argument in a more formal statement, it still refers to a sentence $\delta$ of the type `for each $x$ in $A$ $x$ is false', where this time $A$ is the set of all the statements made by a Cretan and $\delta$ is a member of $A$. Here if  $\delta$ is true then it is false, so we have to conclude that $\delta$ is false, hence there exists $x \in A$ such that $x$ is true. As noticed by Mendelson, it can be unsettling to accept this just because $\delta$ is a member of $A$.\\

We can still use an argument we have shown above with respect to the liar paradox: If $\delta$ is true then for each $x$ in $A$ $x$ is false, so $\delta$ is not in $A$. By contraposition if $\delta$ is in $A$ then $\delta$ is false. And another formulation is the following: $\delta$ is false or $\delta$ is not in $A$.\\

Let $A$ be a set of true/false statements (think to an actual list of statements) and $\delta$ be the statement `for each $x$ in $A$ $x$ is false'. We know from the discussion on the liar paradox that if $A$ has just one element then $\delta$ cannot belong to $A$.\\

In the case of the Cretan paradox we have that $\delta$ could belong to $A$ and there is not a constraint that $A$ has just one element. Is it possible in this case that $\delta$ belongs to $A$? The basic problem is that $\delta$, if it belongs to $A$, makes a reference to itself and this can lead us to suspect that $\delta$ in this case is not something well defined.\\

We could therefore conclude that also in this case it cannot be accepted that $\delta$ belongs to $A$. In this case we could `resolve' the problem by using axioms like
\[ \neg( \in(\forall( x:\psi, \neg (x)), \psi) )  \ , \]

for each expression $\psi$ that represents a set.\\

If instead we accept the possibility that $\delta$ belongs to $A$ it is evident that we must also accept that if $\delta$ belongs to $A$ then it is false, in fact if it were true then it would not belong to $A$.\\

As a conclusion, with respect to paradoxes, we cannot state that our system is designed to prevent for sure every possible form of paradox, for instance it doesn't prevent anyone to conceive something which is unsettling or contradictory. Anyway although I have made some assessments on the matter, I currently have no reason to suppose that the system is subject to some paradox.

\section{Declaration of Generative AI and AI-assisted technologies in the writing process}\label{Ch:decIA}

During the preparation of this manuscript, the author used artificial intelligence tools solely for the purpose of English language editing, stylistic refinement, and grammatical correction of the text. After using this service, the author reviewed and edited the content as needed and takes full responsibility for the framework, mathematical logic, and final content of the publication.

\end{document}